\documentclass[11pt,reqno]{amsart}
\usepackage[T1]{fontenc}
\usepackage[utf8]{inputenc}
\usepackage{lmodern,amsmath,amssymb,amsthm,mathtools,mathrsfs}
\usepackage[a4paper,margin=28mm,headheight=14pt]{geometry}
\usepackage[final,expansion=false]{microtype}
\usepackage{xcolor}
\colorlet{BLACK}{black}
\usepackage[colorlinks=true,linkcolor=blue,citecolor=blue,urlcolor=blue]{hyperref}
\usepackage{bookmark}
\hypersetup{pdfpagemode=UseOutlines}
\bookmarksetup{numbered=true,depth=2,open=true,openlevel=2}
\usepackage{enumitem,needspace}
\setlist[enumerate]{leftmargin=2.2em,itemsep=2pt,topsep=4pt}
\setlist[itemize]{leftmargin=2.2em,itemsep=2pt,topsep=4pt}
\numberwithin{equation}{section}
\allowdisplaybreaks[2]
\newtheorem{theorem}{Theorem}[section]
\newtheorem{proposition}[theorem]{Proposition}
\newtheorem{lemma}[theorem]{Lemma}
\newtheorem{corollary}[theorem]{Corollary}
\theoremstyle{definition}\newtheorem{remark}[theorem]{Remark}
\newtheorem*{pcassumption}{Assumption (PC)}

\makeatletter
\@namedef{ergm@outline@Theorem}{}
\@namedef{ergm@outline@Lemma}{}
\@namedef{ergm@outline@Proposition}{}
\@namedef{ergm@outline@Corollary}{}
\let\ergm@original@begintheorem\@begintheorem
\def\@begintheorem#1#2[#3]{%
  \@ifundefined{ergm@outline@#1}{}{%
    \@ifempty{#3}{%
      \bookmark[level=2,dest={\@currentHref}]{#1 #2}%
    }{%
      \bookmark[level=2,dest={\@currentHref}]{#1 #2 (#3)}%
    }%
  }%
  \ergm@original@begintheorem{#1}{#2}[#3]%
}
\makeatother

\newcommand{\E}{\mathop{\mathbb{E}}\nolimits}\newcommand{\Pp}{\mathop{\mathbb{P}}\nolimits}
\newcommand{\R}{\mathbb R}\newcommand{\one}{\mathbf1}
\newcommand{\dd}{\mathop{}\!\mathrm d}\newcommand{\sig}{\sigma}
\newcommand{\op}{\mathrm{op}}
\newcommand{\dW}{d_{\mathrm W}}\newcommand{\dK}{d_{\mathrm{Kol}}}
\DeclareMathOperator{\inj}{inj}\DeclareMathOperator{\Hom}{Hom}\DeclareMathOperator{\Var}{Var}
\DeclareMathOperator{\Cov}{Cov}\DeclareMathOperator{\cum}{cum}
\DeclareMathOperator{\tr}{tr}\DeclareMathOperator{\Lip}{Lip}
\newcommand{\mc}{m_\star}
\newcommand{\proofstep}[1]{\par\smallskip\noindent\textbf{#1}\enspace\ignorespaces}
\hypersetup{pdftitle={Quantitative non-Gaussian and Gaussian approximation for exponential random graphs in the single-phase critical regime},pdfauthor={Ji-Shuai Miao, Yi-Kun Zhao},pdfsubject={Critical edge counts and quantitative non-Gaussian and Gaussian approximation}}
\title[Non-Gaussian and Gaussian approximation at single-phase criticality]{Quantitative non-Gaussian and Gaussian approximation for exponential random graphs in the single-phase critical regime}
\author[Ji-Shuai Miao and Yi-Kun Zhao]{\NoCaseChange{%
  \begin{tabular}{@{}c@{}}
    {\large Ji-Shuai Miao\textsuperscript{1},\quad
      Yi-Kun Zhao\textsuperscript{1}}\\[0.8ex]
    {\small\itshape \textsuperscript{1}Department of Statistics and Data Science, The Chinese University of Hong Kong}
  \end{tabular}}}
\date{}
\begin{document}
\begin{abstract}
We establish quantitative non-Gaussian and Gaussian approximations for
edge-count fluctuations in dense exponential random graphs at single-phase
criticality of any fixed even order $2\mc\ge4$. The Hamiltonians use
injective homomorphism counts of fixed simple graphs, with nonnegative
leading coefficients for the higher-order terms. In bounded calibrated
windows, edge-density fluctuations on the scale $n^{-1/\mc}$ admit an
explicit non-Gaussian approximation with error $O(n^{-1/\mc})$. For
parameter perturbations by fixed multiples of $n^{-1}$ with a nonzero
first-order drift correction at the critical density, the deterministically
centered edge count admits a Gaussian approximation on the scale
$n^{(3\mc-2)/(2\mc-1)}$, with error $O(n^{-1/(2\mc-1)})$. Both rates
hold in Wasserstein and Kolmogorov distances. The key step is to construct
auxiliary kernels for exterior graphs outside an event of exponentially
small probability. These kernels contract uniformly on the full row state
space and agree with the conditional Gibbs updates on typical row
configurations, yielding the conditional cumulant estimates used in
Stein's method.
\end{abstract}

\maketitle
\begingroup\small\tableofcontents\endgroup
\clearpage

\section{Introduction}\label{sec:introduction}
Dense exponential random graph models are Gibbs measures whose Hamiltonians
are finite linear combinations of suitably scaled subgraph counts.
Their asymptotic behavior is closely related to variational problems
and phase transitions. The graph large deviation principle of Chatterjee and
Varadhan~\cite{CV} and the variational analysis of Chatterjee and
Diaconis~\cite{CD} identify limiting free energies and typical graph structure.
For nonnegative higher-order coefficients, the variational problem reduces
to a one-dimensional maximization problem. Concentration in the cut metric,
however, does not determine edge-count fluctuations: lower-order terms absent from the limiting free
energy can affect their distribution.

Glauber dynamics provide complementary tools. Bhamidi, Bresler and
Sly~\cite{BBS} proved rapid mixing in a high-temperature region and
exponentially slow mixing in a low-temperature region; Reinert and
Ross~\cite{RR} developed Stein comparisons for stationary Glauber laws.
Ganguly and Nam~\cite{GN} established subcritical Poincar\'e inequalities
and Gaussian concentration beyond uniform Dobrushin contraction.
Fang, Liu, Shao and Zhao~\cite{FLSZ} obtained subcritical edge-count central
limit theorems with explicit variance and quantitative Wasserstein and
Kolmogorov bounds, also treating subgraph counts.
Winstein~\cite[Theorem~1.1]{WinsteinCLT} established edge-count central limit
theorems with explicit asymptotic variance and
$O_\varepsilon(n^{-1/2+\varepsilon})$ Wasserstein and Kolmogorov bounds
for every fixed $\varepsilon>0$ throughout the noncritical single-phase
regime, including the supercritical (low-temperature) region. In this regime, the variational function has a unique global
maximizer, which is nondegenerate. The result
concerns the full ERGM measure.
He also treated subgraph counts, vertex degrees and local subgraph counts,
with conditional extensions to noncritical metastable wells at phase
coexistence under additional forest assumptions.
At degenerate criticality, the quadratic restoring force vanishes,
requiring a different fluctuation scale and comparison law.

For the two-star model, Mukherjee and Xu~\cite{MX} determined edge-count
limits in the uniqueness, coexistence and critical regimes. Their
Theorem~1.4(c) gives critical fluctuations on the $n^{3/2}$ scale with both
quadratic and quartic terms in the log-density. Their auxiliary-variable
representation uses the quadratic interaction to obtain conditionally
independent edges, a structure not assumed here.
For the edge--triangle model, Bianchi, Collet and Magnanini~\cite{BCM}
studied fluctuations away from coexistence and the critical endpoint.
Their quartic critical theorem (Theorem~8.8) concerns a mean-field approximation.
They conjectured the corresponding limit for the original graph law
(Conjecture~3.8).
These results motivate a theory of critical fluctuations that retains
the original interactions and their finite-size corrections.

We establish quantitative non-Gaussian and Gaussian approximations for
edge-count fluctuations in exponential random graphs in the critical regime.
The Hamiltonians in \eqref{eq:model} are built from counts of injective
homomorphisms of fixed simple graphs. The leading interaction coefficients
are nonnegative.
The interaction graphs may be disconnected and may contain cycles.
We assume that the variational function $\Phi$ in \eqref{eq:scalar} has a unique global
maximizer $p\in(0,1)$. This is the single-phase condition used here
(cf.\ Winstein~\cite[Section~1.1.2]{WinsteinCLT}). We further assume that,
for a fixed $\mc\ge2$,
\[
 \Phi^{(j)}(p)=0\quad(1\le j<2\mc),\qquad
 a_\star=-\frac{\Phi^{(2\mc)}(p)}{(2\mc-1)!}>0.
\]
Let $E$ denote the edge count and $P=E/N$ denote the edge density, where
$N=\binom n2$.

The polynomial $B(x)$ in \eqref{eq:scalar} depends on the interaction graphs
only through their edge counts. The finite-size correction $\mathcal F_1(x)$
in \eqref{eq:correction} also depends on their numbers of vertices,
adjacent-edge pairs and triangles. Write $b_{r,n}=b_r+u_{r,n}/n$ and
$h_n=h_*+\kappa_n/n$ with bounded corrections as in \eqref{eq:window}, and set
\[
 \mathcal U_n(x)=\sum_r u_{r,n}x^{e_r},\qquad
 C_n(x)=\mathcal F_1(x)+\mathcal U_n(x)+\kappa_n x/2,\qquad
 \beta_{j,n}=n^{1-j/\mc}C_n^{(j)}(p).
\]
The calibrated window \eqref{even:calibration-system} requires the
$\beta_{j,n}$ to be uniformly bounded for $1\le j\le\mc$.

In the calibrated window, Theorem~\ref{thm:main} approximates
$W_n=n^{1/\mc}(P-p)$ by the distribution with density
\[
 \rho_n(z)\propto
 \exp\left\{-\frac{a_\star z^{2\mc}}{2\mc}
                 +\sum_{j=1}^{\mc}\frac{\beta_{j,n}}{j!}z^j\right\}
\]
at rate $O(n^{-1/\mc})$ in the Wasserstein and Kolmogorov distances.
The bound is uniform over parameter sequences satisfying
\eqref{eq:window} and \eqref{even:calibration-system} with prescribed
bounds. The comparison distribution may be asymmetric and may depend on
$n$.

For fixed corrections $u_{r,n}=u_r$ and $\kappa_n=\kappa$ with
$C_n'(p)\ne0$, Theorem~\ref{fixed:main} gives a Gaussian approximation
at rate $O(n^{-1/(2\mc-1)})$ in the same distances. The edge count is
centered at $Np_n$, where the deterministic density $p_n$ differs from
$p$ by order $n^{-1/(2\mc-1)}$, and its fluctuations are of order
$n^{(3\mc-2)/(2\mc-1)}$. This case includes the unperturbed critical
parameters, for which $C_n'(p)=\mathcal F_1'(p)<0$.

In both regimes, the theorems identify the leading means and variances of
the edge count. Corollaries~\ref{even:subgraph-fluctuations}
and \ref{fixed:subgraph-gaussian} extend the approximations and the moment
asymptotics to counts of fixed simple graphs with at least one edge.
Corollary~\ref{even:standardized-counts} gives quantitative bounds after
standardization by the exact means and standard deviations.

In the quartic critical regime, the two theorems cover all fixed correction
coefficients. Indeed, $\beta_{2,n}=C_n''(p)$ is fixed and
$\beta_{1,n}=\sqrt n\,C_n'(p)$. Thus $C_n'(p)=0$ yields a quartic limit
by Theorem~\ref{thm:main}, while $C_n'(p)\ne0$ yields a Gaussian limit by
Theorem~\ref{fixed:main}.

The first calibration equation gives
\[
 h_n=h_*-\frac{2\{\mathcal F_1'(p)+\mathcal U_n'(p)\}}n
                         +2\beta_{1,n}n^{-2+1/\mc}.
\]
This separates the $O(n^{-1})$ shift of the external field from the
critical window of width $n^{-2+1/\mc}$. Proposition~\ref{even:bounded-calibration}
shows that the calibration equations admit uniformly bounded interaction
and field corrections. In the quartic critical regime, calibrating the field
suffices for any bounded interaction corrections. At higher critical
orders, the interaction corrections may also require calibration.

The proofs of Theorems~\ref{thm:main} and~\ref{fixed:main}
follow a common strategy, which we outline below. Both proofs use Stein's
method with an exchangeable pair constructed from a single-edge Gibbs
update. The main inputs are a precise drift expansion, estimates of the
conditional second moment of the increment, and moment bounds at the
relevant fluctuation scales. To obtain these estimates, we first establish
uniform cumulant bounds for the conditional adjacency rows in the critical
regime. The quadratic gap for the
conditional row potential (Proposition~\ref{row:prop:row-gap}) and the
localization estimates allow us to construct auxiliary update kernels
that contract uniformly on the full row state space and agree with
the Gibbs updates on typical row configurations for typical exterior
graphs (Lemma~\ref{lem:local-input}). The resulting estimates under the
conditional Gibbs laws, together with signed subgraph count bounds,
yield the drift expansion for the edge density.

To obtain the required precision in the drift expansion, we exploit
cancellations in the conditional moments of centered row sums. For
example, the conditional third moment of such a sum is $O(n)$, whereas
its third absolute moment can be $O(n^{3/2})$.
Rooted expansions, four-cycle bounds and vertex transpositions control
exterior replacements. Cumulant induction removes logarithmic losses,
and a normalized change of measure handles two-row interactions.
Invariance under a single-edge Gibbs update gives moment bounds at the
relevant fluctuation scales. We then apply Stein's method for exchangeable
pairs. See \cite{CS,EL} for related critical laws.

The paper is organized as follows. Section~\ref{sec:main} introduces
the model and states the main results. Section~\ref{row:sec:conditional-rows}
develops the conditional row estimates. Sections~\ref{sec:counts}--\ref{sec:supports}
establish the subgraph-count and logit expansions needed for the drift
analysis. Section~\ref{sec:drift} states the Stein approximation bound and
combines these estimates to prove the main results. The proofs of
Theorems~\ref{thm:main} and~\ref{fixed:main} are given in
Subsections~\ref{sec:proof-main} and~\ref{sec:proof-fixed-critical-normal},
respectively, and the corollary proofs are collected in
Subsection~\ref{sec:proof-corollaries}.

The appendices contain the proofs of the auxiliary results used in the
main text. Appendix~\ref{row:sec:row-gap} establishes localization estimates
and constructs contractive auxiliary kernels for the conditional row laws.
Appendices~\ref{app:resolvent} and~\ref{ep:row:section} develop the resolvent
estimates and conditional cumulant bounds. Appendix~\ref{app:stein-proofs}
establishes estimates for solutions of the Stein equation and proves the
approximation theorem for exchangeable pairs.
Appendix~\ref{app:calibration-proof} proves the existence of uniformly bounded
solutions to the calibration equations.

\section{Model and main results}\label{sec:main}\label{sec:assembly}
\subsection{The model and homomorphism counts}\label{sec:graph-law}
Let $X$ be a symmetric $\{0,1\}$ adjacency matrix with zero diagonal.
Fix simple graphs $F_1,\ldots,F_k$ without isolated vertices, with
$e_r=e(F_r)\ge2$ and $v_r=v(F_r)$. Define labelled, non-induced counts by
\begin{equation}\label{v8:def:graph-counts}
 \begin{aligned}
 \inj(F,X)&=\sum_{\varphi:V(F)\hookrightarrow[n]}
                \prod_{uv\in E(F)}X_{\varphi(u)\varphi(v)},\\
 \Hom(F,X)&=\sum_{\varphi:V(F)\to[n]}
                \prod_{uv\in E(F)}X_{\varphi(u)\varphi(v)}.
 \end{aligned}
\end{equation}
The same formulas apply to symmetric
zero-diagonal weighted matrices. Put
\begin{equation}\label{v8:def:graph-edge-degree-density}
 E=\sum_{i<j}X_{ij},\qquad D_i=\sum_{j\ne i}X_{ij},\qquad
 N=\binom n2,\qquad P=E/N.
\end{equation}
For the two-star $P_3$ and the triangle $K_3$, the counts are
$\inj(P_3,X)=\sum_iD_i(D_i-1)$ and
$\inj(K_3,X)=6\sum_{i<j<k}X_{ij}X_{ik}X_{jk}$.
For bounded parameter sequences, ERGM assigns probability
\begin{equation}\label{eq:model}
 \begin{aligned}
 \Pp_n(X=x)&=(Z_n^{\mathrm{graph}})^{-1}\exp\{\mathcal H_n(x)\},\\
 \mathcal H_n(X)&=h_nE+\sum_{r=1}^k b_{r,n}n^{2-v_r}\inj(F_r,X).
 \end{aligned}
\end{equation}
Here $\mathcal H_n$ is the Hamiltonian and $Z_n^{\mathrm{graph}}$ is the normalizing constant
over all labelled simple graphs on $[n]$; $n^{2-v_r}$ scales each interaction
to order $n^2$.

\subsection{Variational problem and the critical regime}
Fix $b_r\ge0$, $h_*\in\R$, and an integer $\mc\ge2$, all independent of $n$. Define
\begin{equation}\label{eq:scalar}
 \begin{aligned}
 B(x)&=\sum_r b_rx^{e_r},\qquad
 I(x)=x\log x+(1-x)\log(1-x),\\
 \Phi(x)&=\frac{h_*x}{2}+B(x)-\frac{I(x)}2,
 \qquad 0\log0=0.
 \end{aligned}
\end{equation}
Assume that $p\in(0,1)$ satisfies
\begin{equation}\label{eq:critical}
 \begin{gathered}
 \Phi(p)>\Phi(x)\quad(x\in[0,1]\setminus\{p\}),\qquad
 \Phi^{(j)}(p)=0\quad(1\le j\le2\mc-1),\\
 a_\star=-\frac{\Phi^{(2\mc)}(p)}{(2\mc-1)!}>0.
 \end{gathered}
\end{equation}
The degeneracy order $2\mc$ is distinct from the exterior size $m$;
its parity is explained in Remark~\ref{even:critical-order-parity}. Set
\begin{equation}\label{v8:def:scalar-logistic-field}
 \sig(u)=\frac{e^u}{1+e^u},\qquad
 \ell(x)=h_*+2B'(x),\qquad g(x)=\sig(\ell(x))-x,
\end{equation}
and, for the actual finite-size coefficients, put
\begin{equation}\label{v8:def:finite-size-scalar-logit}
 \ell_n(x)=h_n+2\sum_r b_{r,n}e_rx^{e_r-1}.
\end{equation}

\subsection{Finite-size correction and critical calibration}\label{sec:parameter-ranges}
For a template $F$, let
\[
 a(F)=\sum_{u\in V(F)}\binom{d_F(u)}2,
 \qquad \triangle(F)=\#\{\text{unordered triangles in }F\}.
\]
Thus $a(F)$ counts unordered pairs of adjacent edges. Define the fixed profiles
\begin{equation}\label{eq:profiles}
 \begin{aligned}
 t(x)&=\sum_r b_ra(F_r)x^{e_r-2},\qquad
 \tau(x)=\sum_{r:e_r\ge3}b_r\triangle(F_r)x^{e_r-3},\\
 s(x)&=\sum_r b_r\left[\binom{v_r}2-1\right]x^{e_r},
 \qquad \nu(x)=x(1-x).
 \end{aligned}
\end{equation}
The coefficient in $s(x)$ comes from expanding $(n-2)_{v_r-2}/n^{v_r-2}$
(Proposition~\ref{prop:microscopic}). Since $a(F_r)\le\binom{e_r}2$
and $\Phi''(p)=0$,
\[
 2\nu(p)t(p)\le\nu(p)B''(p)=\frac12.
\]
Choose a fixed interval $J\Subset(0,1)$ containing $p$ on which
$1-2\nu(x)t(x)\ge c_J>0$. On $J$, set
\begin{equation}\label{eq:correction}
 \mathcal F_1(x)=-s(x)-\nu(x)t(x)
 -\frac12\log\{1-2\nu(x)t(x)\}+3\tau(x)^2\nu(x)^3.
\end{equation}
This correction applies at every fixed critical order; its derivatives are taken with
the leading coefficients held fixed.

\begin{pcassumption}\phantomsection\label{ass:PC}
\hypertarget{assumption-PC}{}
\begin{enumerate}[label=\textup{(\roman*)},leftmargin=2.4em,itemsep=3pt]
\item The law is \eqref{eq:model}, with fixed simple templates having
at least two edges and no isolated vertices.
\item The leading coefficients satisfy $b_r\ge0$. The variational function $\Phi$ in
\eqref{eq:scalar} has a unique global maximizer $p\in(0,1)$ (the single-phase
condition). The conditions \eqref{eq:critical} defining the critical regime hold at $p$ for a
fixed integer $\mc\ge2$.
\item For a fixed $K<\infty$,
\begin{equation}\label{eq:parameter-neighborhood}
 |h_n-h_*|+\sum_r|b_{r,n}-b_r|\le K/n.
\end{equation}
\end{enumerate}
\end{pcassumption}

The conditional row estimates below are established under
Assumption~\hyperref[ass:PC]{\textup{(PC)}}. Their constants may depend on
$K$, but are independent of the graph size. Theorem~\ref{thm:main}
additionally requires the calibration condition \eqref{even:calibration-system}.

The bounded neighborhood \eqref{eq:parameter-neighborhood} can equivalently be written as
\begin{equation}\label{eq:window}
 \begin{gathered}
 b_{r,n}=b_r+\frac{u_{r,n}}n,\qquad
 h_n=h_*+\frac{\kappa_n}{n},\qquad
 \mathcal U_n(x)=\sum_r u_{r,n}x^{e_r},\\
 \sup_n\left(\sum_r|u_{r,n}|+|\kappa_n|\right)\le K<\infty.
 \end{gathered}
\end{equation}
For every fixed derivative order $j$,
\begin{equation}\label{even:window-Cr}
 \ell_n-\ell=\frac{\kappa_n+2\mathcal U_n'}n,
 \qquad \|\ell_n-\ell\|_{C^j([0,1])}\le C_j/n.
\end{equation}

For a sequence in \eqref{eq:window}, define its correction and induced coefficients by
\begin{equation}\label{even:calibration}
 \begin{gathered}
 C_n(x)=\mathcal F_1(x)+\mathcal U_n(x)+\kappa_nx/2,\qquad
 \beta_{j,n}=n^{1-j/\mc}C_n^{(j)}(p)\quad(1\le j\le\mc).
 \end{gathered}
\end{equation}
Equivalently,
\begin{equation}\label{even:induced-coefficients}
 \begin{aligned}
 \beta_{1,n}
 &=n^{1-1/\mc}\left\{\mathcal F_1'(p)
       +\sum_{r=1}^k e_ru_{r,n}p^{e_r-1}+\frac{\kappa_n}{2}\right\},\\
 \beta_{j,n}
 &=n^{1-j/\mc}\left\{\mathcal F_1^{(j)}(p)
       +\sum_{r=1}^k u_{r,n}
          \left.\frac{\dd^j}{\dd x^j}x^{e_r}\right|_{x=p}\right\},
       \qquad 2\le j\le\mc.
 \end{aligned}
\end{equation}
The coefficients $\boldsymbol\beta_n$ are induced by the physical parameters.
The non-Gaussian approximation in Theorem~\ref{thm:main} applies to the calibrated
class \eqref{even:calibration-system}. The next proposition shows that these
conditions admit uniformly bounded physical corrections.

\Needspace{6\baselineskip}
\begin{proposition}\label{even:bounded-calibration}
Assume $b_r\ge0$ and \eqref{eq:critical}. The calibration system for the
corrections $u_{1,n},\ldots,u_{k,n},\kappa_n$ is
\begin{equation}\label{even:calibration-system}
 \begin{gathered}
 \begin{aligned}
 \mathcal F_1'(p)+\sum_{r=1}^k e_r u_{r,n}p^{e_r-1}
                       +\frac{\kappa_n}{2}
 &=n^{-1+1/\mc}\beta_{1,n},\\
 \mathcal F_1^{(j)}(p)+\sum_{r=1}^k u_{r,n}
      \left.\frac{\dd^j}{\dd x^j}x^{e_r}\right|_{x=p}
 &=n^{-1+j/\mc}\beta_{j,n},\qquad 2\le j\le\mc,
 \end{aligned}\\
 \sup_n\sum_{j=1}^{\mc}|\beta_{j,n}|\le L<\infty.
 \end{gathered}
\end{equation}
For each coefficient sequence satisfying this bound, the equations admit
a solution with
\begin{equation}\label{even:inverse-calibration-bound}
 \sup_n\left(|\kappa_n|+\sum_{r=1}^k|u_{r,n}|\right)
       \le C(1+L),
\end{equation}
where $C$ depends only on the fixed leading model and critical order.
Thus \eqref{eq:window} and \textup{(PC)} hold for $K\ge C(1+L)$.
If $\boldsymbol\beta_n=0$, the corrections can be chosen independent of $n$.
\end{proposition}
The proof is given in Appendix~\ref{app:calibration-proof}.

\subsection{Non-Gaussian approximation in the calibrated critical window}

We now give a quantitative approximation for the fluctuations of the
edge density defined in \eqref{v8:def:graph-edge-degree-density}. Set

\begin{equation}\label{v8:nearby:main-variable-target}
 W_n=n^{1/\mc}(P-p).
\end{equation}

For each coefficient sequence satisfying \eqref{even:calibration-system},
let $Y_n$ have the target distribution with probability density

\begin{equation}\label{eq:target}
 \begin{gathered}
 V_n(z)=\frac{a_\star}{2\mc}z^{2\mc}
             -\sum_{j=1}^{\mc}\frac{\beta_{j,n}}{j!}z^j,\qquad
 \rho_n(z)=\mathcal Z_n^{-1}e^{-V_n(z)},\\
 \mathcal Z_n=\int_\R e^{-V_n(z)}\,\dd z,
 \qquad Y_n\sim\rho_n.
 \end{gathered}
\end{equation}

\Needspace{6\baselineskip}
We measure approximation errors using the 1-Wasserstein and Kolmogorov
distances. For real random variables with finite first moments, these are
defined by

\begin{equation}\label{v8:def:probability-distances}
 \dW(U,V)=\sup_{\Lip h\le1}|\E h(U)-\E h(V)|,\qquad
 \dK(U,V)=\sup_{x\in\mathbb{R}}|\Pp(U\le x)-\Pp(V\le x)|.
\end{equation}

The following theorem gives a uniform approximation bound over the
calibrated parameter windows.
\begin{theorem}\label{thm:main}
Assume \hyperref[ass:PC]{\textup{(PC)}} and \eqref{even:calibration-system},
with bounds $K,L$. For $W_n$ in \eqref{v8:nearby:main-variable-target} and $Y_n$ with density
\eqref{eq:target},
there exist $C<\infty$ and $n_0$ such that, for $n\ge n_0$,
\begin{equation}\label{eq:main-rate}
 \begin{aligned}
 \E|W_n|^{2\mc}&\le C,\\
 \dW(W_n,Y_n)+\dK(W_n,Y_n)&\le Cn^{-1/\mc}.
 \end{aligned}
\end{equation}
The edge count $E=NP$ satisfies
\begin{equation}\label{eq:main-edge-mean-variance}
 \begin{aligned}
 \E E&=Np+Nn^{-1/\mc}\E Y_n+O(n^{2-2/\mc}),\\
 \Var(E)&=N^2n^{-2/\mc}\Var(Y_n)\{1+O(n^{-1/\mc})\}.
 \end{aligned}
\end{equation}
Here $0<c\le\Var(Y_n)\le C$ uniformly over the calibrated windows.
The constants depend only on $\mc$, the fixed interaction graphs, the leading
parameters, $p$, and the bounds $K,L$. They are uniform over the parameter
sequences satisfying \eqref{eq:window} and \eqref{even:calibration-system}.
Proposition~\ref{even:bounded-calibration} shows that the calibration
equations admit bounded solutions.
\end{theorem}

Theorem~\ref{thm:main} is proved in Section~\ref{sec:proof-main}.
The optimality of the rate $O(n^{-1/\mc})$ remains open.

\subsection{Gaussian approximation under fixed finite-size perturbations}\label{sec:fixed-critical-normal}
Retain \hyperref[ass:PC]{\textup{(PC)}} and consider
\begin{equation}\label{fixed:parameters}
 b_{r,n}=b_r+\frac{u_r}{n},\qquad
 h_n=h_*+\frac{\kappa}{n},\qquad 1\le r\le k,
\end{equation}
where $u_1,\ldots,u_k,\kappa$ are fixed real numbers.
Thus $\mathcal U_n$ and $C_n$ in \eqref{eq:window} and
\eqref{even:calibration} do not depend on $n$. The deterministic drift
correction $A_n$ is defined in \eqref{v8:setup:drift-calibration-coefficients},
and the auxiliary drift $G_n$ is defined in \eqref{fixed:drift-definition}.
Their construction and the existence of the exact center used below are
given in Section~\ref{sec:fixed-center}.

\begin{theorem}
\label{fixed:main}
Assume \hyperref[ass:PC]{\textup{(PC)}}, \eqref{fixed:parameters}
and $C_n'(p)\ne0$.
For each fixed critical order $2\mc\ge4$ and all sufficiently large $n$,
let $p_n\in(0,1)$ be the unique zero of $G_n$ in \eqref{fixed:drift-definition}, and set
$\lambda_n=-G_n'(p_n)>0$, as established in Lemma~\ref{fixed:center}.
Define
\begin{equation}\label{fixed:normalized-variable}
 Z_n=(P-p_n)\sqrt{\frac{N\lambda_n}{\nu(p)}}
     =(E-Np_n)\sqrt{\frac{\lambda_n}{N\nu(p)}}.
\end{equation}
For $G\sim\mathcal N(0,1)$,
\begin{equation}\label{fixed:standard-normal}
 \dW(Z_n,G)+\dK(Z_n,G)\le Cn^{-1/(2\mc-1)}.
\end{equation}
The mean and variance satisfy
\begin{equation}\label{fixed:mean-variance}
 \begin{aligned}
 \E P&=p_n+O(n^{-1})
       =p+\sqrt[2\mc-1]{\frac{C_n'(p)}{a_\star n}}
                       +O(n^{-2/(2\mc-1)}),\\
 \Var(E)&=\frac{N\nu(p)}{\lambda_n}
                  \{1+O(n^{-1/(2\mc-1)})\}\\
 &=\frac{n^{(6\mc-4)/(2\mc-1)}}{4(2\mc-1)a_\star^{1/(2\mc-1)}|C_n'(p)|^{(2\mc-2)/(2\mc-1)}}
                  \{1+O(n^{-1/(2\mc-1)})\}.
 \end{aligned}
\end{equation}
Thus the deterministic center $Np_n$ and squared scale $N\nu(p)/\lambda_n$
in \eqref{fixed:normalized-variable} give the leading mean and variance
of the edge count. The constants depend on the fixed model, $\mc$, the correction
coefficients and $|C_n'(p)|^{-1}$, but not on $n$.
\end{theorem}
Theorem~\ref{fixed:main} is proved in
Section~\ref{sec:proof-fixed-critical-normal}.
We retain the exact deterministic center $p_n$, since its first-order
expansion need not be accurate enough on the fluctuation scale.
These perturbations lie outside the bounded calibrated window, since
$\beta_{1,n}=n^{1-1/\mc}C_n'(p)$ has $|\beta_{1,n}|\to\infty$.
The zero perturbation is included: then $C_n'(p)=\mathcal F_1'(p)<0$
by \eqref{fixed:negative-bound}.

For $\mc=2$, fixed coefficients in \eqref{fixed:parameters} give
$\beta_{1,n}=\sqrt n\,C_n'(p)$ and $\beta_{2,n}=C_n''(p)$.
Consequently, $C_n'(p)=0$ satisfies \eqref{even:calibration-system}
and Theorem~\ref{thm:main} gives a quartic limit; $C_n'(p)\ne0$ is covered
by Theorem~\ref{fixed:main}. This exhausts the fixed correction coefficients
under \hyperref[ass:PC]{\textup{(PC)}} in the quartic critical regime.

\subsection{Corollaries}\label{sec:standardized-counts}
We next state the consequences of Theorems~\ref{thm:main} and~\ref{fixed:main}
for subgraph counts and for counts centered at their expectations and
divided by their standard deviations.
The proofs of Corollaries~\ref{even:subgraph-fluctuations}--\ref{even:standardized-counts}
are given in Section~\ref{sec:proof-corollaries}.

For a fixed graph $H$ and a positive deterministic sequence $a_n$,
we write $r_n=O_H(a_n)$ if $|r_n|\le C_Ha_n$ for all sufficiently large $n$.
Here $C_H$ is independent of $n$, but may depend on $H$ and on the fixed
parameters specified in the relevant result. In subgraph-count estimates,
$C_H$ and $C_{H,q}$ denote constants with these permitted dependencies;
the latter may also depend on the fixed moment order $q$.
These constants may change between occurrences.

\begin{corollary}\label{even:subgraph-fluctuations}
Under the hypotheses of Theorem~\ref{thm:main}, fix a finite simple graph
$H$ with $e(H)\ge1$; disconnected graphs and isolated vertices are allowed.
For $f_H(X)=n^{-v(H)}\inj(H,X)$,
\begin{equation}\label{even:subgraph-rate}
 \begin{aligned}
 \dW\left(\frac{n^{1/\mc}\{f_H(X)-p^{e(H)}\}}
                   {e(H)p^{e(H)-1}},Y_n\right)
 &\le C_Hn^{-1/\mc},\\
 \dK\left(\frac{n^{1/\mc}\{f_H(X)-p^{e(H)}\}}
                   {e(H)p^{e(H)-1}},Y_n\right)
 &\le C_Hn^{-1/(2\mc)}.
 \end{aligned}
\end{equation}
The mean and variance of the injective count satisfy
\begin{equation}\label{even:subgraph-mean-variance}
 \begin{aligned}
 \E\inj(H,X)
 &=n^{v(H)}\bigl[p^{e(H)}
       +e(H)p^{e(H)-1}n^{-1/\mc}\E Y_n+O_H(n^{-2/\mc})\bigr],\\
 \Var(\inj(H,X))
 &=e(H)^2p^{2e(H)-2}n^{2v(H)-2/\mc}\Var(Y_n)\{1+O_H(n^{-1/\mc})\}.
 \end{aligned}
\end{equation}
The constants may depend on $H$ and the bounds in
Theorem~\ref{thm:main}, but not on $n$.
\end{corollary}

The replacement estimate \eqref{eq:count-replacement} transfers the Gaussian
approximation in Theorem~\ref{fixed:main} to subgraph counts.

\begin{corollary}
\label{fixed:subgraph-gaussian}
Under the hypotheses of Theorem~\ref{fixed:main}, fix a finite simple graph
$H$ with $e(H)\ge1$; disconnected graphs and isolated vertices are allowed.
For $G\sim\mathcal N(0,1)$ and all sufficiently large $n$,
\begin{equation}\label{fixed:subgraph-W}
 \dW\left(\sqrt{\frac{N\lambda_n}{\nu(p)}}
             \frac{f_H(X)-p_n^{e(H)}}{e(H)p_n^{e(H)-1}},G\right)
 \le C_Hn^{-1/(2\mc-1)},
\end{equation}
and
\begin{equation}\label{fixed:subgraph-K}
 \dK\left(\sqrt{\frac{N\lambda_n}{\nu(p)}}
             \frac{f_H(X)-p_n^{e(H)}}{e(H)p_n^{e(H)-1}},G\right)
 \le
 \begin{cases}
 C_{H,q}n^{-q/(3(q+1))},&\mc=2,\quad q\ge1\text{ fixed},\\
 C_Hn^{-1/(2\mc-1)},&\mc\ge3.
 \end{cases}
\end{equation}
The corresponding mean and variance satisfy
\begin{equation}\label{fixed:subgraph-variance}
 \begin{aligned}
 \E f_H(X)&=p_n^{e(H)}+O_H(n^{-1})\\
 &=p^{e(H)}+e(H)p^{e(H)-1}
       \sqrt[2\mc-1]{\frac{C_n'(p)}{a_\star n}}
       +O_H(n^{-2/(2\mc-1)}),\\
 \Var(\inj(H,X))
 &=e(H)^2p_n^{2e(H)-2}n^{2v(H)}\frac{\nu(p)}{N\lambda_n}
           \{1+O_H(n^{-1/(2\mc-1)})\}\\
 &=\frac{e(H)^2p^{2e(H)-2}n^{2v(H)-2\mc/(2\mc-1)}}{(2\mc-1)a_\star^{1/(2\mc-1)}|C_n'(p)|^{(2\mc-2)/(2\mc-1)}}
           \{1+O_H(n^{-1/(2\mc-1)})\}.
 \end{aligned}
\end{equation}
Thus the centering and scaling in
\eqref{fixed:subgraph-W}--\eqref{fixed:subgraph-K} use the leading mean
and standard deviation of $f_H(X)$.
The constants may depend on $H$, the fixed model, $\mc$, the
correction coefficients and $|C_n'(p)|^{-1}$, and on $q$ where
indicated, but not on $n$.
\end{corollary}

The mean and variance estimates above allow exact centering and variance
normalization in both parameter regimes.
\begin{corollary}\label{even:standardized-counts}
Fix a finite simple graph $H$ with $e(H)\ge1$, and take $n$ sufficiently large.
Under the hypotheses of Theorem~\ref{thm:main},
\begin{equation}\label{eq:main-standardized-rate}
 \begin{aligned}
 &\dW\left(\frac{E-\E E}{\sqrt{\Var(E)}},
           \frac{Y_n-\E Y_n}{\sqrt{\Var(Y_n)}}\right)\\
 &\quad{}+\dK\left(\frac{E-\E E}{\sqrt{\Var(E)}},
           \frac{Y_n-\E Y_n}{\sqrt{\Var(Y_n)}}\right)
 \le Cn^{-1/\mc},
 \end{aligned}
\end{equation}
and
\begin{equation}\label{even:subgraph-standardized-rate}
 \begin{aligned}
 \dW\left(\frac{\inj(H,X)-\E\inj(H,X)}
                   {\sqrt{\Var(\inj(H,X))}},
           \frac{Y_n-\E Y_n}{\sqrt{\Var(Y_n)}}\right)
 &\le C_Hn^{-1/\mc},\\
 \dK\left(\frac{\inj(H,X)-\E\inj(H,X)}
                   {\sqrt{\Var(\inj(H,X))}},
           \frac{Y_n-\E Y_n}{\sqrt{\Var(Y_n)}}\right)
 &\le C_Hn^{-1/(2\mc)}.
 \end{aligned}
\end{equation}

Under the hypotheses of Theorem~\ref{fixed:main}, for $G\sim\mathcal N(0,1)$,
\begin{equation}\label{fixed:standardized-edge-rate}
 \begin{aligned}
 &\dW\left(\frac{E-\E E}{\sqrt{\Var(E)}},G\right)
 +\dK\left(\frac{E-\E E}{\sqrt{\Var(E)}},G\right)
 \le Cn^{-1/(2\mc-1)},
 \end{aligned}
\end{equation}
and
\begin{equation}\label{fixed:standardized-subgraph-rate}
 \begin{aligned}
 \dW\left(\frac{\inj(H,X)-\E\inj(H,X)}
                   {\sqrt{\Var(\inj(H,X))}},G\right)
 &\le C_Hn^{-1/(2\mc-1)},\\
 \dK\left(\frac{\inj(H,X)-\E\inj(H,X)}
                   {\sqrt{\Var(\inj(H,X))}},G\right)
 &\le
 \begin{cases}
 C_{H,q}n^{-q/(3(q+1))},&\mc=2,\quad q\ge1\text{ fixed},\\
 C_Hn^{-1/(2\mc-1)},&\mc\ge3.
 \end{cases}
 \end{aligned}
\end{equation}
The constants are uniform in $n$; in the Gaussian regime they may also
depend on the fixed correction coefficients and $|C_n'(p)|^{-1}$.
The subgraph bounds may depend on $H$ and on $q$ where indicated.
\end{corollary}

For fixed correction coefficients, these conclusions take a particularly
simple form in the quartic case.

\begin{remark}
\label{rem:quartic-fluctuation-transition}
Assume \hyperref[ass:PC]{\textup{(PC)}} with $\mc=2$ and fix any
correction coefficients in \eqref{fixed:parameters}. The calibration
condition is equivalent to $C_n'(p)=0$, since
$\beta_{1,n}=\sqrt n\,C_n'(p)$ and $\beta_{2,n}=C_n''(p)$ is fixed.
Let $Y$ have density proportional to
$\exp\{-a_\star z^4/4+C_n''(p)z^2/2\}$ and let $G\sim\mathcal N(0,1)$.
Corollary~\ref{even:standardized-counts} gives the limiting laws below,
and \eqref{eq:main-edge-mean-variance} and \eqref{fixed:mean-variance}
give the corresponding fluctuation scales.
\[
 \frac{E-\E E}{\sqrt{\Var(E)}}\ \Rightarrow\
 \begin{cases}
  (Y-\E Y)/\sqrt{\Var(Y)}, & C_n'(p)=0,\\
  G, & C_n'(p)\ne0,
 \end{cases}
 \qquad
 \sqrt{\Var(E)}\asymp
 \begin{cases}
  n^{3/2}, & C_n'(p)=0,\\
  n^{4/3}, & C_n'(p)\ne0.
 \end{cases}
\]
For a fixed simple graph $H$ with at least one edge, its standardized
count has the same two limits, with standard deviations of orders
$n^{v(H)-1/2}$ and $n^{v(H)-2/3}$, respectively, by
\eqref{even:subgraph-mean-variance} and \eqref{fixed:subgraph-variance}.
These cases exhaust the fixed correction coefficients in the quartic critical regime. The calibration hyperplane $C_n'(p)=0$ thus marks a change
in both the fluctuation scale and the standardized limiting law.
Both regimes have the same variational function $\Phi$ and the same
unique global maximizer $p$.
\end{remark}

\begin{remark}\label{even:quartic-recovery}
Given any uniformly bounded sequence of interaction corrections in
\eqref{eq:window}, one can choose a uniformly bounded field correction
so that all coefficients $\beta_{1,n},\ldots,\beta_{\mc,n}$ in
\eqref{even:calibration} remain bounded if and only if
\begin{equation}\label{even:interaction-window-criterion}
 \mathcal F_1^{(j)}(p)+\mathcal U_n^{(j)}(p)
       =O(n^{-1+j/\mc}),\qquad 2\le j\le\mc-1.
\end{equation}
By \eqref{even:induced-coefficients}, the linear field term contributes
only to the first derivative. The coefficient $\beta_{\mc,n}$ is bounded
under \eqref{eq:window}, and $\beta_{1,n}$ is bounded precisely when
\[
 h_n=h_*-\frac{2\{\mathcal F_1'(p)+\mathcal U_n'(p)\}}n
       +2\beta_{1,n}n^{-2+1/\mc},\qquad
       \sup_n|\beta_{1,n}|<\infty.
\]
Condition~\eqref{even:interaction-window-criterion} is vacuous for $\mc=2$;
for $\mc\ge3$, the interaction corrections must also satisfy this condition.

On the fluctuation scale \eqref{v8:nearby:main-variable-target}, Taylor's
formula gives
\begin{equation}\label{even:calibration-scaling}
 n\{C_n(p+n^{-1/\mc}z)-C_n(p)\}
 =\sum_{j=1}^{\mc}\frac{\beta_{j,n}}{j!}z^j
       +O(n^{-1/\mc}|z|^{\mc+1}),
\end{equation}
when the expansion segment lies in $J$. The coefficients on the right
are those appearing in the target density \eqref{eq:target}.

In particular, for $\mc=2$, bounded interaction corrections and a bounded sequence $\theta_n$
give the quartic window
\[
 h_n=h_*-\frac{2\{\mathcal F_1'(p)+\mathcal U_n'(p)\}}n
                 +\frac{\theta_n}{n^{3/2}},
\]
with $a_\star=\alpha_4=-\Phi^{(4)}(p)/6$,
$\beta_{2,n}=\mathcal F_1''(p)+\mathcal U_n''(p)$ and
$\beta_{1,n}=\theta_n/2$.
For general $\mc$, the same field form uses the bounded sequence
$\theta_n=2\beta_{1,n}n^{-1/2+1/\mc}$; the critical order still comes
from \eqref{eq:critical}.
\end{remark}

\begin{remark}\label{even:critical-order-parity}
\hypertarget{even-critical-order-parity-outline}{}
We justify the even order in \eqref{eq:critical} using the variational
function in \eqref{eq:scalar}. Its boundary derivatives satisfy
$\Phi'(0+)=+\infty$ and $\Phi'(1-)=-\infty$, 
so every global maximizer is interior. For even $j>\max_i e_i$,
\[
 \Phi^{(j)}(x)=-\frac{(j-2)!}{2}
     \{x^{1-j}+(1-x)^{1-j}\}<0,\qquad 0<x<1.
\]
Thus the first nonzero derivative at a global maximizer has finite order $j\ge2$.
Applying Taylor's formula on both sides of the maximizer shows that
$j$ is even and $\Phi^{(j)}(p)<0$.
If $\Phi''(p)=0$, then $j=2\mc\ge4$. Subject to their remaining hypotheses,
Theorems~\ref{thm:main} and~\ref{fixed:main} cover every such order.
The nondegenerate quadratic case $j=2$ lies outside the setting of
\eqref{eq:critical}; its relation to our assumptions is discussed in
Remark~\ref{even:quadratic-scope}.
\end{remark}

\begin{remark}\label{even:quadratic-scope}
Theorems~\ref{thm:main} and \ref{fixed:main} admit Gaussian counterparts
for $\mc=1$ when the variational function $\Phi$ has a unique global maximizer $p$ with
$\Phi''(p)<0$. Indeed, \eqref{eq:profiles} and \eqref{eq:scalar} give
\[
 0\le2\nu(p)t(p)\le\nu(p)B''(p)
 =\frac12+\nu(p)\Phi''(p)<\frac12,
\]
which preserves the localized conditional-row contraction gap.
The standardized edge count admits the optimal $O(n^{-1})$ rate in
the Wasserstein and Kolmogorov distances, matching the lower bound imposed
by its lattice spacing. A Gaussian limit also holds for the centered
count of any fixed simple graph with at least one edge. Gaussian fluctuations
in this regime have already been studied in~\cite{FLSZ,WinsteinCLT};
we therefore focus on the degenerate single-phase critical regime, $\mc\ge2$.
\end{remark}

\section{Conditional row laws, moments, and cumulants}\label{row:sec:conditional-rows}\label{sec:conditional-moments}
This section collects the conditional laws, covariance identities, and
moment estimates; Appendices~\ref{row:sec:row-gap}--\ref{ep:row:section}
prove localization, construct the auxiliary kernels, and establish the
cumulant bounds. The ERGM estimates in this section and in
Sections~\ref{sec:counts}--\ref{sec:supports} are proved under
\hyperref[ass:PC]{\textup{(PC)}}.
Their inputs are fixed-tolerance localization
(Lemma~\ref{row:lem:macro}), the quadratic row gap
(Proposition~\ref{row:prop:row-gap}), bounded template derivatives, and
$\|\ell_n-\ell\|_{C^j}=O(n^{-1})$. Under the leading-model assumptions in \textup{(PC)}, the row-gap argument uses
$\Phi'(p)=\Phi''(p)=0$ and uniqueness of the global maximizer.
Counts are centered at $P$ or $P_A$; neither the derivative calibration
\eqref{even:calibration-system} nor critical-scale moments enter these arguments.
Section~\ref{sec:proof-main} obtains the $2\mc$th moment from stationarity.

\subsection{Edge conditionals and moment notation}\label{sec:conditional-notation}
For an unordered edge $e$, let $X^{e,a}$ be the configuration with that
coordinate set to $a\in\{0,1\}$, and write
\begin{equation}\label{audit:physical-edge-conditional}
 \partial_ef(X)=f(X^{e,1})-f(X^{e,0}),\qquad
 L_e=\partial_e\mathcal H_n(X),\qquad
 \Pp_n(X_e=1\mid X_{-e})=\sig(L_e).
\end{equation}
For row coordinates, the same difference operator is denoted by $\nabla_i$.
Mixed differences are taken in distinct coordinates.
For the multiaffine extension of the Hamiltonian, mixed derivatives
in distinct coordinates agree with the corresponding finite differences. Define the centered degrees by
\begin{equation}\label{audit:centered-degree-cancellation}
 z_i=D_i-(n-1)P,\qquad \sum_i z_i=0.
\end{equation}
If vertices $1,2$ are deleted, write $A$ for the induced exterior graph
and $e=X_{12}$ for their mutual edge, and set
\[
 P_A=\binom{n-2}{2}^{-1}\sum_{3\le i<j\le n}X_{ij},\qquad
 x_i=X_{1i},\quad y_i=X_{2i}\quad(i\ge3).
\]
Write $|A|$ for the exterior size. The row laws and their product reference law
are defined in Subsection~\ref{sec:exact-row-laws}.
The fixed point $p_A$ (Lemma~\ref{lem:shape}) is distinct from the observed density $P_A$.
The step graphon associated with $X$ is
\begin{equation}\label{v8:def:graphon}
 \mathcal W_X(u,v)=X_{\lceil nu\rceil,\lceil nv\rceil},
 \qquad (u,v)\in(0,1]^2,
\end{equation}
with value zero when $u=0$ or $v=0$. For an exterior graph $A$, define
$\mathcal W_A$ by the same construction after relabeling its $|A|$ vertices as
$[|A|]$.

We index the model templates by $r$, use $q$ for a moment order, and write
$q_0$ for a fixed number of roots or rows. Cumulant orders are denoted by
$\ell$; two distinct moment orders are written $q_1,q_2$.
For fixed $q\ge1$, the moment norms are
\begin{equation}\label{v8:def:conditional-moment-norm}
 \|Z\|_q=(\E|Z|^q)^{1/q},\qquad
 \|f\|_{L^q(\mu)}=(\E_\mu|f|^q)^{1/q}.
\end{equation}
Unless otherwise specified, unmarked expectations in graph-model estimates
are taken under the original graph law. In particular,
\begin{equation}\label{v8:def:exterior-moment-norm}
 \|\E[Z\mid A,e]\|_q
 =\left(\E_{(A,e)}|\E[Z\mid A,e]|^q\right)^{1/q}.
\end{equation}
Thus the outer expectation averages over the random exterior.
A bound in this norm need not hold uniformly over deterministic exterior
configurations.
For random vectors and matrices, we take the coordinate norm before
the moment norm:
\begin{equation}\label{v8:def:moment-matrix-norms}
 \bigl\|\|R\|_2\bigr\|_q
 =\left[\E\left(\sum_i|R_i|^2\right)^{q/2}\right]^{1/q},\qquad
 \|M\|_F=\left(\sum_{i,j}|M_{ij}|^2\right)^{1/2}.
\end{equation}
The Euclidean operator norm is
\begin{equation}\label{v8:def:operator-norm}
 \|M\|_{\op}=\sup_{\|v\|_2=1}\|Mv\|_2.
\end{equation}
Write $O_{L^q}(a_n)$ for a bound $C_qa_n$ in the indicated norm; moment
bounds written $O(a_n)$ hold at every fixed order. Unless stated otherwise,
constants in graph-model estimates may depend on the templates, the leading
critical point, $K$, and the moment order. In calibrated-window estimates,
they may also depend on the window bound $L$. All such constants are independent
of $n$. Throughout, $n$ is assumed to be sufficiently large.
In subgraph-count estimates, the notation $O_H$ and the constants $C_H$
and $C_{H,q}$ have the meaning specified in Section~\ref{sec:standardized-counts}.

\subsection{Conditional row laws and observables}\label{sec:exact-row-laws}\label{sec:local-input-definitions}
Fix $A=X_{\{3,\ldots,n\}}$, $e=X_{12}$, and
$x_i=X_{1i},y_i=X_{2i}$ for $i\in V(A)$.
All template deletions retain isolated vertices. The single-row Hamiltonian is
\begin{equation}\label{v8:def:single-row-hamiltonian}
\begin{aligned}
 H_A(x)&=h_n\sum_ix_i+\sum_r b_{r,n}n^{2-v_r}
 \sum_{u\in V(F_r)}
 \sum_{\varphi:V(F_r)\setminus\{u\}\hookrightarrow V(A)}\\
 &\quad\prod_{\substack{w\in V(F_r)\\uw\in E(F_r)}}x_{\varphi(w)}
 \prod_{ab\in E(F_r-u)}A_{\varphi(a)\varphi(b)} .
\end{aligned}
\end{equation}
The terms meeting both distinguished vertices contribute
\begin{equation}\label{v8:def:two-row-interaction}
 \begin{aligned}
 V_{A,e}(x,y)
 ={}&\sum_r b_{r,n}n^{2-v_r}
 \sum_{\substack{u,v\in V(F_r)\\u\ne v}}
 \sum_{\varphi:V(F_r)\setminus\{u,v\}\hookrightarrow V(A)}
 e^{\one_{\{uv\in E(F_r)\}}}\\
 &{}\times\prod_{\substack{w\in V(F_r)\setminus\{u,v\}\\uw\in E(F_r)}}x_{\varphi(w)}
 \prod_{\substack{w\in V(F_r)\setminus\{u,v\}\\vw\in E(F_r)}}y_{\varphi(w)}
 \prod_{ab\in E(F_r-\{u,v\})}A_{\varphi(a)\varphi(b)} .
 \end{aligned}
\end{equation}
Root pairs are ordered,
and empty products (including $e^0$) equal one. Partitioning the injective
template maps according to their intersection with $\{1,2\}$ now gives
\begin{equation}\label{row:eq:split}
\begin{aligned}
 \mathcal H_n(X(A,e;x,y))
 &=h_n\left(e+\sum_{i<j}A_{ij}\right)
  +\sum_r b_{r,n}n^{2-v_r}\inj(F_r,A)\\
 &\quad+H_A(x)+H_A(y)+V_{A,e}(x,y).
\end{aligned}
\end{equation}
Define the reference law
\begin{equation}\label{v8:def:reference-row-law}
 \pi_A(x)=\frac{e^{H_A(x)}}{\sum_{u\in\{0,1\}^{V(A)}}e^{H_A(u)}} .
\end{equation}
The corresponding single-coordinate field is independent of $x_i$:
\begin{equation}\label{row:eq:single-row-field}
 \begin{aligned}
 l_{i,A}(x)&=\nabla_iH_A(x)=\partial_iH_A(x)\\
 &=h_n+\sum_r b_{r,n}n^{2-v_r}
 \sum_{\substack{u,v\in V(F_r)\\uv\in E(F_r)}}
 \sum_{\substack{\varphi:V(F_r)\setminus\{u\}\hookrightarrow V(A)\\\varphi(v)=i}}\\
 &\qquad{}\times
 \prod_{\substack{w\in V(F_r)\setminus\{u,v\}\\uw\in E(F_r)}}x_{\varphi(w)}
 \prod_{ab\in E(F_r-u)}A_{\varphi(a)\varphi(b)} .
 \end{aligned}
\end{equation}
Injectivity makes every monomial multiaffine, so coordinate differences
coincide with derivatives. Normalizing the weights of the two possible coordinate values yields
\begin{equation}\label{row:eq:single-site-conditional}
 \pi_A(X_i=1\mid X_{-i})=\sig(l_{i,A}(X)),\qquad
 l_{i,A}(X)=\log\frac{\pi_A(X_i=1\mid X_{-i})}{\pi_A(X_i=0\mid X_{-i})}.
\end{equation}
Extend the field and its logistic transform to the continuous cube.
The densities proportional to
$e^{H_A(x)+H_A(y)+sV_{A,e}(x,y)}$, $0\le s\le1$, interpolate between
$\pi_A\otimes\pi_A$ and the true conditional distribution given $(A,e)$; $\pi_A$ need not
be a marginal of the latter.

For a physical edge, write $L_{uv}=\nabla_{uv}\mathcal H_n$ for its conditional logit. Bounded model parameters give, uniformly in the edge and the state,
\begin{equation}\label{eq:edge-nondegeneracy}
 \kappa\le\Pp_n(X_{uv}=1\mid X_{-uv})=\sig(L_{uv})\le1-\kappa.
\end{equation}
Indeed, each oriented template edge has at most $n^{v_r-2}$ extensions,
which gives the uniform bound $|L_{uv}|\le|h_n|+2\sum_re_r|b_{r,n}|$.

Fix $q_0$ and put $V(A)=[n]\setminus[q_0]$, $m=n-q_0$,
$\mathcal I=[q_0]\times V(A)$, $x_i^b=X_{bi}$, and
$\mathbf e=(X_{bc})_{b<c\le q_0}$.
Keep the ambient normalization $n^{2-v_r}$ in $H_A(x)$.
Let $V_{A,\mathbf e}(\mathbf x)$ collect the template terms meeting at least two
distinguished vertices, and define
\begin{equation}\label{v8:def:conditional-multirow-law}
 \mu_{A,\mathbf e,s}(\mathbf x)
 =\frac{\exp\{\sum_{a=1}^{q_0}H_A(x^a)+sV_{A,\mathbf e}(\mathbf x)\}}
 {\sum_{\mathbf y\in\{0,1\}^{\mathcal I}}
       \exp\{\sum_{a=1}^{q_0}H_A(y^a)+sV_{A,\mathbf e}(\mathbf y)\}},
 \qquad 0\le s\le1.
\end{equation}
Thus $\mu_{A,\mathbf e,0}=\pi_A^{\otimes q_0}$ and
$\mu_{A,\mathbf e,1}$ is the true conditional law; $V_{A,\mathbf e}(\mathbf x)\equiv0$
when $q_0=1$. Random exteriors always have their marginal law under
$\Pp_n$. Write $l_{i,A}(x)=\partial_iH_A(x)$ and
$\partial^\gamma=\prod_j\partial_{x_j}^{\gamma_j}$,
$|\gamma|=\sum_j\gamma_j$.
Partitioning template maps by their distinguished images gives
\begin{equation}\label{row:eq:fixed-rows-split}
 \begin{aligned}
 \mathcal H_n(X)
 ={}&h_n\left(\sum_{\substack{i<j\\i,j\in V(A)}}A_{ij}
                       +\sum_{1\le b<c\le q_0}e_{bc}\right)
       +\sum_{r=1}^k b_{r,n}n^{2-v_r}\inj(F_r,A)\\
 &+\sum_{b=1}^{q_0}H_A(x^b)+V_{A,\mathbf e}(\mathbf x).
 \end{aligned}
\end{equation}
Lemma~\ref{lem:row-deterministic-bounds} gives the deterministic
bounds for $H_A(x)$ and $V_{A,\mathbf e}(\mathbf x)$.

For a function on a binary cube $\{0,1\}^{\mathcal I}$, let $x^{J,b}$
denote the vector obtained by replacing the coordinates in $J$ by $b$, and define
\begin{equation}\label{v8:def:coordinate-differences}
 \nabla_i f(x)=f(x^{i,1})-f(x^{i,0}),\qquad
 \|\nabla_if\|_\infty=\max_x|\nabla_i f(x)|.
\end{equation}
Distinct-coordinate differences commute:
\begin{equation}\label{v8:def:mixed-differences}
 \nabla_Jf=\prod_{i\in J}\nabla_i f
   =\sum_{b\in\{0,1\}^J}(-1)^{|J|-\sum_ib_i}f(x^{J,b}).
\end{equation}
\begin{equation}\label{v8:def:difference-seminorms}
 \|\nabla_Jf\|_\infty=\max_x|\nabla_Jf(x)|,\qquad
 S_j(f)=\sum_{\substack{J\subseteq\mathcal I\\|J|=j}}\|\nabla_Jf\|_\infty .
\end{equation}
Each subset is counted once in these seminorms. For a family of local functions,
define the coordinate supports by
\begin{equation}\label{v8:def:local-support}
 J_\alpha=\{i:\nabla_i\chi_\alpha\not\equiv0\}.
\end{equation}
Bounded overlap means that there are fixed constants $C,d_0,d_1$ such that
\begin{equation}\label{v8:def:bounded-overlap}
 \sup_\alpha\|\chi_\alpha\|_\infty\le C,\qquad
 \sup_\alpha|J_\alpha|\le d_0,\qquad
 \max_i\sum_\alpha\one_{\{i\in J_\alpha\}}\le d_1 .
\end{equation}
For a fixed number $q_0$ of rows, the coordinates are
$x_i^a=X_{ai}$, $\mathcal I=[q_0]\times V(A)$, and the nonconstant
features of the exterior block $i$ are
\begin{equation}\label{v8:def:row-block-features}
 \chi_{i,S}(\mathbf x)=\prod_{a\in S}x_i^a,\qquad
 \varnothing\ne S\subseteq[q_0].
\end{equation}
These features have support size at most $q_0$ and overlap at most
$2^{q_0-1}$. For two rows they are $x_i,y_i,x_iy_i$.
A normalized block polynomial has fixed maximal block degree $d$ and the form
\begin{equation}\label{v8:def:normalized-block-polynomials}
 G=c_\varnothing+\sum_{s=1}^d\sum_{i_1<\cdots<i_s}
 \sum_{\varnothing\ne S_1,\ldots,S_s\subseteq[q_0]}
 c_{i_1,\ldots,i_s;S_1,\ldots,S_s}\prod_{j=1}^s\chi_{i_j,S_j},
\end{equation}
with coefficients fixed by the conditional data and satisfying
\begin{equation}\label{v8:def:block-coefficient-normalization}
 |c_\varnothing|\le C,\qquad
 |c_{i_1,\ldots,i_s;S_1,\ldots,S_s}|\le Cn^{-s}.
\end{equation}
Thus $\|G\|_\infty\le C+C\sum_{s=1}^d(2^{q_0}-1)^s/s!$.
The degree counts exterior blocks rather than individual coordinates.
This normalization bounds the total contribution at each degree without any independence assumption.

\subsection{Localization and covariance bounds}\label{sec:covariance}
A graphon is a symmetric measurable function $W:[0,1]^2\to[0,1]$.
For an integrable signed kernel $D$, use the cut norm
\begin{equation}\label{v8:nearby:app:cut-norm}
 \|D\|_\square
 =\sup_{\substack{B_1,B_2\subseteq[0,1]\\B_1,B_2\text{ measurable}}}
       \left|\int_{B_1\times B_2}D(x,y)\,\dd x\,\dd y\right|.
\end{equation}
The step graphon $\mathcal W_X$ uses equal vertex intervals, including
the zero diagonal cells. For $S\subseteq[n]$, write
$A_S=X_{[n]\setminus S}$ and use $n-|S|$ equal intervals for its graphon.

\begin{lemma}\label{row:lem:macro}
Assume \hyperref[ass:PC]{\textup{(PC)}}.
For every fixed $\delta>0$,
\begin{equation}\label{row:eq:macro}
 \Pp_n(\|\mathcal W_X-p\|_\square>\delta)\le C_\delta e^{-c_\delta n^2}
\end{equation}
for large $n$. Moreover, for every fixed deletion size $d\ge0$,
\begin{equation}\label{audit:macro-deletion-tail}
 \Pp_n\left(
   \max_{\substack{S\subseteq[n]\\|S|=d}}
      \|\mathcal W_{A_S}-p\|_\square>\delta
 \right)
 \le C_{\delta/2}e^{-c_{\delta/2}n^2}
\end{equation}
for all sufficiently large $n>d$, with the threshold also depending on $d$.
The bound is simultaneous over the deleted sets $S$.
\end{lemma}
The proof of Lemma~\ref{row:lem:macro} is given in
Appendix~\ref{app:macro-proof}. Using the good sets $\mathcal A_n$ and
$\Omega_A$ from Lemma~\ref{row:lem:good-exteriors}, set
\begin{equation}\label{v8:nearby:good-multirow-states}
 \Omega(A)=\Omega_A^{q_0}.
\end{equation}
The good sets are independent of $\mathbf e,s$.
Below, constants are uniform in good exteriors $A$, $\mathbf e,s$, row states,
and $n$, but may depend on the fixed model, $q_0,K$, cutoffs,
tolerances, and the indicated orders. The auxiliary kernels contract on
the full cube and agree with the Gibbs updates on typical states.

\begin{lemma}\label{lem:local-input}
Assume \hyperref[ass:PC]{\textup{(PC)}}.
For the conditional data specified above, the good sets and auxiliary
kernels satisfy the following properties.

\textup{(i) Good sets and agreement.}
For some $c,C>0$,
\begin{equation}\label{eq:local-good-exteriors}
 \Pp_n(A\notin\mathcal A_n)\le Ce^{-cn}.
\end{equation}
For each $A\in\mathcal A_n$, $\mathbf e$, and $s$, there are update
probabilities $\widetilde\pi_u^{A,\mathbf e,s}$, $u\in\mathcal I$,
independent of the updated coordinate $x_u$, such that
\begin{equation}\label{eq:local-kernel-agreement}
 \begin{aligned}
 &\sup_{A\in\mathcal A_n}\sup_{\mathbf e}\sup_{0\le s\le1}
       \mu_{A,\mathbf e,s}\bigl(\Omega(A)^c\bigr)\le Ce^{-cn},\\
 &\widetilde\pi_u^{A,\mathbf e,s}(\mathbf x)
   =\mu_{A,\mathbf e,s}(x_u=1\mid x_{-u}=\mathbf x_{-u}),
             \qquad \mathbf x\in\Omega(A).
 \end{aligned}
\end{equation}
The bounds are uniform over the conditional laws on their corresponding good sets.

\textup{(ii) Uniform kernel bounds.}
For each $n\ge n_0$ and $A\in\mathcal A_n$, there is a symmetric
nonnegative $|\mathcal I|\times|\mathcal I|$ matrix $T^{(n,A)}$ such that
\begin{equation}\label{eq:local-influence}
 \begin{gathered}
 \sup_{\mathbf e}\sup_{0\le s\le1}
       \|\nabla_v\widetilde\pi_u^{A,\mathbf e,s}\|_\infty
       \le T_{uv}^{(n,A)},\qquad u,v\in\mathcal I,\\
 T_{uv}^{(n,A)}=T_{vu}^{(n,A)}\ge0,\qquad T_{uu}^{(n,A)}=0,\\
 \sup_{A\in\mathcal A_n}\max_{u,v\in\mathcal I}
                  nT_{uv}^{(n,A)}\le C,\qquad
 \sup_{A\in\mathcal A_n}\max_{u\in\mathcal I}
                  \sum_{v\in\mathcal I}T_{uv}^{(n,A)}\le\rho<1.
 \end{gathered}
\end{equation}
The entries of $T^{(n,A)}$ may depend on $n,A$, but the matrix is common
to all $\mathbf e,s$ and all row states; $C$ and $\rho<1$ are uniform.
All supremum norms in this statement are taken over the full binary cube. Moreover, after increasing
$n_0$ if necessary,
\[
 \frac14\min\{p,1-p\}\le\widetilde\pi_u^{A,\mathbf e,s}(\mathbf x)
                         \le1-\frac14\min\{p,1-p\}.
\]
This positive lower bound depends only on $p$, whereas $n_0$
may depend on all the fixed data listed above.

\hypertarget{local-resolvent-properties}{}
\textup{(iii) Stationarity and resolvent bounds.}
For $n\ge n_0$, $A\in\mathcal A_n$, and fixed $\mathbf e,s$, define
\begin{equation}\label{eq:local-aux-generator}
 \widetilde{\mathcal L}_{A,\mathbf e,s}f
 =\sum_{u\in\mathcal I}
       \{\widetilde\pi_u^{A,\mathbf e,s}(\mathbf x)-x_u\}\nabla_u f .
\end{equation}
It has a unique stationary law
$\widetilde\mu_{A,\mathbf e,s}$, characterized by
$\E_{\widetilde\mu_{A,\mathbf e,s}}
 [\widetilde{\mathcal L}_{A,\mathbf e,s}f]=0$ for every real-valued
function $f$ on the cube. This law may differ from $\mu_{A,\mathbf e,s}$.
Let $\widetilde P_t^{A,\mathbf e,s}
=e^{t\widetilde{\mathcal L}_{A,\mathbf e,s}}$ and define
\begin{equation}\label{eq:local-resolvent-definition}
 \widetilde R_{A,\mathbf e,s}f
 =\int_0^\infty
   (\widetilde P_t^{A,\mathbf e,s}f
          -\E_{\widetilde\mu_{A,\mathbf e,s}}f)\,\dd t .
\end{equation}
The integral converges absolutely in the supremum norm. For every
such $f$ and every fixed integer $j\ge2$,
\begin{equation}\label{eq:local-resolvent-bounds}
 \begin{aligned}
 -\widetilde{\mathcal L}_{A,\mathbf e,s}
       \widetilde R_{A,\mathbf e,s}f
     &=f-\E_{\widetilde\mu_{A,\mathbf e,s}}f,\\
 \|\widetilde R_{A,\mathbf e,s}f\|_\infty
     &\le S_1(\widetilde R_{A,\mathbf e,s}f)\le CS_1(f),\\
 S_j(\widetilde R_{A,\mathbf e,s}f)
     &\le C_j(1+\log n)^{j-1}\sum_{h=1}^jS_h(f).
 \end{aligned}
\end{equation}
The constants are uniform over the specified conditional data and may also depend on the fixed gap $1-\rho$ and the indicated order.
\end{lemma}
The proof is given in Appendix~\ref{app:local-kernel-proof}.

Following F\"ollmer's influence-matrix approach~\cite{Follmer1982}
and the stationary-kernel comparison method~\cite{RvH2014}, we retain the
discrepancy between the true Gibbs updates and the auxiliary kernels,
which agree only on typical states and need not be Gibbs compatible.
The true generator is
\begin{equation}\label{v8:def:true-multirow-generator}
 \mathcal L_{A,\mathbf e,s}f
 =\sum_{(a,i)\in\mathcal I}
   \{\sig(l_{i,A}(x^a)+s\partial_{(a,i)}V_{A,\mathbf e})-x_i^a\}
                                                        \nabla_{(a,i)}f .
\end{equation}
For a general finite Gibbs law $\mu\propto e^H$, define
\begin{equation}\label{v8:nearby:generic-generators-resolvent}
 \mathcal Lf=\sum_i(\sig(\nabla_iH)-x_i)\nabla_i f,\quad
 \widetilde{\mathcal L}f=\sum_i(\widetilde\pi_i-x_i)\nabla_i f,\quad
 \widetilde Rf=\int_0^\infty(e^{t\widetilde{\mathcal L}}f-\E_{\widetilde\mu}f)\,dt .
\end{equation}
Write
\begin{equation}\label{v8:nearby:oscillation}
 \operatorname{osc}(g)=\max_xg(x)-\min_xg(x),\qquad
 S_1(f)=\sum_i\|\nabla_i f\|_\infty .
\end{equation}

\begin{lemma}\label{lem:cov-error}
Let $\mu\propto e^H$ be a Gibbs law on a binary cube
$\{0,1\}^{\mathcal I}$ with $|\mathcal I|\le C_0n$.
Suppose that auxiliary single-coordinate update probabilities
$\widetilde\pi_i$, independent of coordinate $i$, are uniformly
bounded away from zero and one. Assume that they have a symmetric
nonnegative majorant $T$ satisfying
\[
 \|\nabla_j\widetilde\pi_i\|_\infty\le T_{ij},\qquad
 T_{ij}=T_{ji},\quad T_{ii}=0,\quad T_{ij}\le C/n,\quad
 \max_i\sum_jT_{ij}\le\rho<1.
\]
Assume also that $\widetilde\pi_i=\sig(\nabla_iH)$ on a set
$\Omega$ with $\mu(\Omega^c)\le Ce^{-cn}$.
Use the true and auxiliary generators and resolvent in
\eqref{v8:nearby:generic-generators-resolvent}.
For $\mu_{A,\mathbf e,s}$, Lemma~\ref{lem:local-input}
verifies these hypotheses and supplies the resolvent.
For real-valued $f,g$, with coordinate sums over $\mathcal I$,
\begin{equation}\label{eq:cov-error}
 \begin{aligned}
 \Cov_\mu(f,g)
 ={}&\sum_i\E_{\mu}\bigl[\sig(\nabla_iH)
       \{1-\sig(\nabla_iH)\}\nabla_i\widetilde Rf\nabla_ig\bigr]\\
    &+\Cov_\mu((\mathcal L-\widetilde{\mathcal L})\widetilde Rf,g).
 \end{aligned}
\end{equation}
The oscillation and first-difference sum are those in
\eqref{v8:nearby:oscillation}.
The discrepancy in \eqref{eq:cov-error} satisfies
\begin{equation}\label{audit:covariance-discrepancy}
 |\Cov_\mu((\mathcal L-\widetilde{\mathcal L})\widetilde Rf,g)|
 \le C\mu(\Omega^c)S_1(f)\operatorname{osc}(g)
 \le Ce^{-cn}S_1(f)\operatorname{osc}(g).
\end{equation}
The influence-matrix inverse gives
\begin{equation}\label{eq:res-gradient}
 \begin{aligned}
 \|\nabla_j\widetilde Rf\|_\infty
       &\le\sum_k[(I_{|\mathcal I|}-T)^{-1}]_{jk}\|\nabla_k f\|_\infty,\\
 0\le [(I_{|\mathcal I|}-T)^{-1}]_{jk}&\le\mathbf1_{\{j=k\}}+C/n.
 \end{aligned}
\end{equation}
Moreover, $(I_{|\mathcal I|}-T)^{-1}$ is symmetric and
$\max_j\sum_k[(I_{|\mathcal I|}-T)^{-1}]_{jk}\le(1-\rho)^{-1}$. Consequently,
\begin{equation}\label{eq:cov-bound}
 |\Cov_\mu(f,g)|\le C\max_i\|\nabla_i f\|_\infty S_1(g)
       +C e^{-cn}S_1(f)\operatorname{osc}(g).
\end{equation}
All constants depend only on the displayed uniform hypotheses;
in particular, they are uniform over the laws of Lemma~\ref{lem:local-input}.

\end{lemma}

\begin{proof}
The Poisson equation and the true Gibbs conditional identities give
\[
 \begin{gathered}
 \widetilde Rf-\E_\mu[\widetilde Rf\mid X_{-i}]
   =(X_i-\sig(\nabla_iH))\nabla_i\widetilde Rf,\\
 -\E_\mu[(\mathcal L\widetilde Rf)(g-\E_\mu g)]
   =\sum_i\E_\mu[\sig(\nabla_iH)(1-\sig(\nabla_iH))
                      \nabla_i\widetilde Rf\nabla_i g],\\
 \Cov_\mu(f,g)
 =-\E_\mu[(\mathcal L\widetilde Rf)(g-\E_\mu g)]
    +\Cov_\mu((\mathcal L-\widetilde{\mathcal L})\widetilde Rf,g).
 \end{gathered}
\]
This proves \eqref{eq:cov-error}. Agreement on $\Omega$ and
Lemma~\ref{row:lem:resolvent} give
\begin{equation}\label{v8:proof:cov-error:discrepancy}
 |(\mathcal L-\widetilde{\mathcal L})\widetilde Rf|
 \le\one_{\Omega^c}S_1(\widetilde Rf)
 \le C\one_{\Omega^c}S_1(f).
\end{equation}
Since $|g-\E_\mu g|\le\operatorname{osc}(g)$, this yields
\eqref{audit:covariance-discrepancy}.
The same resolvent lemma gives \eqref{eq:res-gradient}; the Neumann
series gives symmetry and row sums at most $(1-\rho)^{-1}$.
Consequently $\max_i\|\nabla_i\widetilde Rf\|_\infty
\le(1-\rho)^{-1}\max_i\|\nabla_i f\|_\infty$.
Substitution into \eqref{eq:cov-error}, with
$\sig(\nabla_iH)(1-\sig(\nabla_iH))\le1/4$,
proves \eqref{eq:cov-bound}.
\end{proof}

\begin{corollary}
\label{cor:disjoint-local-covariance}
Under the hypotheses of Lemma~\ref{lem:cov-error}, suppose that $f,g$
depend on disjoint sets of physical row coordinates. Then
\begin{equation}\label{eq:disjoint-local-covariance}
 |\operatorname{Cov}_\mu(f,g)|
       \le\frac Cn S_1(f)S_1(g).
\end{equation}
In particular, for distinct physical coordinates $i,j$,
\begin{equation}\label{eq:conditional-offdiagonal-covariance}
 |\operatorname{Cov}_\mu(X_i,X_j)|\le C/n.
\end{equation}
Coordinates in different rows at the same exterior vertex are distinct,
and are therefore covered by this assertion. The constants are uniform
over the fixed conditional laws $\mu_{A,\mathbf e,s}$.
\end{corollary}

\begin{proof}
For $j\in\operatorname{supp}(g)$, disjointness and
\eqref{eq:res-gradient} give
$\|\nabla_j\widetilde Rf\|_\infty\le Cn^{-1}S_1(f)$.
Thus \eqref{eq:cov-error} and
$\operatorname{osc}(g)\le S_1(g)$ imply
\[
 |\Cov_\mu(f,g)|
 \le C(n^{-1}+e^{-cn})S_1(f)S_1(g)
 \le Cn^{-1}S_1(f)S_1(g).
\]
Taking $f=X_i$, $g=X_j$ proves the coordinate bound.
\end{proof}

\subsection{Polynomial means and moment bounds}
The disjoint-support covariance bound controls the error in replacing a row polynomial by its value at the marginal means. Fix the full conditional law $\mu_{A,\mathbf e,s}$ on a good exterior
and select two rows $x,y$. Consider the following class of site polynomials of fixed degree:
\begin{equation}\label{v8:nearby:two-row-polynomial-class}
 F_{A,\mathbf e}(x,y)
 =c_\varnothing+\sum_{d=1}^D\sum_{i_1<\cdots<i_d}
 \sum_{\substack{\alpha_j,\beta_j\in\{0,1\}\\\alpha_j+\beta_j\ge1}}
 c_{i_1,\ldots,i_d;\alpha,\beta}\prod_{j=1}^dx_{i_j}^{\alpha_j}y_{i_j}^{\beta_j},
 \quad |c_{i_1,\ldots,i_d;\alpha,\beta}|\le Kn^{-d}.
\end{equation}
Its two marginal mean vectors are
\begin{equation}\label{v8:nearby:selected-row-means}
 \mu_i^x=\E_{\mu_{A,\mathbf e,s}}x_i,\qquad
 \mu_i^y=\E_{\mu_{A,\mathbf e,s}}y_i.
\end{equation}

\begin{lemma}
\label{lem:row-polynomial-mean}
Assume \hyperref[ass:PC]{\textup{(PC)}}.
For a fixed law \eqref{v8:def:conditional-multirow-law} with
$A\in\mathcal A_n$, let $F_{A,\mathbf e}$ and $(\mu^x,\mu^y)$ be as in
\eqref{v8:nearby:two-row-polynomial-class}--\eqref{v8:nearby:selected-row-means}.
The two mean vectors need not agree. We have
\begin{equation}\label{eq:row-polynomial-mean}
 \left|\E_{\mu_{A,\mathbf e,s}} F_{A,\mathbf e}(x,y)
             -F_{A,\mathbf e}(\mu^x,\mu^y)\right|\le CK/n.
\end{equation}
The constant is uniform in the good exterior, internal edges, and
$s$; it may depend on the standing model data, $q_0$, and the fixed degree bounds.
For a polynomial in a single selected row, the corresponding statement is
$|\E_{\mu_{A,\mathbf e,s}} F_{A,\mathbf e}(x)-F_{A,\mathbf e}(\E_{\mu_{A,\mathbf e,s}} x)|\le CK/n$.
This includes the reference law $\pi_A$.
\end{lemma}

\begin{proof}
A monomial is a product of distinct physical coordinates $Z_1,\ldots,Z_h$
(including both coordinates of $x_i y_i$). By
Corollary~\ref{cor:disjoint-local-covariance},
\[
 \E\prod_{a=1}^hZ_a-\prod_{a=1}^h\E Z_a
 =\sum_{a=1}^{h-1}\left(\prod_{b<a}\E Z_b\right)
       \Cov\left(Z_a,\prod_{b>a}Z_b\right)=O(n^{-1}).
\]
There are $O(n^d)$ degree-$d$ monomials, so
\[
 |\E F-F(\mu^x,\mu^y)|
 \le\sum_{d=1}^D Cn^d\cdot Kn^{-d}\cdot n^{-1}\le CK/n.
\]
The same calculation applies to one selected row.
\end{proof}

The difference-moment estimate of Lemma~\ref{lem:conditional-difference-moments}
in Appendix~\ref{app:difference-moments} applies to the normalized block polynomials of
\eqref{v8:def:normalized-block-polynomials}, which can be written as
\begin{equation}\label{v8:nearby:block-polynomial-moment-class}
 \begin{gathered}
 G=c_\varnothing+\sum_{k=1}^d\sum_{i_1<\cdots<i_k}
 \sum_{\varnothing\ne S_1,\ldots,S_k\subseteq[q_0]}
 c_{i_1,\ldots,i_k;S_1,\ldots,S_k}\prod_{j=1}^k\prod_{a\in S_j}x_{i_j}^a,\\
 |c_\varnothing|\le C,\qquad |c_{i_1,\ldots,i_k;S_1,\ldots,S_k}|\le Cn^{-k}.
 \end{gathered}
\end{equation}

\begin{corollary}
\label{cor:row-polynomial-moments}
Assume \hyperref[ass:PC]{\textup{(PC)}}.
Let $\mu_{A,\mathbf e,s}$ be any fixed conditional law from
Lemma~\ref{lem:local-input}. Then the following bounds hold uniformly in $A,\mathbf e,s$.
For every normalized block polynomial $G$ of fixed degree
and every fixed real $q\ge1$, one has
\begin{equation}\label{eq:local-polynomial-moments}
 \begin{gathered}
 \|G-\E_{\mu_{A,\mathbf e,s}}G\|_{L^q(\mu_{A,\mathbf e,s})}\le C_qn^{-1/2},\\
 \Var_{\mu_{A,\mathbf e,s}}(G)\le C/n,\qquad
 \E_{\mu_{A,\mathbf e,s}}|G-\E_{\mu_{A,\mathbf e,s}}G|^4\le Cn^{-2}.
 \end{gathered}
\end{equation}
In particular, each row satisfies
\begin{equation}\label{eq:local-row-moments}
 \max_{1\le a\le q_0}
 \E_{\mu_{A,\mathbf e,s}}\left|
   \sum_{i\in V(A)}(x_i^a-\E_{\mu_{A,\mathbf e,s}}x_i^a)
                \right|^q\le C_qn^{q/2}.
\end{equation}
No logarithmic loss occurs in these bounds. The constants may depend
on the fixed polynomial degrees and coefficient bounds.
\end{corollary}

\begin{proof}
At each block degree $d$, only $O(n^{d-1})$ monomials involve a fixed
coordinate, and their coefficients are $O(n^{-d})$. Summing over
the finitely many degrees gives
\[
 \|\nabla_{(a,i)}G\|_\infty\le C/n,
 \qquad \sum_{(a,i)}\|\nabla_{(a,i)}G\|_\infty^2\le C/n.
\]
Lemma~\ref{lem:conditional-difference-moments} gives
\eqref{eq:local-polynomial-moments}; its cases $q=2,4$ are the
variance and fourth-moment bounds. Applying the same estimate to
$G=n^{-1}\sum_i x_i^a$ gives \eqref{eq:local-row-moments}.
\end{proof}

The same difference-moment estimate also applies to weighted sums of local observables with bounded overlap. Under $\mu=\mu_{A,\mathbf e,s}$, write
\begin{equation}\label{v8:nearby:weighted-local-family}
 J_\alpha=\{u:\nabla_uf_\alpha\not\equiv0\},\qquad
 \|f_\alpha\|_\infty\le M,\quad |J_\alpha|\le d_0,\quad
 \max_u\sum_\alpha\one_{\{u\in J_\alpha\}}\le d_1.
\end{equation}
The three bounds are fixed independently of $n$.

\begin{lemma}
\label{lem:weighted-linear}
Assume \hyperref[ass:PC]{\textup{(PC)}}.
Under a fixed good-exterior law of Lemma~\ref{lem:local-input}, let
$(f_\alpha)_{\alpha\in\Lambda}$ satisfy
\eqref{v8:nearby:weighted-local-family}.
For every deterministic real vector $a=(a_\alpha)$ and fixed real $q\ge1$,
\begin{equation}\label{eq:weighted-local-moments}
 \left\|\sum_{\alpha\in\Lambda}
        a_\alpha(f_\alpha-\E_\mu f_\alpha)\right\|_{L^q(\mu)}
 \le C_q\left(\sum_{\alpha\in\Lambda}a_\alpha^2\right)^{1/2}.
\end{equation}
In particular,
\[
 \|a^{\mathsf T}(X-\E_\mu X)\|_{L^q(\mu)}\le C_q\|a\|_2.
\]
The weights and functions may depend on the conditional data
$(A,\mathbf e,s)$, provided they are fixed under $\mu$. The constants
may depend on $M,d_0,d_1,q$, but not on the weights or these data.
\end{lemma}

\begin{proof}
Cauchy--Schwarz and bounded overlap give
\[
 \|\nabla_i\sum_\alpha a_\alpha f_\alpha\|_\infty^2
 \le4M^2d_1\sum_{\alpha:i\in J_\alpha}a_\alpha^2,
 \qquad \sum_i\|\nabla_i\sum_\alpha a_\alpha f_\alpha\|_\infty^2
 \le4M^2d_0d_1\sum_\alpha a_\alpha^2.
\]
Apply Lemma~\ref{lem:conditional-difference-moments}.
For coordinate features take $M=d_0=d_1=1$.
\end{proof}

\subsection{Row mean equations and contraction}
To sharpen the preliminary cumulant bounds, we first locate the row
means $\bar x^b=\E_{\mu_{A,\mathbf e,s}}X^b$ in a common contraction
neighborhood. The reference field $l_{i,A}(x)=\partial_iH_A(x)$ is used
even when the selected marginal does not have Hamiltonian $H_A(x)$.

\begin{lemma}
\label{ep:row:marginal-mean}
Assume \hyperref[ass:PC]{\textup{(PC)}}.
For a selected row $b$ of \eqref{v8:def:conditional-multirow-law},
with $A\in\mathcal A_n$, the following estimates hold uniformly
as in Lemma~\ref{lem:local-input}. For every fixed real $q\ge1$
and each row coordinate $i$,
\begin{equation}\label{ep:row:field-moments}
 \begin{aligned}
 \|l_{i,A}(X^b)-l_{i,A}(\bar x^b)\|_{L^q(\mu_{A,\mathbf e,s})}
 &\le C_qn^{-1/2},\\
 |\E_{\mu_{A,\mathbf e,s}}\{l_{i,A}(X^b)-l_{i,A}(\bar x^b)\}|&\le Cn^{-1}.
 \end{aligned}
\end{equation}
The mean equation is
\begin{equation}\label{ep:row:marginal-mean-equation}
 \max_i|\bar x_i^b-\sig(l_{i,A}(\bar x^b))|\le C/n.
\end{equation}
The resulting localization and variance estimates are
\begin{align}
 \max_i|\bar x_i^b-p|&\le a+C/n,
                    \label{ep:row:mean-neighborhood}\\
 \max_i|\sig'(l_{i,A}(\bar x^b))
                    -\bar x_i^b(1-\bar x_i^b)|&\le C/n.
                    \label{ep:row:variance-jacobian}
\end{align}
The constants have the same dependence on the fixed data and the same uniformity as in
Lemma~\ref{lem:local-input}; they are also uniform in the selected row.
Here $a$ is the fixed tolerance in \eqref{row:eq:Omega}.
\end{lemma}
The proof is given in Appendix~\ref{app:mean-equation-proof}.

Use the interval $I=[p-\delta,p+\delta]\Subset(0,1)$ fixed in
Appendix~\ref{row:sec:extension}. For $d=n-q_0$, define

\begin{equation}\label{v8:nearby:app:reference-mean-map}
 M_A(y)=(\sig(l_{i,A}(y)))_{i\in V(A)},\qquad y\in I^d.
\end{equation}

\begin{equation}\label{ep:row:jacobian-definition}
 (DM_A(y))_{ij}
 =\sig'(l_{i,A}(y))\partial_jl_{i,A}(y),\qquad (DM_A(y))_{ii}=0.
\end{equation}

Write $H_A''=(\partial_{ij}H_A)$ and
$\|M\|_{\infty\to\infty}=\max_i\sum_j|M_{ij}|$; matrix comparisons
are entrywise.

\begin{lemma}
\label{ep:row:mean-jacobian}
Assume \hyperref[ass:PC]{\textup{(PC)}}.
Fix $q_0\ge1$ and use $I$ and the tolerances fixed before the definition of
$\mathcal A_n$ in Appendix~\ref{row:sec:extension}.
Then, for all sufficiently large $n$ and every
$A\in\mathcal A_n$, each mean vector in
Lemma~\ref{ep:row:marginal-mean} belongs to $I^d$.
Moreover, for this reference mean map, there is a symmetric nonnegative matrix
$C^{(n,A)}$ such that
\begin{equation}\label{ep:row:symmetric-majorant}
 \begin{gathered}
 |DM_A(y)|\le C^{(n,A)},\qquad
 |\operatorname{diag}((\nu(z_i))_i)H_A''(y)|\le C^{(n,A)}
                                    \quad(y,z\in I^d),\\
 C^{(n,A)}_{ij}=C^{(n,A)}_{ji}\ge0,\quad C^{(n,A)}_{ii}=0,
 \quad\max_i\sum_jC^{(n,A)}_{ij}\le\rho=3/4<1.
 \end{gathered}
\end{equation}
The matrix depends only on $n,A$ and the fixed model data;
in particular, it is independent of $\mathbf e,s$ and the selected row.
Absolute values and comparisons are entrywise. The same domination gives
\begin{equation}\label{ep:row:jacobian-gap}
 \sup_{y\in I^d}
 \max\{\|DM_A(y)\|_{\infty\to\infty},\|DM_A(y)\|_{\mathrm{op}}\}
 \le\rho
\end{equation}
and the common inverse bound
\begin{equation}\label{eq:shape-common-inverse}
 \max_{*\in\{\infty\to\infty,\,\mathrm{op}\}}
 \|(I_d-\operatorname{diag}((\nu(z_i))_i)H_A''(y))^{-1}\|_*
 \le(1-\rho)^{-1}\qquad(y,z\in I^d).
\end{equation}
The averaged reference field satisfies
\begin{equation}\label{ep:row:average-field-contraction}
 \left|\sig\!\left(d^{-1}\sum_i l_{i,A}(y)\right)
       -\sig\!\left(d^{-1}\sum_i l_{i,A}(z)\right)\right|
 \le\rho\|y-z\|_\infty\qquad(y,z\in I^d).
\end{equation}
\end{lemma}
The proof is given in Appendix~\ref{app:common-neighborhood-proof}.

\subsection{Uniform cumulants and centered moments}
For bounded variables, define the joint cumulant by
\begin{equation}\label{v8:def:joint-cumulant}
 \cum_\mu(F_1,\ldots,F_\ell)
 =\left.\partial_{t_1}\cdots\partial_{t_\ell}
          \log\E_\mu e^{\sum_jt_jF_j}\right|_{t_1=\cdots=t_\ell=0}.
\end{equation}

In the cumulant statements below, $\mu$ denotes a fixed conditional law
$\mu_{A,\mathbf e,s}$ on a good exterior; the reference row law is included
by setting $s=0$ and restricting the product law to one row.
Index the physical coordinates by $[d]$, with $d=O(n)$, and write $X_j(x)=x_j$.
All orders are fixed independently of the graph size, and all bounds are uniform in the conditional data.

For a finite set $A$ of formal argument positions, let $\mathcal P(A)$
denote its set of partitions. The moment--cumulant identity is
\begin{equation}
\E\prod_{j\in A}Z_j
 =\sum_{\pi\in\mathcal P(A)}
                     \prod_{B\in\pi}\cum(Z_j:j\in B).
                     \label{row:eq:moment-partition}
\end{equation}

\begin{theorem}\label{ep:row:sharp}
Assume \hyperref[ass:PC]{\textup{(PC)}}.
Under a fixed good-exterior law
\eqref{v8:def:conditional-multirow-law} of Lemma~\ref{lem:local-input},
for every local family $(f_\alpha)$ satisfying
\eqref{v8:nearby:weighted-local-family} and fixed $\ell\ge2$,
\begin{equation}\label{ep:row:feature-bound}
 \sup_{\alpha_1}\sum_{\alpha_2,\ldots,\alpha_\ell}
 |\cum_{\mu_{A,\mathbf e,s}}(f_{\alpha_1},\ldots,f_{\alpha_\ell})|\le C_{\ell}.
\end{equation}
The constants are uniform in the good exterior, the internal root-edge
configuration, and $s$. Estimate~\eqref{ep:row:feature-bound} also permits a fixed bounded local
observable as its first argument, with constants depending on its prescribed
supremum-norm and coordinate-support bounds.

Moreover, normalized block polynomials $G_1,\ldots,G_\ell$ of fixed degrees,
as defined in \eqref{v8:nearby:block-polynomial-moment-class}, satisfy
\begin{equation}\label{ep:row:polynomial-bound}
 |\cum_{\mu_{A,\mathbf e,s}}(G_1,\ldots,G_\ell)|\le C_{\ell}n^{1-\ell},\qquad \ell\ge2.
\end{equation}
Finally, let $F$ be a fixed bounded local observable whose supremum norm and
coordinate-support size in \eqref{v8:def:local-support} are bounded
independently of $n$. Then, for $\ell\ge1$,
\begin{equation}\label{ep:row:anchored-polynomial}
 |\cum_{\mu_{A,\mathbf e,s}}(F,G_1,\ldots,G_\ell)|\le C_{\ell}n^{-\ell}.
\end{equation}
\end{theorem}
The proof is given in Appendix~\ref{app:sharp-proof}.

\begin{corollary}
\label{cor:local-feature-sums}
Assume \hyperref[ass:PC]{\textup{(PC)}}.
Let $\mu$ be a conditional law of Theorem~\ref{ep:row:sharp},
or a one-row marginal of such a law, and let $(\chi_\alpha)$ satisfy
\eqref{v8:nearby:weighted-local-family}.
The $\chi_\alpha$ are the original, uncentered features.
For every fixed integer $L\ge1$,
\begin{equation}\label{eq:micro-moment-sum}
 \sum_{\alpha_1,\ldots,\alpha_L}
 \left|\E_\mu\prod_{j=1}^L
       \bigl(\chi_{\alpha_j}-\E_\mu\chi_{\alpha_j}\bigr)\right|
 \le C_Ln^{\lfloor L/2\rfloor}.
\end{equation}
For every fixed integer $s\ge1$ and real coefficients fixed under $\mu$, put
\begin{equation}\label{v8:nearby:centered-feature-polynomial}
 Q_s=\sum_{\alpha_1,\ldots,\alpha_s}
       c_{\alpha_1\ldots\alpha_s}
       \prod_{j=1}^s\bigl(\chi_{\alpha_j}-\E_\mu\chi_{\alpha_j}\bigr),
 \qquad |c_{\alpha_1\ldots\alpha_s}|\le Cn^{-s}.
\end{equation}
Then, for every fixed real $q\ge1$,
\begin{equation}\label{ep:row:homogeneous-moment}
 \|Q_s\|_{L^q(\mu)}\le C_{q,s}n^{-s/2}.
\end{equation}
All sums include repeated labels; tuple restrictions are encoded by zero
coefficients. The constants are uniform as in Theorem~\ref{ep:row:sharp},
with additional dependence on the fixed moment orders, degree and
coefficient bound, and are independent of the tuple restrictions.
The reference row $\pi_A$ and coordinate features $X_i$ are included.
\end{corollary}

\begin{proof}
Discard constant features. Bounded overlap leaves $O(n)$ features, so
Theorem~\ref{ep:row:sharp} gives
$\sum_{\alpha_1,\ldots,\alpha_\ell}|\cum_\mu(\chi_{\alpha_1},\ldots,
\chi_{\alpha_\ell})|\le C_{\ell}n$ for every integer $\ell\ge2$.
Centering removes singleton blocks from \eqref{row:eq:moment-partition}:
\[
 \sum_{\alpha_1,\ldots,\alpha_L}
       \left|\E_\mu\prod_{j=1}^L(\chi_{\alpha_j}-\E_\mu\chi_{\alpha_j})\right|
 \le\sum_{\substack{\pi\in\mathcal P([L])\\|B|\ge2\ (B\in\pi)}}
        C_Ln^{|\pi|}\le C_Ln^{\lfloor L/2\rfloor}.
\]
For $L=1$ the sum is zero. Partitions act on formal positions, so
repeated labels are included.
Choose $2k\ge q$. The coefficient bound gives
\[
 \E_\mu|Q_s|^{2k}
 \le C^{2k}n^{-2ks}
    \sum_{\alpha_{1,1},\ldots,\alpha_{2k,s}}
       \left|\E_\mu\prod_{b=1}^{2k}\prod_{j=1}^s(\chi_{\alpha_{b,j}}-\E_\mu\chi_{\alpha_{b,j}})\right|
 \le C_{k,s}n^{-ks}.
\]
Moment monotonicity proves \eqref{ep:row:homogeneous-moment}.
Tuple restrictions may be discarded after absolute values are taken.
Features of a one-row marginal are functions on the full conditional
cube; $q_0=1$ includes $\pi_A$.
\end{proof}

\section{Graph counts and reduction to the edge density}\label{sec:counts}
We first control exterior changes and whole-row stationarity, then
reduce graph counts to the empirical edge density.
We conclude with bounds on the conditional logit bias and the fixed-point
residuals for the edge density.

\subsection{Row resampling and stationarity}
Use the single-row Hamiltonian $H_A(x)$ defined in
\eqref{v8:def:single-row-hamiltonian}, and fix a good initial exterior $A$.
Let $A^{(1)}$ differ from $A$ at exactly one exterior
edge $\{a,b\}$, and interpolate its value by
\[
 \begin{aligned}
 A^{(u)}_{ab}=A^{(u)}_{ba}
   &=A_{ab}+u(A^{(1)}_{ab}-A_{ab}),\qquad 0\le u\le1,\\
 A^{(u)}_{ij}&=A_{ij}\quad(\{i,j\}\ne\{a,b\}).
 \end{aligned}
\]
The expression $H_{A^{(u)}}$ means evaluation of the same Hamiltonian
polynomial at these exterior entries. Put
\begin{equation}\label{v8:nearby:exterior-edge-interpolation}
 \begin{aligned}
 V(x)&=H_{A^{(1)}}(x)-H_A(x),\\
 H_{A^{(u)}}(x)&=H_A(x)+uV(x),\\
 \pi_{A^{(u)}}(x)&=
 \frac{\exp\{H_A(x)+uV(x)\}}
      {\displaystyle\sum_{y\in\{0,1\}^{V(A)}}\exp\{H_A(y)+uV(y)\}},
 \qquad \pi_{A^{(0)}}=\pi_A.
 \end{aligned}
\end{equation}
The oscillation and first-difference sum are defined in
\eqref{v8:nearby:oscillation}.

\begin{lemma}\label{lem:sensitivity}
Assume \hyperref[ass:PC]{\textup{(PC)}}.
For \eqref{v8:nearby:exterior-edge-interpolation}, uniformly over good
initial exteriors and the changed edge,
\begin{equation}\label{eq:external-energy}
 \operatorname{osc}(V)\le C/n,\qquad S_1(V)\le C/n.
\end{equation}
Moreover, let $f$ be a fixed row test function independent of $u$, and assume
\[
 \max_i\|\nabla_i f\|_\infty\le C_f n^{-2}.
\]
Then
\begin{equation}\label{audit:row-expectation-sensitivity}
 \sup_{0\le u\le1}|\E_{\pi_{A^{(u)}}}f-\E_{\pi_A}f|
 \le CC_fn^{-3}+Ce^{-cn}S_1(f).
\end{equation}
In particular, the same estimate holds for
$|\E_{\pi_{A^{(1)}}}f-\E_{\pi_A}f|$.
Only the initial exterior $A$ is required to be good.
\end{lemma}
\begin{proof}
\proofstep{Affected copies.}
Injectivity gives $H_{A^{(u)}}=H_A+uV$.
A term in $\nabla_iV$ fixes the exterior edge $\{a,b\}$ and the
row edge $\{o,i\}$, and hence fixes three vertices if $i\in\{a,b\}$ and four otherwise.
Multiplying the remaining counts by $O(n^{2-v_r})$ gives
\[
 \|\nabla_iV\|_\infty\le
 \begin{cases}C/n,&i\in\{a,b\},\\ C/n^2,&i\notin\{a,b\},\end{cases}
 \qquad \operatorname{osc}(V)\le S_1(V)\le C/n.
\]
Templates too small to contain the fixed vertices contribute zero;
the field term cancels. This proves \eqref{eq:external-energy}.

\proofstep{Comparison under the initial law.}
The exact tilt identity is
\begin{equation}\label{eq:external-tilt-covariance-comparison}
 \begin{aligned}
 \E_{\pi_{A^{(u)}}}f-\E_{\pi_A}f
 &=\frac{\E_{\pi_A}[f e^{u(V-\min V)}]}
         {\E_{\pi_A}e^{u(V-\min V)}}-\E_{\pi_A}f\\
 &=\frac{\Cov_{\pi_A}(f,e^{u(V-\min V)})}
         {\E_{\pi_A}e^{u(V-\min V)}}.
 \end{aligned}
\end{equation}
Uniformly for $0\le u\le1$,
\[
 \begin{gathered}
 e^{u(V-\min V)}\ge1,\qquad
 \|\nabla_i e^{u(V-\min V)}\|_\infty\le e^{C/n}\|\nabla_iV\|_\infty,\\
 S_1(e^{u(V-\min V)})+\operatorname{osc}(e^{u(V-\min V)})\le C/n.
 \end{gathered}
\]
Only the exterior for $\pi_A$ is required to be good. Applying \eqref{eq:cov-bound} gives
\[
 |\Cov_{\pi_A}(f,e^{u(V-\min V)})|
 \le CC_fn^{-3}+Ce^{-cn}n^{-1}S_1(f).
\]
Division by $\E_{\pi_A}e^{u(V-\min V)}\ge1$ proves
\eqref{audit:row-expectation-sensitivity}.
Any direct exterior dependence of a changing test must be estimated separately.
\end{proof}

We use the preceding stability estimate to control the drift under whole-row
resampling. For a fixed simple graph $H$, retaining all isolated vertices, define
\begin{equation}\label{v8:def:whole-row-generator}
 f_H(X)=n^{-v(H)}\inj(H,X),\qquad
 \mathcal E_i f=\E[f\mid X_{[n]\setminus\{i\}}],\qquad
 \mathcal Gf=\sum_{i=1}^n(\mathcal E_if-f).
\end{equation}
Thus $\mathcal E_i$ resamples the entire incident row from its conditional law.

\begin{proposition}
\label{prop:whole-row}
Assume \hyperref[ass:PC]{\textup{(PC)}}.
For $f_H,\mathcal G$ in \eqref{v8:def:whole-row-generator}
and every fixed real $q\ge1$,
\begin{equation}\label{eq:whole-row}
 \|\mathcal Gf_H\|_q
 =\left(\E|\mathcal Gf_H(X)|^q\right)^{1/q}
 \le C_{H,q}/n.
\end{equation}
\end{proposition}

\begin{proof}
For the true conditional covariance $\Cov_i$ of row $i$,
$\E[h(\mathcal E_i-I)f_H]=-\E\Cov_i(f_H,h)$, hence
\begin{equation}\label{v8:proof:whole-row:stationary-moment}
 \E (\mathcal Gf_H)^{2m}=-\sum_{i=1}^n\E\Cov_i(f_H,(\mathcal Gf_H)^{2m-1}).
\end{equation}

\proofstep{Coordinate differences.}
Let $f_H^{[i]}$ be the contribution of injections whose image contains $i$.
The remaining part is constant in row $i$, so
$(\mathcal E_i-I)f_H=(\mathcal E_i-I)f_H^{[i]}$.
Fixing one vertex or one edge gives, respectively,
\begin{equation}\label{audit:dependency:s5:whole-row-deterministic-envelopes}
 \|f_H^{[i]}\|_\infty\le C_H/n,\qquad
 \max_e\|\partial_e f_H\|_\infty\le C_H/n^2,\qquad
 \|(\mathcal Gf_H)\|_\infty\le\sum_i2\|f_H^{[i]}\|_\infty\le C_H.
\end{equation}
Let $\mathcal D_n$ be the event that all one-row exteriors are good.
Then $\Pp_n(\mathcal D_n^c)\le C e^{-cn}$ by a union bound.
We claim that
\begin{equation}\label{eq:good-oscillation}
 |\partial_e(\mathcal Gf_H)(X)|\le C_H/n^2,\qquad X\in\mathcal D_n.
\end{equation}
For $e=\{a,b\}$ and $i\in\{a,b\}$,
$\partial_e(\mathcal E_i-I)f_H=-\partial_e f_H=O(n^{-2})$.
For $i\notin\{a,b\}$, fixing $i,a,b$ gives
\[
 \|\partial_e f_H^{[i]}\|_\infty\le C_Hn^{-3},\qquad
 \max_{j\ne i}\|\partial_{\{i,j\}}f_H^{[i]}\|_\infty\le C_Hn^{-2},
 \qquad S_1(f_H^{[i]})\le C_H/n.
\]
Starting from the good exterior $A$ of row $i$, let $A'$ be obtained by flipping $e$.
Lemma~\ref{lem:sensitivity} yields
\[
 \begin{aligned}
 |\E_{\pi_{A'}}f_H^{[i]}(A',x)-\E_{\pi_A}f_H^{[i]}(A,x)|
 &\le\|f_H^{[i]}(A',\cdot)-f_H^{[i]}(A,\cdot)\|_\infty\\
 &\quad+|\E_{\pi_{A'}}f_H^{[i]}(A,x)-\E_{\pi_A}f_H^{[i]}(A,x)|\\
 &\le C_Hn^{-3}+C_Hn^{-3}+Ce^{-cn}S_1(f_H^{[i]})
 \le C_H'n^{-3}.
 \end{aligned}
\]
The direct change of $-f_H^{[i]}$ has the same order.
Summing the two incident-row contributions and the $n-2$ remaining
contributions proves \eqref{eq:good-oscillation}; only the starting
exterior must be good.

\proofstep{Recursive bounds for even moments.}
On a good row exterior, \eqref{eq:res-gradient} gives
$\max_{j\ne i}\|\partial_{\{i,j\}}\widetilde R_i f_H\|_\infty
\le C_Hn^{-2}$.
For $X\in\mathcal D_n$, factoring the difference of powers gives
\[
 |\partial_e((\mathcal Gf_H)^{2m-1})(X)|
 \le C_{H,m}n^{-2}
     \bigl[|(\mathcal Gf_H)(X^{e,0})|^{2m-2}+|(\mathcal Gf_H)(X^{e,1})|^{2m-2}\bigr].
\]
The clamped states need not be good. Instead,
\eqref{eq:edge-nondegeneracy} implies, for $b\in\{0,1\}$,
\[
 \begin{aligned}
 \E|(\mathcal Gf_H)(X^{e,b})|^{2m-2}
 &=\E\left[\frac{\one_{\{X_e=b\}}}
                  {\Pp_n(X_e=b\mid X_{-e})}|(\mathcal Gf_H)(X)|^{2m-2}\right]\\
 &\le\kappa^{-1}\E|(\mathcal Gf_H)|^{2m-2}.
 \end{aligned}
\]
For $m=1$, the zeroth powers equal one. The indicator of $\mathcal D_n$
may be dropped from nonnegative terms.
In \eqref{eq:cov-error}, the $n(n-1)$ Dirichlet terms on
$\mathcal D_n$ are therefore bounded in total by
\[
 \begin{aligned}
 &C_{H,m}n(n-1)n^{-2}n^{-2}\E|(\mathcal Gf_H)|^{2m-2}\\
 &\quad\le C_{H,m}n^{-2}\E|(\mathcal Gf_H)|^{2m-2}.
 \end{aligned}
\]
Using \eqref{audit:dependency:s5:whole-row-deterministic-envelopes},
we control three exceptional contributions. On $\mathcal D_n^c$ with
the selected exterior good, the resolvent gradient is still
$O(n^{-2})$ and $|\partial_e(\mathcal Gf_H)^{2m-1}|\le C_{H,m}$, giving
$C_{H,m}\Pp_n(\mathcal D_n^c)$ after summation. For each good exterior,
the retained kernel discrepancy satisfies
\[
 |\Cov_i((\mathcal L_i-\widetilde{\mathcal L}_i)
               \widetilde R_i f_H,(\mathcal Gf_H)^{2m-1})|
 \le Ce^{-cn}S_1(f_H)\operatorname{osc}((\mathcal Gf_H)^{2m-1})
 \le C_{H,m}n^{-1}e^{-cn}.
\]
For a bad exterior,
$|\Cov_i(f_H,(\mathcal Gf_H)^{2m-1})|=|\Cov_i(f_H^{[i]},(\mathcal Gf_H)^{2m-1})|
\le C_{H,m}/n$; its averaged sum is exponentially small.
Thus \eqref{v8:proof:whole-row:stationary-moment} yields
\begin{equation}\label{audit:whole-row-moment-recursion}
 \E|(\mathcal Gf_H)|^{2m}\le C_{H,m}n^{-2}\E|(\mathcal Gf_H)|^{2m-2}
                  +C_{H,m}n^Ce^{-cn}.
\end{equation}
Induction from the zeroth moment gives $\E|(\mathcal Gf_H)|^{2m}\le C_{H,m}n^{-2m}$.
Choose $2m\ge q$ and use moment monotonicity to prove \eqref{eq:whole-row}.
\end{proof}

\subsection{Rooted fluctuations and degree moments}
Let $A=X_{[n]\setminus\mathcal D}$, $|\mathcal D|=O(1)$, $m=|V(A)|$.
For a fixed rooted graph $(H,o)$, with all isolated vertices retained, define
\begin{equation}\label{v8:def:rooted-vector}
 T_i=T_i(H,o;A)
 =n^{-(v(H)-1)}
  \sum_{\substack{\varphi:V(H)\hookrightarrow V(A)\\\varphi(o)=i}}
      \prod_{uv\in E(H)}A_{\varphi(u)\varphi(v)},\qquad
 \bar T=m^{-1}\sum_iT_i .
\end{equation}
The template root $o$ is fixed, and the normalization uses the ambient graph size $n$.
For finitely many fixed rooted pairs $(H_a,o_a)$ and graph-dependent
coefficients that are common to all root images, with
$\max_a\|c_a\|_\infty\le C$, put
\begin{equation}\label{v8:def:rooted-vector-combination}
 U_i=\sum_ac_aT_i(H_a,o_a;A),\qquad\bar U=m^{-1}\sum_iU_i .
\end{equation}
Polynomial edge decorations have the form
\begin{equation}\label{v8:def:rooted-edge-decoration}
 \begin{gathered}
 V_i=n^{-(v(H)-1)}
 \sum_{\substack{\varphi:V(H)\hookrightarrow V(A)\\\varphi(o)=i}}
   \sum_{S\subseteq E(H)}c_S\prod_{uv\in S}A_{\varphi(u)\varphi(v)},\\
 \bar V=m^{-1}\sum_iV_i,\qquad\max_S\|c_S\|_\infty\le C.
 \end{gathered}
\end{equation}
Here the coefficients $c_S$ are again common to all root images and injections.

\begin{lemma}\label{lem:root-exchange}
Assume \hyperref[ass:PC]{\textup{(PC)}}.
Consider the singly rooted counts and fixed-deletion exterior specified above.
Then, for every fixed real $q\ge1$, the following bounds hold
uniformly in $i,j\in V(A)$ and the deleted labels:
\begin{equation}\label{eq:root-deterministic}
 \|T_i-T_j\|_q\le C_{H,q}n^{-1/2},\qquad
 \left\|\left(\sum_{i\in V(A)}|T_i-\overline T|^2\right)^{1/2}
                                      \right\|_q\le C_{H,q}.
\end{equation}
Moreover, the shared-coefficient combinations in
\eqref{v8:def:rooted-vector-combination} satisfy
\begin{equation}\label{eq:root-shared-bounds}
 \|U_i-U_j\|_q\le C_qn^{-1/2},\qquad
 \left\|\left(\sum_{i\in V(A)}|U_i-\overline U|^2\right)^{1/2}
                                      \right\|_q\le C_q.
\end{equation}
The polynomial edge decorations in
\eqref{v8:def:rooted-edge-decoration} also satisfy
\begin{equation}\label{eq:root-decorated-bounds}
 \|V_i-V_j\|_q\le C_qn^{-1/2},\qquad
 \left\|\left(\sum_{i\in V(A)}|V_i-\bar V|^2\right)^{1/2}
                                      \right\|_q\le C_q.
\end{equation}
In particular, $c_S=(-P_A)^{e(H)-|S|}$ replaces each edge by
$A_{uv}-P_A$.
No independence between coefficients and counts is required.
The constants depend on the fixed model, the bound in
\eqref{eq:parameter-neighborhood}, the rooted templates, the numbers of terms
and deleted vertices, the coefficient bound and $q$. They are independent of $n$
and the image labels. Normalization by $m$ gives the same estimates.
\end{lemma}

\begin{proof}
For $i\ne j$ in $V(A)$, put $B=X_{[n]\setminus\{i,j\}}$ and
condition on $\mathcal F_B=\sigma(B,X_{ij})$. Thus
$\mathcal D\subseteq V(B)$.
The transposition $\theta_{ij}$ preserves these data and the conditional law, so
\begin{equation}\label{v8:proof:root-exchange:conditional-centering}
 T_i(H,o;\theta_{ij}A)=T_j(H,o;A),\qquad
 \E[T_i-T_j\mid\mathcal F_B]=0.
\end{equation}
In the two rows, use features $x_u,y_u,x_uy_u$ for
$u\in V(A)\setminus\{i,j\}$.
Fixing $s$ block images leaves $O(n^{v(H)-1-s})$ rooted injections
(one fewer free image if the copy also uses $j$).
Thus $T_i-T_j$ is a normalized block polynomial. Conditional centering,
Corollary~\ref{cor:row-polynomial-moments}, and its bounded envelope on bad
exteriors give
\[
 \E|T_i-T_j|^q\le C_{H,q}n^{-q/2}+Ce^{-cn}.
\]
For $q\ge2$, the identity
$\sum_i|T_i-\bar T|^2=(2m)^{-1}\sum_{i,j}|T_i-T_j|^2$ gives
\[
 \left\|\left(\sum_i|T_i-\bar T|^2\right)^{1/2}\right\|_q^2
 \le(2m)^{-1}\sum_{i,j}\|T_i-T_j\|_q^2\le C_{H,q}.
\]
Use moment monotonicity for $q<2$.
For shared coefficients,
$U_i-U_j=\sum_ac_a(T_i^{(a)}-T_j^{(a)})$ and
$U_i-\bar U=\sum_ac_a(T_i^{(a)}-\bar T^{(a)})$.
The triangle inequalities for absolute values and the Euclidean and $L^q$
norms, together with
$\|c_a\|_\infty\le C$, prove \eqref{eq:root-shared-bounds}.
Every polynomial decoration reduces, by binary edge powers, to
\[
 V_i=\sum_{S\subseteq E(H)}c_ST_i((V(H),S),o;A).
\]
All vertices and the normalization are retained. Applying
\eqref{eq:root-shared-bounds} to this finite linear combination gives
\eqref{eq:root-decorated-bounds}; centered edges have
$c_S=(-P_A)^{e(H)-|S|}$.
Replacing $n$ by $m$ multiplies a count by
$(n/m)^{v(H)-1}=1+O(n^{-1})$ and preserves the bounds.
\end{proof}

For a deterministic deleted set $\mathcal D$ of fixed size $|\mathcal D|=h$,
put $A=X_{[n]\setminus\mathcal D}$, $m=n-h$, and
\begin{equation}\label{v8:nearby:exterior-degrees}
 P_A=\binom m2^{-1}\sum_{i<j}A_{ij},\qquad
 d_{A,i}=\sum_{j\ne i}A_{ij},\qquad z_i^A=d_{A,i}-(m-1)P_A .
\end{equation}

\begin{corollary}
\label{lem:degree-deletion}
Assume \hyperref[ass:PC]{\textup{(PC)}}.
For the degrees and exterior density defined above, the following bounds
hold for every fixed $q\ge1$ under the original graph law and its actual
exterior marginal:
\begin{equation}\label{row:eq:degree-deletion}
 \begin{aligned}
 \|D_i-(n-1)P\|_q+\|D_i-D_j\|_q&\le C_q\sqrt n,\\
 \|P-P_A\|_q&\le C_{q,h}n^{-3/2},
 \end{aligned}
\end{equation}
\begin{equation}\label{eq:degree-vector-moment}
 \max_{i\in V(A)}\|z_i^A\|_q\le C_{q,h}\sqrt n,
 \qquad \left\|\|z^A\|_2\right\|_q\le C_{q,h}n.
\end{equation}
The maximum is taken outside the moment norm. The constants are uniform
over the deleted labels, and $h=0$ recovers the full graph.
\end{corollary}

\begin{proof}
Take $H=K_2$ (a single edge) in Lemma~\ref{lem:root-exchange}. With the normalization
in \eqref{v8:def:rooted-vector},
\[
 T_i=\frac{d_{A,i}}n,\qquad
 \bar T=\frac{(m-1)P_A}{n},\qquad
 z_i^A=n(T_i-\bar T)=\frac n m\sum_{j\in V(A)}(T_i-T_j).
\]
Equation~\eqref{eq:root-deterministic} and the triangle inequality give
\[
 \begin{aligned}
 \|d_{A,i}-d_{A,j}\|_q&=n\|T_i-T_j\|_q\le C_q\sqrt n,\\
 \|z_i^A\|_q&\le\frac n m\sum_{j\in V(A)}\|T_i-T_j\|_q
                       \le C_q\sqrt n,\\
 \bigl\|\|z^A\|_2\bigr\|_q
 &=n\bigl\|\|T-\bar T\one\|_2\bigr\|_q\le C_qn.
 \end{aligned}
\]
These prove \eqref{eq:degree-vector-moment}; taking $\mathcal D=\varnothing$
also gives the first bound in \eqref{row:eq:degree-deletion}.
Counting removed edges exactly gives
\[
 \binom m2(P-P_A)=\sum_{i\in\mathcal D}\bigl[D_i-(n-1)P\bigr]
                 +\binom h2P-\sum_{\substack{i<j\\i,j\in\mathcal D}}X_{ij}.
\]
Consequently, the full-graph degree bound yields
\begin{equation}\label{eq:fixed-deletion-density}
 \|P-P_A\|_q
 \le\binom m2^{-1}
       \left(\sum_{i\in\mathcal D}\|D_i-(n-1)P\|_q+2\binom h2\right)
 \le C_{q,h}n^{-3/2},
\end{equation}
which completes the proof of \eqref{row:eq:degree-deletion}.
\end{proof}

For the exterior estimates below, write
\begin{equation}\label{v8:def:exterior-moment-order}
 R_A=O_{L^q}(n^{-b})
 \quad\Longleftrightarrow\quad(\E_A|R_A|^q)^{1/q}\le C_qn^{-b},
\end{equation}
with the stated inner vector or matrix norm when appropriate.
In this and the following sections, errors with polynomial envelopes on
exponentially unlikely bad exteriors contribute $O(n^{-M})$ for every fixed $M$.

\subsection{Fixed points and row means}
For $m=|A|=n-O(1)$, let $\pi$ be either the reference law $\pi_A$ or the
marginal of a selected row under a full law $\mu_{A,\mathbf e,s}$ of
Lemma~\ref{lem:local-input}. Write $\mu_i=\E_\pi X_i$.
The reference field $l_{i,A}(x)$ defines the same scalar fixed point for
all these choices; $P_A$ always denotes the observed exterior density.

\begin{lemma}\label{lem:shape}
Assume \hyperref[ass:PC]{\textup{(PC)}}.
For the row marginals specified above, $p_A$ is defined by the common
field $l_{i,A}(x)$ and is independent of the selected marginal.
Estimates are uniform in $\mathbf e,s$ and the selected row for fixed $q_0$.
There is a fixed interval
$I=[p-\delta,p+\delta]\subset(0,1)$ such that, for every
$A\in\mathcal A_n$ of Lemma~\ref{lem:local-input}, the equation
\begin{equation}\label{eq:scalar-row}
 p_A=\sig\left(m^{-1}\sum_{i\in V(A)}l_{i,A}(p_A\one)\right)
\end{equation}
has exactly one solution $p_A\in I$, and $\mu_i\in I$ for every $i$.
The row fixed point $p_A$ and the observed exterior density $P_A$ play
different roles. On bad exteriors, set $p_A=p$.
Then, for every fixed real $q\ge1$,
\begin{equation}\label{eq:mean-shape}
 \begin{aligned}
 \left(\E_A\|\mu-p_A\one\|_2^q\right)^{1/q}&\le C_q,\\
 \left(\E_A\left|\sum_{i\in V(A)}(\mu_i-p_A)\right|^q\right)^{1/q}
 &\le C_q,\\
 \|\mu_i-p_A\|_q&\le C_qn^{-1/2}.
 \end{aligned}
\end{equation}
The outer norms use the actual exterior marginal (including $\mathbf e$
when relevant), as in \eqref{v8:def:exterior-moment-norm}; these are not
pointwise exterior bounds. Any uniformly bounded extension on bad
exteriors gives the same conclusions, uniformly over deleted labels
at each fixed deletion size.
The interval $I$, the good exterior set, and the fixed-cube
contraction and inverse bounds are those of
Lemma~\ref{ep:row:mean-jacobian}, in particular
\eqref{ep:row:symmetric-majorant} and \eqref{eq:shape-common-inverse}.
\end{lemma}

\begin{proof}
\proofstep{A common fixed point.}
Equations~\eqref{ep:row:marginal-mean-equation},
\eqref{ep:row:mean-neighborhood}, and \eqref{ep:row:jacobian-gap} give,
on $\mathcal A_n$, a common cube $I^m$ and $\rho=3/4$ such that
\[
 \begin{gathered}
 \|\mu-p\one\|_\infty\le a+C/n,\qquad
 \max_i|\mu_i-\sig(l_{i,A}(\mu))|\le C/n,\\
 \sup_{y\in I^m}\max\{\|DM_A(y)\|_{\infty\to\infty},
                     \|DM_A(y)\|_{\mathrm{op}}\}\le\rho.
 \end{gathered}
\]
For $\Psi_A(u)=\sig(m^{-1}\sum_i l_{i,A}(u\one))$,
\eqref{ep:row:average-field-contraction} gives a Lipschitz constant $\rho$.
Monotonicity of $\sig(u)$, the mean equation, and the prescribed tolerance $a$ give
\[
 \begin{aligned}
 \left|\sig\left(m^{-1}\sum_i l_{i,A}(\mu)\right)-p\right|&\le a+C/n,\\
 |\Psi_A(p)-p|&\le(1+\rho)a+C/n\le(1-\rho)\delta/2,\\
 |\Psi_A(u)-p|&\le\rho|u-p|+(1-\rho)\delta/2\le\delta\quad(u\in I).
 \end{aligned}
\]
Thus $\Psi_A:I\to I$ has a unique fixed point $p_A$ depending only on $A$.

\proofstep{Euclidean mean error.}
Let
\[
 \bar l_A=m^{-1}\sum_i l_{i,A}(p_A\one).
\]
Then $\sum_i\{l_{i,A}(p_A\one)-\bar l_A\}=0$; on $\mathcal A_n$,
\eqref{eq:scalar-row} gives $\sig(\bar l_A)=p_A$.
Setting $x=u\one$ in \eqref{row:eq:single-row-field} and using
\eqref{v8:def:rooted-vector} gives
\[
 l_{i,A}(u\one)
 =h_n+\sum_r b_{r,n}
 \sum_{\substack{a,b\in V(F_r)\\ab\in E(F_r)}}
 u^{d_{F_r}(a)-1}T_i(F_r-a,b;A).
\]
The inner sum runs over ordered adjacent pairs $(a,b)$. All remaining
vertices of $F_r-a$, including isolates, are retained, so its rooted-count
normalization is $n^{-(v(F_r-a)-1)}=n^{2-v_r}$. Collect terms with
$(G,o)\cong(F_r-a,b)$ and $k=d_{F_r}(a)-1$; the coefficient
$c_{G,o,k,n}$ is the sum of the corresponding $b_{r,n}$, counted with
multiplicity. The chain rule gives
\[
 \frac{\mathrm d}{\mathrm du}l_{i,A}(u\one)
 =\sum_{j\in V(A)}\partial_jl_{i,A}(u\one)
 =(H_A''(u\one)\one)_i.
\]
Thus
\[
 \begin{aligned}
 l_{i,A}(u\one)&=h_n+\sum_{(G,o,k)}c_{G,o,k,n}u^kT_i(G,o;A),\\
 (H_A''(u\one)\one)_i
 &=\sum_{(G,o,k):k\ge1}k c_{G,o,k,n}u^{k-1}T_i(G,o;A).
 \end{aligned}
\]
The finite sums have bounded coefficients shared by all image labels,
including at the random value $u=p_A$. Applying \eqref{eq:root-shared-bounds} gives
\begin{equation}\label{audit:row-centered-field-moments}
 \begin{aligned}
 \bigl\|\|l_A(p_A\one)-\bar l_A\one\|_2\bigr\|_q&\le C_q,\\
 \left\|\left\|H_A''(p_A\one)\one
 -\frac{\one^{\mathsf T}H_A''(p_A\one)\one}{m}\one\right\|_2\right\|_q
 &\le C_q.
 \end{aligned}
\end{equation}
For $A\in\mathcal A_n$, the segment from $p_A\one$ to $\mu$ lies in
$I^m$. Using the reference mean map $M_A(y)$ from \eqref{v8:nearby:app:reference-mean-map}, write
\[
 \begin{aligned}
 \mu-p_A\one={}&\{\mu-M_A(\mu)\}+\{M_A(\mu)-M_A(p_A\one)\}\\
              &+\{M_A(p_A\one)-p_A\one\}.
 \end{aligned}
\]
Equations~\eqref{ep:row:marginal-mean-equation},
\eqref{ep:row:jacobian-gap}, and \eqref{eq:scalar-row} give
\begin{equation}\label{audit:mean-map-error-bounds}
 \begin{aligned}
 \|\mu-M_A(\mu)\|_2
 &\le\sqrt m\max_i|\mu_i-\sig(l_{i,A}(\mu))|\le Cn^{-1/2},\\
 \|M_A(\mu)-M_A(p_A\one)\|_2
 &\le\int_0^1\bigl\|DM_A(p_A\one+t(\mu-p_A\one))
                      (\mu-p_A\one)\bigr\|_2\,\dd t\\
 &\le\rho\|\mu-p_A\one\|_2,\\
 \|M_A(p_A\one)-p_A\one\|_2
 &=\left(\sum_i|\sig(l_{i,A}(p_A\one))-\sig(\bar l_A)|^2\right)^{1/2}\\
 &\le\tfrac14\|l_A(p_A\one)-\bar l_A\one\|_2.
 \end{aligned}
\end{equation}
The last line uses $\sig(\bar l_A)=p_A$ and $\|\sig'\|_\infty\le1/4$.
The triangle inequality and \eqref{audit:mean-map-error-bounds} yield,
on $\mathcal A_n$,
\begin{equation}\label{audit:mean-euclidean-contraction}
 (1-\rho)\|\mu-p_A\one\|_2
 \le\tfrac14\|l_A(p_A\one)-\bar l_A\one\|_2+Cn^{-1/2}.
\end{equation}
Since $p_A=p$ on $\mathcal A_n^c$ and $0\le\mu_i\le1$,
\eqref{audit:row-centered-field-moments},
\eqref{audit:mean-euclidean-contraction}, and \eqref{eq:local-good-exteriors}
give the first bound in \eqref{eq:mean-shape}:
\[
 \E_A\|\mu-p_A\one\|_2^q
 \le C_q+Cn^{q/2}e^{-cn}\le C'_q.
\]

\proofstep{The signed sum and coordinate errors.}
Taylor's formula gives
\begin{equation}\label{audit:mean-field-expansion}
 \begin{aligned}
 l_A(\mu)-\bar l_A\one
 &=l_A(p_A\one)-\bar l_A\one+H_A''(p_A\one)(\mu-p_A\one)+R,\\
 R_i&=\int_0^1(1-t)\sum_{j,k}
    \partial_{jk}l_{i,A}(p_A\one+t(\mu-p_A\one))
    (\mu_j-p_A)(\mu_k-p_A)\,\dd t.
 \end{aligned}
\end{equation}
The global bounds \eqref{row:eq:dense-derivatives}, with $m\le n$,
apply to \eqref{audit:mean-field-expansion} on every exterior. Combining
them with \eqref{audit:row-centered-field-moments} and the first bound in
\eqref{eq:mean-shape} at orders $q$ and $2q$ gives
\begin{equation}\label{audit:mean-field-error-bounds}
 \begin{gathered}
 |R_i|\le Cn^{-2}\|\mu-p_A\one\|_1^2,\qquad
 \|R\|_1\le C\|\mu-p_A\one\|_2^2,\\
 \|R\|_2\le Cn^{-1/2}\|\mu-p_A\one\|_2^2,\qquad
 \|H_A''(p_A\one)\|_{\mathrm{op}}\le C,\\
 \bigl\|\|l_A(\mu)-\bar l_A\one\|_2\bigr\|_q\le C_q.
 \end{gathered}
\end{equation}
On $\mathcal A_n$, expanding $\sig(u)$ at $u=\bar l_A$ and using
\eqref{eq:scalar-row} gives
\begin{equation}\label{audit:mean-shift-taylor}
 \begin{aligned}
 \mu_i-p_A={}&\mu_i-\sig(l_{i,A}(\mu))\\
 &+\nu(p_A)\bigl\{l_{i,A}(p_A\one)-\bar l_A
       +(H_A''(p_A\one)(\mu-p_A\one))_i+R_i\bigr\}+r_i,\\
 |r_i|\le{}&\tfrac12\|\sig''\|_\infty|l_{i,A}(\mu)-\bar l_A|^2,
 \quad \sum_i|r_i|\le C\|l_A(\mu)-\bar l_A\one\|_2^2.
 \end{aligned}
\end{equation}
Summing \eqref{audit:mean-shift-taylor}, using
$\sum_i\{l_{i,A}(p_A\one)-\bar l_A\}=0$ and the symmetry of
$H_A''(p_A\one)$, gives
\begin{equation}\label{audit:signed-mean-equation}
 \begin{aligned}
 &\left(1-\frac{\nu(p_A)}m\one^{\mathsf T}H_A''(p_A\one)\one\right)
       \sum_i(\mu_i-p_A)\\
 &\quad=\nu(p_A)\left\{H_A''(p_A\one)\one
       -\frac{\one^{\mathsf T}H_A''(p_A\one)\one}{m}\one\right\}^{\mathsf T}
       (\mu-p_A\one)\\
 &\qquad+\sum_i\{\mu_i-\sig(l_{i,A}(\mu))\}
       +\nu(p_A)\sum_iR_i+\sum_ir_i.
 \end{aligned}
\end{equation}
Since
\[
 \left|\frac{\nu(p_A)}m\one^{\mathsf T}H_A''(p_A\one)\one\right|
 =|\Psi_A'(p_A)|\le\rho,
\]
we obtain from \eqref{audit:signed-mean-equation},
\eqref{ep:row:marginal-mean-equation}, \eqref{audit:mean-field-error-bounds},
\eqref{audit:mean-shift-taylor}, and Cauchy--Schwarz that
\begin{equation}\label{audit:signed-mean-bound}
 \begin{aligned}
 &(1-\rho)\left|\sum_i(\mu_i-p_A)\right|\\
 &\quad\le C\left\|H_A''(p_A\one)\one
        -\frac{\one^{\mathsf T}H_A''(p_A\one)\one}{m}\one\right\|_2
        \|\mu-p_A\one\|_2\\
 &\qquad+C\{1+\|\mu-p_A\one\|_2^2+\|l_A(\mu)-\bar l_A\one\|_2^2\}.
 \end{aligned}
\end{equation}
Applying H\"older's inequality to \eqref{audit:signed-mean-bound}, using
\eqref{audit:row-centered-field-moments}, the first estimate in
\eqref{eq:mean-shape}, and \eqref{audit:mean-field-error-bounds}
at order $2q$, and adding the bad-exterior contribution from
\eqref{eq:local-good-exteriors}, with $|\sum_i(\mu_i-p_A)|\le m$,
gives the second bound in \eqref{eq:mean-shape}:
\[
 \E_A\left|\sum_i(\mu_i-p_A)\right|^q
 \le C_q+Cn^qe^{-cn}\le C'_q.
\]
For a fixed coordinate, \eqref{eq:root-shared-bounds} gives
\begin{equation}\label{audit:mean-coordinate-field-bound}
 \begin{aligned}
 l_{i,A}(p_A\one)-\bar l_A
 &=m^{-1}\sum_j\{l_{i,A}(p_A\one)-l_{j,A}(p_A\one)\},\\
 \|l_{i,A}(p_A\one)-\bar l_A\|_q&\le C_qn^{-1/2}.
 \end{aligned}
\end{equation}
On $\mathcal A_n$, \eqref{row:eq:dense-derivatives} and
\eqref{ep:row:marginal-mean-equation} imply
\[
 |\mu_i-p_A|\le Cn^{-1}\|\mu-p_A\one\|_1
        +\tfrac14|l_{i,A}(p_A\one)-\bar l_A|+C/n.
\]
Using \eqref{audit:mean-coordinate-field-bound}, the first estimate in
\eqref{eq:mean-shape}, and $\|\mu-p_A\one\|_1\le\sqrt m\|\mu-p_A\one\|_2$,
and adding the bad-exterior contribution from \eqref{eq:local-good-exteriors},
with $|\mu_i-p_A|\le1$ on $\mathcal A_n^c$, gives the third bound in
\eqref{eq:mean-shape}:
\[
 \E_A|\mu_i-p_A|^q
 \le C_qn^{-q/2}+Ce^{-cn}\le C'_qn^{-q/2}.
\]
\end{proof}

\subsection{Graph counts at the empirical density}
For $m=|A|$, let $J_m,I_m$ be the all-ones and identity matrices and define
\begin{equation}\label{v8:def:count-centered-adjacency}
 Y=X-P(J_n-I_n),\qquad Y_A=A-P_A(J_m-I_m).
\end{equation}
Both matrices have zero diagonal. Their weighted injective counts use
the formula in \eqref{v8:def:graph-counts}, with the centered matrix
entries in place of the adjacency entries.

\begin{theorem}
\label{thm:counts}
Assume \hyperref[ass:PC]{\textup{(PC)}}.
For every fixed finite simple graph $H$, with its isolated
vertices retained, and every fixed real $q\ge1$, the normalized injective count satisfies
\begin{equation}\label{eq:count-replacement}
 \left\|f_H(X)-P^{e(H)}\right\|_q\le C_{H,q}n^{-1},
 \qquad f_H(X)=n^{-v(H)}\inj(H,X).
\end{equation}
For a deterministic deletion $S\subseteq[n]$ of fixed size $h$,
set $A=X_{[n]\setminus S}$ and $m=n-h$. Then
\begin{equation}\label{audit:exterior-count-replacement}
 \left\|m^{-v(H)}\inj(H,A)-P_A^{e(H)}\right\|_q
 \le C_{H,q,h}n^{-1},
 \qquad
 P_A=\binom m2^{-1}\sum_{i<j:\ i,j\notin S}X_{ij}.
\end{equation}
The constants are uniform over the choice of $S$; the norms are taken
under the marginal exterior law induced by the original graph law.
If $e(H)\ge1$, then the corresponding centered weighted counts also satisfy
\begin{equation}\label{audit:centered-count-replacement}
 \begin{aligned}
 \left\|n^{-v(H)}\inj\bigl(H,X-P(J_n-I_n)\bigr)\right\|_q
       &\le C_{H,q}n^{-1},\\
 \left\|m^{-v(H)}\inj\bigl(H,A-P_A(J_m-I_m)\bigr)\right\|_q
       &\le C_{H,q,h}n^{-1}.
 \end{aligned}
\end{equation}
Here $J_m$ and $I_m$ are the all-ones and identity matrices of order $m$,
respectively. The diagonal remains zero in these weighted counts.
These estimates do not require localization on the critical
$n^{-1/\mc}$ scale.
\end{theorem}

\begin{proof}
\proofstep{1. Compare the row fixed points with $P$.}
Fix $q\ge1$.
For a fixed-deletion exterior $A$, $m_A=n-h_A$, summing the exact row field gives
\begin{equation}\label{eq:averaged-exterior-row-field}
 \frac1{m_A}\sum_i l_{i,A}(u\one)
 =h_n+\frac n{m_A}\sum_{r=1}^k b_{r,n}
       \sum_{v\in V(F_r)}d_{F_r}(v)u^{d_{F_r}(v)-1}
            n^{-(v_r-1)}\inj(F_r-v,A).
\end{equation}
Terms of degree zero are omitted. The factor
$d_{F_r}(v)u^{d_{F_r}(v)-1}$ comes from differentiating and summing
the row monomial, with $n^{2-v_r}/m_A=(n/m_A)n^{-(v_r-1)}$.
For two fixed-deletion exteriors $A,B$ of the same graph,
\begin{equation}\label{audit:dependency:s56:fixed-deletion-count-comparison}
 |n^{-v(G)}\inj(G,A)-n^{-v(G)}\inj(G,B)|
 \le(h_A+h_B)v(G)/n,\qquad
 n/m_A=1+O(n^{-1}).
\end{equation}
Indeed, every removed injection uses a deleted vertex. Applying
\eqref{audit:dependency:s56:fixed-deletion-count-comparison} to
\eqref{eq:averaged-exterior-row-field}, with bounded coefficients and
$\|\sig'\|_\infty\le1/4$, gives
\[
 \sup_{u\in I}\left|
 \sig\!\left(\frac1{m_A}\sum_i l_{i,A}(u\one)\right)
 -\sig\!\left(\frac1{m_B}\sum_i l_{i,B}(u\one)\right)
 \right|\le C/n.
\]
When both exteriors are good, \eqref{eq:scalar-row} and
\eqref{ep:row:average-field-contraction} yield
$(1-\rho)|p_A-p_B|\le C/n$; their bounded extensions therefore satisfy
\begin{equation}\label{eq:row-fixed-point-reference-comparison}
 \E|p_A-p_B|^q\le\left(\frac C{(1-\rho)n}\right)^q+Ce^{-cn},
 \qquad \|p_A-p_B\|_q\le C_qn^{-1}.
\end{equation}
This comparison needs only the two indicated good-exterior events.
For $A_i=X_{[n]\setminus\{i\}}$, the signed-sum estimate in
\eqref{eq:mean-shape} gives
\[
 \mathcal E_iD_i=(n-1)p_{A_i}+\varepsilon_i,
 \qquad \sup_i\|\varepsilon_i\|_q\le C_q.
\]
Since $\sum_iD_i=2NP$, whole-row refreshment gives the exact identity
\begin{equation}\label{audit:dependency:s56:density-generator-identity}
 \mathcal GP=2\left(n^{-1}\sum_i p_{A_i}-P\right)
                  +N^{-1}\sum_i\varepsilon_i,
 \qquad \|N^{-1}\sum_i\varepsilon_i\|_q\le C_q/n.
\end{equation}
Moreover, $f_{K_2}=(n-1)P/n$, so \eqref{eq:whole-row} gives
$\|\mathcal GP\|_q\le C_q/n$.
Using \eqref{eq:row-fixed-point-reference-comparison},
\eqref{audit:dependency:s56:density-generator-identity}, and Minkowski's inequality, we obtain
\begin{equation}\label{eq:scalar-reference-density-comparison}
 \begin{aligned}
 \|p_A-P\|_q
 &\le n^{-1}\sum_i\|p_A-p_{A_i}\|_q
      +\tfrac12\|\mathcal GP\|_q+(2N)^{-1}\|\sum_i\varepsilon_i\|_q\\
 &\le C_qn^{-1}.
 \end{aligned}
\end{equation}
No graph-count replacement has been used in this step.

\proofstep{2. Derive the count recursion.}
For $u\in V(H)$ and $i\in[n]$, set
\[
 T_{i,u}=n^{1-v(H)}
    \sum_{\substack{\varphi:V(H)\hookrightarrow[n]\\\varphi(u)=i}}
       \prod_{ab\in E(H)}X_{\varphi(a)\varphi(b)}.
\]
An injection uses $i$ at most once, and $\sum_iT_{i,u}=nf_H$; thus
\begin{equation}\label{audit:dependency:s56:count-generator-identity}
 \mathcal Gf_H=-v(H)f_H+n^{-1}\sum_i\sum_{u\in V(H)}\mathcal E_iT_{i,u}.
\end{equation}
Condition on $A_i$ and write $\mu_j=\E[X_{ij}\mid A_i]$.
Fixing $d_H(u)$ neighbor images leaves
$O(n^{v(H)-1-d_H(u)})$ extensions, so the row coefficients are
$O(n^{-d_H(u)})$. Equation~\eqref{eq:row-polynomial-mean}, with
\eqref{eq:local-good-exteriors} and bounded bad-exterior envelopes, gives
\begin{equation}\label{audit:dependency:s56:conditional-root-mean}
 \mathcal E_iT_{i,u}=T_{i,u}(\mu;A_i)+O_{L^q}(n^{-1}).
\end{equation}
The error is zero when $d_H(u)=0$.
Expand at $p_{A_i}\one$. For $d_H(u)\ge1$, the linear coefficients are
\[
 a_j=\frac1n\sum_{\substack{v\in V(H)\\uv\in E(H)}}p_{A_i}^{d_H(u)-1}T_j(H-u,v;A_i),
 \qquad \bar a=|A_i|^{-1}\sum_ja_j.
\]
Here $T_j(H-u,v;A_i)$ retains all vertices of $H-u$ and has
normalization $n^{-(v(H)-2)}$. The shared-coefficient bounds
\eqref{eq:root-shared-bounds} and the mean-shape bounds
\eqref{eq:mean-shape} imply
\begin{equation}\label{audit:dependency:s56:count-linear-pairing}
 \begin{gathered}
 |\bar a|\le C_H/n,\qquad
 \bigl\|\|a-\bar a\one\|_2\bigr\|_q\le C_{H,q}/n,\\
 a^{\mathsf T}(\mu-p_{A_i}\one)=\bar a\one^{\mathsf T}(\mu-p_{A_i}\one)
               +(a-\bar a\one)^{\mathsf T}(\mu-p_{A_i}\one),\\
 \|a^{\mathsf T}(\mu-p_{A_i}\one)\|_q
 \le\frac{C_H}{n}\|\one^{\mathsf T}(\mu-p_{A_i}\one)\|_q
 +\bigl\|\|a-\bar a\one\|_2\bigr\|_{2q}
     \bigl\|\|\mu-p_{A_i}\one\|_2\bigr\|_{2q}\le C_{H,q}/n.
 \end{gathered}
\end{equation}
For a degree-$b$ term with $2\le b\le d_H(u)$, fixing its $b$
row sites leaves at most $C_Hn^{v(H)-1-b}$ template injections.
Its coefficient is therefore at most
$C_Hn^{1-v(H)}n^{v(H)-1-b}=C_Hn^{-b}$.
After taking absolute values and dropping the distinct-site restrictions,
the sum of the degree-$b$ terms is at most $C_Hn^{-b}\|\mu-p_{A_i}\one\|_1^b$.
Thus the terms beyond the linear part satisfy
\begin{equation}\label{audit:dependency:s56:count-higher-remainder}
 \begin{aligned}
 &\left\|T_{i,u}(\mu;A_i)-T_{i,u}(p_{A_i}\one;A_i)
                         -a^{\mathsf T}(\mu-p_{A_i}\one)\right\|_q\\
 &\le\sum_{b=2}^{d_H(u)}C_Hn^{-b}
                         \bigl\|\|\mu-p_{A_i}\one\|_1^b\bigr\|_q\\
 &\le\sum_{b=2}^{d_H(u)}C_{H,b}n^{-b/2}
                         \bigl\|\|\mu-p_{A_i}\one\|_2\bigr\|_{bq}^b
 \le C_{H,q}/n,
 \end{aligned}
\end{equation}
where \eqref{eq:mean-shape} bounds the fixed moments of $\|\mu-p_{A_i}\one\|_2$.
The sum is empty for $d_H(u)=1$; for $d_H(u)=0$, the row polynomial is
constant. Combining \eqref{audit:dependency:s56:conditional-root-mean},
\eqref{audit:dependency:s56:count-linear-pairing}, and
\eqref{audit:dependency:s56:count-higher-remainder} gives
\[
 \mathcal E_iT_{i,u}
 =p_{A_i}^{d_H(u)}n^{-(v(H)-1)}\inj(H-u,A_i)+O_{L^q}(n^{-1}).
\]
Equation~\eqref{audit:dependency:s56:fixed-deletion-count-comparison}
replaces the count by $f_{H-u}(X)$;
\eqref{eq:scalar-reference-density-comparison} and
$|x^d-y^d|\le d|x-y|$ replace $p_{A_i}^{d_H(u)}$ by $P^{d_H(u)}$.
Substituting this expansion in
\eqref{audit:dependency:s56:count-generator-identity} and averaging over $i$ proves
\begin{equation}\label{audit:graph-count-recursion}
 \mathcal Gf_H=-v(H)f_H+
       \sum_{u\in V(H)}P^{d_H(u)}f_{H-u}(X)+O_{L^q}(n^{-1}).
\end{equation}

\proofstep{3. Induction, deletions and centering.}
Since $e(H-u)+d_H(u)=e(H)$, the recursion
\eqref{audit:graph-count-recursion} and Minkowski's inequality give
\begin{equation}\label{v8:proof:counts:inductive-reduction}
 \begin{aligned}
 v(H)\|f_H-P^{e(H)}\|_q
 &\le\|\mathcal Gf_H\|_q+
   \sum_{u\in V(H)}\|P^{d_H(u)}(f_{H-u}-P^{e(H-u)})\|_q+C_{H,q}/n.
 \end{aligned}
\end{equation}
The norms are taken under the original graph law or its actual exterior marginals.

For $v$ isolated vertices, $f_H=(n)_v/n^v=1+O(n^{-1})$;
the empty graph has count one.
Since $0\le P\le1$, induction on $v(H)$ gives
\[
 \|P^{d_H(u)}f_{H-u}-P^{e(H)}\|_q
 \le\|f_{H-u}-P^{e(H-u)}\|_q\le C_{H,q}/n.
\]
Insert this estimate and \eqref{eq:whole-row} into
\eqref{v8:proof:counts:inductive-reduction} to obtain
\eqref{eq:count-replacement}.
For $m=n-h$, \eqref{audit:dependency:s56:fixed-deletion-count-comparison}
and the change of normalization give
\[
 \begin{gathered}
 |n^{-v(H)}\inj(H,X)-n^{-v(H)}\inj(H,A)|\le C_{H,h}/n,\\
 |m^{-v(H)}\inj(H,A)-n^{-v(H)}\inj(H,A)|
       \le|1-(m/n)^{v(H)}|\le C_{H,h}/n,\qquad
 |P-P_A|\le C_h/n.
 \end{gathered}
\]
These prove \eqref{audit:exterior-count-replacement}.
Finally, retaining every vertex in each spanning subgraph and applying
\eqref{eq:count-replacement}, expand
\[
 \begin{aligned}
 n^{-v(H)}\inj(H,Y)
 &=\sum_{S\subseteq E(H)}(-P)^{e(H)-|S|}f_{(V(H),S)}(X)\\
 &=P^{e(H)}\sum_{S\subseteq E(H)}(-1)^{e(H)-|S|}
                  +O_{L^q}(n^{-1})=O_{L^q}(n^{-1})
 \end{aligned}
\]
because $e(H)\ge1$. The identical expansion with $m,P_A,A$ and
\eqref{audit:exterior-count-replacement} proves the centered bound for the exterior.
\end{proof}

\subsection{Conditional logits and fixed-point residuals}
Let $X^{e,a}$ be obtained from $X$ by setting $X_e=a$, and define
\begin{equation}\label{v8:setup:physical-edge-logit}
 L_e=\mathcal H_n(X^{e,1})-\mathcal H_n(X^{e,0}),\qquad
 \ell_n(u)=h_n+2\sum_r b_{r,n}e_r u^{e_r-1}.
\end{equation}
We first compute the edge difference defining $L_e$ and then evaluate
the centering field $\ell_n$ at $P$. The following lemma bounds the
conditional logit bias. Section~\ref{sec:drift} derives the averaged
logit correction and higher-moment estimates using the signed bounds
of Section~\ref{sec:supports}.
For the two-row interpolation $\mu_{A,(e),s}$, $0\le s\le1$, with $e=X_{12}$, put
$L_{12}(x,y;A)=\partial_{12}\mathcal H_n(X(A,e;x,y))$.

\begin{lemma}\label{lem:conditional-logit-bias}
Assume \hyperref[ass:PC]{\textup{(PC)}}.
With $\ell_n(x)$ defined in \eqref{v8:setup:physical-edge-logit},
the two-row laws $\mu_{A,(e),s}$ of
\eqref{v8:def:conditional-multirow-law} satisfy, for every fixed moment order $q\ge1$,
\begin{equation}\label{eq:bias}
 \sup_{0\le s\le1}
 \left\|\E_{\mu_{A,(e),s}}L_{12}-\ell_n(P_A)\right\|_q
     \le C_qn^{-1}.
\end{equation}
At $s=1$, the expectation is $\E[L_{12}\mid A,e]$.
The outer norms follow \eqref{v8:def:exterior-moment-norm}
under the actual $(A,e)$ marginal, with the supremum taken outside the norm.
These moment bounds are uniform over the parameter range in
\eqref{eq:parameter-neighborhood}; they are not pointwise exterior estimates.
\end{lemma}

\begin{proof}
Fix $(A,e,s)$. Invariance of $\mu_{A,(e),s}$ under $x\leftrightarrow y$
gives the common row mean
\[
 \mu_i=\E_{\mu_{A,(e),s}}x_i
      =\E_{\mu_{A,(e),s}}y_i,\qquad i\in V(A).
\]
We estimate the three terms in
\begin{equation}\label{v8:proof:conditional-logit-bias:decomposition}
 \begin{aligned}
 \E_{\mu_{A,(e),s}}L_{12}-\ell_n(P_A)
 &=\{\E_{\mu_{A,(e),s}}L_{12}-L_{12}(\mu,\mu;A)\}\\
 &\quad+\{L_{12}(\mu,\mu;A)-L_{12}(p_A\one,p_A\one;A)\}\\
 &\quad+\{L_{12}(p_A\one,p_A\one;A)-\ell_n(P_A)\}.
 \end{aligned}
\end{equation}

\proofstep{Polynomial mean error.}
Differentiating with respect to the template edge mapped to $12$ gives exactly
\begin{equation}\label{eq:conditional-logit-polynomial}
 \begin{aligned}
 L_{12}(x,y;A)-h_n
 ={}&\sum_{r=1}^k b_{r,n}
 \sum_{\substack{(u,v)\in V(F_r)^2\\uv\in E(F_r)}}
 n^{-(v_r-2)}
 \sum_{\psi:V(F_r)\setminus\{u,v\}\hookrightarrow V(A)}\\
 &\quad\cdot\prod_{ab\in E(F_r-\{u,v\})}A_{\psi(a)\psi(b)}
 \prod_{\substack{a\in V(F_r)\setminus\{u,v\}\\ua\in E(F_r)}}x_{\psi(a)}
 \prod_{\substack{b\in V(F_r)\setminus\{u,v\}\\vb\in E(F_r)}}y_{\psi(b)}.
 \end{aligned}
\end{equation}
There are $2e_r$ orientations; $L_{12}$ itself is independent of $e$.
Fixing $d$ exterior sites in this expansion gives
\begin{equation}\label{eq:conditional-logit-coefficients}
 |c_{i_1,\ldots,i_d;S_1,\ldots,S_d}(A)|
 \le\sum_{r:v_r\ge d+2}
 C_{F_r}|b_{r,n}|n^{2-v_r}n^{v_r-d-2}
 \le C_dn^{-d}.
\end{equation}
Thus $L_{12}$ is a normalized block polynomial in $x_i,y_i,x_iy_i$.
Equation~\eqref{eq:row-polynomial-mean}, applied under $\mu_{A,(e),s}$,
gives, on good exteriors,
\begin{equation}\label{audit:dependency:s56:logit-polynomial-mean}
 \left|\E_{\mu_{A,(e),s}}[L_{12}(x,y;A)]
              -L_{12}(\mu,\mu;A)\right|\le C/n.
\end{equation}

\proofstep{Mean shift.}
Set
$a_i=(\partial_{x_i}L_{12}+\partial_{y_i}L_{12})(p_A\one,p_A\one;A)$,
with $\bar a=|A|^{-1}\sum_i a_i$.
Fixing the exterior image $i$ shows that $a_i$ is $n^{-1}$ times
a finite sum of singly rooted counts, with shared bounded coefficients.
Hence \eqref{eq:root-shared-bounds} gives
$|\bar a|\le C/n$ and
$\|\|a-\bar a\one\|_2\|_q\le C_q/n$.
By \eqref{eq:mean-shape},
\begin{equation}\label{v8:proof:conditional-logit-bias:linear-pairing}
 \begin{aligned}
 a^{\mathsf T}(\mu-p_A\one)
 &=\bar a\one^{\mathsf T}(\mu-p_A\one)+(a-\bar a\one)^{\mathsf T}(\mu-p_A\one),\\
 \|a^{\mathsf T}(\mu-p_A\one)\|_q
 &\le\frac Cn\|\one^{\mathsf T}(\mu-p_A\one)\|_q
  +\bigl\|\|a-\bar a\one\|_2\bigr\|_{2q}
     \bigl\|\|\mu-p_A\one\|_2\bigr\|_{2q}\le C_q/n.
 \end{aligned}
\end{equation}
Same-site derivatives satisfy $|\partial_{x_i y_i}L_{12}|\le C/n$;
derivatives involving two distinct sites are $O(n^{-2})$.
Taylor's integral remainder along the segment from
$(p_A\one,p_A\one)$ to $(\mu,\mu)$ therefore obeys
\begin{equation}\label{audit:dependency:s56:logit-taylor-remainder}
 \begin{aligned}
 &\left|L_{12}(\mu,\mu;A)-L_{12}(p_A\one,p_A\one;A)
                  -a^{\mathsf T}(\mu-p_A\one)\right|\\
 &\quad\le\frac Cn\|\mu-p_A\one\|_2^2
          +\frac C{n^2}\|\mu-p_A\one\|_1^2
 \le\frac{C'}n\|\mu-p_A\one\|_2^2.
 \end{aligned}
\end{equation}
Its $L^q$ norm is $O(n^{-1})$ by the first estimate in
\eqref{eq:mean-shape} at order $2q$.
These estimates use H\"older's inequality under the actual joint exterior law, without any independence assumption.

\proofstep{Constant rows.}
At $p_A\one$, \eqref{eq:conditional-logit-polynomial} gives
\[
 L_{12}(p_A\one,p_A\one;A)
 =h_n+\sum_r b_{r,n}
   \sum_{\substack{u,v\in V(F_r)\\uv\in E(F_r)}}
    p_A^{d_{F_r}(u)+d_{F_r}(v)-2}
       n^{-(v_r-2)}\inj(F_r-\{u,v\},A).
\]
Equations~\eqref{eq:scalar-reference-density-comparison},
\eqref{audit:dependency:s56:fixed-deletion-count-comparison}, and
\eqref{audit:exterior-count-replacement} give
\[
 \|p_A-P_A\|_q\le C_q/n,\qquad
 n^{-(v_r-2)}\inj(F_r-\{u,v\},A)
 =P_A^{e_r-d_{F_r}(u)-d_{F_r}(v)+1}+O_{L^q}(n^{-1}).
\]
The exponents add to $e_r-1$, so the $2e_r$ orientations give
\begin{equation}\label{audit:dependency:s56:logit-constant-rows}
 L_{12}(p_A\one,p_A\one;A)=\ell_n(P_A)+O_{L^q}(n^{-1}).
\end{equation}
Combining \eqref{v8:proof:conditional-logit-bias:decomposition},
\eqref{audit:dependency:s56:logit-polynomial-mean},
\eqref{v8:proof:conditional-logit-bias:linear-pairing},
\eqref{audit:dependency:s56:logit-taylor-remainder}, and
\eqref{audit:dependency:s56:logit-constant-rows} proves \eqref{eq:bias}.
All conditional bounds are uniform in $s$; bounded fields and
\eqref{eq:local-good-exteriors} justify averaging under the actual
$(A,e)$ marginal, with the supremum over $s$ outside the norm.
\end{proof}

\begin{corollary}
\label{cor:coarse}
Assume \hyperref[ass:PC]{\textup{(PC)}}.
For every fixed real $q\ge1$, the edge density satisfies
\begin{equation}\label{audit:scalar-defect-full}
 \|\sig(\ell_n(P))-P\|_q+\|g(P)\|_q\le C_q/n.
\end{equation}
If $A$ is obtained by any fixed deterministic vertex deletion and
$P_A$ is its edge density, then
\begin{equation}\label{audit:scalar-defect-exterior}
 \|\sig(\ell_n(P_A))-P_A\|_q+\|g(P_A)\|_q\le C_q/n.
\end{equation}
The norms are defined in
\eqref{v8:def:conditional-moment-norm}--\eqref{v8:def:exterior-moment-norm};
the constants are uniform over deleted labels for each fixed deletion size.
Moreover, there exists $0<\delta_g<\min\{p,1-p\}$, depending only on the
leading model, such that
\begin{equation}\label{eq:scalar-local-cubic-drift}
 g(p+t)=-2\nu(p)a_\star t^{2\mc-1}+O(|t|^{2\mc}),
 \qquad |t|\le\delta_g.
\end{equation}
Here $2\nu(p)a_\star>0$, and the remainder is uniform on this interval,
with its implicit constant depending only on the leading model.
\end{corollary}

\begin{proof}
Delete one vertex and write $m=n-1$.
Since $d_{F_r}(v)-1+e(F_r-v)=e_r-1$ and
$\sum_vd_{F_r}(v)=2e_r$, equations
\eqref{eq:averaged-exterior-row-field}, \eqref{eq:scalar-reference-density-comparison}
and \eqref{audit:exterior-count-replacement} give
\begin{equation}\label{audit:dependency:s56:averaged-field-at-density}
 m^{-1}\sum_i l_{i,A}(p_A\one)=\ell_n(P)+O_{L^q}(n^{-1}).
\end{equation}
Here the exterior counts retain isolated vertices and $n/m=1+O(n^{-1})$.
Equations~\eqref{eq:scalar-row},
\eqref{audit:dependency:s56:averaged-field-at-density}, and
\eqref{eq:scalar-reference-density-comparison}, with
$\|\sig'\|_\infty\le1/4$, yield
\[
 \|\sig(\ell_n(P))-P\|_q
 \le\|P-p_A\|_q+\tfrac14
       \left\|\ell_n(P)-m^{-1}\sum_i l_{i,A}(p_A\one)\right\|_q
       +Ce^{-cn/q}\le C_q/n.
\]
Equation~\eqref{eq:local-good-exteriors} controls the exponential term
from bounded fields on bad exteriors.
The parameter window gives
\begin{equation}\label{eq:scalar-field-window-difference}
 \begin{aligned}
 \ell_n(x)-\ell(x)
 &=(h_n-h_*)+2\sum_r(b_{r,n}-b_r)e_r x^{e_r-1}
   =\frac{\kappa_n+2\mathcal U_n'(x)}n,\\
 \|\ell_n-\ell\|_{C^j([0,1])}&\le C_j/n\qquad(j\ge0\text{ fixed}).
 \end{aligned}
\end{equation}
Combining this defect bound with \eqref{eq:scalar-field-window-difference}
and \eqref{v8:def:scalar-logistic-field} gives $\|g(P)\|_q\le C_q/n$.
Both defect functions have uniformly bounded derivatives; combining this
with \eqref{eq:fixed-deletion-density} proves
\eqref{audit:scalar-defect-exterior}.
Finally, \eqref{eq:scalar}, \eqref{eq:critical}, and
\eqref{v8:def:scalar-logistic-field} give
\[
 \ell(x)-\log\frac{x}{1-x}=2\Phi'(x),\qquad
 \Phi'(p+t)=-a_\star t^{2\mc-1}+O(|t|^{2\mc}).
\]
Since $\sig(\log(x/(1-x)))=x$ and
$\sig'(\log(x/(1-x)))=\nu(x)$, Taylor's formula yields
\[
 g(x)=\sig\left(\log\frac{x}{1-x}+2\Phi'(x)\right)-x
     =2\nu(x)\Phi'(x)+O(\Phi'(x)^2).
\]
Using $\nu(p+t)=\nu(p)+O(|t|)$ and $4\mc-2\ge2\mc$ gives
\[
 \begin{aligned}
 g(p+t)&=-2\nu(p+t)a_\star t^{2\mc-1}
          +O(|t|^{2\mc})+O(|t|^{4\mc-2})\\
       &=-2\nu(p)a_\star t^{2\mc-1}+O(|t|^{2\mc}).
 \end{aligned}
\]
Taking $\delta_g>0$ sufficiently small proves \eqref{eq:scalar-local-cubic-drift}.
\end{proof}

\subsection{Centered matrix squares and triangle moments}
\begin{corollary}\label{cor:matrix-moments}
Assume \hyperref[ass:PC]{\textup{(PC)}}.
For the centered matrices in \eqref{v8:def:count-centered-adjacency},
every deterministic deletion of a fixed number $h$ of vertices satisfies
\begin{equation}\label{eq:four-cycle-square-frobenius}
 \bigl\|\|Y_A^2\|_F\bigr\|_q\le C_{q,h}n^{3/2},\qquad q\ge1 .
\end{equation}
The constants are uniform over the deleted labels. The full graph also
satisfies
\begin{equation}\label{eq:triangle-moment}
 \|\tr Y^3\|_q=\|\inj(K_3,Y)\|_q\le C_qn^2 .
\end{equation}
\end{corollary}
\begin{proof}
In $\tr Y_A^4$, a nonzero noninjective term has $i=k$ or $j=l$.
Writing $C_4$ for the four-cycle, inclusion--exclusion gives
\begin{equation}\label{eq:four-cycle-trace-identity}
 \tr Y_A^4=\inj(C_4,Y_A)
       +2\sum_i\left(\sum_j(Y_A)_{ij}^2\right)^2
       -\sum_{i,j}(Y_A)_{ij}^4.
\end{equation}
Hence $0\le\|Y_A^2\|_F^2=\tr Y_A^4\le|\inj(C_4,Y_A)|+2m^3$.
Equation~\eqref{audit:centered-count-replacement} gives
\begin{equation}\label{eq:four-cycle-trace-moment}
 \|\inj(C_4,Y_A)\|_q+\|\tr Y_A^4\|_q\le C_{q,h}n^3,
 \qquad q\ge1.
\end{equation}
By moment monotonicity and \eqref{eq:four-cycle-trace-moment}, we have
\[
 \begin{aligned}
 \bigl\|\|Y_A^2\|_F\bigr\|_q
 &\le\bigl\|\|Y_A^2\|_F\bigr\|_{\max\{2,q\}}
   =\|\tr Y_A^4\|_{\max\{1,q/2\}}^{1/2}
 \le C_{q,h}n^{3/2},
 \end{aligned}
\]
which proves \eqref{eq:four-cycle-square-frobenius}.
The zero diagonal gives $\tr Y^3=\inj(K_3,Y)$;
\eqref{audit:centered-count-replacement} with $H=K_3$ proves
\eqref{eq:triangle-moment}.
\end{proof}

\section{First-order expansions of rooted counts}\label{sec:refined}
We establish the rooted expansions used in the reference-row analysis.
The coarse mean estimate, a centered-path bound, and a change-of-center
bound provide the required inputs for the rooted Taylor expansions.

\subsection{Row estimates and recentering}
Use the fixed-deletion exterior and degree notation of
\eqref{v8:nearby:exterior-degrees}, and the centered matrix defined in
\eqref{v8:def:count-centered-adjacency}.
For $\pi_A$, write $\mu_A=(\E_{\pi_A}X_i)_{i\in V(A)}$.

\begin{lemma}\label{lem:row-mean-density}
Assume \hyperref[ass:PC]{\textup{(PC)}}.
For $\mu_A=\E_{\pi_A}X$ and every fixed $q\ge1$,
\begin{equation}\label{eq:mean-profiles-a}
 \begin{aligned}
 \bigl\|\|\mu_A-P_A\one\|_2\bigr\|_q
   +\|\one^{\mathsf T}(\mu_A-P_A\one)\|_q&\le C_q,\\
 \|\mu_{A,i}-P_A\|_q&\le C_qn^{-1/2}.
 \end{aligned}
\end{equation}
The same bounds hold for
$(\E_{\mu_{A,\mathbf e,s}}X_{a i})_{i\in V(A)}$ under
\eqref{v8:def:conditional-multirow-law}, uniformly in $a,\mathbf e,s$
for fixed $q_0$; the outer norms use the actual $(A,\mathbf e)$ marginal.
\end{lemma}

\begin{proof}
Equations~\eqref{eq:scalar-reference-density-comparison} and
\eqref{eq:fixed-deletion-density} give
\begin{equation}\label{eq:root-scalar-comparison}
 \|p_A-P_A\|_q\le\|p_A-P\|_q+\|P-P_A\|_q
 \le C_q(n^{-1}+n^{-3/2})\le C'_qn^{-1}.
\end{equation}
Replacing $p_A$ by $P_A$ in \eqref{eq:mean-shape} adds, respectively,
$\sqrt m\|p_A-P_A\|_q$, $m\|p_A-P_A\|_q$, and $\|p_A-P_A\|_q$.
By \eqref{eq:root-scalar-comparison}, these are
$O(n^{-1/2})$, $O(1)$, and $O(n^{-1})$, proving \eqref{eq:mean-profiles-a}.
Since $p_A$ is $A$-measurable, the same argument holds under the actual
$(A,\mathbf e)$ marginal for every selected row.
\end{proof}

The next estimate controls the Hessian acting on the degree correction
in the mean equation. We use the degree vectors and centered matrices from
\eqref{v8:nearby:exterior-degrees} and \eqref{v8:def:count-centered-adjacency};
$h=0$ gives $Y,z$.

\begin{lemma}\label{lem:Yz}
Assume \hyperref[ass:PC]{\textup{(PC)}}.
For the fixed-deletion exterior above and every fixed real $q\ge1$,
the following bounds hold uniformly in the retained and deleted labels:
\begin{equation}\label{eq:Yz}
 \max_{i\in V(A)}\|(Y_Az^A)_i\|_q\le C_qn,\qquad
 \left\|\|Y_Az^A\|_2\right\|_q\le C_qn^{3/2}.
\end{equation}
The maximum in the first bound is outside the $L^q$ norm. Both bounds
also hold with $A=X$, $Y_A=Y$, and $z^A=z$.
\end{lemma}

\begin{proof}
Fix $i\in V(A)$ and put $m=|V(A)|$ and
$B=A-i$. Thus
$V(B)=V(A)\setminus\{i\}$ and $|V(B)|=m-1$.
Condition on $B$ and the internal edges $\mathbf e$ among
$\mathcal D\cup\{i\}$. The full conditional law is
$\mu_{B,\mathbf e,1}$ on $h+1$ rows, including the deleted rows.
Put $N_A=\binom m2$. For $j\in V(B)$, the mean of the selected row
is $\E[X_{ij}\mid B,\mathbf e]$. Exact counting gives
\[
 \begin{gathered}
 N_A(P_A-P_B)=\sum_{j\in V(B)}(X_{ij}-P_B),\qquad
 \sum_{j\in V(B)}z_j^B=0,\\
 z_j^A=z_j^B+(X_{ij}-P_B)-(m-1)(P_A-P_B),\qquad j\in V(B).
 \end{gathered}
\]
Since $(Y_A)_{ii}=0$, the matrix product is
\[
 (Y_Az^A)_i=\sum_{j\in V(A)}(Y_A)_{ij}z_j^A
          =\sum_{j\in V(B)}(X_{ij}-P_A)z_j^A.
\]
For $j\in V(B)$, substitute
$(Y_A)_{ij}=X_{ij}-P_B-(P_A-P_B)$ and split at
$\E[X_{ij}\mid B,\mathbf e]$:
\begin{equation}\label{eq:Yz-full-decomposition}
 \begin{aligned}
 (Y_Az^A)_i
 ={}&\sum_{j\in V(B)}(X_{ij}-\E[X_{ij}\mid B,\mathbf e])z_j^B\\
 &+\sum_{j\in V(B)}(\E[X_{ij}\mid B,\mathbf e]-P_B)z_j^B
       +\sum_{j\in V(B)}(X_{ij}-P_B)^2\\
 &-m(P_A-P_B)\sum_{j\in V(B)}(X_{ij}-P_B)
       +(m-1)^2(P_A-P_B)^2.
 \end{aligned}
\end{equation}
Since $|(P_A-P_B)|\le2/m$,
\begin{equation}\label{eq:Yz-elementary-terms}
 \begin{gathered}
 \sum_{j\in V(B)}(X_{ij}-P_B)^2\le m-1,\\
 \left|m(P_A-P_B)\sum_{j\in V(B)}(X_{ij}-P_B)\right|\le2(m-1),\qquad
 (m-1)^2(P_A-P_B)^2\le4.
 \end{gathered}
\end{equation}
Equation~\eqref{eq:degree-vector-moment} gives
$\|\|z^B\|_2\|_q\le C_{q,h}n$ for every fixed $q\ge1$.
On good $B$, \eqref{eq:weighted-local-moments}, with zero weights on the other
rows, bounds the conditional norm of the first sum by $C_q\|z^B\|_2$.
Its bad-exterior envelope is $m^2$, so \eqref{eq:local-good-exteriors} gives
\begin{equation}\label{eq:Yz-centered-weighted-term}
 \begin{aligned}
 &\E\left|\sum_{j\in V(B)}(X_{ij}-\E[X_{ij}\mid B,\mathbf e])z_j^B\right|^q\\
 &\quad\le C_q\E\|z^B\|_2^q+m^{2q}\Pp_n(B\notin\mathcal A_n)
 \le C_{q,h}n^q.
 \end{aligned}
\end{equation}
The Cauchy--Schwarz and H\"older inequalities, together with the selected-row bound
\eqref{eq:mean-profiles-a} under the same actual $(B,\mathbf e)$ marginal
and the degree bound \eqref{eq:degree-vector-moment}, give
\begin{equation}\label{eq:Yz-mean-weighted-term}
 \begin{aligned}
 &\left\|\sum_{j\in V(B)}(\E[X_{ij}\mid B,\mathbf e]-P_B)z_j^B\right\|_q\\
 &\le\bigl\|\|(\E[X_{ij}\mid B,\mathbf e]-P_B)_{j\in V(B)}\|_2\bigr\|_{2q}
       \bigl\|\|z^B\|_2\bigr\|_{2q}\le C_{q,h}n.
 \end{aligned}
\end{equation}
Combining \eqref{eq:Yz-full-decomposition}, \eqref{eq:Yz-elementary-terms},
\eqref{eq:Yz-centered-weighted-term}, and \eqref{eq:Yz-mean-weighted-term}
proves $\|(Y_Az^A)_i\|_q\le C_{q,h}n$ uniformly in $i\in V(A)$.
For $q_* =\max\{q,2\}$, Minkowski's inequality and moment monotonicity give,
after summing over the original vertex set $V(A)$,
\[
 \bigl\|\|Y_Az^A\|_2\bigr\|_q^2
 \le\sum_{i\in V(A)}\|(Y_Az^A)_i\|_{q_*}^2\le C_{q,h}n^3.
\]
The deletion $B=A-i$ is made separately for each $i$ in this final sum.
All estimates use the actual exterior marginal; $h=0$ includes the full graph.
\end{proof}

\begin{lemma}
\label{lem:row-center-moments}
Assume \hyperref[ass:PC]{\textup{(PC)}}.
For $\mu=\mu_{A,\mathbf e,s}$ in \eqref{v8:def:conditional-multirow-law},
take $A\in\mathcal A_n$, fixed $q_0,d$, and a selected row $X$.
Let $u$ and $c_{i_1\ldots i_d}$ be fixed under this law, with
$|c_{i_1\ldots i_d}|\le Kn^{-d}$.
The selected-row marginal need not have Hamiltonian $H_A(x)$.
For every fixed $q\ge1$,
\begin{equation}\label{eq:row-center-moments}
 \left\|\sum_{i_1,\ldots,i_d\ {\rm distinct}}
       c_{i_1\ldots i_d}(A,\mathbf e,s)
                      \prod_{a=1}^d(X_{i_a}-u_{i_a})\right\|_{L^q(\mu_{A,\mathbf e,s})}
 \le C_{d,q}K n^{-d/2}(1+\|\E_{\mu_{A,\mathbf e,s}} X-u\|_2)^d.
\end{equation}
This is a conditional bound for fixed $(A,\mathbf e,s)$;
subsequent applications average its right-hand side over $(A,\mathbf e)$.
The one-row reference law is included by taking $q_0=1$.
\end{lemma}

\begin{proof}
Write $X_i-u_i=\xi_i+\delta_i$, where
$\xi_i=X_i-\E_\mu X_i$ and $\delta_i=\E_\mu X_i-u_i$.
For $S\subseteq[d]$, $b=|S|$, sum over the indices outside $S$ first.
The resulting centered degree-$b$ polynomial has coefficients
at most $Kn^{-d}\|\delta\|_1^{d-b}$.
Equation~\eqref{ep:row:homogeneous-moment} gives
\begin{equation}\label{audit:dependency:s56:recentered-monomial-moment}
 \begin{aligned}
 &\left\|\sum_{i_1,\ldots,i_d\ {\rm distinct}}
 c_{i_1\ldots i_d}(A,\mathbf e,s)
 \prod_{a\in S}\xi_{i_a}\prod_{a\notin S}\delta_{i_a}\right\|_{L^q(\mu)}\\
 &\le C_{d,q}Kn^{-d+b/2}\|\delta\|_1^{d-b}
 \le C_{d,q}Kn^{-d/2}\|\delta\|_2^{d-b}.
 \end{aligned}
\end{equation}
For $b=0$, this is the deterministic coefficient bound.
Since all coefficients and centers are fixed under $\mu$, summing
\eqref{audit:dependency:s56:recentered-monomial-moment} over $S$ yields
\begin{align*}
 &\left\|\sum_{i_1,\ldots,i_d\ {\rm distinct}}
 c_{i_1\ldots i_d}(A,\mathbf e,s)
                   \prod_{a=1}^d(X_{i_a}-u_{i_a})\right\|_{L^q(\mu)}\\*
 &\quad\le C_{d,q}Kn^{-d/2}\sum_{b=0}^d\binom db\|\delta\|_2^{d-b}
 =C_{d,q}Kn^{-d/2}(1+\|\delta\|_2)^d.\qedhere
\end{align*}
\end{proof}

\subsection{Singly and doubly rooted count expansions}
For an ordered pair of distinct template roots $(u,v)$, the full count
includes the edge between the roots when that edge is present:
\begin{equation}\label{eq:full-two-root-count}
 \begin{aligned}
 T^{\mathrm{full}}_{ij}(H;u,v;X)
 &=n^{2-v(H)}
  \sum_{\substack{\varphi:V(H)\hookrightarrow[n]\\\varphi(u)=i,\ \varphi(v)=j}}
       \prod_{ab\in E(H)}X_{\varphi(a)\varphi(b)}\\
 &=X_{ij}^{a_H(u,v)}T_{H,(u,v);(i,j)}(X),\qquad
 a_H(u,v)=\one_{\{uv\in E(H)\}} .
 \end{aligned}
\end{equation}
The rooted count on the right omits the edge between the roots, in accordance
with \eqref{row:eq:rooted-count}; the prefactor restores that edge in the full count.
For an exterior $A$, restrict injections to $V(A)$ and retain the normalization
based on the ambient size $n$. For a fixed deletion $A=X_{[n]\setminus\mathcal D}$,
$|\mathcal D|=h$, write $P_A$ for its density and
$z_i^A=d_{A,i}-(n-h-1)P_A$. The singly rooted count
$T_i(H,o;A)$ is defined in \eqref{v8:def:rooted-vector}; put $d_o=d_H(o)$.

\begin{lemma}\label{lem:root-linear}
Assume \hyperref[ass:PC]{\textup{(PC)}}.
Consider the singly rooted count and fixed-deletion exterior specified above.
Then, for every fixed real $q\ge1$, the following expansion holds
uniformly over the retained root images:
\begin{equation}\label{eq:one-root}
 T_i(H,o;A)=P_A^{e(H)}
              +d_oP_A^{e(H)-1}\frac{z_i^A}{n}+O_{L^q}(n^{-1}).
\end{equation}
Omit the linear term when $d_o=0$. For fixed $H$, the constants are uniform
in the designated root $o$, while $d_o$ retains its dependence on the root.
The norms are taken under the induced marginal of the original graph, uniformly in
the deleted labels and root images; $h=0$ includes the full graph.
The coordinate estimate also gives
\begin{equation}\label{eq:one-root-vector-remainder}
 \left\|\left(\sum_{i\in V(A)}
  \left|T_i(H,o;A)-P_A^{e(H)}
       -d_oP_A^{e(H)-1}z_i^A/n\right|^2\right)^{1/2}\right\|_q
       \le C_{H,q,h}n^{-1/2},
\end{equation}
again omitting the linear term when $d_o=0$.
\end{lemma}

\begin{proof}
\proofstep{Conditional expansion.}
Fix $i\in V(A)$ and put
\[
 B=A-i,\qquad
 m=|V(B)|=n-h-1.
\]
Condition on $\mathcal F_B=\sigma(B,
X_{ab}:a<b,\ a,b\in\mathcal D\cup\{i\})$.
The full conditional law has $h+1$ rows. Write
$X_i=(X_{ij})_{j\in V(B)}$ and
$\mu=(\E[X_{ij}\mid\mathcal F_B])_{j\in V(B)}$.
Let $T_i(y;B)$ be the row polynomial, so that
$T_i(H,o;A)=T_i(X_i;B)$.
All sums over $j$ and coordinate sets $J$ in this proof range over $V(B)$.
Fixing $s$ neighbor images gives
$\sup_{y\in[0,1]^{V(B)}}|\partial_JT_i(y;B)|\le C_Hn^{-s}$ for $|J|=s$.
Equations~\eqref{audit:exterior-count-replacement} and
\eqref{eq:mean-profiles-a} imply
\begin{equation}\label{eq:l82-exterior-count}
 n^{-v(K)}\inj(K,B)=P_B^{e(K)}+O_{L^q}(n^{-1})
 \qquad(K\text{ fixed, including isolates}),
\end{equation}
\begin{equation}\label{eq:l82-mean-at-density}
 \bigl\|\|\mu-P_B\one\|_2\bigr\|_q\le C_q.
\end{equation}
If $d_o=0$, \eqref{eq:l82-exterior-count} and density replacement prove
the claim directly. For $d_o\ge1$, multiaffinity gives
\begin{equation}\label{audit:dependency:s56:rooted-row-polynomial-expansion}
 T_i(X_i;B)=T_i(P_B\one;B)+\sum_ja_j(X_{ij}-P_B)+\sum_{s=2}^{d_o}R_s,
\end{equation}
where $a_j=\partial_jT_i(P_B\one;B)$ and
$R_s=\sum_{|J|=s}\partial_JT_i(P_B\one;B)\prod_{j\in J}(X_{ij}-P_B)$.

\proofstep{Constant and linear coefficients.}
Retaining the isolated vertices of $H-o$, we have
\begin{equation}\label{audit:dependency:s56:rooted-constant-linear-coefficients}
 \begin{aligned}
 T_i(P_B\one;B)&=P_B^{d_o}n^{-v(H-o)}\inj(H-o,B)
                         =P_B^{e(H)}+O_{L^q}(n^{-1}),\\
 a_j&=\frac{P_B^{d_o-1}}n\sum_{\substack{u\in V(H)\\ou\in E(H)}}T_j(H-o,u;B),\\
 \bar a:=m^{-1}\sum_ja_j
 &=\frac{d_oP_B^{d_o-1}}m n^{-v(H-o)}\inj(H-o,B).
 \end{aligned}
\end{equation}
Equations~\eqref{eq:root-shared-bounds}, \eqref{eq:l82-exterior-count},
and \eqref{audit:dependency:s56:rooted-constant-linear-coefficients},
with $m=n-O(1)$, give
\[
 \bigl\|\|a-\bar a\one\|_2\bigr\|_q\le C_{H,q}/n,\qquad
 \left\|\bar a-\frac{d_oP_B^{e(H)-1}}n\right\|_q
 \le C_{H,q}(mn)^{-1}+C_H|m^{-1}-n^{-1}|\le C_{H,q}n^{-2}.
\]
Thus the coefficient error $\varepsilon_j=a_j-d_oP_B^{e(H)-1}/n$ satisfies
\begin{equation}\label{audit:root-linear-derivative-error}
 \bigl\|\|\varepsilon\|_2\bigr\|_q
 \le\bigl\|\|a-\bar a\one\|_2\bigr\|_q
    +\sqrt m\left\|\bar a-d_oP_B^{e(H)-1}/n\right\|_q
 \le C_{H,q}/n.
\end{equation}
Split $\sum_j\varepsilon_j(X_{ij}-P_B)$ at $\mu_j$.
The coefficients are $\mathcal F_B$-measurable, so
equations~\eqref{eq:weighted-local-moments} and
\eqref{audit:root-linear-derivative-error}, with zero weights on the other rows, give
\begin{equation}\label{audit:root-linear-centered-pairing}
 \left\|\sum_j\varepsilon_j(X_{ij}-\mu_j)\right\|_q
 \le C_q\bigl\|\|\varepsilon\|_2\bigr\|_q+Ce^{-cn/q}
 \le C_{H,q}n^{-1}.
\end{equation}
Here $|\varepsilon_j|\le C_H/n$ supplies a bounded bad-exterior envelope.
For the mean shift, H\"older's inequality, \eqref{audit:root-linear-derivative-error},
and \eqref{eq:l82-mean-at-density} give
\begin{equation}\label{audit:root-linear-mean-pairing}
 \left\|\sum_j\varepsilon_j(\mu_j-P_B)\right\|_q
 \le\bigl\|\|\varepsilon\|_2\bigr\|_{2q}
        \bigl\|\|\mu-P_B\one\|_2\bigr\|_{2q}
 \le C_{H,q}n^{-1}.
\end{equation}

\proofstep{Higher terms and the final center.}
Represent each set $J$ by its increasing tuple, with zero coefficients
elsewhere. Equations~\eqref{eq:row-center-moments},
\eqref{eq:l82-mean-at-density}, and \eqref{eq:local-good-exteriors},
with the coefficient bound $Cn^{-s}$ and the bounded bad-exterior envelope of $R_s$, yield
\begin{equation}\label{audit:root-linear-higher-remainder}
 \begin{aligned}
 \|R_s\|_q
 &\le C_{H,q,s}n^{-s/2}
      \bigl(1+\bigl\|\|\mu-P_B\one\|_2\bigr\|_{sq}\bigr)^s
          +Ce^{-cn/q}\le C_{H,q,s}n^{-s/2},\\
 \left\|\sum_{s=2}^{d_o}R_s\right\|_q&\le C_{H,q}n^{-1}.
 \end{aligned}
\end{equation}
The sum is empty if $d_o=1$. Combining
\eqref{audit:dependency:s56:rooted-row-polynomial-expansion},
\eqref{audit:dependency:s56:rooted-constant-linear-coefficients},
\eqref{audit:root-linear-centered-pairing},
\eqref{audit:root-linear-mean-pairing}, and
\eqref{audit:root-linear-higher-remainder} gives
\begin{equation}\label{audit:dependency:s56:rooted-expansion-before-recentering}
 T_i(H,o;A)=P_B^{e(H)}+
      \frac{d_oP_B^{e(H)-1}}n\sum_j(X_{ij}-P_B)+O_{L^q}(n^{-1}).
\end{equation}
Exact counting gives $\sum_j(X_{ij}-P_B)=z_i^A+m(P_A-P_B)$.
Equations~\eqref{eq:fixed-deletion-density} and
\eqref{eq:degree-vector-moment} give
$\|P_A-P_B\|_q\le C_qn^{-3/2}$ and $\|z_i^A\|_q\le C_q\sqrt n$, hence
\begin{equation}\label{audit:root-linear-density-replacement}
 \begin{aligned}
 \|P_B^{e(H)}-P_A^{e(H)}\|_q&\le C_{H,q}n^{-3/2},\\
 \left\|\frac{md_o}{n}P_B^{e(H)-1}(P_A-P_B)\right\|_q
 &\le C_{H,q}n^{-3/2},\\
 \frac{d_o}{n}\|(P_B^{e(H)-1}-P_A^{e(H)-1})z_i^A\|_q
 &\le\frac{C_H}{n}\|P_B-P_A\|_{2q}\|z_i^A\|_{2q}
 \le C_{H,q}n^{-2}.
 \end{aligned}
\end{equation}
The last term vanishes for $e(H)=1$.
Equations~\eqref{audit:dependency:s56:rooted-expansion-before-recentering}
and \eqref{audit:root-linear-density-replacement} prove \eqref{eq:one-root};
the case $d_o=0$ uses no negative powers.
For its remainder $R_i$ and $q_* =\max\{q,2\}$, we have
\[
 \left\|\left(\sum_i|R_i|^2\right)^{1/2}\right\|_q^2
 \le\sum_i\|R_i\|_{q_*}^2\le C_{H,q,h}n^{-1},
\]
which proves \eqref{eq:one-root-vector-remainder}.
\end{proof}

For the same fixed-deletion exterior $A$ and two ordered roots
$(u,v)$, use the full count \eqref{eq:full-two-root-count},
including every remaining isolated vertex.

\begin{lemma}\label{lem:two-root-counts}
Assume \hyperref[ass:PC]{\textup{(PC)}}.
Consider the two-root count and fixed-deletion exterior specified above.
Then, for every fixed real $q\ge1$,
\begin{equation}\label{eq:two-root}
 T^{\mathrm{full}}_{ij}(H;u,v;A)
 =A_{ij}^{a_H(u,v)}P_A^{e(H)-a_H(u,v)}+O_{L^q}(n^{-1/2}),
 \qquad i\ne j.
\end{equation}
The error also satisfies
\begin{equation}\label{eq:two-root-matrix-remainder}
 \left\|\left(\sum_{i\ne j}
 |T^{\mathrm{full}}_{ij}(H;u,v;A)-A_{ij}^{a_H(u,v)}P_A^{e(H)-a_H(u,v)}|^2
                              \right)^{1/2}\right\|_q
 \le C_{H,q,h}n^{1/2}.
\end{equation}
All norms are taken under the induced marginal of the original graph law.
The constants are uniform over the deleted labels and root images, with maxima over the
latter taken outside the norms. Taking $h=0$ includes the full graph.
\end{lemma}

\begin{proof}
Fix $(u,v)$ and distinct retained images $(i,j)$.
Put
\[
 B=A-\{i,j\},\qquad
 m=|V(B)|=n-h-2.
\]
Condition on $B$ and the internal edges $\mathbf e$ among
$\mathcal D\cup\{i,j\}$.
The full conditional distribution is $\mu=\mu_{B,\mathbf e,1}$ on $h+2$ rows.
Write $e=A_{ij}=X_{ij}$, $x_k=X_{ik}$ and $y_k=X_{jk}$ for $k\in V(B)$.
For $x,y\in[0,1]^{V(B)}$, define
\[
 \begin{aligned}
 T_{H,(u,v);(i,j)}(x,y;B)&=n^{2-v(H)}
 \sum_{\psi:V(H-\{u,v\})\hookrightarrow V(B)}
       \prod_{bc\in E(H-\{u,v\})}B_{\psi(b)\psi(c)}\\
 &\qquad\times\prod_{\substack{w\in V(H)\setminus\{u,v\}\\uw\in E(H)}}x_{\psi(w)}
                  \prod_{\substack{w\in V(H)\setminus\{u,v\}\\vw\in E(H)}}y_{\psi(w)}.
 \end{aligned}
\]
Then $T^{\mathrm{full}}_{ij}(H;u,v;A)=e^{a_H(u,v)} T_{H,(u,v);(i,j)}(x,y;B)$, with $e^0=1$.
The means $\mu^x=\E_\mu x$ and $\mu^y=\E_\mu y$ need not agree,
because the other internal edges may break row symmetry.
We use the exact decomposition
\begin{equation}\label{eq:root-two-error-decomposition}
 \begin{aligned}
 &T_{H,(u,v);(i,j)}(x,y;B)-P_B^{e(H)-a_H(u,v)}\\
 &=\{T_{H,(u,v);(i,j)}(P_B\one,P_B\one;B)-P_B^{e(H)-a_H(u,v)}\}\\
 &\quad+\{T_{H,(u,v);(i,j)}(\mu^x,\mu^y;B)
                     -T_{H,(u,v);(i,j)}(P_B\one,P_B\one;B)\}\\
 &\quad+\{T_{H,(u,v);(i,j)}(x,y;B)
                     -T_{H,(u,v);(i,j)}(\mu^x,\mu^y;B)\}.
 \end{aligned}
\end{equation}
We claim, under the actual marginal of $(B,\mathbf e)$, that
\begin{equation}\label{eq:root-two-scalar-error}
 \|T_{H,(u,v);(i,j)}(P_B\one,P_B\one;B)-P_B^{e(H)-a_H(u,v)}\|_q\le C_{H,q}n^{-1},
\end{equation}
\begin{equation}\label{eq:root-two-mean-error}
 \|T_{H,(u,v);(i,j)}(\mu^x,\mu^y;B)-T_{H,(u,v);(i,j)}(P_B\one,P_B\one;B)\|_q\le C_{H,q}n^{-1/2},
\end{equation}
\begin{equation}\label{eq:l82-two-random-error}
 \|T_{H,(u,v);(i,j)}(x,y;B)-T_{H,(u,v);(i,j)}(\mu^x,\mu^y;B)\|_q\le C_{H,q}n^{-1/2}.
\end{equation}

\proofstep{The three errors.}
Since $e(H-\{u,v\})=e(H)-d_H(u)-d_H(v)+a_H(u,v)$,
\eqref{eq:l82-exterior-count}, applied to the fixed-deletion exterior $B$, gives
\[
 \begin{aligned}
 T_{H,(u,v);(i,j)}(P_B\one,P_B\one;B)
 &=P_B^{d_H(u)+d_H(v)-2a_H(u,v)}
       n^{-(v(H)-2)}\inj(H-\{u,v\},B)\\
 &=P_B^{e(H)-a_H(u,v)}+O_{L^q}(n^{-1}).
 \end{aligned}
\]
Fixing one exterior image gives
$\sup_{x,y}(|\partial_{x_a}T_{H,(u,v);(i,j)}|+|\partial_{y_a}T_{H,(u,v);(i,j)}|)\le C_H/n$.
Integrating along the mean shifts and using \eqref{eq:mean-profiles-a} yields
\[
 \begin{aligned}
 &\|T_{H,(u,v);(i,j)}(\mu^x,\mu^y;B)-T_{H,(u,v);(i,j)}(P_B\one,P_B\one;B)\|_q\\
 &\quad\le\frac{C_H\sqrt m}{n}
  \left(\bigl\|\|\mu^x-P_B\one\|_2\bigr\|_q
       +\bigl\|\|\mu^y-P_B\one\|_2\bigr\|_q\right)
 \le C_{H,q}n^{-1/2}.
 \end{aligned}
\]
Fixing $s$ exterior sites gives coefficients $O(n^{-s})$ in
$x_a,y_a,x_ay_a$. Equations~\eqref{eq:local-polynomial-moments}
and \eqref{eq:row-polynomial-mean} therefore give, on good exteriors,
\[
 \begin{aligned}
 &\|T_{H,(u,v);(i,j)}(x,y;B)
              -T_{H,(u,v);(i,j)}(\mu^x,\mu^y;B)\|_{L^q(\mu)}\\
 &\quad\le\|T_{H,(u,v);(i,j)}(x,y;B)
              -\E_\mu T_{H,(u,v);(i,j)}(x,y;B)\|_{L^q(\mu)}\\
 &\qquad+|\E_\mu T_{H,(u,v);(i,j)}(x,y;B)
              -T_{H,(u,v);(i,j)}(\mu^x,\mu^y;B)|\\
 &\quad\le C_{H,q}n^{-1/2}+C_H/n.
 \end{aligned}
\]
The second term includes the $x_ay_a$ covariance correction.
The uniform bound of two and \eqref{eq:local-good-exteriors}
give \eqref{eq:l82-two-random-error} after averaging.
If there are no row factors, the last two errors vanish; if $v(H)=2$,
$T_{H,(u,v);(i,j)}(x,y;B)=1$ and all three vanish. All density exponents are nonnegative.

\proofstep{Restore the center and sum.}
Equations~\eqref{eq:root-two-error-decomposition}--\eqref{eq:l82-two-random-error}
give $T_{H,(u,v);(i,j)}(x,y;B)=P_B^{e(H)-a_H(u,v)}+O_{L^q}(n^{-1/2})$.
Replacing $P_B^{e(H)-a_H(u,v)}$ by $P_A^{e(H)-a_H(u,v)}$ costs $O_{L^q}(n^{-3/2})$
by \eqref{eq:fixed-deletion-density}; the cost is zero if $e(H)-a_H(u,v)=0$.
Multiplication by $e^{a_H(u,v)}\le1$ proves \eqref{eq:two-root}.
For its remainder $R_{ij}$ and $q_* =\max\{q,2\}$, we have
\[
 \left\|\left(\sum_{i\ne j}|R_{ij}|^2\right)^{1/2}\right\|_q
 \le\left(\sum_{i\ne j}\|R_{ij}\|_{q_*}^2\right)^{1/2}
 \le C_{H,q,h}n^{1/2},
\]
which proves \eqref{eq:two-root-matrix-remainder}.
\end{proof}

\section{Reference-row mean and covariance expansions}\label{sec:reference-row-profiles}
Rooted expansions yield the reference Hessian, mean and covariance profiles.

Using the fixed $\delta_0$ chosen in Appendix~\ref{row:sec:extension}, take a
smooth function $\chi_0:[0,1]\to[0,1]$ equal to one on
$|x-p|\le\delta_0$ and supported in $|x-p|<2\delta_0$. Define
\begin{equation}\label{eq:logit-profiles}
 v_0(x)=
 \begin{cases}
 \nu(x)\left(1-\chi_0(x)+
           \dfrac{\chi_0(x)}{1-2\nu(x)t(x)}\right),&x\in J,\\
 \nu(x),&x\notin J
 \end{cases}.
\end{equation}
The cutoff makes $v_0$ smooth and bounded on $[0,1]$.
This profile describes the reference covariance with the row total.
We use $v_0(x)=\nu(x)/(1-2\nu(x)t(x))$ only when
$|x-p|\le\delta_0$.

\subsection{Reference Hessians and covariance equations}
Unless stated otherwise, row moments are taken under $\pi_A$:
\begin{equation}\label{v8:def:reference-row-data}
 \mu_{A,i}=\E_{\pi_A}X_i,\quad \mu_A=(\mu_{A,i})_i,\qquad
 \Sigma_{A,ij}=\Cov_{\pi_A}(X_i,X_j),\quad\Sigma_A=(\Sigma_{A,ij})_{i,j}.
\end{equation}
We use the exterior density and degree vector from
\eqref{v8:nearby:exterior-degrees} and the centered matrix from
\eqref{v8:def:count-centered-adjacency}.
We work with the good exterior sets $\mathcal A_n$ of
Lemma~\ref{row:lem:good-exteriors}, whose construction includes
$|P_A-p|\le\delta$ for the fixed neighborhood selected in
Appendix~\ref{row:sec:extension}.
Lemmas~\ref{lem:shape} and~\ref{ep:row:mean-jacobian} then give
\begin{equation}\label{eq:row-profile-good-exterior-event}
 A\in\mathcal A_n\Longrightarrow P_A\one,\mu_A\in[p-\delta,p+\delta]^m,
 \qquad \Pp_n(A\notin\mathcal A_n)\le Ce^{-cn}.
\end{equation}
Intermediate logit and inverse estimates are restricted to $\mathcal A_n$.
The fixed global profile extensions and polynomial envelopes ensure that its
complement contributes $O(n^{-M})$ to every fixed moment, for every fixed $M$.
Combining \eqref{v11:eq:tolerance-contraction} and
\eqref{row:eq:good-properties} with the definition
\eqref{eq:logit-profiles} and $0\le\nu(x)\le1/4$ gives
\begin{equation}\label{eq:row-profile-row-sum-gap}
 2t(P_A)\nu(P_A)\le3/4,\qquad
 1-2t(P_A)\nu(P_A)\ge1/4,\qquad |v_0(P_A)|\le1
 \quad(A\in\mathcal A_n).
\end{equation}
Applying Theorem~\ref{ep:row:sharp}, \eqref{ep:row:feature-bound}, to the
coordinate family at order two, and then using symmetry, gives
\begin{equation}\label{eq:reference-covariance-operator}
 \|\Sigma_A\|_{\infty\to\infty}
   =\max_i\sum_j|\Sigma_{A,ij}|\le C,\qquad \|\Sigma_A\|_{\op}\le C .
\end{equation}
For the covariance expansion, we also need a variance estimate that
preserves cancellation when summed over the coordinates:
\begin{equation}\label{eq:row-profile-variance-input}
 \bigl\|\|(\nu(\mu_{A,i})-\nu(P_A))_i\|_2\bigr\|_q
 +\left\|\sum_i(\nu(\mu_{A,i})-\nu(P_A))\right\|_q\le C_q.
\end{equation}
This follows from Lemma~\ref{lem:row-mean-density},
\eqref{eq:mean-profiles-a}, and the exact identities
\[
 \begin{aligned}
 \nu(\mu_{A,i})-\nu(P_A)
   &=(1-2P_A)(\mu_{A,i}-P_A)-(\mu_{A,i}-P_A)^2,\\
 \sum_i\{\nu(\mu_{A,i})-\nu(P_A)\}
   &=(1-2P_A)\one^{\mathsf T}(\mu_A-P_A\one)
                                  -\|\mu_A-P_A\one\|_2^2.
 \end{aligned}
\]
Define the reference Hessian by
\begin{equation}\label{v8:def:reference-row-hessian}
 [H_A''(x)]_{ij}=\partial_j\partial_iH_A(x)=\partial_jl_{i,A}(x).
\end{equation}
We define the Hessian remainders by the exact differences
$R_{A,\rm const}=H_A''(P_A\one)-2t(P_A)J_m/n-6\tau(P_A)Y_A/n$ and
$R_{A,\rm mean}=H_A''(\mu_A)-2t(P_A)J_m/n-6\tau(P_A)Y_A/n$.

\begin{lemma}\label{lem:row-hessian-profile}
Assume \hyperref[ass:PC]{\textup{(PC)}}.
With the data \eqref{v8:def:reference-row-data} and the remainders defined
above, for every fixed $q\ge1$,
\[
 \begin{aligned}
 H_A''(P_A\one)&=2t(P_A)J_m/n+6\tau(P_A)Y_A/n+R_{A,\mathrm{const}},\\
 H_A''(\mu_A)&=2t(P_A)J_m/n+6\tau(P_A)Y_A/n+R_{A,\mathrm{mean}},\\
 \bigl\|\|R_{A,\mathrm{const}}\|_F\bigr\|_q
   +\bigl\|\|R_{A,\mathrm{mean}}\|_F\bigr\|_q&\le C_qn^{-1/2}.
 \end{aligned}
\]
The mean-row Hessian also has the signed row-sum bound
\begin{equation}\label{eq:row-profile-hessian-row-sums}
 \begin{aligned}
 &\bigl\|\|H_A''(\mu_A)\one-2t(P_A)\one\|_2\bigr\|_q\\
 &\quad+\|\one^{\mathsf T}(H_A''(\mu_A)\one-2t(P_A)\one)\|_q
      \le C_q.
 \end{aligned}
\end{equation}
\end{lemma}

\begin{proof}
Work on $\mathcal A_n$; all exterior estimates below hold at every fixed
moment order. For $i\ne j$, differentiation and Lemma~\ref{lem:two-root-counts} give
\[
\begin{aligned}
[H_A''(u\one)]_{ij}
 &=\frac1n\sum_r b_{r,n}\sum_{\substack{o\in V(F_r)\\d(o)\ge2}}
   \sum_{\substack{v,w\in V(F_r)\\ov,ow\in E(F_r),\ v\ne w}}
       u^{d(o)-2}T^{\mathrm{full}}_{ij}(F_r-o;v,w;A),\\
T^{\mathrm{full}}_{ij}(F_r-o;v,w;A)
 &=A_{ij}^{a}P_A^{e_r-d(o)-a}+O_{L^q}(n^{-1/2}),
 \qquad a=\one_{\{vw\in E(F_r)\}}.
\end{aligned}
\]
The rooted count retains isolates and has normalization $n^{3-v_r}$.
Only nonempty summands are evaluated. Since
\[
 \sum_{o\in V(F_r)}d(o)(d(o)-1)=2a(F_r),\qquad
 \sum_{o\in V(F_r)}\sum_{\substack{v,w\in V(F_r)\\ov,ow\in E(F_r),\ v\ne w}}
       \one_{\{vw\in E(F_r)\}}=6\triangle(F_r),
\]
the definitions of $t,\tau$ and $b_{r,n}-b_r=O(n^{-1})$ yield
\[
 [H_A''(P_A\one)]_{ij}
 =\frac{2t(P_A)+6\tau(P_A)(A_{ij}-P_A)}n
       +O_{L^q}(n^{-3/2}).
\]
The Hessian diagonal is zero. Minkowski's inequality at order $q\vee2$ gives
\begin{equation}\label{eq:row-profile-constant-hessian}
 H_A''(P_A\one)=2t(P_A)J_m/n+6\tau(P_A)Y_A/n+R_{A,\rm const},
 \qquad\bigl\|\|R_{A,\rm const}\|_F\bigr\|_q\le C_qn^{-1/2},
\end{equation}
because the off-diagonal and diagonal errors contribute
$O(mn^{-3/2})$ and $O(\sqrt m/n)$.

Equations~\eqref{row:eq:dense-derivatives} and
\eqref{eq:mean-profiles-a} imply
\[
\begin{aligned}
 [H_A''(\mu_A)-H_A''(P_A\one)]_{ij}
 &=\int_0^1\sum_k\partial_{ijk}H_A(P_A\one+s(\mu_A-P_A\one))(\mu_{A,k}-P_A)\,\dd s,\\
 \|H_A''(\mu_A)-H_A''(P_A\one)\|_F
 &\le Cmn^{-2}\|\mu_A-P_A\one\|_1\le Cn^{-1/2}\|\mu_A-P_A\one\|_2.
\end{aligned}
\]
Combining this estimate with \eqref{eq:row-profile-constant-hessian}
and the Hessian bound from \eqref{row:eq:dense-derivatives} proves
\begin{equation}\label{eq:K-profile}
\begin{gathered}
 H_A''(\mu_A)=2t(P_A)J_m/n+6\tau(P_A)Y_A/n+R_{A,\rm mean},\\
 \bigl\|\|R_{A,\rm mean}\|_F\bigr\|_q\le C_qn^{-1/2},
 \qquad \|H_A''(\mu_A)\|_{\op}\le C.
\end{gathered}
\end{equation}
Corollary~\ref{lem:degree-deletion} gives
\begin{equation}\label{eq:row-profile-exterior-degree-vector}
 \bigl\|\|z^A\|_2\bigr\|_q\le C_qn.
\end{equation}
Since $Y_A\one=z^A$ and $m=n-O(1)$,
\[
 H_A''(\mu_A)\one-2t(P_A)\one
 =2t(P_A)(m/n-1)\one+6\tau(P_A)z^A/n+R_{A,\rm mean}\one.
\]
Equations~\eqref{eq:K-profile} and
\eqref{eq:row-profile-exterior-degree-vector} therefore give the Euclidean
bound in \eqref{eq:row-profile-hessian-row-sums}.

For the signed total, sum before estimating.
Equations~\eqref{eq:one-root} and \eqref{eq:parameter-neighborhood} give
\begin{equation}\label{audit:dependency:s79:hessian-signed-average}
\begin{aligned}
 \frac1m\one^{\mathsf T}H_A''(P_A\one)\one
 &=\frac nm\sum_r b_{r,n}\sum_{o\in V(F_r)}d(o)(d(o)-1)
 P_A^{d(o)-2}n^{-(v_r-1)}\inj(F_r-o,A)\\
 &=2t(P_A)+O_{L^q}(n^{-1}).
\end{aligned}
\end{equation}
Writing $c_j=\sum_{i,l}\partial_{ilj}H_A(P_A\one)$, Taylor's formula gives
\begin{equation}\label{audit:dependency:s79:hessian-signed-taylor}
 \begin{aligned}
 &\left|\one^{\mathsf T}\{H_A''(\mu_A)-H_A''(P_A\one)\}\one
                                -c^{\mathsf T}(\mu_A-P_A\one)\right|\\
 &\quad\le Cn^{-1}\|\mu_A-P_A\one\|_1^2\le C\|\mu_A-P_A\one\|_2^2.
 \end{aligned}
\end{equation}
Here $|\sum_{i,l}\partial_{iljk}H_A(x)|\le C/n$, and we have the exact identity
\[
 c_j=\sum_r b_{r,n}
 \sum_{\substack{(o,v):ov\in E(F_r)\\d(o)\ge3}}
 (d(o)-1)(d(o)-2)P_A^{d(o)-3}T_j(F_r-o,v;A).
\]
For $\bar c=m^{-1}\sum_jc_j$, \eqref{eq:root-shared-bounds} gives
$|\bar c|\le C$ and $\|\|c-\bar c\one\|_2\|_q\le C_q$. Thus
\begin{equation}\label{audit:dependency:s79:hessian-linear-pairing}
 \|c^{\mathsf T}(\mu_A-P_A\one)\|_q
 \le C\|\one^{\mathsf T}(\mu_A-P_A\one)\|_q
 +\bigl\|\|c-\bar c\one\|_2\bigr\|_{2q}
  \bigl\|\|\mu_A-P_A\one\|_2\bigr\|_{2q}\le C_q.
\end{equation}
Combining
\eqref{audit:dependency:s79:hessian-signed-average},
\eqref{audit:dependency:s79:hessian-signed-taylor}, and
\eqref{audit:dependency:s79:hessian-linear-pairing} proves the signed bound
in \eqref{eq:row-profile-hessian-row-sums}.
Polynomial envelopes and $\Pp_n(\mathcal A_n^c)\le Ce^{-cn}$ remove
the good-event restriction.
\end{proof}

We next relate the reference covariance to the Hessian. An absolute
row-sum bound on the residual controls the error after inversion.

\begin{lemma}\label{lem:row-covariance-equation}
Assume \hyperref[ass:PC]{\textup{(PC)}}.
For every $A\in\mathcal A_n$ as in
\eqref{eq:row-profile-good-exterior-event}, the following residual bound holds:
\begin{equation}\label{eq:row-profile-covariance-residual}
 \begin{aligned}
 \mathcal R_{\Sigma,A}
 &=\left[I_m-\operatorname{diag}((\nu(\mu_{A,a}))_a)H_A''(\mu_A)\right]\Sigma_A
       -\operatorname{diag}((\nu(\mu_{A,a}))_a),\\
 \|\mathcal R_{\Sigma,A}\|_{\infty\to\infty}&\le C/n.
 \end{aligned}
\end{equation}
The constant is uniform over good exteriors; $\Sigma_A$ is the covariance
matrix under $\pi_A$, as defined in \eqref{v8:def:reference-row-data}.
\end{lemma}

\begin{proof}
Fix $A\in\mathcal A_n$ and let $X\sim\pi_A$, with $\pi_A$ defined in
\eqref{v8:def:reference-row-law}. The exterior $A$ is held fixed:
all expectations, covariances and variances below are under $\pi_A$,
with no outer averaging under $\mathbb P_n$. Put $\xi=X-\mu_A$.
Equation~\eqref{row:eq:single-site-conditional} and conditioning on the
other row coordinates give
\[
 \begin{gathered}
 \pi_A(X_i=1\mid X_{-i})
   =\E_{\pi_A}[X_i\mid X_{-i}]=\sig(l_{i,A}(X)),\\
 \E_{\pi_A}\sig(l_{i,A}(X))=\mu_{A,i},\qquad
 \Cov_{\pi_A}(\sig(l_{i,A}(X)),X_i)
   =\Var_{\pi_A}(\sig(l_{i,A}(X))),
 \end{gathered}
\]
and hence
\begin{equation}\label{audit:reference-covariance-conditioning}
 \Sigma_{A,ij}=\Cov_{\pi_A}(\sig(l_{i,A}(X)),X_j)
 +\one_{\{i=j\}}\{\nu(\mu_{A,i})-\Var_{\pi_A}(\sig(l_{i,A}(X)))\}.
\end{equation}
Equation~\eqref{row:eq:dense-derivatives} gives the finite expansion
\[
 \begin{aligned}
 l_{i,A}(X)-l_{i,A}(\mu_A)&=(H_A''(\mu_A)\xi)_i+q_i,\\
 q_i&=\sum_{d=2}^D\sum_{\mathbf k}c^{(i,d)}_{\mathbf k}
                         \prod_{a=1}^d\xi_{k_a},
 \qquad |c^{(i,d)}_{\mathbf k}|\le C_dn^{-d}.
 \end{aligned}
\]
Taylor's formula and \eqref{ep:row:field-moments} imply
\begin{equation}\label{audit:reference-covariance-logistic-expansion}
 \begin{aligned}
 \sig(l_{i,A}(X))
 &=\sig(l_{i,A}(\mu_A))
   +\sig'(l_{i,A}(\mu_A))(H_A''(\mu_A)\xi)_i
   +\mathcal Q_i+\mathcal T_i,\\
 \E_{\pi_A}|\mathcal T_i|
 &\le C\E_{\pi_A}|l_{i,A}(X)-l_{i,A}(\mu_A)|^4\le Cn^{-2},
 \end{aligned}
\end{equation}
where
\[
 \begin{aligned}
 \mathcal Q_i
 &=\sig'(l_{i,A}(\mu_A))q_i
   +\tfrac12\sig''(l_{i,A}(\mu_A))
                         \{l_{i,A}(X)-l_{i,A}(\mu_A)\}^2\\
 &\quad+\tfrac16\sig^{(3)}(l_{i,A}(\mu_A))
                         \{l_{i,A}(X)-l_{i,A}(\mu_A)\}^3,\\
 \mathcal T_i
 &=\frac{\{l_{i,A}(X)-l_{i,A}(\mu_A)\}^4}{6}
   \int_0^1(1-s)^3
     \sig^{(4)}\!\left(l_{i,A}(\mu_A)
                   +s\{l_{i,A}(X)-l_{i,A}(\mu_A)\}\right)\,\dd s.
 \end{aligned}
\]
At formal degree $d\in[2,3D]$, the coefficients of $\mathcal Q_i$
are $O(n^{-d})$, including terms with repeated centered positions.
The coefficient bounds and \eqref{eq:micro-moment-sum} give
\begin{equation}\label{audit:dependency:s79:covariance-nonlinear-error}
\begin{aligned}
 \sum_j|\Cov_{\pi_A}(\mathcal Q_i,X_j)|
 &\le\sum_{d=2}^{3D}C_dn^{-d}
       \sum_{j,k_1,\ldots,k_d}\left|\E_{\pi_A}\xi_j\prod_{a=1}^d\xi_{k_a}\right|\\
 &\le\sum_{d=2}^{3D}C_dn^{-d+\lfloor(d+1)/2\rfloor}\le C/n,\\
 \sum_j|\Cov_{\pi_A}(\mathcal T_i,X_j)|&\le m\E_{\pi_A}|\mathcal T_i|\le C/n.
\end{aligned}
\end{equation}
Equations~\eqref{ep:row:marginal-mean-equation} and
\eqref{ep:row:field-moments} also give
\begin{equation}\label{audit:dependency:s79:covariance-local-corrections}
 \begin{aligned}
 |\sig'(l_{i,A}(\mu_A))-\nu(\mu_{A,i})|
 &\le|\sig(l_{i,A}(\mu_A))-\mu_{A,i}|\le C/n,\\
 \Var_{\pi_A}(\sig(l_{i,A}(X)))
 &\le\tfrac1{16}\E_{\pi_A}|l_{i,A}(X)-l_{i,A}(\mu_A)|^2\le C/n.
 \end{aligned}
\end{equation}
Substituting \eqref{audit:reference-covariance-logistic-expansion} into
\eqref{audit:reference-covariance-conditioning} identifies the residual:
\[
 \begin{aligned}
 \mathcal R_{\Sigma,A,ij}
 &=(\sig'(l_{i,A}(\mu_A))-\nu(\mu_{A,i}))
                                      (H_A''(\mu_A)\Sigma_A)_{ij}\\
 &\quad+\Cov_{\pi_A}(\mathcal Q_i+\mathcal T_i,X_j)
                          -\one_{\{i=j\}}\Var_{\pi_A}(\sig(l_{i,A}(X))).
 \end{aligned}
\]
By \eqref{row:eq:dense-derivatives} and
\eqref{eq:reference-covariance-operator},
$\|H_A''(\mu_A)\|_{\infty\to\infty}
 +\|\Sigma_A\|_{\infty\to\infty}\le C$; hence
\eqref{audit:dependency:s79:covariance-nonlinear-error} and
\eqref{audit:dependency:s79:covariance-local-corrections} give
\[
 \begin{aligned}
 \max_i\sum_j|\mathcal R_{\Sigma,A,ij}|
 &\le C/n+\max_i\sum_j|\Cov_{\pi_A}(\mathcal Q_i+\mathcal T_i,X_j)|\\
 &\quad+\max_i\Var_{\pi_A}(\sig(l_{i,A}(X)))\le C/n.
 \end{aligned}
\]
All constants are uniform over $A\in\mathcal A_n$.
\end{proof}

Lemma~\ref{ep:row:mean-jacobian}, \eqref{eq:shape-common-inverse}, with
$y=z=\mu_A$, gives the common inverse bound
\begin{equation}\label{eq:covariance-inverse-norms}
 \begin{aligned}
 &\|\left[I_m-\operatorname{diag}((\nu(\mu_{A,a}))_a)H_A''(\mu_A)\right]^{-1}
                                             \|_{\infty\to\infty}\\
 &\quad+\|\left[I_m-\operatorname{diag}((\nu(\mu_{A,a}))_a)H_A''(\mu_A)\right]^{-1}
                                             \|_{\op}\le2/(1-\rho).
 \end{aligned}
\end{equation}
In using this inverse, we retain the displayed order of the matrix
factors; no commutation is needed.

\subsection{First-order reference-row mean expansion}
The preceding Hessian estimate and the centered-path bound identify the mean correction.

\begin{lemma}
\label{lem:row-profiles}
Assume \hyperref[ass:PC]{\textup{(PC)}}.\label{lem:row-mean-profile}
With $\mu_A$ as in \eqref{v8:def:reference-row-data}, for every fixed
$q\ge1$,
\begin{equation}\label{eq:mean-profiles-b}
 \left\|\left\|\mu_A-P_A\one
       -\frac{2\nu(P_A)t(P_A)}n z^A\right\|_2\right\|_q
       \le C_qn^{-1/2}.
\end{equation}
This refined expansion is asserted only for the reference law $\pi_A$.
\end{lemma}

\begin{proof}
Work on $\mathcal A_n$ under $\pi_A$.
The exact rooted expansion is
\[
 l_{i,A}(u\one)=h_n+\sum_r b_{r,n}
       \sum_{(o,v):ov\in E(F_r)}u^{d(o)-1}T_i(F_r-o,v;A).
\]
Equation~\eqref{eq:one-root} and the ordered-edge identities give
\[
\begin{aligned}
 T_i(F_r-o,v;A)
 &=P_A^{e_r-d(o)}
 +(d(v)-1)P_A^{e_r-d(o)-1}z_i^A/n+O_{L^q}(n^{-1}),\\
 \sum_{(o,v)}1&=2e_r,\qquad
 \sum_{(o,v)}(d(v)-1)=2a(F_r),\\
 l_{i,A}(P_A\one)&=\ell_n(P_A)+2t(P_A)z_i^A/n+O_{L^q}(n^{-1}).
\end{aligned}
\]
Terms with zero coefficient are omitted; replacing $b_{r,n}$ by
$b_r$ in the degree term costs $O(n^{-2}\|z_i^A\|_q)=O(n^{-3/2})$.
By \eqref{audit:scalar-defect-exterior}, the inverse logistic map
$I'(x)=\log(x/(1-x))$ is Lipschitz on the fixed interior interval.
Therefore, with
\[
 E_{A,i}=l_{i,A}(P_A\one)-I'(P_A)-2t(P_A)z_i^A/n,
\]
we have $\|E_{A,i}\|_q\le C_q/n$.

Taylor's formula, the dense derivative bounds and
\eqref{eq:mean-profiles-a} give
\[
\begin{aligned}
 l_A(\mu_A)-l_A(P_A\one)&=H_A''(P_A\one)(\mu_A-P_A\one)+r_A^{\rm field},\\
 r_{A,i}^{\rm field}
 &=\int_0^1(1-s)\sum_{j,k}\partial_{jk}l_{i,A}(P_A\one+s(\mu_A-P_A\one))
                                   (\mu_{A,j}-P_A)(\mu_{A,k}-P_A)\,\dd s,\\
 |r_{A,i}^{\rm field}|&\le Cn^{-2}\|\mu_A-P_A\one\|_1^2
                         \le Cn^{-1}\|\mu_A-P_A\one\|_2^2,\\
 &\|\{l_{i,A}(\mu_A)-I'(P_A)\}\|_q
    \le Cn^{-1}\|z_i^A\|_q+\|E_{A,i}\|_q\\
 &\qquad+Cn^{-1}\bigl\|\|\mu_A-P_A\one\|_1\bigr\|_q
    \le C_qn^{-1/2}.
\end{aligned}
\]
The logistic expansion and \eqref{ep:row:marginal-mean-equation} give
\[
\begin{aligned}
 \sig(l_{i,A}(\mu_A))-P_A&=\nu(P_A)\{l_{i,A}(\mu_A)-I'(P_A)\}+r_{A,i}^{\rm logistic},\\
 r_{A,i}^{\rm logistic}
 &=\{l_{i,A}(\mu_A)-I'(P_A)\}^2
        \int_0^1(1-s)\sig''(I'(P_A)+s\{l_{i,A}(\mu_A)-I'(P_A)\})\,\dd s,\\
 \|r_{A,i}^{\rm logistic}\|_q&\le C_q/n,\\
 |\mu_{A,i}-\sig(l_{i,A}(\mu_A))|&\le C/n.
\end{aligned}
\]
The preceding coordinatewise bounds and Minkowski's inequality give
\begin{equation}\label{audit:dependency:s79:mean-profile-error-budget}
 \begin{aligned}
 &\bigl\|\|E_A\|_2\bigr\|_q
   +\bigl\|\|r_A^{\rm field}\|_2\bigr\|_q
   +\bigl\|\|r_A^{\rm logistic}\|_2\bigr\|_q\\
 &\qquad+\bigl\|\|\mu_A-\bigl(\sig(l_{i,A}(\mu_A))\bigr)_{i\in V(A)}\|_2\bigr\|_q
      \le C_qn^{-1/2}.
 \end{aligned}
\end{equation}
Adding the mean-equation error
$\mu_A-\bigl(\sig(l_{i,A}(\mu_A))\bigr)_{i\in V(A)}$
to the combined Taylor expansions gives
\[
\begin{aligned}
 \mu_A-P_A\one
 & =\nu(P_A)H_A''(P_A\one)(\mu_A-P_A\one)
      +\frac{2\nu(P_A)t(P_A)}n z^A\\
 &\quad+\nu(P_A)E_A+\nu(P_A)r_A^{\rm field}+r_A^{\rm logistic}\\
 &\quad+\mu_A-\bigl(\sig(l_{i,A}(\mu_A))\bigr)_{i\in V(A)}.
\end{aligned}
\]
Rearranging this identity around the first-order degree correction
$2\nu(P_A)t(P_A)z^A/n$ yields
\begin{equation}\label{audit:reference-mean-linear-equation}
 \begin{aligned}
 &\bigl[I_m-\nu(P_A)H_A''(P_A\one)\bigr]
       \left[\mu_A-P_A\one-\frac{2\nu(P_A)t(P_A)}n z^A\right]\\
 &\quad=\frac{2\nu(P_A)^2t(P_A)}n H_A''(P_A\one)z^A+\nu(P_A)E_A\\
 &\qquad+\nu(P_A)r_A^{\rm field}+r_A^{\rm logistic}
       +\mu_A-\bigl(\sig(l_{i,A}(\mu_A))\bigr)_{i\in V(A)}.
 \end{aligned}
\end{equation}
Since $J_mz^A=0$, \eqref{eq:row-profile-constant-hessian} and
\eqref{eq:Yz}, together with \eqref{eq:degree-vector-moment}, yield
\begin{equation}\label{audit:dependency:s79:mean-profile-forcing}
\begin{aligned}
 \frac{2\nu(P_A)t(P_A)}n H_A''(P_A\one)z^A
 &=2\nu(P_A)t(P_A)
       \{6\tau(P_A)Y_Az^A/n^2+R_{A,\rm const}z^A/n\},\\
 \left\|\left\|\frac{2\nu(P_A)t(P_A)}n
                             H_A''(P_A\one)z^A\right\|_2\right\|_q
 &\le Cn^{-2}\bigl\|\|Y_Az^A\|_2\bigr\|_q\\
 &\quad+Cn^{-1}\bigl\|\|R_{A,\rm const}\|_F\bigr\|_{2q}
             \bigl\|\|z^A\|_2\bigr\|_{2q}
 \le C_qn^{-1/2}.
\end{aligned}
\end{equation}
Apply \eqref{eq:shape-common-inverse}, with $y=z=P_A\one$, to
\eqref{audit:reference-mean-linear-equation}. Equations
\eqref{audit:dependency:s79:mean-profile-error-budget} and
\eqref{audit:dependency:s79:mean-profile-forcing} then prove
\eqref{eq:mean-profiles-b} on $\mathcal A_n$.
Polynomial envelopes remove this restriction.
\end{proof}

\subsection{Reference covariance profiles}
The coarse mean profile and the off-diagonal covariance bound already
suffice to control the following quadratic covariance sums.

\begin{lemma}\label{lem:row-covariance-squares}
Assume \hyperref[ass:PC]{\textup{(PC)}}.
With $\Sigma_A$ as in \eqref{v8:def:reference-row-data}, for every fixed
$q\ge1$,
\begin{equation}\label{eq:cov-profiles-b}
 \left\|\sum_{i,j}\Sigma_{A,ij}^2-m\nu(P_A)^2\right\|_q
       =\|\tr\Sigma_A^2-m\nu(P_A)^2\|_q\le C_q,
\end{equation}
and the vector of row sums of these squares satisfies
\begin{equation}\label{eq:covariance-square-row-sums}
 \left\|\left\|\left(\sum_j\Sigma_{A,ij}^2-\nu(P_A)^2\right)_i
                      \right\|_2\right\|_q\le C_q.
\end{equation}
\end{lemma}

\begin{proof}
Work on $\mathcal A_n$.
Equations~\eqref{eq:disjoint-local-covariance} and
\eqref{eq:row-profile-variance-input} give
\[
 \begin{gathered}
 |\Sigma_{A,ij}|\le C/n\quad(i\ne j),\\
 \bigl\|\|(\nu(\mu_{A,i})-\nu(P_A))_i\|_2\bigr\|_q
 +\left\|\sum_i\{\nu(\mu_{A,i})-\nu(P_A)\}\right\|_q\le C_q.
 \end{gathered}
\]
Since $\Sigma_{A,ii}=\nu(\mu_{A,i})$,
\[
 \begin{aligned}
 \tr\Sigma_A^2-m\nu(P_A)^2
 &=2\nu(P_A)\sum_i\{\nu(\mu_{A,i})-\nu(P_A)\}\\
 &\quad+\sum_i\{\nu(\mu_{A,i})-\nu(P_A)\}^2
             +\sum_{i\ne j}\Sigma_{A,ij}^2
 =O_{L^q}(1).
 \end{aligned}
\]
Also $\sum_{j\ne i}\Sigma_{A,ij}^2\le C/n$ and
$|\nu(\mu_{A,i})^2-\nu(P_A)^2|\le\tfrac12|\mu_{A,i}-P_A|$, whence
\[
 \left\|\left(\sum_j\Sigma_{A,ij}^2-\nu(P_A)^2\right)_i\right\|_2
 \le\tfrac12\|\mu_A-P_A\one\|_2+C\sqrt m/n=O_{L^q}(1).
\]
Equations~\eqref{eq:mean-profiles-a} and
\eqref{eq:row-profile-good-exterior-event} complete the proof.
\end{proof}

Here $(\Sigma_A\one)_i=\Cov_{\pi_A}(X_i,\sum_jX_j)$ and
$\one^{\mathsf T}\Sigma_A\one=\Var_{\pi_A}(\sum_jX_j)$.

\begin{lemma}\label{lem:row-total-covariance}
Assume \hyperref[ass:PC]{\textup{(PC)}}.
With $\Sigma_A$ as in \eqref{v8:def:reference-row-data}, for every fixed
$q\ge1$,
\begin{equation}\label{eq:cov-profiles-a}
 \begin{aligned}
 \bigl\|\|\Sigma_A\one-v_0(P_A)\one\|_2\bigr\|_q&\le C_q,\\
 \|\one^{\mathsf T}\Sigma_A\one-mv_0(P_A)\|_q&\le C_q.
 \end{aligned}
\end{equation}
\end{lemma}

\begin{proof}
Work on $\mathcal A_n$. Equations
\eqref{eq:row-profile-covariance-residual},
\eqref{eq:row-profile-variance-input}, and
\eqref{eq:row-profile-hessian-row-sums} give
\begin{equation}\label{audit:reference-total-covariance-inputs}
 \begin{aligned}
 &\left[I_m-\operatorname{diag}((\nu(\mu_{A,i}))_i)H_A''(\mu_A)\right]
                                  \Sigma_A
       =\operatorname{diag}((\nu(\mu_{A,i}))_i)+\mathcal R_{\Sigma,A},\\
 &\|\mathcal R_{\Sigma,A}\one\|_2\le C\sqrt m/n,
       \qquad |\one^{\mathsf T}\mathcal R_{\Sigma,A}\one|\le C,\\
 &\bigl\|\|(\nu(\mu_{A,i})-\nu(P_A))_i\|_2\bigr\|_q
                 +\|\sum_i\bigl[\nu(\mu_{A,i})-\nu(P_A)\bigr]\|_q\le C_q,\\
 &\bigl\|\|H_A''(\mu_A)\one-2t(P_A)\one\|_2\bigr\|_q
   +\|\one^{\mathsf T}\bigl[H_A''(\mu_A)\one-2t(P_A)\one\bigr]\|_q\le C_q.
 \end{aligned}
\end{equation}
Multiplying the first identity in \eqref{audit:reference-total-covariance-inputs}
by $\one$ and using
$v_0(P_A)=\nu(P_A)+2t(P_A)\nu(P_A)v_0(P_A)$ from
\eqref{eq:logit-profiles} gives the exact equation
\begin{equation}\label{audit:reference-total-covariance-vector}
 \begin{aligned}
 &\left[I_m-\operatorname{diag}((\nu(\mu_{A,i}))_i)H_A''(\mu_A)\right]
                       \bigl[\Sigma_A\one-v_0(P_A)\one\bigr]\\
 &\quad=\bigl[1+2t(P_A)v_0(P_A)\bigr](\nu(\mu_{A,i})-\nu(P_A))_i+\mathcal R_{\Sigma,A}\one\\
 &\qquad+v_0(P_A)\operatorname{diag}((\nu(\mu_{A,i}))_i)
                           \bigl[H_A''(\mu_A)\one-2t(P_A)\one\bigr].
 \end{aligned}
\end{equation}
Applying \eqref{eq:covariance-inverse-norms} to
\eqref{audit:reference-total-covariance-vector} and using
\eqref{audit:reference-total-covariance-inputs} and
\eqref{eq:row-profile-row-sum-gap}, with
$\|\operatorname{diag}((\nu(\mu_{A,i}))_i)\|_{\op}\le1/4$, gives
\begin{equation}\label{audit:dependency:s79:total-covariance-vector-moment}
 \bigl\|\|\Sigma_A\one-v_0(P_A)\one\|_2\bigr\|_q\le C_q.
\end{equation}

For the signed total, put
\begin{equation}\label{audit:reference-total-covariance-test-decomposition}
 \begin{aligned}
 c&=H_A''(\mu_A)(\nu(\mu_{A,i}))_i
                              -2t(P_A)\nu(P_A)\one\\
  &=H_A''(\mu_A)(\nu(\mu_{A,i})-\nu(P_A))_i
       +\nu(P_A)\bigl[H_A''(\mu_A)\one-2t(P_A)\one\bigr].
 \end{aligned}
\end{equation}
Using the symmetry of $H_A''(\mu_A)$, \eqref{eq:K-profile} and
\eqref{audit:reference-total-covariance-inputs} in
\eqref{audit:reference-total-covariance-test-decomposition} gives
\begin{equation}\label{audit:dependency:s79:total-covariance-test-moments}
 \begin{aligned}
 \bigl\|\|c\|_2\bigr\|_q&\le C_q,\\
 \one^{\mathsf T}c
 &=\nu(P_A)\one^{\mathsf T}\bigl[H_A''(\mu_A)\one-2t(P_A)\one\bigr]
          +2t(P_A)\sum_i\bigl[\nu(\mu_{A,i})-\nu(P_A)\bigr]\\
 &\quad+\bigl[H_A''(\mu_A)\one-2t(P_A)\one\bigr]^{\mathsf T}(\nu(\mu_{A,i})-\nu(P_A))_i,\\
 \|\one^{\mathsf T}c\|_q
 &\le C\|\one^{\mathsf T}\bigl[H_A''(\mu_A)\one-2t(P_A)\one\bigr]\|_q
                  +C\|\sum_i\bigl[\nu(\mu_{A,i})-\nu(P_A)\bigr]\|_q\\
 &\quad+\bigl\|\|H_A''(\mu_A)\one-2t(P_A)\one\|_2\bigr\|_{2q}
          \bigl\|\|(\nu(\mu_{A,i})-\nu(P_A))_i\|_2\bigr\|_{2q}\le C_q.
 \end{aligned}
\end{equation}
Summing the first identity in \eqref{audit:reference-total-covariance-inputs}
and using symmetry and
\eqref{audit:reference-total-covariance-test-decomposition} gives
\begin{equation}\label{audit:reference-total-covariance-summed-equation}
 \begin{aligned}
 \one^{\mathsf T}\Sigma_A\one
 &=\sum_i\nu(\mu_{A,i})
        +2t(P_A)\nu(P_A)\one^{\mathsf T}\Sigma_A\one
        +v_0(P_A)\one^{\mathsf T}c\\
 &\quad+c^{\mathsf T}\bigl[\Sigma_A\one-v_0(P_A)\one\bigr]
                   +\one^{\mathsf T}\mathcal R_{\Sigma,A}\one.
 \end{aligned}
\end{equation}
By \eqref{eq:logit-profiles},
$\bigl[1-2t(P_A)\nu(P_A)\bigr]mv_0(P_A)=m\nu(P_A)$.
Subtracting this identity from
\eqref{audit:reference-total-covariance-summed-equation} gives
\begin{equation}\label{audit:reference-total-covariance-signed}
 \begin{aligned}
 &\bigl[1-2t(P_A)\nu(P_A)\bigr]
                  \bigl[\one^{\mathsf T}\Sigma_A\one-mv_0(P_A)\bigr]\\
 &\quad=\sum_i\bigl[\nu(\mu_{A,i})-\nu(P_A)\bigr]+v_0(P_A)\one^{\mathsf T}c
         +c^{\mathsf T}\bigl[\Sigma_A\one-v_0(P_A)\one\bigr]
         +\one^{\mathsf T}\mathcal R_{\Sigma,A}\one.
 \end{aligned}
\end{equation}
By H\"older's inequality, \eqref{audit:reference-total-covariance-inputs},
\eqref{audit:dependency:s79:total-covariance-vector-moment}, and
\eqref{audit:dependency:s79:total-covariance-test-moments}, the right side of
\eqref{audit:reference-total-covariance-signed} is $O_{L^q}(1)$.
The gap in \eqref{eq:row-profile-row-sum-gap} gives the bound for the
variance of the row total.
Polynomial envelopes and \eqref{eq:row-profile-good-exterior-event}
extend both estimates to all exteriors.
\end{proof}

We next expand the off-diagonal entries. For $i\ne j$, the remainder
$R_{A,ij}$ is defined by the exact difference in the following profile.

\begin{lemma}\label{lem:row-entry-covariance}
Assume \hyperref[ass:PC]{\textup{(PC)}}.
With $\Sigma_A$ as in \eqref{v8:def:reference-row-data}, for distinct
$i,j\in V(A)$ and every fixed exterior moment order $q\ge1$,
\begin{equation}\label{eq:cov-entry-profile}
 \begin{aligned}
 \Sigma_{A,ij}
 &=\frac{v_0(P_A)-\nu(P_A)}n
   +\frac{6\nu(P_A)^2\tau(P_A)}n(A_{ij}-P_A)+R_{A,ij},\\
 \left\|\left(\sum_{i\ne j}|R_{A,ij}|^2\right)^{1/2}\right\|_q
 &\le C_qn^{-1/2}.
 \end{aligned}
\end{equation}
The off-diagonal Frobenius norm is taken before the exterior moment.
\end{lemma}

\begin{proof}
Work first on $\mathcal A_n$, with exterior moment norms restricted to
this event until the final estimate. Equations
\eqref{eq:row-profile-covariance-residual} and \eqref{eq:K-profile} give
\[
\begin{aligned}
 \Sigma_A-\operatorname{diag}((\nu(\mu_{A,i}))_i)
 &=\operatorname{diag}((\nu(\mu_{A,i}))_i)
      H_A''(\mu_A)\Sigma_A+\mathcal R_{\Sigma,A},\\
 H_A''(\mu_A)
 &=\frac{2t(P_A)}nJ_m+\frac{6\tau(P_A)}nY_A+R_{A,\rm mean}.
\end{aligned}
\]
Substitution and subtraction of the leading terms yield
\begin{equation}\label{audit:reference-covariance-error-decomposition}
\begin{aligned}
 &\Sigma_A-\operatorname{diag}((\nu(\mu_{A,i}))_i)
 -\frac{2t(P_A)\nu(P_A)v_0(P_A)}nJ_m
 -\frac{6\tau(P_A)\nu(P_A)^2}nY_A\\
 &=\frac{2t(P_A)}n
   \bigl[\operatorname{diag}((\nu(\mu_{A,i}))_i)J_m\Sigma_A
                                  -\nu(P_A)v_0(P_A)J_m\bigr]\\
 &\quad+\frac{6\tau(P_A)}n
   \bigl[\operatorname{diag}((\nu(\mu_{A,i}))_i)Y_A\Sigma_A
                                  -\nu(P_A)^2Y_A\bigr]\\
 &\quad+\operatorname{diag}((\nu(\mu_{A,i}))_i)
                         R_{A,\rm mean}\Sigma_A+\mathcal R_{\Sigma,A}.
\end{aligned}
\end{equation}
For the last two terms, \eqref{eq:row-profile-covariance-residual},
\eqref{eq:K-profile}, and
\eqref{eq:reference-covariance-operator}, with
$\|\operatorname{diag}((\nu(\mu_{A,i}))_i)\|_{\op}\le1/4$, give
\begin{equation}\label{audit:reference-covariance-direct-remainders}
\begin{aligned}
 \|\mathcal R_{\Sigma,A}\|_F
 &\le\sqrt m\,\|\mathcal R_{\Sigma,A}\|_{\infty\to\infty}
       \le Cn^{-1/2},\\
 \bigl\|\|\operatorname{diag}((\nu(\mu_{A,i}))_i)
                  R_{A,\rm mean}\Sigma_A\|_F\bigr\|_q
 &\le C\bigl\|\|R_{A,\rm mean}\|_F\bigr\|_q\le C_qn^{-1/2}.
\end{aligned}
\end{equation}
For the first two terms in
\eqref{audit:reference-covariance-error-decomposition}, it remains to prove
\begin{equation}\label{audit:reference-covariance-product-bounds}
\begin{aligned}
 &\bigl\|\|\operatorname{diag}((\nu(\mu_{A,i}))_i)Y_A\Sigma_A
                                       -\nu(P_A)^2Y_A\|_F\bigr\|_q\\
 &\quad+\bigl\|\|\operatorname{diag}((\nu(\mu_{A,i}))_i)J_m\Sigma_A
                                  -\nu(P_A)v_0(P_A)J_m\|_F\bigr\|_q
 \le C_q\sqrt n.
\end{aligned}
\end{equation}

For the $Y_A$ bound in \eqref{audit:reference-covariance-product-bounds},
use $Y_AJ_m=z^A\one^{\mathsf T}$ and $\|Y_A\|_{\op}\le m$ in
\eqref{eq:K-profile} to obtain
\begin{equation}\label{audit:reference-covariance-y-hessian}
\begin{aligned}
 Y_AH_A''(\mu_A)
 &=\frac{2t(P_A)}n z^A\one^{\mathsf T}
   +\frac{6\tau(P_A)}nY_A^2+Y_AR_{A,\rm mean},\\
 \bigl\|\|Y_AH_A''(\mu_A)\|_F\bigr\|_q
 &\le\frac{C\sqrt m}{n}\bigl\|\|z^A\|_2\bigr\|_q
   +\frac Cn\bigl\|\|Y_A^2\|_F\bigr\|_q
   +m\bigl\|\|R_{A,\rm mean}\|_F\bigr\|_q\\
 &\le C_q\sqrt n,
\end{aligned}
\end{equation}
where the last line uses \eqref{eq:row-profile-exterior-degree-vector},
\eqref{eq:four-cycle-square-frobenius}, and \eqref{eq:K-profile}.
Next, the symmetry of $Y_A$, the bound $|Y_{A,ij}|\le1$, and
\eqref{eq:row-profile-variance-input} give
\begin{equation}\label{audit:reference-covariance-diagonal-products}
\begin{aligned}
 &\|Y_A\operatorname{diag}((\nu(\mu_{A,i})-\nu(P_A))_i)\|_F^2
  =\|\operatorname{diag}((\nu(\mu_{A,i})-\nu(P_A))_i)Y_A\|_F^2\\
 &\quad=\sum_j\bigl[\nu(\mu_{A,j})-\nu(P_A)\bigr]^2\left(\sum_iY_{A,ij}^2\right)
  \le m\|(\nu(\mu_{A,i})-\nu(P_A))_i\|_2^2,\\
 &\bigl\|\|Y_A\operatorname{diag}((\nu(\mu_{A,i})-\nu(P_A))_i)\|_F\bigr\|_q
  +\bigl\|\|\operatorname{diag}((\nu(\mu_{A,i})-\nu(P_A))_i)Y_A\|_F\bigr\|_q
  \le C_q\sqrt n.
\end{aligned}
\end{equation}
Multiplying \eqref{eq:row-profile-covariance-residual} by $Y_A$ on the
left and splitting the diagonal matrix into $\nu(P_A)I_m$ and
$\operatorname{diag}((\nu(\mu_{A,i})-\nu(P_A))_i)$ gives
\[
\begin{aligned}
 &Y_A\bigl[\Sigma_A-\operatorname{diag}((\nu(\mu_{A,i}))_i)\bigr]
 =\nu(P_A)Y_AH_A''(\mu_A)\Sigma_A\\
 &\quad+Y_A\operatorname{diag}((\nu(\mu_{A,i})-\nu(P_A))_i)
                      H_A''(\mu_A)\Sigma_A+Y_A\mathcal R_{\Sigma,A}.
\end{aligned}
\]
The operator bounds \eqref{eq:reference-covariance-operator} and
\eqref{eq:K-profile}, together with
\eqref{audit:reference-covariance-y-hessian},
\eqref{audit:reference-covariance-diagonal-products}, and
\eqref{audit:reference-covariance-direct-remainders}, therefore yield
\begin{equation}\label{audit:dependency:s79:centered-covariance-product}
 \bigl\|\|Y_A\bigl[\Sigma_A-\operatorname{diag}((\nu(\mu_{A,i}))_i)\bigr]\|_F\bigr\|_q
 \le C_q\sqrt n.
\end{equation}
Adding and subtracting the diagonal covariance term and using the same
decomposition of the diagonal factors in their original order gives
\[
\begin{aligned}
 &\operatorname{diag}((\nu(\mu_{A,i}))_i)Y_A\Sigma_A-\nu(P_A)^2Y_A\\
 &=\operatorname{diag}((\nu(\mu_{A,i}))_i)Y_A
                 \bigl[\Sigma_A-\operatorname{diag}((\nu(\mu_{A,i}))_i)\bigr]\\
 &\quad+\operatorname{diag}((\nu(\mu_{A,i})-\nu(P_A))_i)
                 Y_A\operatorname{diag}((\nu(\mu_{A,i}))_i)\\
 &\quad+\nu(P_A)Y_A
                 \operatorname{diag}((\nu(\mu_{A,i})-\nu(P_A))_i).
\end{aligned}
\]
By \eqref{audit:dependency:s79:centered-covariance-product},
\eqref{audit:reference-covariance-diagonal-products}, and
\eqref{eq:row-profile-variance-input},
\[
\begin{aligned}
 &\bigl\|\|\operatorname{diag}((\nu(\mu_{A,i}))_i)Y_A\Sigma_A
                                  -\nu(P_A)^2Y_A\|_F\bigr\|_q\\
 &\quad\le\tfrac14
   \bigl\|\|Y_A\bigl[\Sigma_A-\operatorname{diag}((\nu(\mu_{A,i}))_i)\bigr]\|_F\bigr\|_q
   +C\sqrt m\,\bigl\|\|(\nu(\mu_{A,i})-\nu(P_A))_i\|_2\bigr\|_q
 \le C_q\sqrt n.
\end{aligned}
\]
This proves the first bound in
\eqref{audit:reference-covariance-product-bounds}.

For the $J_m$ bound, use $J_m=\one\one^{\mathsf T}$ and the symmetry of
$\Sigma_A$ to obtain
\[
\begin{aligned}
 &\operatorname{diag}((\nu(\mu_{A,i}))_i)J_m\Sigma_A
                                  -\nu(P_A)v_0(P_A)J_m\\
 &=\bigl[(\nu(\mu_{A,i}))_i-\nu(P_A)\one\bigr](\Sigma_A\one)^{\mathsf T}
  +\nu(P_A)\one\bigl[\Sigma_A\one-v_0(P_A)\one\bigr]^{\mathsf T}.
\end{aligned}
\]
Since $\|\Sigma_A\one\|_2\le C\sqrt m$ by
\eqref{eq:reference-covariance-operator}, the identity
$\|ab^{\mathsf T}\|_F=\|a\|_2\|b\|_2$, together with
\eqref{eq:row-profile-variance-input} and \eqref{eq:cov-profiles-a}, gives
\[
\begin{aligned}
 &\bigl\|\|\operatorname{diag}((\nu(\mu_{A,i}))_i)J_m\Sigma_A
                                  -\nu(P_A)v_0(P_A)J_m\|_F\bigr\|_q\\
 &\quad\le C\sqrt m\,\bigl\|\|(\nu(\mu_{A,i})-\nu(P_A))_i\|_2\bigr\|_q
   +C\sqrt m\,\bigl\|\|\Sigma_A\one-v_0(P_A)\one\|_2\bigr\|_q
 \le C_q\sqrt n.
\end{aligned}
\]
This completes \eqref{audit:reference-covariance-product-bounds}.

We now return to \eqref{audit:reference-covariance-error-decomposition}.
By \eqref{eq:logit-profiles},
$2t(P_A)\nu(P_A)v_0(P_A)=v_0(P_A)-\nu(P_A)$ on $\mathcal A_n$.
Since $J_{m,ij}=1$ and $Y_{A,ij}=A_{ij}-P_A$ for $i\ne j$,
the off-diagonal entries on the left of
\eqref{audit:reference-covariance-error-decomposition} are exactly
$R_{A,ij}$. On the complement, bounded covariances and the globally
bounded profiles give $(\sum_{i\ne j}|R_{A,ij}|^2)^{1/2}\le Cn$.
Thus \eqref{audit:reference-covariance-direct-remainders},
\eqref{audit:reference-covariance-product-bounds}, and
\eqref{eq:row-profile-good-exterior-event} give
\[
 \left\|\left(\sum_{i\ne j}|R_{A,ij}|^2\right)^{1/2}\right\|_q
 \le C_q(n^{-1}\sqrt n+n^{-1/2})
       +Cn\,\Pp_n(A\notin\mathcal A_n)^{1/q}
 \le C_q'n^{-1/2}.
\]
This proves \eqref{eq:cov-entry-profile}.
\end{proof}

\section{Two-row logit expansions and pairings}\label{sec:two-row}
We expand the physical-edge logit
\eqref{v8:setup:physical-edge-logit} under the two-row reference law
and transfer its variance and common-neighbor pairings to the true law.
\subsection{Two-row laws and common-neighbor moments}
Fix $A=X_{\{3,\ldots,n\}}$, $e=X_{12}$ and $m=n-2$.
Write $x_i=X_{1i}$, $y_i=X_{2i}$ and use the reference data
$\mu_{A,i}=\E_{\pi_A}X_i$, $\Sigma_{A,ij}=\Cov_{\pi_A}(X_i,X_j)$
from Section~\ref{sec:reference-row-profiles}. The interpolation
$\mu_{A,(e),s}$ has endpoints
$\mu_{A,(e),0}=\pi_A\otimes\pi_A$ and
$\mu_{A,(e),1}=\mathcal L((x,y)\mid A,e)$.
Its density relative to the reference product law is
\begin{equation}\label{v8:setup:two-row-interpolation-density}
 \frac{\dd\mu_{A,(e),s}}{\dd(\pi_A\otimes\pi_A)}
 =\frac{\exp(s V_{A,e})}{\E_{\pi_A\otimes\pi_A}\exp(s V_{A,e})}.
\end{equation}
For fixed $(A,e)$ and observables independent of $s$, differentiation gives
\begin{align}
 \frac{\dd}{\dd s}\E_{\mu_{A,(e),s}}F
 &=\Cov_{\mu_{A,(e),s}}(F,V_{A,e}),\label{eq:two-row-expectation-derivative}\\
 \frac{\dd}{\dd s}\Cov_{\mu_{A,(e),s}}(F,G)
 &=\cum_{\mu_{A,(e),s}}(F,G,V_{A,e}),\label{eq:two-row-covariance-derivative}\\
 \Cov_{\mu_{A,(e),s}}(FG,V)&=\cum_{\mu_{A,(e),s}}(F,G,V)\notag\\
 &\quad+\E_{\mu_{A,(e),s}}F\,\Cov_{\mu_{A,(e),s}}(G,V)\notag\\
 &\quad+\E_{\mu_{A,(e),s}}G\,\Cov_{\mu_{A,(e),s}}(F,V).
 \label{eq:two-row-product-covariance}
\end{align}
The last identity follows from the third-cumulant formula and holds
under every interpolating law, including when $F=G$.
To apply the conditional cumulant bounds, note that both $L_{12}$
and $V_{A,e}$ are normalized block polynomials: fixing $d$ exterior
sites leaves $O(n^{v_r-2-d})$ template completions and hence a
coefficient of order $O(n^{-d})$. A prescribed root edge contributes only the bounded factor
$e$. The same property holds for normalized row totals and $B_{12}/n$,
where
\begin{equation}\label{v8:setup:common-neighbor-polynomial}
 B_{ij}(u)=\sum_{k\notin\{i,j\}}(X_{ik}-u)(X_{jk}-u).
\end{equation}
Theorem~\ref{ep:row:sharp}, \eqref{ep:row:polynomial-bound}, therefore
bounds the corresponding covariances by $C/n$ and third cumulants by
$C/n^2$, uniformly in $e,s$ on the common event $\mathcal A_n$ of
\eqref{eq:row-profile-good-exterior-event}.

Outer $\|\cdot\|_q$ norms use the actual marginal of $(A,e)$,
as defined in \eqref{v8:setup:two-row-observables}; the probability measure
is always indicated for inner row norms.
Intermediate estimates are restricted to $\mathcal A_n$, and products
use sufficiently high fixed moments. The bounds $|L_{12}|\le C$,
$|B_{12}|\le n$ and $\Pp_n(\mathcal A_n^c)\le Ce^{-cn}$ make
exceptional-exterior errors $O(n^{-M})$ for every fixed $M$.
\begin{equation}\label{v8:setup:two-row-observables}
\begin{gathered}
 P_A=\binom{n-2}{2}^{-1}\sum_{i<j:\ i,j\in V(A)}A_{ij},
 \qquad B_{12}(u)=\sum_i(x_i-u)(y_i-u),\\
 \qquad \|R\|_q=(\E_{(A,e)}|R|^q)^{1/q}.
\end{gathered}
\end{equation}

\begin{lemma}\label{lem:common-neighbor-moments}
Assume \hyperref[ass:PC]{\textup{(PC)}}.
For every fixed $q\ge1$, the following bounds hold uniformly
for $0\le s\le1$:
\begin{equation}\label{eq:revised-common-neighbor-moments}
 \begin{aligned}
 \|\E_{\mu_{A,(e),s}}B_{12}(P_A)\|_q&\le C_q,\\
 \left(\E_{(A,e)}\E_{\mu_{A,(e),s}}|B_{12}(P_A)|^q\right)^{1/q}
       &\le C_q\sqrt n.
 \end{aligned}
\end{equation}
In particular, under the original graph law,
\begin{equation}\label{eq:B-bounds}
 \|\E[B_{12}(P_A)\mid A,e]\|_q
       +n^{-1/2}\|B_{12}(P_A)\|_{L^q(\Pp_n)}\le C_q.
\end{equation}
\end{lemma}

\begin{proof}
We bound the conditional mean and then the centered moments.
Root interchange gives $a_i:=\E_{\mu_{A,(e),s}}x_i-P_A
=\E_{\mu_{A,(e),s}}y_i-P_A$. Hence
\[
\begin{aligned}
 \E_{\mu_{A,(e),s}}B_{12}(P_A)
 &=\sum_i\Cov_{\mu_{A,(e),s}}(x_i,y_i)+\|a\|_2^2,\\
 \left\|\E_{\mu_{A,(e),s}}B_{12}(P_A)\right\|_q
 &\le Cm/n+\bigl\|\|a\|_2\bigr\|_{2q}^2\le C_q.
\end{aligned}
\]
Here \eqref{eq:disjoint-local-covariance} bounds the distinct-edge
covariances, and \eqref{eq:mean-profiles-a} bounds $a$.
For the centered moments, the features $(x_i-P_A)(y_i-P_A)$ have
disjoint coordinate-pair supports, so \eqref{eq:weighted-local-moments} gives, on good exteriors,
\[
 \|B_{12}(P_A)-\E_{\mu_{A,(e),s}}B_{12}(P_A)\|_{L^q(\mu_{A,(e),s})}
 \le C_q\sqrt m.
\]
Combining the mean and centered-moment bounds by the mixed-norm triangle
inequality proves \eqref{eq:revised-common-neighbor-moments};
$|B_{12}|\le n$ handles bad exteriors. Taking $s=1$ gives
\eqref{eq:B-bounds}.
\end{proof}

\subsection{Single-site coefficients and reference logit expansion}
In the following expansion, derivatives of $L_{12}(x,y;A)$ are taken
before substituting the common row argument. Root interchange gives
$L_{12}(x,y;A)=L_{12}(y,x;A)$, so the two linear derivatives agree at
$x=y=\mu_A$.
Define the exterior coefficient errors by
\begin{equation}\label{v8:setup:single-site-logit-errors}
\begin{aligned}
 e_{A,i}^{(1)}
 &=\left.\partial_{x_i}L_{12}\right|_{x=y=\mu_A}
                  -2t(P_A)/n
  =\left.\partial_{y_i}L_{12}\right|_{x=y=\mu_A}-2t(P_A)/n,\\
 e_{A,i}^{(12)}
 &=\left.\partial_{x_i}\partial_{y_i}L_{12}\right|_{x=y=\mu_A}
                  -6\tau(P_A)/n .
\end{aligned}
\end{equation}

\begin{lemma}\label{lem:logit-coefficients}
Assume \hyperref[ass:PC]{\textup{(PC)}}.
For the coefficient errors in \eqref{v8:setup:single-site-logit-errors}
and every fixed exterior moment order $q\ge1$,
\begin{equation}\label{eq:revised-logit-coefficient-errors}
 \sum_{\alpha\in\{1,12\}}
 \left(\bigl\|\|e_A^{(\alpha)}\|_2\bigr\|_q+
       \left\|\sum_i e_{A,i}^{(\alpha)}\right\|_q\right)
       \le C_qn^{-1}.
\end{equation}
Both norms are taken under the exterior law; in particular, the
coordinate sum is estimated before taking its absolute value.
\end{lemma}

\begin{proof}
\proofstep{1. Exact coefficients.}
For an ordered template edge $(u,v)$, put $K=F_r-\{u,v\}$,
retaining isolates. Then
\[
 \partial_{12}\inj(F_r,X)
 =\sum_{\substack{u,v\in V(F_r)\\uv\in E(F_r)}}\sum_{\varphi:V(K)\hookrightarrow V(A)}
 \prod_{ab\in E(K)}A_{\varphi(a)\varphi(b)}
 \prod_{\substack{a\in V(K)\\ua\in E(F_r)}}x_{\varphi(a)}
 \prod_{\substack{b\in V(K)\\vb\in E(F_r)}}y_{\varphi(b)}.
\]
A derivative at site $i$ fixes one image, giving
$n^{2-v_r}n^{v_r-3}=n^{-1}$; two derivatives at that site
must select the same vertex by injectivity. Thus
\begin{equation}\label{audit:single-site-logit-coefficient-formulas}
\begin{aligned}
 &b_{r,n}n^{2-v_r}
 \left.\partial_{x_i}\partial_{12}\inj(F_r,X)\right|_{x=y=s\one}\\
 &=\frac{b_{r,n}}n\sum_{\substack{u,v\in V(F_r)\\uv\in E(F_r)}}
 \sum_{\substack{w\in V(K)\\uw\in E(F_r)}}
 s^{d(u)+d(v)-3}T_i(K,w;A),\\
 &b_{r,n}n^{2-v_r}
 \left.\partial_{x_i}\partial_{y_i}\partial_{12}\inj(F_r,X)
                                      \right|_{x=y=s\one}\\
 &=\frac{b_{r,n}}n\sum_{\substack{u,v\in V(F_r)\\uv\in E(F_r)}}
 \sum_{\substack{w\in V(K)\\uw,vw\in E(F_r)}}s^{d(u)+d(v)-4}T_i(K,w;A).
\end{aligned}
\end{equation}
Interchanging $(u,v,x,y)$ proves equality of the two linear coefficients.
All nonempty summands have nonnegative exponents, and
\[
 \begin{gathered}
 e(K)=e_r-d(u)-d(v)+1,\qquad
 \sum_{\substack{u,v\in V(F_r)\\uv\in E(F_r)}}(d(u)-1)=2a(F_r),\\
 \sum_{\substack{u,v,w\in V(F_r)\\uv,uw,vw\in E(F_r)}}1=6\triangle(F_r).
 \end{gathered}
\]
Consequently, replacing $T_i(K,w;A)$ by $P_A^{e(K)}$ and
$b_{r,n}$ by $b_r$ gives $2t(P_A)/n$ for the first derivative
and $6\tau(P_A)/n$ for the mixed derivative.

\proofstep{2. Constant-row errors.}
For the rooted-count errors in
\eqref{audit:single-site-logit-coefficient-formulas}, take
$T=(T_i(K,w;A))_i$, with $\bar T=m^{-1}\sum_iT_i$.
Equations~\eqref{audit:exterior-count-replacement} and
\eqref{eq:root-deterministic} give
\[
\begin{aligned}
 \frac1n\sum_iT_i&=n^{-v(K)}\inj(K,A)
                         =P_A^{e(K)}+O_{L^q}(n^{-1}),\\
 \bigl\|\|(T_i-P_A^{e(K)})_i\|_2\bigr\|_q
 &\le\bigl\|\|T-\bar T\one\|_2\bigr\|_q
       +\sqrt m\|\bar T-P_A^{e(K)}\|_q\le C_q.
\end{aligned}
\]
Using \eqref{audit:single-site-logit-coefficient-formulas}, $m=n-2$,
and $b_{r,n}-b_r=O(n^{-1})$ gives
\begin{equation}\label{audit:dependency:s79:logit-constant-row-errors}
\begin{aligned}
 &\left\|\left\|\left(
 \partial_{x_i}L_{12}(P_A\one,P_A\one;A)-\frac{2t(P_A)}n
                              \right)_i\right\|_2\right\|_q\\
 &\quad+\left\|\sum_i\left\{
 \partial_{x_i}L_{12}(P_A\one,P_A\one;A)-\frac{2t(P_A)}n
                         \right\}\right\|_q\le C_q/n,\\
 &\left\|\left\|\left(
 \partial_{x_i}\partial_{y_i}L_{12}(P_A\one,P_A\one;A)
                   -\frac{6\tau(P_A)}n\right)_i\right\|_2\right\|_q\\
 &\quad+\left\|\sum_i\left\{
 \partial_{x_i}\partial_{y_i}L_{12}(P_A\one,P_A\one;A)
                   -\frac{6\tau(P_A)}n\right\}\right\|_q\le C_q/n.
\end{aligned}
\end{equation}

\proofstep{3. Shift to $\mu_A$.}
Counting distinct prescribed sites gives the following bounds on the whole row cube:
\begin{equation}\label{audit:dependency:s79:logit-common-row-derivatives}
\begin{aligned}
 &\left|\partial_{s_j}\{\partial_{x_i}L_{12}(\mathbf s,\mathbf s;A)\}\right|
 +\left|\partial_{s_j}\{\partial_{x_i}\partial_{y_i}
                         L_{12}(\mathbf s,\mathbf s;A)\}\right|
 \le Cn^{-|\{i,j\}|},\\
 &\left|\partial_{s_j}\partial_{s_k}
                \{\partial_{x_i}L_{12}(\mathbf s,\mathbf s;A)\}\right|
 +\left|\partial_{s_j}\partial_{s_k}
                \{\partial_{x_i}\partial_{y_i}L_{12}(\mathbf s,\mathbf s;A)\}\right|
 \le Cn^{-|\{i,j,k\}|}.
\end{aligned}
\end{equation}
Repeated common-row derivatives are retained: $x_jy_j=s_j^2$, and
\[
\begin{aligned}
 \partial_{s_j}\{\partial_{x_i}L_{12}(\mathbf s,\mathbf s;A)\}
 &=\left.(\partial_{x_j}+\partial_{y_j})\partial_{x_i}L_{12}
                                                    \right|_{x=y=\mathbf s},\\
 \partial_{s_i}\{\partial_{x_i}L_{12}(\mathbf s,\mathbf s;A)\}
 &=\partial_{x_i}\partial_{y_i}L_{12}(\mathbf s,\mathbf s;A),\\
 \partial_{s_i}\{\partial_{x_i}\partial_{y_i}L_{12}(\mathbf s,\mathbf s;A)\}
 &=0.
\end{aligned}
\]
For the Euclidean shift, integrate
\eqref{audit:dependency:s79:logit-common-row-derivatives}
along the segment to obtain
\[
\begin{aligned}
 &|\partial_{x_i}L_{12}(\mu_A,\mu_A;A)
                  -\partial_{x_i}L_{12}(P_A\one,P_A\one;A)|\\
 &\quad+|\partial_{x_i}\partial_{y_i}L_{12}(\mu_A,\mu_A;A)
                  -\partial_{x_i}\partial_{y_i}L_{12}(P_A\one,P_A\one;A)|\\
 &\quad\le Cn^{-1}|\mu_{A,i}-P_A|
                    +Cn^{-2}\|\mu_A-P_A\one\|_1,
\end{aligned}
\]
and therefore
\begin{equation}\label{audit:dependency:s79:logit-vector-shift}
\begin{aligned}
 &\left\|\left(
 \partial_{x_i}L_{12}(\mu_A,\mu_A;A)
                  -\partial_{x_i}L_{12}(P_A\one,P_A\one;A)
                              \right)_i\right\|_2\\
 &\quad+\left\|\left(
 \partial_{x_i}\partial_{y_i}L_{12}(\mu_A,\mu_A;A)
                  -\partial_{x_i}\partial_{y_i}L_{12}(P_A\one,P_A\one;A)
                              \right)_i\right\|_2\\
 &\quad\le Cn^{-1}\|\mu_A-P_A\one\|_2
                    +Cn^{-3/2}\|\mu_A-P_A\one\|_1
 \le Cn^{-1}\|\mu_A-P_A\one\|_2=O_{L^q}(n^{-1}).
\end{aligned}
\end{equation}
For the signed shift, summing first removes the selected-image constraint.
The two derivative families give
\begin{align*}
 &\left.\partial_{s_j}\sum_i\partial_{x_i}L_{12}(\mathbf s,\mathbf s;A)
                                      \right|_{\mathbf s=P_A\one}\\*
 &\quad=\frac1n\sum_r b_{r,n}\sum_{\substack{u,v\in V(F_r)\\uv\in E(F_r)}}
 \sum_{\substack{w\in V(K)\\uw\in E(F_r)}}P_A^{d(u)+d(v)-4}\\*
 &\quad\times\left\{
 \sum_{\substack{a\in V(K)\setminus\{w\}\\ua\in E(F_r)}}T_j(K,a;A)
 +\sum_{\substack{b\in V(K)\\vb\in E(F_r)}}T_j(K,b;A)\right\},\\
 &\left.\partial_{s_j}\sum_i\partial_{x_i}\partial_{y_i}
                         L_{12}(\mathbf s,\mathbf s;A)
                                      \right|_{\mathbf s=P_A\one}\\*
 &\quad=\frac1n\sum_r b_{r,n}\sum_{\substack{u,v\in V(F_r)\\uv\in E(F_r)}}
 \sum_{\substack{w\in V(K)\\uw,vw\in E(F_r)}}P_A^{d(u)+d(v)-5}\\*
 &\quad\times\left\{
 \sum_{\substack{a\in V(K)\setminus\{w\}\\ua\in E(F_r)}}T_j(K,a;A)
 +\sum_{\substack{b\in V(K)\setminus\{w\}\\vb\in E(F_r)}}T_j(K,b;A)\right\}.
\end{align*}
The two sums in each brace retain the multiplicity of any common
vertex. Terms with no surviving row factor are omitted before evaluating
the power of $P_A$.
Each rooted summand satisfies
\begin{equation}\label{audit:dependency:s79:logit-rooted-linear-pairing}
\begin{aligned}
 \left|\frac1n\sum_jT_j(\mu_{A,j}-P_A)\right|
 &\le\frac{|\bar T|}{n}\left|\sum_j(\mu_{A,j}-P_A)\right|\\
 &\quad+\frac1n\|T-\bar T\one\|_2\|\mu_A-P_A\one\|_2
 =O_{L^q}(n^{-1})
\end{aligned}
\end{equation}
by Lemmas~\ref{lem:root-exchange} and~\ref{lem:row-mean-density}.
For the quadratic Taylor errors, the distinct-site derivative bounds imply
\[
\begin{aligned}
 &\left|\partial_{s_j}\partial_{s_k}
                 \sum_i\partial_{x_i}L_{12}(\mathbf s,\mathbf s;A)\right|\\
 &\quad+\left|\partial_{s_j}\partial_{s_k}
                 \sum_i\partial_{x_i}\partial_{y_i}
                                   L_{12}(\mathbf s,\mathbf s;A)\right|
 \le\begin{cases}C/n,&j=k,\\C/n^2,&j\ne k.\end{cases}
\end{aligned}
\]
Taylor's formula therefore gives
\begin{equation}\label{audit:dependency:s79:logit-signed-shift}
\begin{aligned}
 &\left|\sum_i\{\partial_{x_i}L_{12}(\mu_A,\mu_A;A)
               -\partial_{x_i}L_{12}(P_A\one,P_A\one;A)\}\right|\\
 &\quad\le\left|\sum_j
  \left.\partial_{s_j}\sum_i\partial_{x_i}L_{12}(\mathbf s,\mathbf s;A)
                                  \right|_{\mathbf s=P_A\one}
                       (\mu_{A,j}-P_A)\right|
      +Cn^{-1}\|\mu_A-P_A\one\|_2^2,\\
 &\left|\sum_i\{\partial_{x_i}\partial_{y_i}L_{12}(\mu_A,\mu_A;A)
               -\partial_{x_i}\partial_{y_i}L_{12}(P_A\one,P_A\one;A)\}\right|\\
 &\quad\le\left|\sum_j
  \left.\partial_{s_j}\sum_i\partial_{x_i}\partial_{y_i}
                         L_{12}(\mathbf s,\mathbf s;A)
                                  \right|_{\mathbf s=P_A\one}
                       (\mu_{A,j}-P_A)\right|
      +Cn^{-1}\|\mu_A-P_A\one\|_2^2.
\end{aligned}
\end{equation}
Equations~\eqref{audit:dependency:s79:logit-rooted-linear-pairing} and
\eqref{eq:mean-profiles-a} bound both right-hand sides of
\eqref{audit:dependency:s79:logit-signed-shift} in exterior $L^q$ by $C_q/n$.
Combining \eqref{audit:dependency:s79:logit-constant-row-errors},
\eqref{audit:dependency:s79:logit-vector-shift}, and
\eqref{audit:dependency:s79:logit-signed-shift} proves
\eqref{eq:revised-logit-coefficient-errors}.
\end{proof}

Under $\pi_A\otimes\pi_A$, set
\begin{equation}\label{v8:setup:reference-centered-sums}
 S_{A,x}=\sum_i(x_i-\mu_{A,i}),\quad
 S_{A,y}=\sum_i(y_i-\mu_{A,i}),\quad
 B_{A,0}=\sum_i(x_i-\mu_{A,i})(y_i-\mu_{A,i}),
\end{equation}
and define the centered leading term
\begin{equation}\label{eq:revised-logit-leading-term}
 L_{A,0}=\frac{2t(P_A)}n(S_{A,x}+S_{A,y})
                   +\frac{6\tau(P_A)}nB_{A,0}.
\end{equation}

\begin{lemma}\label{lem:reference-logit-expansion}
Assume \hyperref[ass:PC]{\textup{(PC)}}.

Under $\pi_A\otimes\pi_A$, for every fixed $q\ge1$,
\begin{equation}\label{audit:reference-logit-approximation}
 \left\|\|L_{12}-\ell_n(P_A)-L_{A,0}
                         \|_{L^q(\pi_A\otimes\pi_A)}\right\|_q\le C_qn^{-1}.
\end{equation}
As stipulated above, the exterior norm on the left is restricted to
$\mathcal A_n$.
\end{lemma}

\begin{proof}
We expand about $\mu_A$ and bound the constant, single-site, and
higher-site errors. Outer $\|\cdot\|_q$ norms use the actual marginal
of $(A,e)$ in \eqref{v8:setup:two-row-observables}, with the stated
restriction to $\mathcal A_n$. For fixed $A$, under
$\pi_A\otimes\pi_A$, set
\[
 \psi_{i,1}=x_i-\mu_{A,i},\quad \psi_{i,2}=y_i-\mu_{A,i},\quad
 \psi_{i,12}=\psi_{i,1}\psi_{i,2},\qquad \E_{\pi_A\otimes\pi_A}\psi_{i,\alpha}=0.
\]
Row independence centers the bilinear feature; all supports and overlaps
are bounded. Expanding each template product and grouping by distinct
selected sites gives exactly
\begin{equation}\label{audit:reference-logit-exact-expansion}
 L_{12}-\ell_n(P_A)=L_{A,0}
 +\{L_{12}(\mu_A,\mu_A;A)-\ell_n(P_A)\}
 +R_{A,1}+\sum_{d=2}^D H_{A,d},
\end{equation}
where $D=\max_r(v_r-2)$ and
\begin{equation}\label{audit:reference-logit-remainder-polynomials}
 R_{A,1}=\sum_i\{e_{A,i}^{(1)}(\psi_{i,1}+\psi_{i,2})
                         +e_{A,i}^{(12)}\psi_{i,12}\},\qquad
 H_{A,d}=\sum_{\beta_1,\ldots,\beta_d}
        c_{A,\beta_1\ldots\beta_d}\prod_{a=1}^d\psi_{\beta_a}.
\end{equation}
Here $\beta=(i,\alpha)$, and $H_{A,d}$ has $d$ distinct sites,
counting a bilinear feature once. Fixing these sites leaves
$O(n^{v_r-2-d})$ completions, so $|c_{A,\boldsymbol\beta}|\le Cn^{-d}$.
For the random remainders, Corollary~\ref{cor:local-feature-sums} and
Lemmas~\ref{lem:weighted-linear} and~\ref{lem:logit-coefficients} give
\begin{equation}\label{audit:dependency:s79:reference-logit-remainder-moments}
\begin{aligned}
 \|H_{A,d}\|_{L^q(\pi_A\otimes\pi_A)}&\le C_{q,d}n^{-d/2},\\
 \|R_{A,1}\|_{L^q(\pi_A\otimes\pi_A)}
 &\le C_q(2\|e_A^{(1)}\|_2+\|e_A^{(12)}\|_2),\qquad \E_{\pi_A\otimes\pi_A}R_{A,1}=0,\\
 \left\|\|R_{A,1}\|_{L^q(\pi_A\otimes\pi_A)}\right\|_q&\le C_q/n.
\end{aligned}
\end{equation}
For the constant error in \eqref{audit:reference-logit-exact-expansion},
Lemma~\ref{lem:conditional-logit-bias} supplies
\begin{equation}\label{eq:revised-logit-interpolated-bias}
 \|\one_{\mathcal A_n}\E_{\mu_{A,(e),s}}[L_{12}-\ell_n(P_A)]\|_q
 \le C_qn^{-1},\qquad 0\le s\le1.
\end{equation}
At $s=0$, $\mu_{A,(e),0}=\pi_A\otimes\pi_A$. Taking
$\E_{\pi_A\otimes\pi_A}$ in
\eqref{audit:reference-logit-exact-expansion} isolates this error. Using
\eqref{audit:dependency:s79:reference-logit-remainder-moments} and
\eqref{eq:revised-logit-interpolated-bias} gives
\begin{equation}\label{audit:dependency:s79:reference-logit-constant-bias}
\begin{aligned}
 &L_{12}(\mu_A,\mu_A;A)-\ell_n(P_A)
   =\E_{\pi_A\otimes\pi_A}[L_{12}-\ell_n(P_A)]
       -\sum_{d=2}^D\E_{\pi_A\otimes\pi_A}H_{A,d},\\
 &\|\one_{\mathcal A_n}\{L_{12}(\mu_A,\mu_A;A)-\ell_n(P_A)\}\|_q
 \le C_q/n+\sum_{d=2}^D C_{q,d}n^{-d/2}\le C_q/n.
\end{aligned}
\end{equation}
The mixed-norm triangle inequality in
\eqref{audit:reference-logit-exact-expansion}, using
\eqref{audit:dependency:s79:reference-logit-remainder-moments} and
\eqref{audit:dependency:s79:reference-logit-constant-bias}, proves
\eqref{audit:reference-logit-approximation}.
\end{proof}

\subsection{Conditional variance and common-neighbor pairings}
Row independence gives
$\E_{\pi_A\otimes\pi_A}(\psi_{i,1}\psi_{j,1})=\E_{\pi_A\otimes\pi_A}(\psi_{i,2}\psi_{j,2})=\Sigma_{A,ij}$
and $\E_{\pi_A\otimes\pi_A}(\psi_{i,12}\psi_{j,12})=\Sigma_{A,ij}^2$; every mixed-type
pairing vanishes. Thus, for vectors fixed by the exterior,
\begin{equation}\label{eq:two-row-site-pairing}
\begin{split}
 &\E_{\pi_A\otimes\pi_A}
 \left[\sum_i(a_i\psi_{i,1}+b_i\psi_{i,2}+c_i\psi_{i,12})
       \sum_j(\widetilde a_j\psi_{j,1}
             +\widetilde b_j\psi_{j,2}+\widetilde c_j\psi_{j,12})\right]\\
 &\qquad=a^{\mathsf T}\Sigma_A\widetilde a+
          b^{\mathsf T}\Sigma_A\widetilde b+
          \sum_{i,j}c_i\widetilde c_j\Sigma_{A,ij}^2 .
\end{split}
\end{equation}
For the higher-site terms in \eqref{audit:reference-logit-remainder-polynomials}
of Lemma~\ref{lem:reference-logit-expansion},
Corollary~\ref{cor:local-feature-sums}, \eqref{eq:micro-moment-sum},
also gives, for $d\ge2$,
\begin{equation}\label{eq:two-row-higher-site-signed-pairing}
\begin{aligned}
 &\left|\E_{\pi_A\otimes\pi_A}
       \left[H_{A,d}\sum_i\psi_{i,\alpha}\right]\right|\\
 &\le Cn^{-d}\sum_{\beta_0,\ldots,\beta_d}
          \left|\E_{\pi_A\otimes\pi_A}\prod_{a=0}^d\psi_{\beta_a}\right|\\
 &\le C_dn^{-d+\lfloor(d+1)/2\rfloor}\le C_dn^{-1}.
\end{aligned}
\end{equation}
Here cancellation is retained within each expectation before taking
absolute values. In particular, the three-feature sum is $O(n)$.
Using the fixed profiles in \eqref{eq:profiles} and
\eqref{eq:logit-profiles}, define the smooth logit variance profile
\begin{equation}\label{v8:setup:two-row-variance-profile}
 w_2(u)=8t(u)^2v_0(u)+36\tau(u)^2\nu(u)^2,\qquad u\in[0,1].
\end{equation}

\begin{lemma}\label{lem:logit-pairings}
Assume \hyperref[ass:PC]{\textup{(PC)}}.\label{lem:logit-variance}
For every fixed $q\ge1$,
\begin{equation}\label{eq:variance}
 \left\|\Var(L_{12}\mid A,e)-\frac{w_2(P_A)}n\right\|_q
       \le C_qn^{-2}.
\end{equation}
The conditional variance is under $\mu_{A,(e),1}$ and the exterior
norm is under the actual marginal of $(A,e)$.
\end{lemma}

\begin{proof}
First estimate the variance under $\pi_A\otimes\pi_A$. Write
$R_A=R_{A,1}+\sum_{d=2}^D H_{A,d}$ as in
\eqref{audit:reference-logit-exact-expansion}. For its variance, use
\eqref{audit:dependency:s79:reference-logit-remainder-moments}: its mixed norm
$\left\|\|R_A\|_{L^2(\pi_A\otimes\pi_A)}\right\|_{2q}$
is $O(n^{-1})$, so
\begin{equation}\label{audit:dependency:s79:reference-variance-remainder}
 \|\Var_{\pi_A\otimes\pi_A}(R_A)\|_q\le C_qn^{-2}.
\end{equation}
The exact variance decomposition is
\begin{equation}\label{audit:dependency:s79:reference-variance-decomposition}
 \begin{aligned}
 \Var_{\pi_A\otimes\pi_A}(L_{12})
 &=\E_{\pi_A\otimes\pi_A}L_{A,0}^2
       +2\E_{\pi_A\otimes\pi_A}(L_{A,0}R_{A,1})\\
 &\quad+2\sum_{d=2}^D\E_{\pi_A\otimes\pi_A}(L_{A,0}H_{A,d})
       +\Var_{\pi_A\otimes\pi_A}(R_A).
 \end{aligned}
\end{equation}
For the leading variance and higher-site pairings, row independence,
\eqref{eq:cov-profiles-a}, and \eqref{eq:cov-profiles-b} give
\begin{equation}\label{audit:dependency:s79:reference-variance-leading-bounds}
\begin{aligned}
 \E_{\pi_A\otimes\pi_A} L_{A,0}^2
 &=\frac{8t(P_A)^2}{n^2}\one^{\mathsf T}\Sigma_A\one
   +\frac{36\tau(P_A)^2}{n^2}\tr\Sigma_A^2
 =\frac{w_2(P_A)}n+O_{L^q}(n^{-2}),\\
 |\E_{\pi_A\otimes\pi_A}(L_{A,0}H_{A,d})|
 &\le\frac Cn\sum_{\alpha\in\{1,2,12\}}
    \left|\E_{\pi_A\otimes\pi_A}\left[H_{A,d}\sum_i\psi_{i,\alpha}\right]\right|
 \le C_dn^{-2}.
\end{aligned}
\end{equation}
The first line uses $m=n-2$ and the vanishing mixed-type pairings;
the second uses \eqref{eq:two-row-higher-site-signed-pairing}.
For the single-site pairing, \eqref{eq:two-row-site-pairing} gives
\[
 \begin{aligned}
 \E_{\pi_A\otimes\pi_A}(L_{A,0}R_{A,1})
 &=\frac{4t(P_A)}n(\Sigma_A\one)^{\mathsf T}e_A^{(1)}\\
 &\quad+\frac{6\tau(P_A)}n
                \sum_{i\in V(A)} e_{A,i}^{(12)}\left(\sum_{j\in V(A)}\Sigma_{A,ij}^2\right).
 \end{aligned}
\]
To estimate this pairing, retain the signed coefficient sums:
\begin{equation}\label{audit:dependency:s79:reference-variance-error-pairings}
\begin{aligned}
 &(\Sigma_A\one)^{\mathsf T}e_A^{(1)} =v_0(P_A)\sum_i e_{A,i}^{(1)}
       +(\Sigma_A\one-v_0(P_A)\one)^{\mathsf T}e_A^{(1)},\\
 &\sum_{i\in V(A)} e_{A,i}^{(12)}\left(\sum_{j\in V(A)}\Sigma_{A,ij}^2\right) =\nu(P_A)^2\sum_i e_{A,i}^{(12)}
\\
 &\quad+((\sum_j\Sigma_{A,ij}^2-\nu(P_A)^2)_i)^{\mathsf T}e_A^{(12)},\\
 &\|(\Sigma_A\one)^{\mathsf T}e_A^{(1)}\|_q\\
 &\le C\left\|\sum_i e_{A,i}^{(1)}\right\|_q
   +\bigl\|\|\Sigma_A\one-v_0(P_A)\one\|_2\bigr\|_{2q}
       \bigl\|\|e_A^{(1)}\|_2\bigr\|_{2q}\le C_q/n,\\
 &\left\|\sum_{i\in V(A)} e_{A,i}^{(12)}\left(\sum_{j\in V(A)}\Sigma_{A,ij}^2\right)\right\|_q\\
 &\le C\left\|\sum_i e_{A,i}^{(12)}\right\|_q\\
 &\quad+\bigl\|\|(\sum_j\Sigma_{A,ij}^2-\nu(P_A)^2)_i\|_2\bigr\|_{2q}
       \bigl\|\|e_A^{(12)}\|_2\bigr\|_{2q}\le C_q/n.
\end{aligned}
\end{equation}
Equations~\eqref{eq:revised-logit-coefficient-errors},
\eqref{eq:cov-profiles-a}, and \eqref{eq:covariance-square-row-sums}
justify \eqref{audit:dependency:s79:reference-variance-error-pairings};
the additional $1/n$ gives
\[
 \left\|\E_{\pi_A\otimes\pi_A}(L_{A,0}R_{A,1})\right\|_q\le C_qn^{-2}.
\]
Substituting this bound, \eqref{audit:dependency:s79:reference-variance-remainder},
and \eqref{audit:dependency:s79:reference-variance-leading-bounds} into
\eqref{audit:dependency:s79:reference-variance-decomposition} proves
\begin{equation}\label{audit:reference-logit-variance}
 \left\|\Var_{\pi_A\otimes\pi_A}(L_{12})-w_2(P_A)/n\right\|_q
 \le C_qn^{-2}.
\end{equation}
To transfer \eqref{audit:reference-logit-variance}, use
\eqref{eq:two-row-covariance-derivative} and
\eqref{ep:row:polynomial-bound} on $\mathcal A_n$:
\[
 \Var_{\mu_{A,(e),1}}(L_{12})-\Var_{\pi_A\otimes\pi_A}(L_{12})
 =\int_0^1\cum_{\mu_{A,(e),s}}(L_{12},L_{12},V_{A,e})\,\dd s
 =O(n^{-2}).
\]
Together with \eqref{audit:reference-logit-variance} and the
exceptional-exterior bound, this proves \eqref{eq:variance} under the
actual law, uniformly in $e,s$.
\end{proof}

The reference expansion identifies the pairing, which interpolation
then transfers to the true conditional law.
\begin{lemma}
\label{lem:common-neighbor-pairings}
Assume \hyperref[ass:PC]{\textup{(PC)}}.
Under $\mu_{A,(e),1}$ in \eqref{v8:setup:two-row-interpolation-density},
for every fixed $q\ge1$,
\begin{equation}\label{eq:FB}
 \left\|\E[(L_{12}-\ell_n(P_A))B_{12}(P_A)\mid A,e]
       -6\tau(P_A)\nu(P_A)^2\right\|_q\le C_qn^{-1}.
\end{equation}
\end{lemma}

\begin{proof}
First work under $\pi_A\otimes\pi_A$ and recenter:
\[
 \begin{aligned}
 B_{12}(P_A)
 &=B_{A,0}+\sum_i(\mu_{A,i}-P_A)(\psi_{i,1}+\psi_{i,2})
                       +\|\mu_A-P_A\one\|_2^2,\\
 \E_{\pi_A\otimes\pi_A} B_{12}(P_A)&=\|\mu_A-P_A\one\|_2^2.
 \end{aligned}
\]
Lemmas~\ref{lem:weighted-linear} and~\ref{lem:row-mean-density} give,
for every fixed $q_1,q_2\ge1$,
\begin{equation}\label{audit:dependency:s79:common-neighbor-recentering-moment}
 \left\|\|B_{12}(P_A)-B_{A,0}\|_{L^{q_2}(\pi_A\otimes\pi_A)}\right\|_{q_1}
 \le C_{q_2}\bigl\|\|\mu_A-P_A\one\|_2\bigr\|_{q_1}
       +\bigl\|\|\mu_A-P_A\one\|_2\bigr\|_{2q_1}^2\le C_{q_1,q_2}.
\end{equation}
Equation~\eqref{eq:two-row-site-pairing} gives the leading contribution
\begin{equation}\label{audit:dependency:s79:linear-common-neighbor-leading-pairing}
 \begin{aligned}
 \E_{\pi_A\otimes\pi_A}(L_{A,0}B_{12}(P_A))
 &=\frac{4t(P_A)}n(\mu_A-P_A\one)^{\mathsf T}\Sigma_A\one
   +\frac{6\tau(P_A)}n\tr\Sigma_A^2\\
 &=6\tau(P_A)\nu(P_A)^2+O_{L^q}(n^{-1}).
 \end{aligned}
\end{equation}
because $\tr\Sigma_A^2=m\nu(P_A)^2+O_{L^q}(1)$, $m=n-2$, and
\[
\begin{aligned}
 &(\mu_A-P_A\one)^{\mathsf T}\Sigma_A\one\\
 &\quad=v_0(P_A)\sum_i(\mu_{A,i}-P_A)+
           (\mu_A-P_A\one)^{\mathsf T}(\Sigma_A\one-v_0(P_A)\one),\\
 &\|(\mu_A-P_A\one)^{\mathsf T}\Sigma_A\one\|_q\\
 &\quad\le C\left\|\sum_i(\mu_{A,i}-P_A)\right\|_q
 +\bigl\|\|\mu_A-P_A\one\|_2\bigr\|_{2q}
     \bigl\|\|\Sigma_A\one-v_0(P_A)\one\|_2\bigr\|_{2q}\le C_q.
\end{aligned}
\]
By \eqref{audit:dependency:s79:reference-logit-constant-bias} and
\eqref{eq:mean-profiles-a}, the constant error in
\eqref{audit:reference-logit-exact-expansion} contributes
\[
 \begin{aligned}
 &\|\{L_{12}(\mu_A,\mu_A;A)-\ell_n(P_A)\}\E_{\pi_A\otimes\pi_A} B_{12}(P_A)\|_q\\
 &\quad\le\|L_{12}(\mu_A,\mu_A;A)-\ell_n(P_A)\|_{2q}
              \bigl\|\|\mu_A-P_A\one\|_2\bigr\|_{4q}^2
 \le C_q/n.
 \end{aligned}
\]
For the single-site error, \eqref{audit:dependency:s79:reference-variance-error-pairings},
\eqref{eq:revised-logit-coefficient-errors}, \eqref{eq:mean-profiles-a}, and
\eqref{eq:reference-covariance-operator} give
\[
\begin{aligned}
 \E_{\pi_A\otimes\pi_A}(R_{A,1}B_{12}(P_A))&=\sum_{i\in V(A)} e_{A,i}^{(12)}\left(\sum_{j\in V(A)}\Sigma_{A,ij}^2\right)
                   +2(e_A^{(1)})^{\mathsf T}\Sigma_A(\mu_A-P_A\one),\\
 \|\E_{\pi_A\otimes\pi_A}(R_{A,1}B_{12}(P_A))\|_q
 &\le C_q/n+C\bigl\|\|e_A^{(1)}\|_2\bigr\|_{2q}
             \bigl\|\|\mu_A-P_A\one\|_2\bigr\|_{2q}\le C_q/n.
\end{aligned}
\]
For the higher-site errors, separate the $B_{A,0}$ pairing from the
recentering term. For each $d\ge2$,
\eqref{eq:two-row-higher-site-signed-pairing},
\eqref{audit:dependency:s79:reference-logit-remainder-moments}, and
\eqref{audit:dependency:s79:common-neighbor-recentering-moment}, with $q_1=2q$ and $q_2=2$, give
\[
\begin{aligned}
 |\E_{\pi_A\otimes\pi_A}(H_{A,d}B_{A,0})|&\le C_dn^{-d+\lfloor(d+1)/2\rfloor}\le C_d/n,\\
 \|\E_{\pi_A\otimes\pi_A}[H_{A,d}(B_{12}(P_A)-B_{A,0})]\|_q
 &\le\left\|\|H_{A,d}\|_{L^2(\pi_A\otimes\pi_A)}\right\|_{2q}
       \left\|\|B_{12}(P_A)-B_{A,0}\|_{L^2(\pi_A\otimes\pi_A)}\right\|_{2q}
  \\
  &\le C_{q,d}n^{-d/2}\le C_{q,d}/n.
\end{aligned}
\]
Combining these error bounds with
\eqref{audit:dependency:s79:linear-common-neighbor-leading-pairing} in
\eqref{audit:reference-logit-exact-expansion} proves
the reference-law counterpart of \eqref{eq:FB}.

To transfer this reference estimate, use
\eqref{eq:two-row-expectation-derivative} and
\eqref{eq:two-row-product-covariance}:
\begin{equation}\label{audit:dependency:s79:linear-pairing-transfer-derivative}
\begin{aligned}
 &\frac{\dd}{\dd s}
       \E_{\mu_{A,(e),s}}[\{L_{12}-\ell_n(P_A)\}B_{12}(P_A)/n]\\
 &=\cum_{\mu_{A,(e),s}}(L_{12},B_{12}(P_A)/n,V_{A,e})\\
 &\quad+\E_{\mu_{A,(e),s}}[L_{12}-\ell_n(P_A)]\,
                       \Cov_{\mu_{A,(e),s}}(B_{12}(P_A)/n,V_{A,e})\\
 &\quad+\E_{\mu_{A,(e),s}}[B_{12}(P_A)/n]\,\Cov_{\mu_{A,(e),s}}(L_{12},V_{A,e}).
\end{aligned}
\end{equation}
On $\mathcal A_n$, Theorem~\ref{ep:row:sharp} bounds the cumulant by
$Cn^{-2}$ and both covariances by $Cn^{-1}$, whereas
\eqref{eq:revised-logit-interpolated-bias} and
Lemma~\ref{lem:common-neighbor-moments} give
\[
 \|\E_{\mu_{A,(e),s}}[L_{12}-\ell_n(P_A)]\|_q
       +\|\E_{\mu_{A,(e),s}}[B_{12}(P_A)/n]\|_q\le C_q/n
\]
uniformly in $s$. Substitution into
\eqref{audit:dependency:s79:linear-pairing-transfer-derivative} bounds the
derivative by $O_{L^q}(n^{-2})$; integration gives
\[
\begin{aligned}
 &\left\|\E_{\mu_{A,(e),1}}[\{L_{12}-\ell_n(P_A)\}B_{12}(P_A)]\right.\\
 &\quad\left.
       -\E_{\pi_A\otimes\pi_A}[\{L_{12}-\ell_n(P_A)\}B_{12}(P_A)]
                                                      \right\|_q\\
 &\quad\le n\int_0^1 C_qn^{-2}\,\dd s=C_q/n.
\end{aligned}
\]
Combining this transfer bound with the reference estimate and the
exceptional-exterior bound proves \eqref{eq:FB}.
\end{proof}

We also need the quadratic pairing to be $O_{L^q}(n^{-1})$.
The moment--cumulant expansion gives this bound directly.
The centering field $\ell_n$ from \eqref{v8:setup:physical-edge-logit} is
\begin{equation}\label{v8:setup:two-row-centering-field}
 \ell_n(u)=h_n+2\sum_{r=1}^k b_{r,n}e_ru^{e_r-1}.
\end{equation}

\begin{lemma}
\label{lem:quadratic-common-neighbor-pairing}
Assume \hyperref[ass:PC]{\textup{(PC)}}.
Under $\mu_{A,(e),1}$ in \eqref{v8:setup:two-row-interpolation-density},
for every fixed $q\ge1$,
\begin{equation}\label{eq:F2B}
 \left\|\E[(L_{12}-\ell_n(P_A))^2B_{12}(P_A)\mid A,e]\right\|_q
       \le C_qn^{-1}.
\end{equation}
\end{lemma}

\begin{proof}
Under $\mu_{A,(e),1}$, apply the moment--cumulant identity:
\begin{equation}\label{audit:dependency:s79:quadratic-pairing-identity}
\begin{aligned}
 &\E_{\mu_{A,(e),1}}[\{L_{12}-\ell_n(P_A)\}^2B_{12}(P_A)]\\
 &=\cum_{\mu_{A,(e),1}}(L_{12},L_{12},B_{12}(P_A))
       +\E_{\mu_{A,(e),1}} B_{12}(P_A)\,\Var_{\mu_{A,(e),1}}(L_{12})\\
 &\quad+2\{\E_{\mu_{A,(e),1}}[L_{12}]-\ell_n(P_A)\}
                    \Cov_{\mu_{A,(e),1}}(L_{12},B_{12}(P_A))\\
 &\quad+\{\E_{\mu_{A,(e),1}}[L_{12}]-\ell_n(P_A)\}^2\E_{\mu_{A,(e),1}} B_{12}(P_A).
\end{aligned}
\end{equation}
On good exteriors, \eqref{ep:row:polynomial-bound}, applied to $L_{12}$
and $B_{12}(P_A)/n$, gives
\[
\begin{aligned}
 &|\cum_{\mu_{A,(e),1}}(L_{12},L_{12},B_{12}(P_A))|+\Var_{\mu_{A,(e),1}}(L_{12})\le C/n,\\
 &|\Cov_{\mu_{A,(e),1}}(L_{12},B_{12}(P_A))|\le C.
\end{aligned}
\]
Lemmas~\ref{lem:conditional-logit-bias} and
\ref{lem:common-neighbor-moments} give
\[
 \|\E_{\mu_{A,(e),1}}[L_{12}]-\ell_n(P_A)\|_q\le C_q/n,
 \qquad \|\E_{\mu_{A,(e),1}} B_{12}(P_A)\|_q\le C_q
\]
at every fixed order by \eqref{eq:bias} and
\eqref{eq:revised-common-neighbor-moments}. Substitution into
\eqref{audit:dependency:s79:quadratic-pairing-identity} and H\"older's inequality yield
\[
 \|\E_{\mu_{A,(e),1}}[\{L_{12}-\ell_n(P_A)\}^2B_{12}(P_A)]\|_q
 \le C_q(n^{-1}+n^{-1}+n^{-1}+n^{-2})\le C_q/n.
\]
The exceptional-exterior bound then proves \eqref{eq:F2B} under the
full marginal law.
\end{proof}

\section{Expectation estimates for centered subgraph counts}\label{sec:supports}
For a fixed simple graph $H$, write
\begin{equation}\label{v8:setup:signed-centered-count}
 \inj(H,Y)=\sum_{\varphi:V(H)\hookrightarrow[n]}
           \prod_{uv\in E(H)}Y_{\varphi(u)\varphi(v)},
 \qquad Y=X-P(J_n-I_n).
\end{equation}
We group centered products by subgraph support, as in the graph-factor
method of Janson~\cite{Janson1994}, with the empirical center $P$
and the critical Gibbs law. Degree-square and triangle terms are
identified explicitly; conditional row estimates and two-row pairings
control the remaining supports. Absolute values are taken after expectation,
and isolated vertices retain their exact falling-factorial factors.
Unless a conditional or reference measure is displayed, expectations and
random-variable $L^q$ norms are under $\Pp_n$. For quantities depending
only on an exterior $A$, write $\E_A$ for its induced marginal expectation;
their outer moment norms use the same marginal.
For Lipschitz tests on $[0,1]$, write
\begin{equation}\label{v8:setup:signed-test-norms}
 \|\phi\|_\infty=\sup_{x\in[0,1]}|\phi(x)|,\qquad
 \Lip\phi=\sup_{x\ne y}\frac{|\phi(x)-\phi(y)|}{|x-y|}.
\end{equation}

\subsection{Signed degree-square and triangle approximations}
Define the degree fluctuation by
\begin{equation}\label{v8:setup:signed-degree-square}
 z_i=D_i-(n-1)P,\qquad \sum_i z_i=0.
\end{equation}

\begin{lemma}\label{lem:low-supports}
Assume \hyperref[ass:PC]{\textup{(PC)}}.
Let $v_0(x)$ be the fixed profile in \eqref{eq:logit-profiles}, and let
$\phi:[0,1]\to\mathbb R$ be any Lipschitz function. Then
\begin{equation}\label{eq:V-weak}
 |\E[\phi(P)\{\frac1{n^2}\sum_{i=1}^n z_i^2-v_0(P)\}]|
 \le Cn^{-1}\|\phi\|_\infty+Cn^{-3/2}\Lip\phi.
\end{equation}
In particular, $|\E\frac1{n^2}\sum_{i=1}^n z_i^2-\E v_0(P)|\le C/n$.
This is a signed estimate, not an $L^1$ approximation.
\end{lemma}

\begin{proof}
Reduce the degree square to a two-row variance. Since $\sum_i z_i=0$,
\begin{equation}\label{audit:dependency:s79:degree-square-exchangeability}
 \sum_{i,j}(D_i-D_j)^2=2n\sum_i z_i^2,
 \qquad
 \E\left[\phi(P)n^{-2}\sum_i z_i^2\right]
 =\frac{n-1}{2n^2}\E[\phi(P)(D_1-D_2)^2].
\end{equation}
For $A=X_{\{3,\ldots,n\}}$ and $e=X_{12}$, root interchange gives
$\E_{\mu_{A,(e),s}}(D_1-D_2)=0$, since $D_1-D_2=\sum_i(x_i-y_i)$. Thus
\[
 \begin{aligned}
 \frac{\dd}{\dd s}\E_{\mu_{A,(e),s}}(D_1-D_2)^2
 &=\Cov_{\mu_{A,(e),s}}((D_1-D_2)^2,V_{A,e})\\
 &=n^2\cum_{\mu_{A,(e),s}}((D_1-D_2)/n,(D_1-D_2)/n,V_{A,e})=O(1)
 \end{aligned}
\]
on good exteriors by Theorem~\ref{ep:row:sharp}.
At $s=0$, row independence and Lemma~\ref{lem:row-total-covariance}
give $\E_{\pi_A\otimes\pi_A}(D_1-D_2)^2=2\one^{\mathsf T}\Sigma_A\one
=2(n-2)v_0(P_A)+O_{L^q}(1)$. Integration yields
\begin{equation}\label{audit:dependency:s79:degree-square-conditional-variance}
 \|\E[(D_1-D_2)^2\mid A,e]-2(n-2)v_0(P_A)\|_q\le C_q.
\end{equation}
The envelope $(D_1-D_2)^2\le n^2$ handles bad exteriors.
Corollary~\ref{lem:degree-deletion} gives
$\|(D_1-D_2)\|_4\le C\sqrt n$ and $\|P-P_A\|_2\le Cn^{-3/2}$, hence
\[
\begin{aligned}
 &\left|\frac{n-1}{2n^2}\E[\phi(P)(D_1-D_2)^2]
       -\frac{(n-1)(n-2)}{n^2}\E[\phi(P_A)v_0(P_A)]\right|\\
 &\quad\le\frac Cn\Lip\phi\,\|P-P_A\|_2\|(D_1-D_2)\|_4^2
             +\frac Cn\|\phi\|_\infty
 \le Cn^{-3/2}\Lip\phi+Cn^{-1}\|\phi\|_\infty.
\end{aligned}
\]
Boundedness of $v_0(x),v_0'(x)$ and
$|(n-1)(n-2)/n^2-1|\le3/n$ imply
\[
\begin{aligned}
 &\left|\frac{(n-1)(n-2)}{n^2}\E[\phi(P_A)v_0(P_A)]
                      -\E[\phi(P)v_0(P)]\right|\\
 &\quad\le Cn^{-1}\|\phi\|_\infty
       +C(\|\phi\|_\infty+\Lip\phi)\|P_A-P\|_1.
\end{aligned}
\]
Combining these bounds with
\eqref{audit:dependency:s79:degree-square-exchangeability},
\eqref{audit:dependency:s79:degree-square-conditional-variance}, and
\eqref{row:eq:degree-deletion} proves \eqref{eq:V-weak}.
\end{proof}

The centered triangle statistic is
\begin{equation}\label{v8:setup:signed-triangle-statistic}
 \tr Y^3=\sum_{i,j,k\ \mathrm{distinct}}
            (X_{ij}-P)(X_{jk}-P)(X_{ki}-P).
\end{equation}

\begin{lemma}\label{lem:triangle-profile}
Assume \hyperref[ass:PC]{\textup{(PC)}}.
For every Lipschitz $\phi:[0,1]\to\mathbb R$,
\begin{equation}\label{eq:T-weak}
 \left|\E\left[\phi(P)
       \{n^{-2}\tr Y^3-6\tau(P)\nu(P)^3\}\right]\right|
 \le Cn^{-1}\|\phi\|_\infty+Cn^{-3/2}\Lip\phi.
\end{equation}
In particular,
$|n^{-2}\E\tr Y^3-6\E[\tau(P)\nu(P)^3]|\le C/n$.
This identifies the leading signed term, whereas
\eqref{eq:triangle-moment} of Corollary~\ref{cor:matrix-moments}
controls its absolute moments. The constant
is independent of $n$ and $\phi$.
\end{lemma}

\begin{proof}
For an edge $e=\{a,b\}$ with two-root exterior $A$, write
$B_e(u)=\sum_{k\notin\{a,b\}}(X_{ak}-u)(X_{bk}-u)$. Then
\[
 \tr Y^3=2\sum_e(X_e-P)B_e(P).
\]
\proofstep{1. Restore the exterior center.}
Equations~\eqref{row:eq:degree-deletion}, \eqref{eq:degree-vector-moment},
and \eqref{eq:B-bounds} give
\[
 \|P-P_A\|_q\le C_qn^{-3/2},\qquad
 \|z_a\|_q+\|z_b\|_q+\|B_e(P_A)\|_q\le C_q\sqrt n.
\]
Since $B_e'(P)=-z_a-z_b+2(X_e-P)$ and $B_e''=2(n-2)$,
\[
\begin{aligned}
 \|B_e'(P_A)\|_q&\le C_q\sqrt n,\\
 B_e(P)-B_e(P_A)&=(P-P_A)B_e'(P_A)+(n-2)(P-P_A)^2,\\
 \|B_e(P)-B_e(P_A)\|_q&\le C_qn^{-1}.
\end{aligned}
\]
Thus the change of center costs
\begin{equation}\label{audit:dependency:s79:triangle-center-restoration}
\begin{aligned}
 &\left|\E[\phi(P)(X_e-P)B_e(P)]-\E[\phi(P)(X_e-P_A)B_e(P_A)]\right|\\
 &\quad\le\|\phi\|_\infty
       \{\|B_e(P)-B_e(P_A)\|_1+\|P-P_A\|_2\|B_e(P_A)\|_2\}
 \le Cn^{-1}\|\phi\|_\infty.
\end{aligned}
\end{equation}

\proofstep{2. Replace the edge and freeze the test.}
Conditioning on $X_{-e}$ gives the exact Bernoulli identity
\[
\begin{aligned}
 &\E[(X_e-\sig(L_e))\phi(P)\mid X_{-e}]\\
 &\quad=\sig(L_e)(1-\sig(L_e))
       \{\phi((E-X_e+1)/N)-\phi((E-X_e)/N)\}.
\end{aligned}
\]
Because $B_e(P_A)$ is $X_{-e}$-measurable,
\begin{equation}\label{audit:dependency:s79:triangle-edge-replacement}
 |\E[(X_e-\sig(L_e))\phi(P)B_e(P_A)]|
 \le\frac{\Lip\phi}{4N}\|B_e(P_A)\|_1\le Cn^{-3/2}\Lip\phi.
\end{equation}
Equations~\eqref{eq:local-polynomial-moments}, \eqref{eq:bias}, and
\eqref{audit:scalar-defect-exterior} give
\begin{equation}\label{audit:dependency:s79:triangle-logit-moments}
\begin{gathered}
 \|L_e-\ell_n(P_A)\|_q\le C_qn^{-1/2},\qquad
 \|\sig(\ell_n(P_A))-P_A\|_q\le C_q/n,\\
 \|\sig(L_e)-P_A\|_q\le C_qn^{-1/2}.
\end{gathered}
\end{equation}
Thus freezing the test costs at most
\begin{equation}\label{audit:dependency:s79:triangle-test-freezing}
\begin{aligned}
 &|\E[(\phi(P)-\phi(P_A))(\sig(L_e)-P_A)B_e(P_A)]|\\
 &\quad\le\Lip\phi\,\|P-P_A\|_2\|\sig(L_e)-P_A\|_4\|B_e(P_A)\|_4
 \le Cn^{-3/2}\Lip\phi.
\end{aligned}
\end{equation}

\proofstep{3. Evaluate the conditional pairing.}
Taylor's formula at $\ell_n(P_A)$ gives
\[
\begin{aligned}
 \sig(L_e)-P_A
 &=\sig(\ell_n(P_A))-P_A
    +\sig'(\ell_n(P_A))\{L_e-\ell_n(P_A)\}\\
 &\quad+\tfrac12\sig''(\ell_n(P_A))\{L_e-\ell_n(P_A)\}^2
    +O(|L_e-\ell_n(P_A)|^3).
\end{aligned}
\]
Equations~\eqref{eq:B-bounds}, \eqref{eq:FB}, and \eqref{eq:F2B},
together with \eqref{audit:dependency:s79:triangle-logit-moments}, give
\begin{equation}\label{audit:dependency:s79:triangle-conditional-pairing-inputs}
\begin{aligned}
 \|\E[B_e(P_A)\mid A,X_e]\|_q&\le C_q,\\
 \E[\{L_e-\ell_n(P_A)\}B_e(P_A)\mid A,X_e]&=6\tau(P_A)\nu(P_A)^2+O_{L^q}(n^{-1}),\\
 \|\E[\{L_e-\ell_n(P_A)\}^2B_e(P_A)\mid A,X_e]\|_q&\le C_q/n,\\
 \E_{\Pp_n}|L_e-\ell_n(P_A)|^3|B_e(P_A)|&\le\|L_e-\ell_n(P_A)\|_6^3\|B_e(P_A)\|_2\le C/n.
\end{aligned}
\end{equation}
Moreover, $\sig'(u)=\nu(\sig(u))$ and
\eqref{audit:dependency:s79:triangle-logit-moments} give
\begin{equation}\label{audit:dependency:s79:triangle-logistic-derivative-comparison}
\begin{aligned}
 |\sig'(\ell_n(P_A))-\nu(P_A)|
 &=|\sig(\ell_n(P_A))-P_A|\,|1-\sig(\ell_n(P_A))-P_A|\\
 &\le|\sig(\ell_n(P_A))-P_A|,\\
 \|\sig'(\ell_n(P_A))-\nu(P_A)\|_q&\le C_q/n.
\end{aligned}
\end{equation}
Equations~\eqref{audit:dependency:s79:triangle-conditional-pairing-inputs} and
\eqref{audit:dependency:s79:triangle-logit-moments} bound the constant and
quadratic terms and the cubic Taylor remainder by
$O(n^{-1}\|\phi\|_\infty)$. For the linear term,
\eqref{eq:FB} and
\eqref{audit:dependency:s79:triangle-logistic-derivative-comparison} give
\begin{equation}\label{audit:dependency:s79:triangle-conditional-pairing}
 \begin{aligned}
 &\E[\phi(P_A)\sig'(\ell_n(P_A))
          \{L_e-\ell_n(P_A)\}B_e(P_A)]\\
 &\quad=\E[6\phi(P_A)\tau(P_A)\nu(P_A)^3]
                 +O(n^{-1}\|\phi\|_\infty).
 \end{aligned}
\end{equation}
The probability bound \eqref{eq:local-good-exteriors} and the polynomial
envelopes account for bad exteriors in these moment estimates.

\proofstep{4. Sum over edges.}
Combining \eqref{audit:dependency:s79:triangle-center-restoration},
\eqref{audit:dependency:s79:triangle-edge-replacement},
\eqref{audit:dependency:s79:triangle-test-freezing}, and
\eqref{audit:dependency:s79:triangle-conditional-pairing} gives, for each $e$,
\[
\begin{aligned}
 &\left|\E[\phi(P)(X_e-P)B_e(P)]
            -\E[6\phi(P_A)\tau(P_A)\nu(P_A)^3]\right|\\
 &\quad\le C\{n^{-1}\|\phi\|_\infty+n^{-3/2}\Lip\phi\}.
\end{aligned}
\]
Replacing $P_A$ by $P$ in the tested smooth profile costs at most
$Cn^{-3/2}(\|\phi\|_\infty+\Lip\phi)$.
Multiply by $2/n^2$, sum over the $N$ edges, and use
$2N/n^2=1-1/n$ to obtain \eqref{eq:T-weak}.
\end{proof}

\subsection{Rooted reduction and branching graphs}
Fix a simple graph $H$, a designated vertex $u$ of degree $d\ge2$,
and $v=v(H)$. With $A=X_{[n]\setminus\{o\}}$, define
\begin{equation}\label{v8:setup:rooted-centered-count}
 T_o(H,u;t)=\sum_{\substack{\varphi:V(H)\hookrightarrow[n]\\\varphi(u)=o}}
            \prod_{ab\in E(H)}(X_{\varphi(a)\varphi(b)}-t).
\end{equation}

\begin{lemma}
\label{lem:centered-root-moments}
Assume \hyperref[ass:PC]{\textup{(PC)}}.
For the fixed simple graph $H$, vertex $u$ of degree $d\ge2$, and count
\eqref{v8:setup:rooted-centered-count}, uniformly in $o$,
\begin{equation}\label{eq:high-support-common-moments}
 \|T_o(H,u;P_A)\|_2\le C_Hn^{v-1-d/2},\qquad
 \|\partial_tT_o(H,u;P_A)\|_2\le C_Hn^{v-1-(d-1)/2}.
\end{equation}
For every Lipschitz $\phi:[0,1]\to\mathbb R$, these moments give
\begin{equation}\label{eq:high-support-common-transfer}
 \begin{aligned}
 |\E[\phi(P)T_o(H,u;P)]|
 &\le|\E[\phi(P_A)T_o(H,u;P_A)]|\\
 &\quad+C_H\{n^{v-3}\|\phi\|_\infty+n^{v-7/2}\Lip\phi\}.
 \end{aligned}
\end{equation}
All norms are taken under the original graph law, with exterior
quantities integrated against its induced marginal. The graph $H$ need not be connected.
\end{lemma}

\begin{proof}
For fixed $A=X_{[n]\setminus\{o\}}$, the conditional law of
$(X_{oi})_{i\in V(A)}$ under $\Pp_n$ is $\pi_A$.
Fixing the $d$ neighbor images leaves at most $n^{v-1-d}$ completions, so
$n^{-(v-1)}T_o(H,u;P_A)$ is a degree-$d$ row polynomial with
coefficients bounded by $n^{-d}$ and equal to zero on repeated neighbor images.
Apply Lemmas~\ref{lem:row-center-moments} and~\ref{lem:row-mean-density}
under $\pi_A$, then average over the actual marginal of $A$ to obtain
\[
 \|\one_{\{A\in\mathcal A_n\}}T_o(H,u;P_A)\|_{L^2(\Pp_n)}
 \le C_Hn^{v-1-d/2}
       \{\E_A(1+\|\mu_A-P_A\one\|_2)^{2d}\}^{1/2}
 \le C_Hn^{v-1-d/2}.
\]
The complement contributes at most $n^{v-1}Ce^{-cn/2}$.
Edge deletion retains all $v$ vertices and lowers the root degree by
at most one, so the same argument applied to
\[
 \partial_tT_o(H,u;t)=-\sum_{f\in E(H)}T_o(H-f,u;t)
\]
proves both bounds in \eqref{eq:high-support-common-moments}.
For the change of center, use $|\partial_t^2T_o(H,u;t)|\le C_Hn^{v-1}$
for $t\in[0,1]$ and Taylor's formula:
\[
 \begin{aligned}
 &T_o(H,u;P)-T_o(H,u;P_A)\\
 &\quad=(P-P_A)\partial_tT_o(H,u;P_A)\\
 &\qquad+(P-P_A)^2\int_0^1(1-s)\partial_t^2T_o(H,u;P_A+s(P-P_A))\,\dd s.
 \end{aligned}
\]
By \eqref{row:eq:degree-deletion},
\eqref{eq:high-support-common-moments}, and H\"older's inequality,
\begin{equation}\label{audit:rooted-center-restoration-errors}
\begin{aligned}
 \E|T_o(H,u;P)-T_o(H,u;P_A)|
 &\le C_H\{n^{-3/2}n^{v-1-(d-1)/2}+n^{v-1}n^{-3}\}
 \le C_Hn^{v-3},\\
 \E|(\phi(P)-\phi(P_A))T_o(H,u;P_A)|
 &\le C_H\Lip\phi\,n^{-3/2}n^{v-1-d/2}
 \le C_Hn^{v-7/2}\Lip\phi.
\end{aligned}
\end{equation}
Here only $d\ge2$ is used. Add and subtract
$\phi(P)T_o(H,u;P_A)$ and apply
\eqref{audit:rooted-center-restoration-errors} to obtain
\eqref{eq:high-support-common-transfer}.
\end{proof}

A root of degree three requires a signed bound for the product of
three degree fluctuations. We establish this bound before treating
graphs with a branching vertex.
Fix $L<\infty$ independently of $n$. For fixed distinct vertices $i,j,k$
and arbitrary deterministic edge-label sequences $\mathcal E_0,\mathcal E_1$
with $|\mathcal E_0|+|\mathcal E_1|\le L$ (overlaps and repetitions allowed),
define, for a graph $A$ containing their endpoints,
\begin{equation}\label{v8:setup:decorated-degree-products}
 B_0(A;t)=\prod_{e\in\mathcal E_0}A_e
                   \prod_{e\in\mathcal E_1}(A_e-t).
\end{equation}
For a graph $A$ obtained by a fixed vertex deletion, put
\begin{equation}\label{v8:setup:exterior-centered-degrees}
 P_A=\binom{|A|}{2}^{-1}\sum_{a<b}A_{ab},\quad
 d_{A,i}=\sum_{j\ne i}A_{ij},\quad
 z_i^A=d_{A,i}-(|A|-1)P_A .
\end{equation}

\begin{lemma}\label{lem:three-degrees}
Assume \hyperref[ass:PC]{\textup{(PC)}}.
For fixed distinct $i,j,k$, the decoration $B_0$ in
\eqref{v8:setup:decorated-degree-products}, and every Lipschitz
$\phi:[0,1]\to\mathbb R$,
\begin{equation}\label{eq:three-degrees}
 |\E\phi(P)z_iz_jz_kB_0(X;P)|
 \le C\bigl(n\|\phi\|_\infty+\Lip\phi\bigr).
\end{equation}
For $A$ obtained by deleting $h$ fixed vertices and distinct retained vertices
$i,j,k$, with all decoration endpoints in $V(A)$, the same bound holds
under the actual marginal exterior law,
with edge factors and degrees centered at $P_A$:
\begin{equation}\label{audit:decorated-exterior-degree-bound}
 \left|\E\phi(P_A)z_i^Az_j^Az_k^A B_0(A;P_A)\right|
 \le C\bigl(n\|\phi\|_\infty+\Lip\phi\bigr)
\end{equation}
In both bounds, $C$ may depend on $L$, $h$ and the fixed data in
\textup{(PC)}, but not on $n$ or the vertex and edge labels;
take $h=0$ for the full graph.
\end{lemma}

\begin{proof}
Freeze the test and decoration on a smaller exterior.
Use $A=X_{[n]\setminus\mathcal D}$ with $|\mathcal D|=h$ fixed;
$\mathcal D=\varnothing$ gives the full graph.
Put $m=|V(A)|=n-|\mathcal D|$ and $N_A=\binom m2$.
Set
\[
 T=\{i,j,k\}\cup
 \{v\in V(A):v\in e\text{ for some }e\in\mathcal E_0\cup\mathcal E_1\},
 \qquad B:=A-T.
\]
Thus $|T|\le3+2L$ and $B=X_{[n]\setminus(\mathcal D\cup T)}$.
Condition on
\[
 \mathcal F_B
 =\sigma(B,X_{ab}:a<b,\ a,b\in\mathcal D\cup T).
\]
Given $\mathcal F_B$, all decoration edges are frozen. Conditional
expectations, covariances and cumulants below are taken under $\Pp_n$
given $\mathcal F_B$, the true interacting law on the rows indexed
by $\mathcal D\cup T$.
For $a\in T$ set
$S_a=\sum_{t\in V(B)}(X_{at}-\E[X_{at}\mid\mathcal F_B])$.
The exact count
\[
 N_AP_A=\binom{|V(B)|}{2}P_B+
       \sum_{a\in T}\sum_{t\in V(B)}X_{at}
       +\sum_{\{a,b\}\subseteq T}X_{ab}
\]
and $(m-1)/N_A=2/m$ give
\[
 P_A-\E[P_A\mid\mathcal F_B]=N_A^{-1}\sum_{b\in T}S_b,\qquad
 z_a^A-\E[z_a^A\mid\mathcal F_B]=S_a-\frac2m\sum_{b\in T}S_b.
\]
Lemma~\ref{lem:row-mean-density} gives
$\|\sum_{t\in V(B)}\E[X_{at}\mid\mathcal F_B]-|V(B)|P_B\|_q\le C_q$, and
\[
\begin{aligned}
 N_A\{\E[P_A\mid\mathcal F_B]-P_B\}
 &=\sum_{a\in T}\left\{\sum_{t\in V(B)}\E[X_{at}\mid\mathcal F_B]-|V(B)|P_B\right\}\\
 &\quad+\sum_{\{a,b\}\subseteq T}(X_{ab}-P_B),\\
 m_a:=\E[z_a^A\mid\mathcal F_B]
 &=\sum_{t\in V(B)}\E[X_{at}\mid\mathcal F_B]-|V(B)|P_B\\
 &\quad+\sum_{b\in T\setminus\{a\}}(X_{ab}-P_B)
       -(m-1)\{\E[P_A\mid\mathcal F_B]-P_B\}.
\end{aligned}
\]
Thus, under the actual marginal of $\mathcal F_B$,
\begin{equation}\label{audit:dependency:s79:decorated-degree-conditional-means}
 \|\E[P_A\mid\mathcal F_B]-P_B\|_q\le C_qn^{-2},\qquad
 \max_{a\in T}\|m_a\|_q\le C_q.
\end{equation}

On $B\in\mathcal A_n$, \eqref{ep:row:feature-bound} gives, for $\ell=2,3$,
\[
 |\cum(S_{a_1},\ldots,S_{a_{\ell}}\mid\mathcal F_B)|
 \le\sum_{t_1,\ldots,t_{\ell}\in V(B)}
       |\cum(X_{a_1t_1},\ldots,X_{a_{\ell}t_{\ell}}\mid\mathcal F_B)|
 \le C_{\ell}n.
\]
Since $\sum_{b\in T}|\one_{\{b=a\}}-2/m|\le3$, the same bound
holds for conditional degree covariances and third cumulants.
The exact identity
\[
 \E[z_i^Az_j^Az_k^A\mid\mathcal F_B]
 =\cum(z_i^A,z_j^A,z_k^A\mid\mathcal F_B)
  +\sum_{\rm cyc}m_i\Cov(z_j^A,z_k^A\mid\mathcal F_B)+m_im_jm_k
\]
therefore gives
\begin{equation}\label{audit:dependency:s79:decorated-degree-frozen-bound}
 |\E[\phi(P_B)B_0(A;P_B)z_i^Az_j^Az_k^A]|
 \le Cn\|\phi\|_\infty.
\end{equation}
Here the multiplier is $\mathcal F_B$-measurable; H\"older's inequality uses the
bounded moments of $m_a$ from
\eqref{audit:dependency:s79:decorated-degree-conditional-means}, and the event
$B\notin\mathcal A_n$ contributes at most
$Cn^3e^{-cn}\|\phi\|_\infty$.

For the restoration, Corollary~\ref{lem:degree-deletion} gives
$\|P_A-P_B\|_2\le Cn^{-3/2}$ and
$\|z_i^Az_j^Az_k^A\|_2\le Cn^{3/2}$.
The latter also follows from
$z_a^A-z_a=\sum_{v\in\mathcal D}(P-X_{av})+(m-1)(P-P_A)$
and the sixth degree moments. Telescoping the decoration factors,
with repetitions, gives
\[
 |B_0(A;P_A)-B_0(A;P_B)|\le|\mathcal E_1||P_A-P_B|.
\]
Replacing both arguments costs at most
\[
 (\Lip\phi+|\mathcal E_1|\|\phi\|_\infty)
 \|P_A-P_B\|_2\|z_i^Az_j^Az_k^A\|_2
 \le C(\|\phi\|_\infty+\Lip\phi).
\]
This bound and \eqref{audit:dependency:s79:decorated-degree-frozen-bound}
prove both assertions.
\end{proof}

\begin{lemma}
\label{lem:branching-support}
Assume \hyperref[ass:PC]{\textup{(PC)}}.
Let $H$ be a fixed connected simple graph with maximum degree at least
three, and let $\phi:[0,1]\to\mathbb R$ be any Lipschitz function. Then
\begin{equation}\label{eq:branching-support}
 |\E[\phi(P)\inj(H,Y)]|
 \le C_H\{n^{v(H)-2}\|\phi\|_\infty
                       +n^{v(H)-5/2}\Lip\phi\}.
\end{equation}
\end{lemma}

\begin{proof}
Put $v=v(H)$ and choose a maximum-degree root $u$, of degree $d$.
For $d\ge4$, Lemma~\ref{lem:centered-root-moments} gives directly
\[
 \E|\inj(H,Y)|
 \le C_Hn\{n^{v-1-d/2}+n^{v-3}\}\le C_Hn^{v-2}.
\]
Suppose $d=3$. Order the neighbors as $u_1,u_2,u_3$, fix a root
image $o$ and an exterior $A$, and set
\[
 C_{ijk}=\sum_{\substack{\psi:V(H)\setminus\{u\}\hookrightarrow V(A)\\
          (\psi(u_1),\psi(u_2),\psi(u_3))=(i,j,k)}}
       \prod_{ab\in E(H-u)}(A_{\psi(a)\psi(b)}-P_A).
\]
Then $|C_{ijk}|\le n^{v-4}$, with $C_{ijk}=0$ on coincidences.
For fixed $A$, decompose the conditional mean under $\pi_A$ as
\[
\begin{aligned}
 \E_{\pi_A}[T_o(H,u;P_A)]
 &=\sum_{i,j,k}C_{ijk}
            (\mu_{A,i}-P_A)(\mu_{A,j}-P_A)(\mu_{A,k}-P_A)\\
 &\quad+\sum_{i,j,k}C_{ijk}(\mu_{A,i}-P_A)\Sigma_{A,jk}
        +\sum_{i,j,k}C_{ijk}(\mu_{A,j}-P_A)\Sigma_{A,ik}\\
 &\quad+\sum_{i,j,k}C_{ijk}(\mu_{A,k}-P_A)\Sigma_{A,ij}
        +\sum_{i,j,k}C_{ijk}\cum_{\pi_A}(X_{oi},X_{oj},X_{ok}).
\end{aligned}
\]
All terms on the right are $A$-measurable; their outer expectations
and random $L^q$ norms use the actual marginal of $A$.
For the cumulant term, Theorem~\ref{ep:row:sharp} gives, on good exteriors,
\[
 n^{v-4}\sum_{i,j,k}|\cum_{\pi_A}(X_{oi},X_{oj},X_{ok})|
 \le C_Hn^{v-3}.
\]

For the mean--covariance term with $(\mu_{A,i}-P_A)\Sigma_{A,jk}$,
\eqref{eq:cov-entry-profile} gives, for $j\ne k$,
\[
 \Sigma_{A,jk}
 =\frac{v_0(P_A)-\nu(P_A)}n
  +\frac{6\nu(P_A)^2\tau(P_A)}n(Y_A)_{jk}+R_{A,jk}.
\]
The first two terms contribute bounded profile multiples of
$n^{-1}\sum_i(\mu_{A,i}-P_A)G_i$, corresponding, respectively, to
\[
 G_i=\sum_{j,k}C_{ijk}
 \qquad\text{and}\qquad
 G_i=\sum_{j,k}C_{ijk}(Y_A)_{jk}.
\]
The first choice corresponds to the rooted graph $(H-u,u_1)$.
For the second, the rooted graph is $((H-u)+u_2u_3,u_1)$ if
$u_2u_3\notin E(H-u)$. If $u_2u_3\in E(H-u)$, the identity
\[
 (Y_A)_{jk}^2=\nu(P_A)+(1-2P_A)(Y_A)_{jk}
\]
expresses the second choice as a linear combination of the counts with
rooted graphs
\[
 ((H-u)-\{u_2u_3\},u_1)
 \qquad\text{and}\qquad (H-u,u_1),
\]
with coefficients $\nu(P_A)$ and $1-2P_A$, respectively.
Here $(H-u)-\{u_2u_3\}$ deletes only the edge $u_2u_3$;
all $v-1$ vertices are retained, including isolates.
Lemma~\ref{lem:root-exchange} applies to these three graph types and their
bounded-coefficient combinations, giving, for $\bar G=m^{-1}\sum_iG_i$,
\[
 |\bar G|\le C_Hn^{v-2},\qquad
 \bigl\|\|G-\bar G\one\|_2\bigr\|_2\le C_Hn^{v-2}.
\]
Thus \eqref{eq:mean-profiles-a}, \eqref{eq:cov-entry-profile}, and
the Cauchy--Schwarz inequality yield
\begin{equation}\label{audit:dependency:s79:branching-mixed-terms}
\begin{aligned}
 &\frac1n\E_A\left|\sum_i(\mu_{A,i}-P_A)G_i\right|\\
 &\quad\le\frac1n\|\bar G\|_2\left\|\sum_i(\mu_{A,i}-P_A)\right\|_2\\
 &\qquad+\frac1n\bigl\|\|\mu_A-P_A\one\|_2\bigr\|_2
       \bigl\|\|G-\bar G\one\|_2\bigr\|_2
 \le C_Hn^{v-3},\\
 &\E_A\left|\sum_{i,j,k}C_{ijk}(\mu_{A,i}-P_A)R_{A,jk}\right|\\
 &\quad\le n^{v-5/2}\bigl\|\|\mu_A-P_A\one\|_2\bigr\|_2
       \bigl\|\|R_A\|_{F,\mathrm{off}}\bigr\|_2
 \le C_Hn^{v-3}.
\end{aligned}
\end{equation}
Here we used
$\|\mu_A-P_A\one\|_1\le\sqrt n\|\mu_A-P_A\one\|_2$ and
$\sum_{j\ne k}|R_{A,jk}|\le n\|R_A\|_{F,\mathrm{off}}$.
The other mixed terms use $u_2$ and $u_3$ as
carrier roots; no symmetry of $C$ is assumed.

For the three-mean term, Lemma~\ref{lem:row-mean-profile} gives
\[
 \left\|\left\|\mu_A-P_A\one-\frac{2\nu(P_A)t(P_A)}n z^A
                 \right\|_1\right\|_3=O(1),
 \qquad\bigl\|\|\mu_A-P_A\one\|_1\bigr\|_3=O(\sqrt n).
\]
Successive replacement and H\"older's inequality give
\begin{equation}\label{audit:dependency:s79:branching-mean-replacement-error}
\begin{aligned}
 &C_Hn^{v-4}
       \left\|\left\|\mu_A-P_A\one-\frac{2\nu(P_A)t(P_A)}n z^A
                                               \right\|_1\right\|_3\\
 &\quad\times\left(
       \bigl\|\|\mu_A-P_A\one\|_1\bigr\|_3
       +\left\|\left\|\mu_A-P_A\one-\frac{2\nu(P_A)t(P_A)}n z^A
                                               \right\|_1\right\|_3
                                                    \right)^2\\
 &\quad\le C_Hn^{v-3}.
\end{aligned}
\end{equation}
The remaining tested term equals
\[
 n^{-3}\sum_{\psi:V(H-u)\hookrightarrow V(A)}
 \E_A\left[\phi(P_A)\{2\nu(P_A)t(P_A)\}^3
       \prod_{a=1}^3z_{\psi(u_a)}^A
       \prod_{bc\in E(H-u)}(A_{\psi(b)\psi(c)}-P_A)\right].
\]
The three degree labels are distinct, so
\eqref{audit:decorated-exterior-degree-bound} applies with multiplier
$x\mapsto\phi(x)\{2\nu(x)t(x)\}^3$. Sum over at most $n^{v-1}$ injections.
Combining this estimate with \eqref{ep:row:feature-bound},
\eqref{audit:dependency:s79:branching-mixed-terms}, and
\eqref{audit:dependency:s79:branching-mean-replacement-error} gives
\[
 |\E[\phi(P_A)T_o(H,u;P_A)]|
 \le C_H\{n^{v-3}\|\phi\|_\infty+n^{v-4}\Lip\phi\}.
\]
Polynomial envelopes handle bad exteriors. Restore the center and test by
\eqref{eq:high-support-common-transfer} and sum over $o$ to obtain
\eqref{eq:branching-support}.
\end{proof}

\subsection{Signed bounds for paths and cycles}
For $v\ge4$, use the labelled path and cycle
\begin{equation}\label{v8:setup:path-cycle-graphs}
\begin{gathered}
 V(P_v)=V(C_v)=[v],\quad
 E(P_v)=\{\{j,j+1\}:1\le j<v\},\\
 E(C_v)=E(P_v)\cup\{\{v,1\}\}.
\end{gathered}
\end{equation}

\begin{lemma}
\label{lem:path-cycle-support}
Assume \hyperref[ass:PC]{\textup{(PC)}}.
Let $H=P_v$ or $H=C_v$ with $v\ge4$. For every Lipschitz
$\phi:[0,1]\to\mathbb R$, we have
\begin{equation}\label{eq:path-cycle-support}
 |\E[\phi(P)\inj(H,Y)]|
 \le C_H\{n^{v-2}\|\phi\|_\infty+n^{v-5/2}\Lip\phi\}.
\end{equation}
\end{lemma}

\begin{proof}
Root $P_v$ at its second vertex, taking its first
and third vertices as the ordered neighbors; root $C_v$ anywhere. Write these vertices as
$u,u_1,u_2$, fix a root image $o$, and put $A=X_{[n]\setminus\{o\}}$.
With $Y_A=A-P_A(J-I)$, define $C_{ii}=0$ and
\[
 C_{ij}=\sum_{\substack{\psi:V(H)\setminus\{u\}\hookrightarrow V(A)\\
                (\psi(u_1),\psi(u_2))=(i,j)}}
       \prod_{ab\in E(H-u)}(Y_A)_{\psi(a)\psi(b)}.
\]
For paths, $u_1$ is isolated in $H-u$, so $C$ need not be symmetric.
For fixed $A$, taking the expectation under $\pi_A$ gives
\[
 \E_{\pi_A}[T_o(H,u;P_A)]
 =\sum_{i\ne j}C_{ij}\Sigma_{A,ij}+(\mu_A-P_A\one)^{\mathsf T}C(\mu_A-P_A\one).
\]
Every term on the right is $A$-measurable. The outer expectations
$\E_A$, random moment norms and $O_{L^q}$ bounds below use the actual
marginal of $A$. Each term will be $O(n^{v-3})$ in absolute expectation.

\proofstep{1. Covariance term.}
The carriers $H-u$ and $(H-u)+u_1u_2$ are respectively
$(K_1\sqcup P_{v-2},P_{v-1})$ for paths and
$(P_{v-1},C_{v-1})$ for cycles.
Summing over the distinct images $i,j\in V(A)$ removes the two root
constraints, giving
\[
 \begin{aligned}
 \sum_{i\ne j}C_{ij}&=\inj(H-u,Y_A),\\
 \sum_{i\ne j}C_{ij}(Y_A)_{ij}
   &=\inj((H-u)+u_1u_2,Y_A).
 \end{aligned}
\]
The extra factor $(Y_A)_{ij}$ corresponds to the added edge $u_1u_2$.
Both graphs have $v-1$ vertices and at least one edge, so
Theorem~\ref{thm:counts}, \eqref{audit:centered-count-replacement}, yields
\begin{equation}\label{audit:dependency:s79:path-cycle-centered-count-moments}
 \E_A|\inj(H-u,Y_A)|
 +\E_A|\inj((H-u)+u_1u_2,Y_A)|
 \le C_Hn^{v-2}.
\end{equation}
For each injection $\psi$, expand the centered edge factors as
\[
 \prod_{ab\in E(H-u)}(A_{\psi(a)\psi(b)}-P_A)
 =\sum_{S\subseteq E(H-u)}(-P_A)^{e(H-u)-|S|}
       \prod_{ab\in S}A_{\psi(a)\psi(b)}.
\]
The ordered roots $u_1,u_2$ are nonadjacent in $H-u$, and hence in
every spanning subgraph with edge set $S$.
Lemma~\ref{lem:two-root-counts}, \eqref{eq:two-root}, therefore gives
$P_A^{|S|}$ as the leading normalized two-root count, with no root-edge
factor $A_{ij}$. Since $e(H-u)\ge1$, these leading terms cancel:
\[
 \sum_{S\subseteq E(H-u)}(-P_A)^{e(H-u)-|S|}P_A^{|S|}
 =(P_A-P_A)^{e(H-u)}=0.
\]
Consequently,
\[
\begin{aligned}
 n^{-(v-3)}C_{ij}
 &=\sum_{S\subseteq E(H-u)}(-P_A)^{e(H-u)-|S|}
       \{P_A^{|S|}+O_{L^q}(n^{-1/2})\}
 =O_{L^q}(n^{-1/2}),\\
 \bigl\|\|C\|_F\bigr\|_q
 &\le\left(\sum_{i\ne j}\|C_{ij}\|_{q\vee2}^2\right)^{1/2}
 \le C_{H,q}n^{v-5/2}.
\end{aligned}
\]
Substituting \eqref{eq:cov-entry-profile} and using the two count
identities gives
\[
 \begin{aligned}
 \sum_{i\ne j}C_{ij}\Sigma_{A,ij}
 &=\frac{v_0(P_A)-\nu(P_A)}n\inj(H-u,Y_A)\\
 &\quad+\frac{6\nu(P_A)^2\tau(P_A)}n
            \inj((H-u)+u_1u_2,Y_A)
       +\sum_{i\ne j}C_{ij}R_{A,ij}.
 \end{aligned}
\]
Thus \eqref{audit:dependency:s79:path-cycle-centered-count-moments},
the preceding Frobenius bound, and the remainder estimate in
\eqref{eq:cov-entry-profile} yield
\begin{equation}\label{audit:dependency:s79:path-cycle-covariance-term}
\begin{aligned}
 \E_A\left|\sum_{i\ne j}C_{ij}\Sigma_{A,ij}\right|
 &\le\frac Cn\{\E_A|\inj(H-u,Y_A)|
        +\E_A|\inj((H-u)+u_1u_2,Y_A)|\}\\
 &\quad+\bigl\|\|C\|_F\bigr\|_2
             \bigl\|\|R_A\|_{F,\mathrm{off}}\bigr\|_2
 \le C_Hn^{v-3}.
\end{aligned}
\end{equation}

\proofstep{2. Quadratic mean term.}
Replace injective embeddings by unrestricted walks, whose coefficients
are $Y_A^{v-2}$ for cycles and
$\one(Y_A^{v-3}\one)^{\mathsf T}$ for paths.
Off-diagonal collision errors are at most $C_Hn^{v-4}$ per entry,
and diagonal errors are at most $n^{v-3}$. Thus
\begin{equation}\label{audit:dependency:s79:path-cycle-collision-bounds}
\begin{aligned}
 \max_{\ast\in\{\infty\to\infty,1\to1\}}
 \|C-Y_A^{v-2}\|_\ast&\le C_Hn^{v-3}&&(H=C_v),\\
 \max_{\ast\in\{\infty\to\infty,1\to1\}}
 \|C-\one(Y_A^{v-3}\one)^{\mathsf T}\|_\ast
 &\le C_Hn^{v-3}&&(H=P_v).
\end{aligned}
\end{equation}
Both matrix norms give a quadratic-form error of at most
$C_Hn^{v-3}\|\mu_A-P_A\one\|_2^2$, without requiring symmetry of $C$.
For cycles, \eqref{eq:mean-profiles-b} and \eqref{eq:Yz}, with $\|Y_A\|_{\op}\le n$, give
\begin{equation}\label{audit:dependency:s79:path-cycle-cycle-quadratic-form}
\begin{aligned}
 &\bigl\|\|Y_A(\mu_A-P_A\one)\|_2\bigr\|_2\\
 &\le Cn^{-1}\bigl\|\|Y_Az^A\|_2\bigr\|_2\\
 &\quad+n\left\|\left\|\mu_A-P_A\one-\frac{2\nu(P_A)t(P_A)}n z^A
                                                \right\|_2\right\|_2
 \le C\sqrt n,\\
 &\E_A|(\mu_A-P_A\one)^{\mathsf T}Y_A^{v-2}(\mu_A-P_A\one)|\\
 &\le n^{v-4}\E_A\|Y_A(\mu_A-P_A\one)\|_2^2\le C_Hn^{v-3}.
\end{aligned}
\end{equation}
For paths,
\begin{equation}\label{audit:dependency:s79:path-cycle-path-quadratic-form}
\begin{aligned}
 &\bigl\|\|Y_A^{v-3}\one\|_2\bigr\|_q
 =\bigl\|\|Y_A^{v-4}z^A\|_2\bigr\|_q\le C_{H,q}n^{v-3},\\
 &\E_A|((\mu_A-P_A\one)^{\mathsf T}\one)(Y_A^{v-3}\one)^{\mathsf T}(\mu_A-P_A\one)|\\
 &\quad\le\left\|\sum_i(\mu_{A,i}-P_A)\right\|_2
       \bigl\|\|\mu_A-P_A\one\|_2\bigr\|_4
       \bigl\|\|Y_A^{v-3}\one\|_2\bigr\|_4
 \le C_Hn^{v-3}.
\end{aligned}
\end{equation}

Equations~\eqref{audit:dependency:s79:path-cycle-cycle-quadratic-form} and
\eqref{audit:dependency:s79:path-cycle-path-quadratic-form}, with
\eqref{audit:dependency:s79:path-cycle-collision-bounds} and \eqref{eq:mean-profiles-a}, bound the quadratic mean
term. Combine this with \eqref{audit:dependency:s79:path-cycle-covariance-term},
restore the center and test by \eqref{eq:high-support-common-transfer}, and
sum over the $n$ root images to prove \eqref{eq:path-cycle-support}.
\end{proof}

\subsection{General signed count bounds}
For disconnected graphs we control collisions between component
embeddings, using \eqref{v8:setup:signed-centered-count} and
\eqref{v8:setup:signed-test-norms}. In particular,
\begin{equation}\label{v8:setup:general-centered-count-notation}
 \inj(H,Y)=\sum_{\varphi:V(H)\hookrightarrow[n]}
                 \prod_{ab\in E(H)}(X_{\varphi(a)\varphi(b)}-P).
\end{equation}
\begin{proposition}\label{prop:high-support}
Assume \hyperref[ass:PC]{\textup{(PC)}}.
Let $H$ be a fixed simple graph without isolated vertices, with at least
three edges and not equal to $K_3$, and let $\phi:[0,1]\to\mathbb R$ be
any Lipschitz function. Then
\begin{equation}\label{eq:high-support}
 |\E\phi(P)\inj(H,Y)|
 \le C_H\bigl(n^{v(H)-2}\|\phi\|_\infty
                    +n^{v(H)-5/2}\Lip\phi\bigr).
\end{equation}
The excluded triangle is treated in Lemma~\ref{lem:triangle-profile},
\eqref{eq:T-weak}; its absolute moments satisfy
\eqref{eq:triangle-moment} of Corollary~\ref{cor:matrix-moments}.
The constant $C_H$ in \eqref{eq:high-support} depends on $H$,
the fixed model, and the parameter bound $K$, but not on $n$ or $\phi$.
\end{proposition}

\begin{proof}
By \eqref{eq:branching-support} and \eqref{eq:path-cycle-support}, only the
disconnected case remains. Let $H$ have components $H_1,\ldots,H_c$,
$c\ge2$, and write $v_a=v(H_a)$, $v=\sum_a v_a$.
Every component has at least one edge, so Theorem~\ref{thm:counts} gives
$\|\inj(H_a,Y)\|_2\le C_Hn^{v_a-1}$ and
$|\inj(H_a,Y)|\le n^{v_a}$. Hence
\begin{equation}\label{eq:disconnected-component-product}
 \left\|\prod_{a=1}^c\inj(H_a,Y)\right\|_1
 \le n^{v-v_1-v_2}\|\inj(H_1,Y)\|_2\|\inj(H_2,Y)\|_2
 \le C_Hn^{v-2}.
\end{equation}
Separate the collisions in this product by the number of distinct images:
\begin{equation}\label{eq:disconnected-collision-expansion}
\begin{aligned}
 \prod_{a=1}^c\inj(H_a,Y)
 &=\inj(H,Y)+\sum_{a<b}
    \sum_{\substack{u\in V(H_a)\\w\in V(H_b)}}
      \inj(H/\{u=w\},Y)+R_H,\\
 |R_H|&\le\sum_{j=1}^{v-2}\binom{n}{j}j^v\le C_Hn^{v-2}.
\end{aligned}
\end{equation}
A map with $v-1$ images identifies one cross-component pair with
multiplicity one.
The remainder $R_H$ collects maps that are injective within each component
and have at most $v-2$ distinct images; choosing $j$ images and assigning
the $v$ labelled vertices gives at most $\binom{n}{j}j^v$ maps.
Since every edge product has absolute value at most one and $v$ is fixed,
summing over $1\le j\le v-2$ gives the bound for $R_H$ in
\eqref{eq:disconnected-collision-expansion}.
Each single identification gives a simple
quotient with $v-1$ vertices and at least one edge, since it creates neither loops
nor repeated edges. Hence
\[
 \|\inj(H/\{u=w\},Y)\|_1\le C_Hn^{v-2}
\]
by \eqref{audit:centered-count-replacement}. There are only $C_H$
such pairs. Combining \eqref{eq:disconnected-component-product} and
\eqref{eq:disconnected-collision-expansion} consequently gives
\[
 \|\inj(H,Y)\|_1
 \le\left\|\prod_a\inj(H_a,Y)\right\|_1
   +\sum_{a<b}\sum_{\substack{u\in V(H_a)\\w\in V(H_b)}}
       \|\inj(H/\{u=w\},Y)\|_1+\|R_H\|_1
 \le C_Hn^{v-2}.
\]
Multiplication by $\|\phi\|_\infty$ proves the disconnected case.
\end{proof}

\begin{remark}\label{rem:disconnected-support-scope}
Proposition~\ref{prop:high-support} combines
Lemmas~\ref{lem:branching-support} and~\ref{lem:path-cycle-support}
with the disconnected estimate above.
For disconnected $H$ without isolated vertices, the restriction
$e(H)\ge3$ is unnecessary: the proof gives
\[
 \E|\inj(H,Y)|\le C_Hn^{v(H)-2},
\]
including $H=2K_2$, so the tested bound requires only
$\|\phi\|_\infty$ for bounded $\phi$.
The restriction in the unified statement excludes the adjacent two-edge
path $P_3$, whose degree-square contribution is treated by
Lemma~\ref{lem:low-supports}.
\end{remark}

\section{Proofs of the main results}\label{sec:drift}
The signed support estimates and two-row pairings identify the
microscopic drift. Stationarity supplies the required moment bounds;
the exchangeable-pair theorem then yields Theorems~\ref{thm:main} and
\ref{fixed:main}.

\subsection{Stein approximation bounds}\label{sec:stein}

The approximations in Theorems~\ref{thm:main} and~\ref{fixed:main}
require bounds on the weak regression error, the conditional second
moment of the increment, and the higher-order increment remainder.
We first state the exchangeable-pair bound of Theorem~\ref{prop:Stein},
which combines these quantities. The subsequent subsections verify its
hypotheses for the edge-count updates in each parameter regime.

We give distributional approximation bounds for laws with polynomial
potentials of even degree and positive leading coefficient, including
Gaussian laws in the quadratic case. The results are stated for general
exchangeable pairs and will be applied to the edge-count approximations
below; see \cite[Chapter~2]{CGS2011} and
\cite{Stein86,EL,CS} for background.
The notation in this subsection is local. Fix an integer $\mc\ge1$ and a
compact parameter set $\mathcal K\Subset(0,\infty)\times\R^{2\mc-1}$.
For $\gamma=(a,c_1,\ldots,c_{2\mc-1})\in \mathcal K$, set
\begin{equation}\label{stein:quartic-potential}
 \begin{aligned}
 V_\gamma(w)&=\frac a{2\mc}w^{2\mc}
 +\sum_{j=1}^{2\mc-1}c_jw^j,\\
 \rho_\gamma(w)&=\mathcal Z_\gamma^{-1}e^{-V_\gamma(w)},\qquad
 \mathcal Z_\gamma=\int_\R e^{-V_\gamma(t)}\,\dd t.
 \end{aligned}
\end{equation}
Write $V=V_\gamma$ and $\rho=\rho_\gamma$, and let $Y$ have density $\rho$;
the fixed parameter $\gamma$ is suppressed in the notation. Constants $c_{\mathcal K},C_{\mathcal K}>0$ depend only on $\mathcal K$ and may
change between occurrences; convexity of $V$ is not required. For Theorem~\ref{thm:main}, take
$a=a_\star$, $c_j=-\beta_{j,n}/j!$ for $j\le\mc$, and $c_j=0$ otherwise.

Let $X$ be a random element on a probability space and $W$ a real-valued
measurable function of $X$. Suppose $(W,W')$ is exchangeable. For a deterministic $\Lambda>0$, set
\begin{equation}\label{v8:nearby:app:exchangeable-pair-data}
 \Delta=W'-W,\qquad
 D=\frac{\E[\Delta^2\mid X]}{2\Lambda}.
\end{equation}

\begin{theorem}
\label{prop:Stein}
Let $Y$ have density \eqref{stein:quartic-potential}, with $\gamma\in \mathcal K$.
Suppose $|\Delta|\le\delta\le1$ almost surely and
$\E|W|^{2\mc-1}<\infty$.
Define the bounded-Lipschitz weak regression error by
\begin{equation}\label{stein:weak-regression}
 \begin{aligned}
 \mathcal R
 &:=\sup_{\|f\|_\infty+\Lip f\le1}
 \left|\E\left[
   \left\{\frac{\E[\Delta\mid X]}\Lambda+V'(W)\right\}
           f(W)\right]\right|\\
 &=\frac1\Lambda\sup_{\|f\|_\infty+\Lip f\le1}
       |\E[\Delta f(W)]+\Lambda\E[V'(W)f(W)]|<\infty.
 \end{aligned}
\end{equation}
Then
\begin{equation}\label{eq:Stein-bound}
 \begin{aligned}
 &\dW(W,Y)+\dK(W,Y)\\
 &\quad\le C_{\mathcal K}\left\{
       \E|1-D|+\mathcal R
       +\frac1\Lambda\E[|\Delta|^3(1+|W|^{2\mc-1})]
       +\delta\E|V'(W)|\right\}.
 \end{aligned}
\end{equation}
For $\mc\ge2$, \eqref{eq:Stein-bound} also holds with
$\sum_{q=1}^2|\E W^q-\E Y^q|$ added to its left-hand side.
The last term is needed only for the Kolmogorov estimate.
No uniform moment bound of order $2\mc$ is required.
\end{theorem}

The Stein solution estimates and the proof of
Theorem~\ref{prop:Stein} are given in Appendix~\ref{app:stein-proofs}.

\subsection{Averaged logit expansions and moments}
For $L_e$ and $\ell_n$ defined in \eqref{v8:setup:physical-edge-logit},
define the averaged centered logit moments by
\begin{equation}\label{eq:edge-centered-logit}
 M_j=N^{-1}\sum_e[L_e-\ell_n(P)]^j,\qquad j=1,2,3,4.
\end{equation}
With the fixed profiles from \eqref{eq:profiles} and
\eqref{eq:logit-profiles}, define
\begin{equation}\label{eq:averaged-logit-correction}
 B_{\mathrm{logit}}(x)
 =2t'(x)\{v_0(x)-\nu(x)\}-2s'(x)
     +12\tau(x)\tau'(x)\nu(x)^3,\qquad x\in[0,1].
\end{equation}
The three terms arise from pairs of adjacent edges, the normalization of
injective counts, and triangles, respectively.
Proposition~\ref{prop:microscopic} identifies $B_{\mathrm{logit}}(P)/n$
as the leading term of $M_1$ in expectations tested against Lipschitz
functions of $P$; this correction subsequently enters the drift expansion
for the edge density.

We first identify the mean correction;
Proposition~\ref{prop:high-support} controls the remaining supports.
\begin{proposition}
\label{prop:microscopic}
Assume \hyperref[ass:PC]{\textup{(PC)}}.
For $M_1$ in \eqref{eq:edge-centered-logit} and every Lipschitz
$\phi:[0,1]\to\mathbb R$,
\begin{equation}\label{eq:micro-mean-input}
 \left|\E\left[\phi(P)
              \{M_1-B_{\mathrm{logit}}(P)/n\}\right]\right|
 \le C\{n^{-2}\|\phi\|_\infty+n^{-5/2}\Lip\phi\}.
\end{equation}
The correction profile $B_{\mathrm{logit}}(x)$ is defined in
\eqref{eq:averaged-logit-correction}.
We first compute the edge difference defining $L_e$ and then evaluate
the centering field $\ell_n$ at $P$.
\end{proposition}

\begin{proof}
\proofstep{Support expansion.}
We separate the empty support, degree-square and triangle terms from
the remaining centered supports. For $S\subseteq E(F_r)$, let $F_r[S]$ consist of its edges and endpoints,
with the empty count equal to one. Differentiation precedes centering:
\[
 \sum_e\partial_e\inj(F_r,X)=\sum_{f\in E(F_r)}\inj(F_r-f,X),
 \qquad
 \prod_{g\ne f}X_{\varphi(g)}
 =\sum_{S\subseteq E(F_r)\setminus\{f\}}
       P^{e_r-1-|S|}\prod_{g\in S}Y_{\varphi(g)}.
\]
Edge deletion retains all $v_r$ vertices. Each proper $S$ has
$e_r-|S|$ possible deleted edges and $(n-h)_{v_r-h}$ extensions,
where $h=v(F_r[S])$. Therefore
\begin{equation}\label{audit:averaged-logit-support-expansion}
\begin{aligned}
 N^{-1}\sum_eL_e-h_n
 &=\sum_r\sum_{S\subsetneq E(F_r)}
                 c_{r,S,n}(P)\inj(F_r[S],Y),\\
 c_{r,S,n}(x)
 &=\frac{b_{r,n}n^{2-v_r}}N(e_r-|S|)
      x^{e_r-1-|S|}(n-h)_{v_r-h}.
\end{aligned}
\end{equation}
The exponents are nonnegative, so the following bounds hold uniformly on $[0,1]$:
\begin{equation}\label{audit:dependency:s1011:support-coefficient-bounds}
 \frac{n^{2-v_r}}N(n-h)_{v_r-h}
   =2n^{-h}(1+O(n^{-1})),\qquad
 \|c_{r,S,n}\|_{C^1}\le Cn^{-h}.
\end{equation}
Write $s_n(x),t_n(x),\tau_n(x)$ for the profiles using $b_{r,n}$. The empty
support gives
\begin{equation}\label{audit:dependency:s1011:empty-support}
\begin{gathered}
 \frac{n^{2-v_r}(n)_{v_r}}N
 =2-\frac2n\{\tbinom{v_r}2-1\}+O(n^{-2}),\\
 \sum_r c_{r,\varnothing,n}(P)
 =\ell_n(P)-h_n-\frac{2s_n'(P)}n+O(n^{-2}).
\end{gathered}
\end{equation}

\proofstep{Low supports.}
Since $\sum_{j\ne i}Y_{ij}=z_i$, $\sum_i z_i=0$ and
$Y_{ij}^2=(1-2P)Y_{ij}+\nu(P)$, splitting the square of the total edge sum gives
\[
 0=\Bigl(\sum_{i\ne j}Y_{ij}\Bigr)^2
   =\inj(2K_2,Y)+4\inj(P_3,Y)+2\sum_{i\ne j}Y_{ij}^2.
\]
Consequently,
\begin{equation}\label{audit:low-support-count-identities}
\begin{gathered}
 \inj(K_2,Y)=0,\qquad \sum_{i\ne j}Y_{ij}^2=n(n-1)\nu(P),\\
 \inj(P_3,Y)=\sum_i z_i^2-n(n-1)\nu(P),\\
 \inj(2K_2,Y)=-4\sum_i z_i^2+2n(n-1)\nu(P),\qquad
 \inj(K_3,Y)=\tr Y^3 .
\end{gathered}
\end{equation}
Adjacent pairs and triangles contribute $2t_n'(P)n^{-3}(1+O(n^{-1}))$
and $2\tau_n'(P)n^{-3}(1+O(n^{-1}))$, respectively, because
\[
 t_n'(x)=\sum_r b_{r,n}a(F_r)(e_r-2)x^{e_r-3},\qquad
 \tau_n'(x)=\sum_r b_{r,n}\triangle(F_r)(e_r-3)x^{e_r-4},
\]
with zero derivatives omitted. Combining
\eqref{audit:averaged-logit-support-expansion},
\eqref{audit:dependency:s1011:empty-support}, and
\eqref{audit:low-support-count-identities} gives
\begin{equation}\label{audit:averaged-logit-low-support}
 M_1=\frac{2t_n'(P)}n\{n^{-2}\sum_i z_i^2-\nu(P)\}
       -\frac{2s_n'(P)}n+\frac{2\tau_n'(P)}{n^3}\tr Y^3+R_n.
\end{equation}
By \eqref{audit:dependency:s1011:empty-support} and
\eqref{audit:low-support-count-identities}, the empty-support remainder
and the replacement of $n(n-1)$ by $n^2$ contribute $O(n^{-2})$;
the $O(n^{-4})$ adjacent-pair coefficient error
and the disjoint-pair coefficient contribute
$O(n^{-4}\sum_i z_i^2+n^{-2})$; and the $O(n^{-4})$ triangle coefficient
error contributes $O(n^{-4}|\tr Y^3|)$.
Equations~\eqref{eq:degree-vector-moment}, with no deletion, and
\eqref{eq:triangle-moment} give
\[
 \E\left[n^{-2}\left(1+n^{-2}\sum_i z_i^2\right)
                 +n^{-4}|\tr Y^3|\right]\le Cn^{-2}.
\]

\proofstep{Tested remainder and coefficient.}
All other supports have at least three edges and are not triangles.
By \eqref{audit:dependency:s1011:support-coefficient-bounds}, their coefficients satisfy
\[
 \|c_{r,S,n}\phi\|_\infty\le Cn^{-h}\|\phi\|_\infty,\qquad
 \Lip(c_{r,S,n}\phi)\le Cn^{-h}(\|\phi\|_\infty+\Lip\phi).
\]
Applying Proposition~\ref{prop:high-support} to each of the finitely
many supports, then adding the absolute errors above, gives
\begin{equation}\label{audit:averaged-logit-tested-remainder}
 |\E[\phi(P)R_n]|
 \le C\{n^{-2}\|\phi\|_\infty+n^{-5/2}\Lip\phi\}.
\end{equation}
In \eqref{audit:averaged-logit-low-support}, apply
\eqref{eq:V-weak} and \eqref{eq:T-weak} with tests
$2t_n'\phi$ and $2\tau_n'\phi$. The uniform $C^1$ bounds on $t_n'$
and $\tau_n'$, together with the outer $n^{-1}$, preserve the error in
\eqref{audit:averaged-logit-tested-remainder}. Replacing the coefficients
by their leading values costs $O(n^{-2}\|\phi\|_\infty)$ and gives
\[
 2t'(P)\{v_0(P)-\nu(P)\}-2s'(P)
       +12\tau'(P)\tau(P)\nu(P)^3=B_{\mathrm{logit}}(P).
\]
Together with \eqref{audit:averaged-logit-tested-remainder} and the
definition \eqref{eq:averaged-logit-correction}, this proves \eqref{eq:micro-mean-input}.
\end{proof}

The variance profile and conditional bias now give the averaged
second and third moments; the conditional moment bounds control
the fourth moment.

\begin{lemma}
\label{lem:logit-higher-moments}
Assume \hyperref[ass:PC]{\textup{(PC)}}.
For $M_j$ in \eqref{eq:edge-centered-logit} and every Lipschitz
$\phi:[0,1]\to\mathbb R$,
\begin{equation}\label{eq:micro-higher-inputs}
 \begin{aligned}
 &|\E[\phi(P)\{M_2-w_2(P)/n\}]|+|\E[\phi(P)M_3]|\\
 &\qquad\le C\{n^{-2}\|\phi\|_\infty+n^{-5/2}\Lip\phi\},
 \qquad \E M_4\le Cn^{-2}.
 \end{aligned}
\end{equation}
The third-moment bound is signed, and $w_2(x)$ is defined in
\eqref{v8:setup:two-row-variance-profile}.
\end{lemma}

\begin{proof}
By exchangeability, it suffices to consider edge $\{1,2\}$.
Set $A=X_{\{3,\ldots,n\}}$, $e=X_{12}$ and $\mathcal F=\sigma(A,e)$;
the conditional row law is $\mu_{A,(e),1}$ in
\eqref{v8:setup:two-row-interpolation-density}.
All moment norms below use $\Pp_n$, or equivalently the actual $(A,e)$
marginal for $\mathcal F$-measurable quantities. We use the decomposition
\begin{equation}\label{audit:dependency:s1011:logit-decomposition}
 \begin{aligned}
 L_{12}-\ell_n(P)
 &=\bigl[L_{12}-\E[L_{12}\mid\mathcal F]\bigr]
   +\bigl[\E[L_{12}\mid\mathcal F]-\ell_n(P_A)\bigr]\\
 &\quad+\ell_n(P_A)-\ell_n(P).
 \end{aligned}
\end{equation}
Equations~\eqref{eq:local-polynomial-moments}, \eqref{eq:bias}, and
\eqref{row:eq:degree-deletion}, followed by Minkowski's inequality, give,
for every fixed $q\ge1$,
\begin{equation}\label{audit:centered-logit-decomposition-moments}
 \begin{gathered}
 \bigl\|L_{12}-\E[L_{12}\mid\mathcal F]\bigr\|_q\le C_qn^{-1/2},
 \qquad
 \bigl\|\E[L_{12}\mid\mathcal F]-\ell_n(P_A)\bigr\|_q\le C_qn^{-1},\\
 \|\ell_n(P_A)-\ell_n(P)\|_q+\|P-P_A\|_q\le C_qn^{-3/2},\\
 \|L_{12}-\ell_n(P_A)\|_q
 \le\bigl\|L_{12}-\E[L_{12}\mid\mathcal F]\bigr\|_q
    +\bigl\|\E[L_{12}\mid\mathcal F]-\ell_n(P_A)\bigr\|_q
 \le C_qn^{-1/2}.
 \end{gathered}
\end{equation}
Conditional centering, \eqref{ep:row:polynomial-bound}, and
\eqref{eq:variance} also give
\begin{equation}\label{audit:dependency:s1011:conditional-logit-moments}
 \begin{gathered}
 \E\bigl[L_{12}-\E[L_{12}\mid\mathcal F]\mid\mathcal F\bigr]=0,
 \qquad
 \left\|\Var(L_{12}\mid\mathcal F)-\frac{w_2(P_A)}n\right\|_q
       \le C_qn^{-2},\\
 \left\|\E\bigl[(L_{12}-\E[L_{12}\mid\mathcal F])^3
                         \mid\mathcal F\bigr]\right\|_q\le C_qn^{-2}.
 \end{gathered}
\end{equation}
Polynomial envelopes account for exceptional exteriors in these bounds.

We first evaluate the conditional moments centered at $\ell_n(P_A)$.
For the second moment, \eqref{audit:centered-logit-decomposition-moments}
and \eqref{audit:dependency:s1011:conditional-logit-moments} yield
\begin{equation}\label{audit:dependency:s1011:logit-square-errors}
 \begin{aligned}
 \E\bigl[(L_{12}-\ell_n(P_A))^2\mid\mathcal F\bigr]
 &=\Var(L_{12}\mid\mathcal F)
       +\bigl[\E[L_{12}\mid\mathcal F]-\ell_n(P_A)\bigr]^2,\\
 \left\|\E\bigl[(L_{12}-\ell_n(P_A))^2\mid\mathcal F\bigr]
                  -\frac{w_2(P_A)}n\right\|_q
 &\le C_qn^{-2}.
 \end{aligned}
\end{equation}
For the signed third moment, conditional centering gives
\begin{equation}\label{audit:dependency:s1011:logit-signed-cube}
 \begin{aligned}
 &\E\bigl[(L_{12}-\ell_n(P_A))^3\mid\mathcal F\bigr]\\
 &\quad=\E\bigl[(L_{12}-\E[L_{12}\mid\mathcal F])^3\mid\mathcal F\bigr]
       +3\bigl[\E[L_{12}\mid\mathcal F]-\ell_n(P_A)\bigr]
                       \Var(L_{12}\mid\mathcal F)\\
 &\qquad+\bigl[\E[L_{12}\mid\mathcal F]-\ell_n(P_A)\bigr]^3,\\
 &\left\|\E\bigl[(L_{12}-\ell_n(P_A))^3\mid\mathcal F\bigr]\right\|_q
 \le C_q(n^{-2}+n^{-1}\cdot n^{-1}+n^{-3})\le C_q'n^{-2}.
 \end{aligned}
\end{equation}
The last bound follows from H\"older's inequality and
\eqref{audit:centered-logit-decomposition-moments}--\eqref{audit:dependency:s1011:conditional-logit-moments}.

We next freeze the test at $P_A$ and restore the center $\ell_n(P)$.
For $j\in\{2,3\}$, measurability gives
\[
 \E[\phi(P_A)(L_{12}-\ell_n(P_A))^j]
 =\E\Bigl[\phi(P_A)\E[(L_{12}-\ell_n(P_A))^j\mid\mathcal F]\Bigr].
\]
H\"older's inequality, factorization of the difference of powers, and
\eqref{audit:centered-logit-decomposition-moments} yield
\begin{equation}\label{audit:dependency:s1011:logit-cube-cross-errors}
 \begin{aligned}
 &\left|\E[(\phi(P)-\phi(P_A))(L_{12}-\ell_n(P_A))^j]\right|\\
 &\quad\le\Lip\phi\,\|P-P_A\|_2\|L_{12}-\ell_n(P_A)\|_{2j}^j
       \le Cn^{-(j+3)/2}\Lip\phi,\\
 &\left|\E\Bigl[\phi(P)\bigl\{(L_{12}-\ell_n(P))^j
                                  -(L_{12}-\ell_n(P_A))^j\bigr\}\Bigr]\right|\\
 &\quad\le C\|\phi\|_\infty\|\ell_n(P_A)-\ell_n(P)\|_j
       \bigl[\|L_{12}-\ell_n(P_A)\|_j
                      +\|\ell_n(P_A)-\ell_n(P)\|_j\bigr]^{j-1}\\
 &\quad\le Cn^{-(j+2)/2}\|\phi\|_\infty.
 \end{aligned}
\end{equation}
For the second moment, the bounded $C^1$ norm of $w_2$ further gives
\[
 \frac1n\left|\E[\phi(P_A)w_2(P_A)-\phi(P)w_2(P)]\right|
 \le Cn^{-5/2}(\|\phi\|_\infty+\Lip\phi).
\]
Combining this estimate with
\eqref{audit:dependency:s1011:logit-square-errors} and
\eqref{audit:dependency:s1011:logit-cube-cross-errors}, with $j=2$,
proves the required second-moment bound.
Equations~\eqref{audit:dependency:s1011:logit-signed-cube} and
\eqref{audit:dependency:s1011:logit-cube-cross-errors}, with $j=3$,
give the signed third-moment bound, since $n^{-3}\le n^{-5/2}$.

Finally, the three terms in
\eqref{audit:dependency:s1011:logit-decomposition} satisfy
\[
 \begin{aligned}
 \E|L_{12}-\ell_n(P)|^4
 &\le27\Bigl[\E\bigl|L_{12}-\E[L_{12}\mid\mathcal F]\bigr|^4
       +\E\bigl|\E[L_{12}\mid\mathcal F]-\ell_n(P_A)\bigr|^4\\
 &\quad+\E|\ell_n(P_A)-\ell_n(P)|^4\Bigr]\le Cn^{-2}
 \end{aligned}
\]
by \eqref{audit:centered-logit-decomposition-moments}.
Averaging the three single-edge estimates proves \eqref{eq:micro-higher-inputs}.
\end{proof}

\subsection{Drift calibration and expansion}\label{sec:final-proof}
Combine the microscopic mean and variance profiles into
\begin{equation}\label{eq:a-function}
 \mathfrak a(x)=\sig'(\ell(x))B_{\mathrm{logit}}(x)
                     +\tfrac12\sig''(\ell(x))w_2(x).
\end{equation}
For the bounded parameter window \eqref{eq:window}, put
\begin{equation}\label{v8:setup:drift-calibration-coefficients}
 A_n(x)=\mathfrak a(x)+\sig'(\ell(x))
                    \{\kappa_n+2\mathcal U_n'(x)\}.
\end{equation}
Only the derivatives of $\mathcal F_1$ at $p$ and local values of $C_n$
in \eqref{even:calibration} are used; no extension outside $J$ is needed.
The global profiles \eqref{v8:setup:two-row-variance-profile} and
\eqref{eq:averaged-logit-correction} define $A_n$ independently on $[0,1]$.

\begin{lemma}
\label{lem:drift-calibration}
Let $\Phi$ be as in \eqref{eq:scalar}. Assume $b_r\ge0$, the conditions \eqref{eq:critical} defining the critical regime, and the bounded window \eqref{eq:window}. There is a fixed
$\delta>0$ such that, for $|x-p|\le2\delta$,
\begin{equation}\label{eq:F-identity}
 \mathcal F_1'(x)=\tfrac12B_{\mathrm{logit}}(x)
                         +\tfrac14\nu'(x)w_2(x),
\end{equation}
and $A_n(x)$ admits the exact factorization
\begin{equation}\label{eq:H-calibration}
 A_n(x)=\sig'(\ell(x))\{2C_n'(x)-g(x)w_2(x)\}.
\end{equation}
With $a_\star=-\Phi^{(2\mc)}(p)/(2\mc-1)!$ and
$|u|\le\delta$, one has
\begin{equation}\label{even:local-expansions}
 \begin{aligned}
 g(p+u)&=-2\nu(p)a_\star u^{2\mc-1}+O(|u|^{2\mc}),\\
 ug(p+u)&\le-c|u|^{2\mc},\\
 |\sig'(\ell(p+u))-\nu(p)|&\le C|u|,
 \end{aligned}
\end{equation}
The estimates in \eqref{even:local-expansions} depend only on the leading
parameters and therefore hold throughout the bounded window. Under the
additional calibration
\eqref{even:calibration-system},
\begin{equation}\label{eq:correction-taylor}
 C_n'(p+u)=\sum_{j=1}^{\mc}\frac{\beta_{j,n}}{(j-1)!}
               n^{-1+j/\mc}u^{j-1}+O(|u|^{\mc}).
\end{equation}
The remainder is uniform over bounded calibrated windows.
\end{lemma}

\begin{proof}
Choose $0<2\delta\le\delta_0$ so that $[p-2\delta,p+2\delta]\subset J$.
On this interval, $1-2\nu(x)t(x)\ge c_J>0$ and
$v_0(x)=\nu(x)/(1-2\nu(x)t(x))$.
Differentiating \eqref{eq:correction} and using
\eqref{eq:logit-profiles}, \eqref{v8:setup:two-row-variance-profile}, and
\eqref{eq:averaged-logit-correction} gives
\[
 \begin{aligned}
 &\frac{\dd}{\dd x}\{-\nu(x)t(x)-\tfrac12\log(1-2\nu(x)t(x))\}\\
 &\qquad=2t(x)^2v_0(x)\nu'(x)+t'(x)\{v_0(x)-\nu(x)\},\\
 &\frac{\dd}{\dd x}\{3\tau(x)^2\nu(x)^3\}\\
 &\qquad=6\tau(x)\tau'(x)\nu(x)^3+9\tau(x)^2\nu(x)^2\nu'(x),\\
 &\mathcal F_1'(x)=-s'(x)+t'(x)\{v_0(x)-\nu(x)\}
                     +6\tau(x)\tau'(x)\nu(x)^3\\
 &\quad+\nu'(x)\{2t(x)^2v_0(x)+9\tau(x)^2\nu(x)^2\}\\
 &\qquad=\tfrac12B_{\mathrm{logit}}(x)+\tfrac14\nu'(x)w_2(x).
 \end{aligned}
\]
This proves \eqref{eq:F-identity}. Substitute it into
\eqref{eq:a-function} and \eqref{v8:setup:drift-calibration-coefficients}.
The definition \eqref{even:calibration} and
$\sig''(\ell(x))/\sig'(\ell(x))=1-2\sig(\ell(x))=\nu'(x)-2g(x)$ give
\[
 \begin{aligned}
 A_n(x)&=\sig'(\ell(x))\bigl\{B_{\mathrm{logit}}(x)
                   +\tfrac12(\nu'(x)-2g(x))w_2(x)\\
 &\quad+\kappa_n+2\mathcal U_n'(x)\bigr\}\\
 &=\sig'(\ell(x))\{2C_n'(x)-g(x)w_2(x)\},
 \end{aligned}
\]
proving \eqref{eq:H-calibration} without derivative calibration.

Next, \eqref{eq:scalar}, \eqref{eq:critical}, and
\eqref{v8:def:scalar-logistic-field} give $\ell(x)=I'(x)+2\Phi'(x)$,
$\sig(I'(x))=x$, and
\[
 \Phi'(p+u)=-a_\star u^{2\mc-1}+O(|u|^{2\mc}),\qquad
 g(x)=2\nu(x)\Phi'(x)+O(\Phi'(x)^2).
\]
Reducing $\delta$ if necessary, the preceding expansion and smoothness give
\[
 \begin{aligned}
 g(p+u)&=-2\nu(p)a_\star u^{2\mc-1}+O(|u|^{2\mc}),\\
 ug(p+u)&=-2\nu(p)a_\star|u|^{2\mc}
                +O(|u|^{2\mc+1})\le-c|u|^{2\mc},\\
 |\sig'(\ell(p+u))-\nu(p)|&\le C|u|.
 \end{aligned}
\]
Finally, \eqref{eq:window} and \eqref{eq:correction} imply
$\sup_{n,|x-p|\le2\delta}|C_n^{(\mc+1)}(x)|<\infty$.
Taylor's formula and
$C_n^{(j)}(p)=n^{-1+j/\mc}\beta_{j,n}$ from \eqref{even:calibration} prove
\eqref{eq:correction-taylor}.
\end{proof}

The average conditional edge drift is
\begin{equation}\label{eq:drift-definition}
 \mathfrak d_n=N^{-1}\sum_e
   \{\Pp_n(X_e=1\mid X_{-e})-X_e\}
       =N^{-1}\sum_e\sig(L_e)-P.
\end{equation}

\begin{proposition}
\label{prop:weak-drift}
Assume \hyperref[ass:PC]{\textup{(PC)}}. For $A_n(x)$ defined in
\eqref{v8:setup:drift-calibration-coefficients} and every Lipschitz
$\phi:[0,1]\to\R$,
\begin{equation}\label{eq:weak-drift}
 \left|\E\left[\phi(P)
       \{\mathfrak d_n-g(P)-A_n(P)/n\}\right]\right|
 \le C\{n^{-2}\|\phi\|_\infty+n^{-5/2}\Lip\phi\}.
\end{equation}
The absolute value is taken after integrating the drift discrepancy
against the test function. The estimate is established under
\hyperref[ass:PC]{\textup{(PC)}} and does not require the calibration condition
\eqref{even:calibration-system} imposed in Theorem~\ref{thm:main}.
\end{proposition}

\begin{proof}
By \eqref{eq:scalar-field-window-difference}, for every fixed nonnegative integer $j$,
\[
 \ell_n(x)=\ell(x)+\frac{\kappa_n+2\mathcal U_n'(x)}n,
 \qquad\|\kappa_n+2\mathcal U_n'\|_{C^j}\le C_j.
\]
Taylor expansion at $\ell_n(P)$ gives
\begin{equation}\label{audit:dependency:s1011:weak-drift-taylor}
 \mathfrak d_n=\sig(\ell_n(P))-P
   +\sum_{j=1}^3\frac{\sig^{(j)}(\ell_n(P))}{j!}M_j+R_4,
 \qquad |R_4|\le\frac{\|\sig^{(4)}\|_\infty}{24}M_4.
\end{equation}
For $1\le j\le3$, the tests $x\mapsto\phi(x)\sig^{(j)}(\ell_n(x))$ satisfy
\[
 \begin{aligned}
 \|\phi(\cdot)\sig^{(j)}(\ell_n(\cdot))\|_\infty
     &\le C\|\phi\|_\infty,\\
 \Lip\{\phi(\cdot)\sig^{(j)}(\ell_n(\cdot))\}
     &\le C(\|\phi\|_\infty+\Lip\phi).
 \end{aligned}
\]
Applying \eqref{eq:micro-mean-input} and \eqref{eq:micro-higher-inputs}
to \eqref{audit:dependency:s1011:weak-drift-taylor} with these tests gives
\begin{equation}\label{audit:dependency:s1011:weak-drift-retained-profile}
 \begin{aligned}
 &\left|\E\left[\phi(P)\left\{\mathfrak d_n-\sig(\ell_n(P))+P
 -\frac{\sig'(\ell_n(P))B_{\mathrm{logit}}(P)
       +\tfrac12\sig''(\ell_n(P))w_2(P)}n\right\}\right]\right|\\
 &\qquad\le C\{n^{-2}\|\phi\|_\infty+n^{-5/2}\Lip\phi\}.
 \end{aligned}
\end{equation}
Uniformly in $x\in[0,1]$, bounded derivatives yield
\begin{equation}\label{audit:dependency:s1011:weak-drift-profile-expansion}
 \begin{gathered}
 \sig(\ell_n(x))-x
   =g(x)+n^{-1}\sig'(\ell(x))
                    \{\kappa_n+2\mathcal U_n'(x)\}+O(n^{-2}),\\
 n^{-1}\{\sig'(\ell_n(x))B_{\mathrm{logit}}(x)
       +\tfrac12\sig''(\ell_n(x))w_2(x)\}
   =n^{-1}\mathfrak a(x)+O(n^{-2}).
 \end{gathered}
\end{equation}
By \eqref{v8:setup:drift-calibration-coefficients}, the sum in
\eqref{audit:dependency:s1011:weak-drift-profile-expansion}
is $g(x)+A_n(x)/n+O(n^{-2})$. Substitution into
\eqref{audit:dependency:s1011:weak-drift-retained-profile} proves
\eqref{eq:weak-drift}.
\end{proof}

\subsection{Proof of the critical edge-density approximation}\label{sec:proof-main}
We first use stationarity to control moments at the critical scale.
We then identify the scaled drift and verify the remaining
hypotheses of the exchangeable-pair bound. Recall
\begin{equation}\label{v8:setup:main-proof-target}
 \begin{gathered}
 W_n=n^{1/\mc}(P-p),\qquad
 a_\star=-\frac{\Phi^{(2\mc)}(p)}{(2\mc-1)!}>0,\\
 V_n'(w)=a_\star w^{2\mc-1}
              -\sum_{j=1}^{\mc}\frac{\beta_{j,n}}{(j-1)!}w^{j-1}.
 \end{gathered}
\end{equation}

\begin{proof}[Proof of Theorem~\ref{thm:main}]
\proofstep{1. Critical moments.}
Equation~\eqref{even:local-expansions} gives
\begin{equation}\label{v8:proof:main:restoring-input}
 (P-p)g(P)\le-c|P-p|^{2\mc}\qquad(|P-p|\le\delta).
\end{equation}
Equation~\eqref{row:eq:macro} and
$|P-p|\le\|\mathcal W_X-p\|_\square+n^{-1}$ give
the bounds
\[
 \Pp_n(|P-p|>\delta)\le Ce^{-cn^2},\qquad
 c\E|P-p|^{2\mc}\le-\E[(P-p)g(P)]+Ce^{-cn^2}.
\]
Given $X\sim\Pp_n$, choose $e$ uniformly from $\binom{[n]}2$,
independently of $X$, and resample that edge from its true conditional law:
\[
 \begin{gathered}
 X'_e\mid(X,e)\sim\operatorname{Bernoulli}(\sig(L_e)),\qquad
 X'_f=X_f\quad(f\ne e),\\
 P':=\frac1N\sum_{f\in\binom{[n]}2}X'_f
     =P+\frac{X'_e-X_e}{N}.
 \end{gathered}
\]
Conditional on $X$, probabilities and expectations are taken over
edge selection and resampling.
Reversibility makes $(X,X')$ exchangeable, so
\[
 \begin{aligned}
 \E[(P'-P)(P'-p)]
 &=-\E[(P'-P)(P-p)],\\
\E(P'-P)^2&=-2\E[(P-p)(P'-P)].
 \end{aligned}
\]
Conditioning on $X$ therefore gives
\begin{equation}\label{even:quadratic-stationarity}
 \begin{gathered}
 \E[P'-P\mid X]=\mathfrak d_n/N,\qquad |P'-P|\le N^{-1},\\
 0=\frac2N\E_{\Pp_n}[(P-p)\mathfrak d_n]+\E(P'-P)^2,\qquad
 0\le-\E_{\Pp_n}[(P-p)\mathfrak d_n]\le(2N)^{-1}.
 \end{gathered}
\end{equation}
Test \eqref{eq:weak-drift} with $x-p$. On $|P-p|\le\delta$,
\eqref{eq:H-calibration} and \eqref{eq:correction-taylor} imply
\begin{equation}\label{audit:dependency:s1011:main-calibrated-correction}
 \frac{|(P-p)A_n(P)|}{n}
 \le C\sum_{j=1}^{\mc}n^{-2+j/\mc}|P-p|^j
       +Cn^{-1}|P-p|^{\mc+1}+Cn^{-1}|P-p|^{2\mc}.
\end{equation}
Outside this event, the uniform bound on $A_n(P)$ absorbs the exponentially
small contribution. Combining \eqref{v8:proof:main:restoring-input},
\eqref{even:quadratic-stationarity}, and \eqref{eq:weak-drift}
with \eqref{audit:dependency:s1011:main-calibrated-correction},
and multiplying by $n^2$, gives
\begin{equation}\label{audit:dependency:s1011:main-moment-bootstrap}
 c\E_{\Pp_n}|W_n|^{2\mc}
 \le C+C\sum_{j=1}^{\mc}\E_{\Pp_n}|W_n|^j
       +Cn^{-1/\mc}\E_{\Pp_n}|W_n|^{\mc+1}
       +Cn^{-1}\E_{\Pp_n}|W_n|^{2\mc}.
\end{equation}
Since $\mc+1<2\mc$, Young's inequality
$C|w|^q\le\varepsilon|w|^{2\mc}+C_{\varepsilon,q}$ for $q<2\mc$
absorbs the lower powers in \eqref{audit:dependency:s1011:main-moment-bootstrap};
taking $n$ sufficiently large absorbs its last term. H\"older's inequality then yields
\begin{equation}\label{eq:critical-moment}
\begin{gathered}
 \E_{\Pp_n}|P-p|^{2\mc}\le Cn^{-2},\qquad
 \E_{\Pp_n}|W_n|^{2\mc}\le C,\\
 \E_{\Pp_n}|P-p|^q\le C_qn^{-q/\mc}\quad(0\le q\le2\mc).
\end{gathered}
\end{equation}

\proofstep{2. Scaled drift.}
For bounded Lipschitz $f$, test \eqref{eq:weak-drift} with
$\phi(x)=f(n^{1/\mc}(x-p))$. Its error after multiplication by
$n^{2-1/\mc}$ is
$C\{n^{-1/\mc}\|f\|_\infty+n^{-1/2}\Lip f\}$.
On $|P-p|\le\delta$, \eqref{even:local-expansions} and
\eqref{eq:correction-taylor} give
\[
 \begin{aligned}
 n^{2-1/\mc}g(P)
   &=-2\nu(p)a_\star W_n^{2\mc-1}
         +O(n^{-1/\mc}|W_n|^{2\mc}),\\
 n^{1-1/\mc}C_n'(P)
   &=\sum_{j=1}^{\mc}\frac{\beta_{j,n}}{(j-1)!}W_n^{j-1}+O(n^{-1/\mc}|W_n|^\mc),\\
 |\sig'(\ell(P))-\nu(p)|&\le Cn^{-1/\mc}|W_n|,\qquad
 n^{1-1/\mc}|g(P)w_2(P)|\le Cn^{-1}|W_n|^{2\mc-1}.
 \end{aligned}
\]
Substituting these estimates into \eqref{eq:H-calibration} and
\eqref{v8:setup:main-proof-target} yields
\begin{equation}\label{audit:dependency:s1011:scaled-drift-error}
 \begin{aligned}
 &\left|n^{2-1/\mc}\{g(P)+A_n(P)/n\}+2\nu(p)V_n'(W_n)\right|\\
 &\quad\le Cn^{-1/\mc}
       \left(|W_n|^{2\mc}+|W_n|^\mc+
                     \sum_{j=1}^{\mc}|W_n|^j\right)
         +Cn^{-1}|W_n|^{2\mc-1}.
 \end{aligned}
\end{equation}
Equation~\eqref{eq:critical-moment} bounds the expectation of
\eqref{audit:dependency:s1011:scaled-drift-error}.
The exceptional event contributes exponentially little, since the
integrands grow polynomially in $n$. Combining \eqref{eq:weak-drift} and
\eqref{audit:dependency:s1011:scaled-drift-error}, with
$n^{-1/2}\le n^{-1/\mc}$, gives
\begin{equation}\label{even:scaled-regression}
 \left|\E\left[
  \{n^{2-1/\mc}\mathfrak d_n+2\nu(p)V_n'(W_n)\}f(W_n)\right]\right|
 \le Cn^{-1/\mc}(\|f\|_\infty+\Lip f).
\end{equation}

\proofstep{3. Exchangeable-pair errors.}
It remains to control the variance and increment errors in
Theorem~\ref{prop:Stein}. The heat-bath update is reversible, so
$(W_n,W_n')$ is exchangeable. Define
\begin{equation}\label{even:exchangeable-scales}
 \Delta=W_n'-W_n,\qquad
 \delta_n=\frac{n^{1/\mc}}N,\qquad
 \Lambda_n=\frac{2\nu(p)}{Nn^{2-2/\mc}}.
\end{equation}
Then
$\Delta\in\{-\delta_n,0,\delta_n\}$ and
$\E[\Delta\mid X]/\Lambda_n=n^{2-1/\mc}\mathfrak d_n/(2\nu(p))$,
so \eqref{even:scaled-regression} gives the weak regression bound
required in \eqref{stein:weak-regression}:
\[
 \sup_{\|f\|_\infty+\Lip f\le1}
 \left|\E_{\Pp_n}\left[
 \left\{\frac{\E[\Delta\mid X]}{\Lambda_n}+V_n'(W_n)\right\}
 f(W_n)\right]\right|\le Cn^{-1/\mc}.
\]
Since a flip has squared increment $\delta_n^2$,
\begin{equation}\label{even:conditional-variance-normalization}
 D:=\frac{\E[\Delta^2\mid X]}{2\Lambda_n}
       =\frac{n^2}{4\nu(p)N}\Pp(P'\ne P\mid X),\qquad 0\le D\le \nu(p)^{-1}.
\end{equation}
Equations~\eqref{eq:micro-higher-inputs} and \eqref{eq:critical-moment} imply
\begin{equation}\label{audit:dependency:s1011:averaged-flip-probability-error}
 \begin{aligned}
 \frac1N\sum_e\E|\sig(L_e)-p|
 &\le C(\E M_4)^{1/4}+C/n+C\E|P-p|\\
 &\le C(n^{-1/2}+n^{-1/\mc})\le Cn^{-1/\mc}.
 \end{aligned}
\end{equation}
Here Jensen's and H\"older's inequalities bound the averaged centered logit by
$(\E M_4)^{1/4}$, while
$|\sig(\ell_n(P))-p|\le C|P-p|+C/n$.
The exact identities
\[
 \begin{aligned}
 \Pp(P'\ne P\mid X)&=\frac1N\sum_e\{\sig(L_e)(1-X_e)+(1-\sig(L_e))X_e\},\\
 \Pp(P'\ne P\mid X)-2\nu(p)
   &=(1-2p)(P-p)+\frac1N\sum_e(\sig(L_e)-p)(1-2X_e),\\
 D-1&=\frac{n^2}{4\nu(p)N}(\Pp(P'\ne P\mid X)-2\nu(p))+\frac1{n-1}
 \end{aligned}
\]
together with \eqref{eq:critical-moment} and
\eqref{audit:dependency:s1011:averaged-flip-probability-error}, give
\begin{equation}\label{audit:dependency:s1011:conditional-variance-error}
 \E|1-D|\le Cn^{-1/\mc}.
\end{equation}
Finally, $|\Delta|^3=\delta_n\Delta^2$,
\eqref{even:exchangeable-scales}, \eqref{even:conditional-variance-normalization},
and \eqref{eq:critical-moment} yield
\begin{equation}\label{audit:dependency:s1011:exchangeable-increment-errors}
 \begin{aligned}
 \Lambda_n^{-1}\E[|\Delta|^3(1+|W_n|^{2\mc-1})]
   &=2\delta_n\E[D(1+|W_n|^{2\mc-1})]
     \le Cn^{-2+1/\mc},\\
 \delta_n\E|V_n'(W_n)|&\le Cn^{-2+1/\mc}.
 \end{aligned}
\end{equation}
The target parameters
$a=a_\star$, $c_{j,n}=-\beta_{j,n}/j!$ for $j\le\mc$ and
$c_{j,n}=0$ otherwise lie in a fixed compact subset of
$(0,\infty)\times\R^{2\mc-1}$.
Equations~\eqref{eq:critical-moment}, \eqref{even:scaled-regression},
\eqref{audit:dependency:s1011:conditional-variance-error}, and
\eqref{audit:dependency:s1011:exchangeable-increment-errors}
therefore verify the hypotheses of Theorem~\ref{prop:Stein}, with
$W=W_n$ and $\Lambda=\Lambda_n$.
Applying \eqref{eq:Stein-bound} and its moment extension for $\mc\ge2$,
with target moments evaluated under the density $\rho_n$ in \eqref{eq:target}, gives
\begin{equation}\label{eq:main-moment-comparison}
 \begin{aligned}
 \dW(W_n,Y_n)+\dK(W_n,Y_n)
 &\le C\{n^{-1/\mc}+n^{-2+1/\mc}\}\le Cn^{-1/\mc},\\
 \sum_{q=1}^2|\E W_n^q-\E_{\rho_n}Y_n^q|
 &\le Cn^{-1/\mc}.
 \end{aligned}
\end{equation}
Together with \eqref{eq:critical-moment}, this proves \eqref{eq:main-rate}.

By \eqref{stein:quartic-common-bounds},
\begin{equation}\label{eq:main-target-moments}
 \E Y_n^2\le C,\qquad
 \Var(Y_n)\ge c\int_{-1}^1(w-\E Y_n)^2\,\dd w
             \ge\frac{2c}{3}>0.
\end{equation}
Thus \eqref{eq:main-moment-comparison} also gives $|\E W_n|+|\E Y_n|\le C$ and
\begin{equation}\label{eq:main-variance-comparison}
 \begin{aligned}
 &|\Var(W_n)-\Var(Y_n)|\\
 \le{}&|\E W_n^2-\E Y_n^2|
       +|\E W_n-\E Y_n|(|\E W_n|+|\E Y_n|)
 \le Cn^{-1/\mc}.
 \end{aligned}
\end{equation}
Finally, $E=Np+Nn^{-1/\mc}W_n$ and
\eqref{eq:main-moment-comparison}--\eqref{eq:main-variance-comparison} yield
\[
 \begin{aligned}
 \E E&=Np+Nn^{-1/\mc}\E Y_n+O(Nn^{-2/\mc}),\\
 \Var(E)&=N^2n^{-2/\mc}
              \{\Var(Y_n)+O(n^{-1/\mc})\}\\
 &=N^2n^{-2/\mc}\Var(Y_n)\{1+O(n^{-1/\mc})\}.
 \end{aligned}
\]
Using $N=O(n^2)$ gives \eqref{eq:main-edge-mean-variance}.
\end{proof}

\subsection{Center, moments and variance under fixed perturbations}\label{sec:fixed-center}
We work under the hypotheses of Theorem~\ref{fixed:main}.
Proposition~\ref{prop:weak-drift} applies to \eqref{fixed:parameters}
without the calibration condition \eqref{even:calibration-system}.
The functions $C_n$ and $A_n$ are independent of $n$, and
\eqref{eq:H-calibration} gives $A_n(p)=2\nu(p)C_n'(p)\ne0$.

Choose a sufficiently small fixed $\delta>0$ and
$\chi\in C^\infty([0,1];[0,1])$ such that
\begin{equation}\label{fixed:cutoff}
 2\delta<\min\{p,1-p\},\qquad
 \chi(x)=
 \begin{cases}
 1,&|x-p|\le\delta,\\
 0,&|x-p|\ge2\delta.
 \end{cases}
\end{equation}
Using $A_n$ in \eqref{v8:setup:drift-calibration-coefficients}, define, for $x\in[0,1]$,
\begin{equation}\label{fixed:drift-definition}
 \begin{aligned}
 G_n(x)=\frac{A_n(x)}n
 +\int_p^x\left[\chi(u)g'(u)
                  -(1-\chi(u))(u-p)^{2\mc-2}\right]\,\dd u.
 \end{aligned}
\end{equation}
Since $g(p)=0$, \eqref{fixed:cutoff} and
\eqref{fixed:drift-definition} give
\begin{equation}\label{fixed:local-identification}
 G_n(x)=g(x)+\frac{A_n(x)}n
 \qquad(|x-p|\le\delta).
\end{equation}

\begin{lemma}\label{fixed:center}
One has $\mathcal F_1'(p)<0$. For all sufficiently large $n$,
$G_n$ in \eqref{fixed:drift-definition} has a unique zero $p_n\in(0,1)$, with
\begin{equation}\label{fixed:center-slope}
 \begin{aligned}
 p_n&=p+\sqrt[2\mc-1]{\frac{C_n'(p)}{a_\star n}}
                              +O(n^{-2/(2\mc-1)}),\\
 \lambda_n:=-G_n'(p_n)
 &=2\nu(p)(2\mc-1)a_\star^{1/(2\mc-1)}|C_n'(p)|^{(2\mc-2)/(2\mc-1)}\\
 &\quad\times n^{-(2\mc-2)/(2\mc-1)}\{1+O(n^{-1/(2\mc-1)})\}>0.
 \end{aligned}
\end{equation}
For $2\le q\le2\mc$,
\begin{equation}\label{fixed:moment-bounds}
 \E|P-p_n|^q\le C_qn^{-(2\mc-2+q)/(2\mc-1)},\qquad
 |\E P-p_n|\le Cn^{-1}.
\end{equation}
Moreover,
\begin{equation}\label{fixed:variance-asymptotics}
 \Var(P)=\frac{\nu(p)}{N\lambda_n}
                   \{1+O(n^{-1/(2\mc-1)})\}.
\end{equation}
\end{lemma}
\begin{proof}
\proofstep{1. Stable center.}
By \eqref{eq:critical}, $b_r\ge0$, and the simple-graph counting bounds,
\begin{equation}\label{fixed:critical-identities}
 \begin{gathered}
 B''(p)=\frac1{2\nu(p)},\qquad
 B'''(p)=\frac{2p-1}{2\nu(p)^2}\ge0,\qquad p\ge\tfrac12,\\
 a(F)\le\binom{e(F)}2,\qquad
 3\triangle(F)\le a(F),\qquad
 \binom{v(F)}2-1\ge e(F)-1.
 \end{gathered}
\end{equation}
Applying these inequalities termwise in \eqref{eq:profiles} gives
\begin{equation}\label{fixed:profile-bounds}
 \begin{gathered}
 t(p)\le\tfrac12B''(p),\quad t'(p)\le\tfrac12B'''(p),\quad
 s'(p)\ge pB''(p),\\
 \tau(p)\le\frac{B''(p)}{6p},\qquad
 \tau'(p)\le\frac{B'''(p)}{6p}.
 \end{gathered}
\end{equation}
Equations~\eqref{fixed:critical-identities} and \eqref{fixed:profile-bounds},
together with \eqref{eq:logit-profiles} and
\eqref{v8:setup:two-row-variance-profile}, give
$0\le v_0(p)-\nu(p)\le\nu(p)$, $w_2(p)\ge0$ and $\nu'(p)\le0$.
Substituting \eqref{eq:averaged-logit-correction} into \eqref{eq:F-identity}
and applying \eqref{fixed:profile-bounds} yields
\begin{equation}\label{fixed:negative-bound}
 \begin{aligned}
 \mathcal F_1'(p)
 &=-s'(p)+t'(p)\{v_0(p)-\nu(p)\}
      +6\tau(p)\tau'(p)\nu(p)^3+\tfrac14\nu'(p)w_2(p)\\
 &\le-s'(p)+\nu(p)t'(p)+6\tau(p)\tau'(p)\nu(p)^3\\
 &\le-\frac p{2\nu(p)}+\frac{\nu(p)}2B'''(p)
       +\frac{\nu(p)^3}{6p^2}B''(p)B'''(p)\\
 &=-\frac1{4\nu(p)}+\frac{2p-1}{24p^2}\le-\frac{23}{24}<0.
 \end{aligned}
\end{equation}
Here $2\nu(p)t(p)\le1/2$, $\nu(p)\le1/4$ and $2p-1\le p^2$.
Thus zero corrections satisfy the hypothesis $C_n'(p)\ne0$.
For general fixed corrections, Taylor's formula, \eqref{eq:H-calibration},
\eqref{even:local-expansions} and
\eqref{v8:setup:drift-calibration-coefficients} give,
after reducing $\delta$ if necessary, uniformly for $|u|\le2\delta$,
\begin{equation}\label{fixed:local-drift}
 \begin{aligned}
 A_n(p+u)&=2\nu(p)C_n'(p)+O(|u|),\qquad
          \max_{0\le j\le3}\|A_n^{(j)}\|_\infty\le C,\\
 g'(p+u)&=-2\nu(p)(2\mc-1)a_\star u^{2\mc-2}+O(|u|^{2\mc-1}),\\
 |g^{(j)}(p+u)|&\le C|u|^{2\mc-1-j},\qquad j=2,3.
 \end{aligned}
\end{equation}
Thus $\delta$ in \eqref{fixed:cutoff} can be chosen so that
$-g'(x)\asymp|x-p|^{2\mc-2}$ on $|x-p|\le2\delta$.
Differentiating \eqref{fixed:drift-definition}, using
\eqref{fixed:local-drift} and the cutoff \eqref{fixed:cutoff}, gives, for
$x\in[0,1]$,
\begin{equation}\label{fixed:extension-bounds}
 \begin{gathered}
 -G_n'(x)=\chi(x)(-g'(x))+(1-\chi(x))(x-p)^{2\mc-2}
                                      -\frac{A_n'(x)}n,\\
 c|x-p|^{2\mc-2}-C/n\le-G_n'(x)
                     \le C|x-p|^{2\mc-2}+C/n,\\
 |G_n^{(j)}(x)|\le C\{|x-p|^{2\mc-1-j}+n^{-1}\},\qquad j=2,3.
 \end{gathered}
\end{equation}
Substituting \eqref{even:local-expansions} and \eqref{fixed:local-drift}
into \eqref{fixed:local-identification} and its derivative gives, for each
fixed $R<\infty$ and all sufficiently large $n$, uniformly for $|u|\le R$,
\begin{equation}\label{fixed:scaled-drift}
 \begin{aligned}
 &nG_n(p+n^{-1/(2\mc-1)}u) =2\nu(p)[C_n'(p)-a_\star u^{2\mc-1}]+O_R(n^{-1/(2\mc-1)}),\\
 &n^{(2\mc-2)/(2\mc-1)}G_n'(p+n^{-1/(2\mc-1)}u)\\
 &=-2\nu(p)(2\mc-1)a_\star u^{2\mc-2}+O_R(n^{-1/(2\mc-1)}).
 \end{aligned}
\end{equation}
The limiting polynomial has the simple zero
$\sqrt[2\mc-1]{C_n'(p)/a_\star}\ne0$. Combining
\eqref{fixed:scaled-drift} with the local inverse function theorem and the
mean value theorem \cite[Theorems~9.24 and~5.10]{Rudin1976} gives a locally
unique zero $p_n$ of $G_n$ satisfying
\begin{equation}\label{fixed:center-location-error}
 \left|n^{1/(2\mc-1)}(p_n-p)
       -\sqrt[2\mc-1]{\frac{C_n'(p)}{a_\star}}\right|
       \le Cn^{-1/(2\mc-1)}.
\end{equation}
Equation~\eqref{fixed:center-location-error} gives the expansion of $p_n$
in \eqref{fixed:center-slope}. Substitution into the derivative estimate
in \eqref{fixed:scaled-drift} gives the expansion of $\lambda_n$.

Since $2\mc-2$ is even, for all $u,v\in\R$,
\[
 \begin{aligned}
 &2(u-v)(u^{2\mc-1}-v^{2\mc-1})
       -(u-v)^2(u^{2\mc-2}+v^{2\mc-2})\\
 &\qquad=(u^2-v^2)(u^{2\mc-2}-v^{2\mc-2})\ge0.
 \end{aligned}
\]
Together with
$|u-v|^{2\mc-2}\le2^{2\mc-3}(|u|^{2\mc-2}+|v|^{2\mc-2})$,
this gives
\begin{equation}\label{fixed:odd-power-monotonicity}
 (u-v)(u^{2\mc-1}-v^{2\mc-1})
 \ge c(u-v)^2\bigl(|v|^{2\mc-2}+|u-v|^{2\mc-2}\bigr).
\end{equation}
By \eqref{fixed:extension-bounds},
\[
 \frac{\dd}{\dd x}\left[
  -G_n(x)-\frac{c}{2\mc-1}(x-p)^{2\mc-1}+\frac Cn x
 \right]\ge0,\qquad x\in[0,1].
\]
Moreover, \eqref{fixed:center-location-error} gives
\[
 |p_n-p|^{2\mc-2}\asymp n^{-(2\mc-2)/(2\mc-1)},\qquad
 n^{-1}=o\bigl(n^{-(2\mc-2)/(2\mc-1)}\bigr).
\]
For $p_n+y\in[0,1]$, monotonicity gives a nonnegative product of the
increment from $p_n$ to $p_n+y$ with $y$.
Using $G_n(p_n)=0$ and \eqref{fixed:odd-power-monotonicity}, we obtain,
for all sufficiently large $n$,
\begin{equation}\label{fixed:coercivity}
 \begin{aligned}
 -yG_n(p_n+y)
 &\ge\frac{cy}{2\mc-1}
       \bigl[(p_n-p+y)^{2\mc-1}-(p_n-p)^{2\mc-1}\bigr]
       -\frac Cn y^2\\
 &\ge y^2\bigl[c_1(|p_n-p|^{2\mc-2}+|y|^{2\mc-2})-C/n\bigr]\\
 &\ge cy^2\{n^{-(2\mc-2)/(2\mc-1)}+|y|^{2\mc-2}\}.
 \end{aligned}
\end{equation}
In particular, \eqref{fixed:coercivity} proves that $p_n$ is the only zero
on $[0,1]$.

\proofstep{2. Moments and mean.}
By \eqref{fixed:local-identification}, boundedness of the profiles, and
Lemma~\ref{row:lem:macro}, \eqref{row:eq:macro}, every Lipschitz
$\phi:[0,1]\to\R$ satisfies
\begin{equation}\label{fixed:drift-replacement}
 \begin{aligned}
 &\left|\E\left[\phi(P)\left\{G_n(P)-g(P)-\frac{A_n(P)}n\right\}\right]\right|\\
 &\qquad\le C\|\phi\|_\infty\Pp_n(|P-p|>\delta)
       \le Ce^{-cn^2}\|\phi\|_\infty.
 \end{aligned}
\end{equation}
Both \eqref{eq:weak-drift} and \eqref{fixed:drift-replacement} apply to
$n$-dependent tests with the indicated supremum and Lipschitz norms.
For the single-edge Gibbs update in \eqref{even:quadratic-stationarity},
exchangeability gives
\begin{equation}\label{fixed:update}
 \begin{gathered}
 \E[P'-P\mid X]=\mathfrak d_n/N,\qquad
 \E[(P'-P)^2\mid X]=N^{-2}\Pp(P'\ne P\mid X),\\
 0=\E[(P'-p_n)^2-(P-p_n)^2]
   =\frac2N\E[(P-p_n)\mathfrak d_n]+\E(P'-P)^2,\\
 0\le-\E[(P-p_n)\mathfrak d_n]
 =\frac1{2N}\E\Pp(P'\ne P\mid X)\le\frac1{2N}.
 \end{gathered}
\end{equation}
For $\phi(x)=x-p_n$, one has $\|\phi\|_\infty\le1$ and $\Lip\phi=1$.
Applying \eqref{fixed:coercivity}, then \eqref{eq:weak-drift} and
\eqref{fixed:drift-replacement}, and finally \eqref{fixed:update}, gives
\begin{equation}\label{fixed:moment-start}
 \begin{aligned}
 n^{-(2\mc-2)/(2\mc-1)}\E(P-p_n)^2+\E|P-p_n|^{2\mc}
 &\le-C\E[(P-p_n)G_n(P)]\\
 &\le-C\E[(P-p_n)\mathfrak d_n]+Cn^{-2}\\
 &\le Cn^{-2}.
 \end{aligned}
\end{equation}
H\"older's inequality applied to \eqref{fixed:moment-start} gives, for
$2\le q\le2\mc$,
\begin{equation}\label{fixed:moment-interpolation}
 \begin{aligned}
 \E|P-p_n|^q
 &\le\{\E(P-p_n)^2\}^{(2\mc-q)/(2\mc-2)}
      \{\E|P-p_n|^{2\mc}\}^{(q-2)/(2\mc-2)}\\
 &\le C_qn^{-(2\mc-2+q)/(2\mc-1)}.
 \end{aligned}
\end{equation}
For the mean, Taylor's formula at $p_n$, using
\eqref{fixed:center-slope} and \eqref{fixed:extension-bounds}, followed by
\eqref{fixed:moment-interpolation}, gives
\begin{equation}\label{fixed:linearization}
 \begin{aligned}
 |G_n(p_n+y)+\lambda_ny|
 &\le Cy^2\int_0^1|G_n''(p_n+ty)|\,\dd t\\
 &\le C\{n^{-(2\mc-3)/(2\mc-1)}y^2+|y|^{2\mc-1}\},\\
 \E|G_n(P)+\lambda_n(P-p_n)|
 &\le Cn^{-(2\mc-3)/(2\mc-1)}\E(P-p_n)^2
                         +C\E|P-p_n|^{2\mc-1}\\
 &\le Cn^{-(4\mc-3)/(2\mc-1)}.
 \end{aligned}
\end{equation}
By \eqref{fixed:update}, $\E\mathfrak d_n=0$. Testing
\eqref{eq:weak-drift} with $\phi=1$, then using
\eqref{fixed:drift-replacement}, \eqref{fixed:linearization}, and
\eqref{fixed:center-slope}, gives
\begin{equation}\label{fixed:mean-center-error}
 \begin{aligned}
 \lambda_n|\E P-p_n|
 &\le|\E G_n(P)|+\E|G_n(P)+\lambda_n(P-p_n)|\\
 &\le C(n^{-2}+n^{-(4\mc-3)/(2\mc-1)}),\qquad
 |\E P-p_n|\le Cn^{-1}.
 \end{aligned}
\end{equation}

This proves \eqref{fixed:moment-bounds}.

\proofstep{3. Variance.}
Equations~\eqref{eq:critical}, \eqref{v8:def:scalar-logistic-field}
and \eqref{even:window-Cr} give
$|\sig(\ell_n(p))-p|\le C/n$.
Conditioning on the updated edge, then applying
\eqref{eq:micro-higher-inputs} with $\phi=1$,
\eqref{fixed:center-location-error}, and \eqref{fixed:moment-interpolation}, gives
\begin{equation}\label{fixed:flip-error}
 \begin{aligned}
 \Pp(P'\ne P\mid X)-2\nu(p)
 &=(1-2p)(P-p)+\frac1N\sum_e(\sig(L_e)-p)(1-2X_e),\\
 \frac1N\sum_e\E|\sig(L_e)-p|
 &\le C\{(\E M_2)^{1/2}+\E|P-p|+n^{-1}\}\\
 &\le Cn^{-1/(2\mc-1)},\\
 \E|\Pp(P'\ne P\mid X)-2\nu(p)|
 &\le |1-2p|\E|P-p|+\frac1N\sum_e\E|\sig(L_e)-p|\\
 &\le Cn^{-1/(2\mc-1)}.
 \end{aligned}
\end{equation}
For $x\in[0,1]$, set
\begin{equation}\label{fixed:adapted-test}
 \psi_n(x)=-\frac{(x-p_n)^2}{G_n(x)}\quad(x\ne p_n),\qquad
 \psi_n(p_n)=0.
\end{equation}
Since $G_n'(p_n)=-\lambda_n<0$, $\psi_n$ extends smoothly to $p_n$.
Differentiating $\psi_n(x)G_n(x)=-(x-p_n)^2$ twice gives, for $x\ne p_n$,
\[
 \begin{aligned}
 \psi_n'(x)
 &=-\frac{2(x-p_n)}{G_n(x)}
   +\frac{(x-p_n)^2G_n'(x)}{G_n(x)^2},\\
 \psi_n''(x)
 &=\frac{(x-p_n)^2G_n''(x)}{G_n(x)^2}
   -\frac{2[G_n(x)-(x-p_n)G_n'(x)]^2}{G_n(x)^3}.
 \end{aligned}
\]
Taylor's formula at $x$, using $G_n(p_n)=0$, gives
\[
 G_n(x)-(x-p_n)G_n'(x)
 =-\frac{(x-p_n)^2}{2}G_n''(\xi_x),
\]
where $\xi_x$ lies between $p_n$ and $x$.
These identities, \eqref{fixed:extension-bounds},
\eqref{fixed:center-location-error}, and \eqref{fixed:coercivity},
together with continuity at $p_n$, give
\begin{equation}\label{fixed:adapted-bounds}
 \begin{gathered}
 \psi_n(x)G_n(x)=-(x-p_n)^2,\qquad \psi_n'(p_n)=\lambda_n^{-1},\\
 \begin{aligned}
 \|\psi_n^{(j)}\|_\infty
 &\le C\sup_{x\in[0,1]}
 \bigl(n^{-1/(2\mc-1)}+|x-p_n|\bigr)^{-(2\mc+j-3)}\\
 &\le Cn^{(2\mc+j-3)/(2\mc-1)},\qquad j=0,1,2.
 \end{aligned}
 \end{gathered}
\end{equation}
Exchangeability, \eqref{fixed:update}, and Taylor's formula give the
first expansion below; \eqref{eq:weak-drift}, \eqref{fixed:drift-replacement},
and \eqref{fixed:adapted-bounds} give the second estimate:
\begin{equation}\label{fixed:adapted-stationarity}
 \begin{aligned}
 \E[\psi_n(P)\mathfrak d_n]
 &=-\frac N2\E[(P'-P)\{\psi_n(P')-\psi_n(P)\}]\\
 &=-\frac1{2N}\E[\Pp(P'\ne P\mid X)\psi_n'(P)]
             +O(N^{-2}\|\psi_n''\|_\infty),\\[3pt]
 |\E[\psi_n(P)\{\mathfrak d_n-G_n(P)\}]|
 &\le C\{(n^{-2}+e^{-cn^2})\|\psi_n\|_\infty
                       +n^{-5/2}\|\psi_n'\|_\infty\}\\
 &\le C\{n^{-\frac{2\mc+1}{2\mc-1}}
                       +n^{-5/2+\frac{2\mc-2}{2\mc-1}}\}
   \le Cn^{-\frac{2\mc+1}{2\mc-1}}.
 \end{aligned}
\end{equation}
By \eqref{fixed:adapted-bounds},
$N^{-2}\|\psi_n''\|_\infty=O(n^{-3})$; hence
\eqref{fixed:adapted-stationarity} gives
\begin{equation}\label{fixed:variance-identity}
 \E(P-p_n)^2
 =\frac1{2N}\E[\Pp(P'\ne P\mid X)\psi_n'(P)]
                       +O(n^{-(2\mc+1)/(2\mc-1)}).
\end{equation}
From \eqref{fixed:adapted-bounds} and \eqref{fixed:moment-interpolation},
followed by \eqref{fixed:flip-error}, we obtain
\begin{equation}\label{fixed:variance-coefficient-errors}
 \begin{aligned}
 |\E\psi_n'(P)-\lambda_n^{-1}|
 &\le\|\psi_n''\|_\infty\E|P-p_n|\le Cn^{(\mc-1)/(2\mc-1)},\\
 |\E[(\Pp(P'\ne P\mid X)-2\nu(p))\psi_n'(P)]|
 &\le Cn^{-1/(2\mc-1)}\|\psi_n'\|_\infty\le Cn^{(2\mc-3)/(2\mc-1)}.
 \end{aligned}
\end{equation}
Equation~\eqref{fixed:variance-identity}, with
\eqref{fixed:variance-coefficient-errors}, followed by
\eqref{fixed:mean-center-error} and \eqref{fixed:center-slope}, gives
\begin{equation}\label{fixed:variance-expansion}
 \begin{aligned}
 \E(P-p_n)^2
 &=\frac{\nu(p)}{N\lambda_n}+O(n^{-(2\mc+1)/(2\mc-1)}),\\
 \Var(P)
 &=\E(P-p_n)^2-(\E P-p_n)^2
 =\frac{\nu(p)}{N\lambda_n}\{1+O(n^{-1/(2\mc-1)})\}.
 \end{aligned}
\end{equation}
This proves \eqref{fixed:variance-asymptotics}.

\end{proof}

\subsection{Proof of Gaussian approximation under fixed perturbations}\label{sec:proof-fixed-critical-normal}
\begin{proof}[Proof of Theorem~\ref{fixed:main}]
\proofstep{1. Gaussian approximation.}
We apply Theorem~\ref{prop:Stein}, \eqref{eq:Stein-bound}, with
$V(w)=w^2/2$, so that the target distribution is standard normal.
Use the single-edge update in \eqref{fixed:update}.
With $Z_n$ from \eqref{fixed:normalized-variable}, set
\begin{equation}\label{fixed:normalized-increment}
 Z_n'=(P'-p_n)\sqrt{\frac{N\lambda_n}{\nu(p)}},\qquad
 \Delta=Z_n'-Z_n.
\end{equation}
Then $(Z_n,Z_n')$ is exchangeable and $\E Z_n^2\le C$ by
\eqref{fixed:center-slope} and \eqref{fixed:moment-start}.
For bounded Lipschitz $f$, take
$\phi(x)=f((x-p_n)\sqrt{N\lambda_n/\nu(p)})$. Then
\[
 \|\phi\|_\infty\le\|f\|_\infty,\qquad
 \Lip\phi\le\sqrt{N\lambda_n/\nu(p)}\,\Lip f.
\]
Equations~\eqref{eq:weak-drift}, \eqref{fixed:drift-replacement}, and
\eqref{fixed:linearization} give the weak regression estimate
\begin{equation}\label{fixed:weak-normal-regression}
 \begin{aligned}
 &\left|\E\left[\left\{\sqrt{\frac N{\nu(p)\lambda_n}}
                         \mathfrak d_n+Z_n\right\}f(Z_n)\right]\right|\\
 \le{}&\sqrt{\frac N{\nu(p)\lambda_n}}
 \Bigl[Cn^{-2}\|f\|_\infty
       +Cn^{-5/2}\sqrt{N\lambda_n/\nu(p)}\,\Lip f
       +\|f\|_\infty\E|G_n(P)+\lambda_n(P-p_n)|\Bigr]\\
 \le{}&C\{[n^{-(\mc-1)/(2\mc-1)}+n^{-\mc/(2\mc-1)}]\|f\|_\infty
                               +n^{-1/2}\Lip f\}\\
 \le{}&Cn^{-1/(2\mc-1)}(\|f\|_\infty+\Lip f).
 \end{aligned}
\end{equation}
In \eqref{eq:Stein-bound}, take $W=Z_n$,
$\Lambda=\lambda_n/N$ and $\delta=\sqrt{\lambda_n/(N\nu(p))}$;
\eqref{fixed:center-slope} ensures $\delta\le1$ for all sufficiently large $n$.
The update identities \eqref{even:quadratic-stationarity} and
\eqref{fixed:update}, together with \eqref{fixed:flip-error}, give
\begin{equation}\label{fixed:clock-variance}
 \begin{gathered}
 \frac{\E[\Delta\mid X]}\Lambda
 =\sqrt{\frac N{\nu(p)\lambda_n}}\mathfrak d_n,\qquad
 D:=\frac{\E[\Delta^2\mid X]}{2\Lambda}
     =\frac{\Pp(P'\ne P\mid X)}{2\nu(p)},\\
 \E|1-D|\le Cn^{-1/(2\mc-1)},\qquad
 0\le D\le[2\nu(p)]^{-1},\qquad |\Delta|\le\delta.
 \end{gathered}
\end{equation}
Equation~\eqref{fixed:clock-variance}, together with $\E Z_n^2\le C$ from
\eqref{fixed:center-slope} and \eqref{fixed:moment-start}, gives
\[
 \begin{gathered}
 \frac1\Lambda\E[|\Delta|^3(1+|Z_n|)]
 \le2\delta\E[D(1+|Z_n|)]\le C\delta,\qquad
 \delta\E|Z_n|\le\delta(\E Z_n^2)^{1/2}\le C\delta.
 \end{gathered}
\]
Thus \eqref{eq:Stein-bound},
\eqref{fixed:weak-normal-regression}, and \eqref{fixed:clock-variance} give
\begin{equation}\label{fixed:reference-bound}
 \begin{aligned}
 \dW(Z_n,G)+\dK(Z_n,G)
 &\le C\left\{n^{-1/(2\mc-1)}
                  +\sqrt{\frac{\lambda_n}{N\nu(p)}}\right\}\\
 &\le Cn^{-1/(2\mc-1)}.
 \end{aligned}
\end{equation}
Equation~\eqref{fixed:reference-bound} proves \eqref{fixed:standard-normal}.

\proofstep{2. Mean and variance.}
By \eqref{fixed:moment-bounds} and \eqref{fixed:variance-asymptotics},
\[
 \E P=p_n+O(n^{-1}),\qquad
 \Var(E)=\frac{N\nu(p)}{\lambda_n}
                      \{1+O(n^{-1/(2\mc-1)})\}.
\]
Substituting \eqref{fixed:center-slope} and
$N=n^2\{1+O(n^{-1})\}/2$ gives \eqref{fixed:mean-variance}.
\end{proof}

\subsection{Proofs of the subgraph and standardized-count corollaries}
\label{sec:proof-corollaries}
We now deduce the subgraph approximations from the edge-density theorems
and the replacement estimate \eqref{eq:count-replacement}, and then pass
to standardization by the exact means and variances.
In each proof below, $H$ is fixed, and the constants $C_H$ and $C_{H,q}$
have the parameter dependence specified in the corresponding corollary
in Section~\ref{sec:standardized-counts}.

\begin{proof}[Proof of Corollary~\ref{even:subgraph-fluctuations}]
We separate the graph-count error from the Taylor remainder of $x^{e(H)}$.
\[
 \begin{aligned}
 &f_H(X)-p^{e(H)}-e(H)p^{e(H)-1}(P-p)\\
 ={}&\bigl[f_H(X)-P^{e(H)}\bigr]
   +\bigl[P^{e(H)}-p^{e(H)}-e(H)p^{e(H)-1}(P-p)\bigr].
 \end{aligned}
\]
Since $x\mapsto x^{e(H)}$ has bounded second derivative on $[0,1]$,
Taylor's formula bounds the second bracket in absolute value by
$C_H(P-p)^2$. Hence
\[
 \left|f_H(X)-p^{e(H)}-e(H)p^{e(H)-1}(P-p)\right|
 \le |f_H(X)-P^{e(H)}|+C_H(P-p)^2.
\]
Multiplying by $n^{1/\mc}/\{e(H)p^{e(H)-1}\}$ and taking $L^j$ norms,
we obtain from \eqref{eq:count-replacement} and the moment bound in
\eqref{eq:main-rate}
\begin{equation}\label{audit:subgraph-linear-coupling}
 \begin{aligned}
 &\left\|\frac{n^{1/\mc}\{f_H(X)-p^{e(H)}\}}
                   {e(H)p^{e(H)-1}}-W_n\right\|_j\\
 &\quad\le C_Hn^{-1+1/\mc}+C_Hn^{1/\mc}\|P-p\|_{2j}^2
 \le C_Hn^{-1/\mc},\qquad j=1,2.
 \end{aligned}
\end{equation}
Coupling, \eqref{audit:subgraph-linear-coupling}, and
\eqref{eq:main-rate} prove the Wasserstein bound in
\eqref{even:subgraph-rate}.
For the Kolmogorov distance, compare distribution functions at $t-\delta$
and $t+\delta$. For $t\in\R$ and $\delta>0$,
\[
 \begin{aligned}
 &\Pp(W_n\le t-\delta)
 -\Pp\left(\left|\frac{n^{1/\mc}\{f_H(X)-p^{e(H)}\}}
                   {e(H)p^{e(H)-1}}-W_n\right|>\delta\right)\\
 \le{}&\Pp\left(\frac{n^{1/\mc}\{f_H(X)-p^{e(H)}\}}
                   {e(H)p^{e(H)-1}}\le t\right)\\
 \le{}&\Pp(W_n\le t+\delta)
 +\Pp\left(\left|\frac{n^{1/\mc}\{f_H(X)-p^{e(H)}\}}
                   {e(H)p^{e(H)-1}}-W_n\right|>\delta\right).
 \end{aligned}
\]
By \eqref{stein:quartic-common-bounds},
\[
 \sup_{t\in\R}\Pp(t-\delta<Y_n\le t+\delta)
 \le2\delta\|\rho_n\|_\infty\le C\delta.
\]
Taking the supremum over $t$ and applying Markov's inequality,
\eqref{eq:main-rate}, and \eqref{audit:subgraph-linear-coupling}, with
$\delta=n^{-1/(2\mc)}$, yields
\[
 \begin{aligned}
 &\dK\left(\frac{n^{1/\mc}\{f_H(X)-p^{e(H)}\}}
                   {e(H)p^{e(H)-1}},Y_n\right)\\
 \le{}&\dK(W_n,Y_n)+C\delta
 +\Pp\left(\left|\frac{n^{1/\mc}\{f_H(X)-p^{e(H)}\}}
                   {e(H)p^{e(H)-1}}-W_n\right|>\delta\right)\\
 \le{}&\dK(W_n,Y_n)+C\delta
       +\delta^{-1}\E\left|\frac{n^{1/\mc}\{f_H(X)-p^{e(H)}\}}
                   {e(H)p^{e(H)-1}}-W_n\right|\\
 \le{}&C_H\bigl(n^{-1/\mc}+\delta+\delta^{-1}n^{-1/\mc}\bigr)
 \le C_Hn^{-1/(2\mc)}.
 \end{aligned}
\]

The $j=1$ bound in \eqref{audit:subgraph-linear-coupling} controls the mean.
Centering the $j=2$ bound and applying the reverse triangle inequality in
$L^2$, together with \eqref{eq:main-edge-mean-variance} and
\eqref{stein:quartic-common-bounds}, gives
\begin{equation}\label{even:subgraph-moment-comparison}
 \begin{aligned}
 &\left|\frac{n^{1/\mc}\{\E f_H(X)-p^{e(H)}\}}
                   {e(H)p^{e(H)-1}}-\E Y_n\right|\\
 &\quad{}+\left|\frac{n^{2/\mc}\Var(f_H(X))}
                   {e(H)^2p^{2e(H)-2}}-\Var(Y_n)\right|
 \le C_Hn^{-1/\mc}.
 \end{aligned}
\end{equation}
Using $\inj(H,X)=n^{v(H)}f_H(X)$ and $\Var(Y_n)\ge c>0$ proves
\eqref{even:subgraph-mean-variance}.
\end{proof}

\begin{proof}[Proof of Corollary~\ref{fixed:subgraph-gaussian}]
We use the same decomposition, now centered at $p_n$:
\[
 \begin{aligned}
 &f_H(X)-p_n^{e(H)}-e(H)p_n^{e(H)-1}(P-p_n)\\
 ={}&\bigl[f_H(X)-P^{e(H)}\bigr]
   +\bigl[P^{e(H)}-p_n^{e(H)}-e(H)p_n^{e(H)-1}(P-p_n)\bigr].
 \end{aligned}
\]
Taylor's formula on $[0,1]$ bounds the second bracket in absolute value by
$C_H(P-p_n)^2$. The triangle inequality therefore gives
\begin{equation}\label{fixed:subgraph-centered-error}
 |f_H(X)-p_n^{e(H)}-e(H)p_n^{e(H)-1}(P-p_n)|
 \le |f_H(X)-P^{e(H)}|+C_H(P-p_n)^2.
\end{equation}
Since $p_n\to p\in(0,1)$ and
$\sqrt{N\lambda_n/\nu(p)}\asymp n^{\mc/(2\mc-1)}$,
\eqref{fixed:subgraph-centered-error}, \eqref{eq:count-replacement},
\eqref{fixed:center-slope}, and \eqref{fixed:moment-bounds} imply
\begin{equation}\label{fixed:subgraph-L2}
 \begin{aligned}
 &\left\|\sqrt{\frac{N\lambda_n}{\nu(p)}}
       \frac{f_H(X)-p_n^{e(H)}}{e(H)p_n^{e(H)-1}}-Z_n\right\|_j\\
 &\quad\le C_Hn^{\mc/(2\mc-1)}
                 \{n^{-1}+\|P-p_n\|_{2j}^2\}
 \le C_Hn^{-1/(2\mc-1)},\qquad j=1,2.
 \end{aligned}
\end{equation}
Coupling, \eqref{fixed:subgraph-L2}, and \eqref{fixed:standard-normal}
prove \eqref{fixed:subgraph-W}.
For fixed $q\ge1$ and $\delta>0$, apply the smoothing comparison in
Corollary~\ref{even:subgraph-fluctuations} and Markov's inequality separately
to the graph-count error and the quadratic Taylor remainder in
\eqref{fixed:subgraph-centered-error}. Using \eqref{eq:count-replacement},
\eqref{fixed:moment-bounds}, and \eqref{fixed:standard-normal} gives
\begin{equation}\label{fixed:subgraph-smoothing}
 \begin{aligned}
 &\dK\left(\sqrt{\frac{N\lambda_n}{\nu(p)}}
                   \frac{f_H(X)-p_n^{e(H)}}{e(H)p_n^{e(H)-1}},G\right)\\
 &\quad\le\dK(Z_n,G)+C_H\delta
       +C_{H,q}\delta^{-q}n^{q\mc/(2\mc-1)}
                         \E|f_H(X)-P^{e(H)}|^q\\
 &\qquad{}+C_H\delta^{-1}n^{\mc/(2\mc-1)}\E(P-p_n)^2\\
 &\quad\le C_H(n^{-1/(2\mc-1)}+\delta)
        +C_{H,q}n^{-q(\mc-1)/(2\mc-1)}\delta^{-q}
        +C_Hn^{-\mc/(2\mc-1)}\delta^{-1}.
 \end{aligned}
\end{equation}
For $\mc\ge3$, take $q=1$ and $\delta=n^{-1/(2\mc-1)}$; for $\mc=2$,
take $\delta=n^{-q/(3(q+1))}$. Substitution in
\eqref{fixed:subgraph-smoothing} proves \eqref{fixed:subgraph-K}.

Taking expectations in \eqref{fixed:subgraph-centered-error} controls the
mean. Centering its remainder and applying the reverse triangle inequality
in $L^2$ controls the standard deviation. Together with
\eqref{eq:count-replacement}, \eqref{fixed:moment-bounds}, and
\eqref{fixed:subgraph-L2}, this gives
\begin{equation}\label{fixed:subgraph-variance-transfer}
 \begin{gathered}
 |\E f_H(X)-p_n^{e(H)}|
 \le C_H\bigl(|\E P-p_n|+\E|f_H(X)-P^{e(H)}|
                         +\E(P-p_n)^2\bigr)\le C_Hn^{-1},\\
 \begin{aligned}
 &\left|\sqrt{\Var(f_H(X))}
       -e(H)p_n^{e(H)-1}\sqrt{\Var(P)}\right|\\
 \le{}&\|f_H(X)-p_n^{e(H)}
                   -e(H)p_n^{e(H)-1}(P-p_n)\|_2
 \le C_Hn^{-(\mc+1)/(2\mc-1)}.
 \end{aligned}
 \end{gathered}
\end{equation}
Using $\inj(H,X)=n^{v(H)}f_H(X)$, \eqref{fixed:center-slope}, and
\eqref{fixed:variance-asymptotics} now proves \eqref{fixed:subgraph-variance}.
\end{proof}

\begin{proof}[Proof of Corollary~\ref{even:standardized-counts}]
We use the same comparison for all four approximations. Choose
\[
 \begin{aligned}
 Y=Y_n:\quad &U=W_n\quad\text{or}\quad
       \frac{n^{1/\mc}\{f_H(X)-p^{e(H)}\}}{e(H)p^{e(H)-1}},\\
 Y=G:\quad &U=Z_n\quad\text{or}\quad
       \sqrt{\frac{N\lambda_n}{\nu(p)}}
       \frac{f_H(X)-p_n^{e(H)}}{e(H)p_n^{e(H)-1}}.
 \end{aligned}
\]
Equations~\eqref{eq:main-edge-mean-variance} and
\eqref{even:subgraph-moment-comparison} give the first case below;
\eqref{fixed:center-slope}, \eqref{fixed:mean-variance}, and
\eqref{fixed:subgraph-variance} give the second:
\begin{equation}\label{eq:main-standardization-errors}
 |\E U-\E Y|+|\Var(U)-\Var(Y)|
 \le\begin{cases}
 Cn^{-1/\mc},&Y=Y_n,\\
 Cn^{-1/(2\mc-1)},&Y=G.
 \end{cases}
\end{equation}
In the Gaussian case, the mean contribution is at most
$Cn^{\mc/(2\mc-1)}n^{-1}\le Cn^{-1/(2\mc-1)}$ since $\mc\ge2$.
These variance estimates, \eqref{stein:quartic-common-bounds}, and the
standard normal density also give
\[
 c\le\Var(U),\Var(Y)\le C,\qquad
 \E|Y|+\sup_{w\in\R}(1+|w|)\rho(w)\le C,
\]
where $\rho$ is the density of $Y$.

Match the mean and variance of $Y$ to those of $U$ by setting
\[
 a=\sqrt{\frac{\Var(U)}{\Var(Y)}},\qquad b=\E U-a\E Y.
\]
The preceding moment and variance bounds imply
\[
 |a-1|+|\log a|+|b|
 \le C\bigl[|\E U-\E Y|+|\Var(U)-\Var(Y)|\bigr].
\]
Coupling and the identity
\[
 \Pp(aY\le t)-\Pp(Y\le t)
 =-\int_1^a\frac{t}{s^2}\rho(t/s)\,\dd s
\]
give
\begin{equation}\label{eq:standardized-affine-comparison}
 \begin{aligned}
 \dW(aY+b,Y)&\le |a-1|\E|Y|+|b|,\\
 \dK(aY+b,Y)&\le\dK(aY,Y)+\dK(Y+b,Y)\\
 &\le\|\rho\|_\infty|b|
           +\sup_{w\in\R}|w|\rho(w)\,|\log a|.
 \end{aligned}
\end{equation}
Since $(aY+b-\E U)/\sqrt{\Var(U)}=(Y-\E Y)/\sqrt{\Var(Y)}$,
the scaling property of $\dW$, the affine invariance of $\dK$, and the
triangle inequality yield, for either $d\in\{\dW,\dK\}$,
\begin{equation}\label{eq:standardized-distance-transfer}
 \begin{aligned}
 d\left(\frac{U-\E U}{\sqrt{\Var(U)}},
         \frac{Y-\E Y}{\sqrt{\Var(Y)}}\right)
 &\le C d(U,aY+b)\\
 &\le C\bigl[d(U,Y)+d(aY+b,Y)\bigr]\\
 &\le C\bigl[d(U,Y)+|\E U-\E Y|+|\Var(U)-\Var(Y)|\bigr].
 \end{aligned}
\end{equation}
Each $U$ is a positive affine function of the corresponding edge or
subgraph count, so standardizing $U$ gives the variable in the statement.
Combining \eqref{eq:standardized-distance-transfer} and
\eqref{eq:main-standardization-errors} with \eqref{eq:main-rate},
\eqref{even:subgraph-rate}, \eqref{fixed:standard-normal}, and
\eqref{fixed:subgraph-W}--\eqref{fixed:subgraph-K} proves the claimed rates.
When $\mc=2$, the additional error $n^{-1/3}$ is bounded by
$n^{-q/(3(q+1))}$ for every fixed $q\ge1$.
\end{proof}

\appendix

\section{Graph localization and contractive row kernels}\label{row:sec:row-gap}\label{app:counting-smoothness}
We prove the localization used in Section~\ref{row:sec:conditional-rows}
and construct the auxiliary updates. The deterministic counting bounds
come first; the tolerances are fixed before the good exterior sets are chosen.

\subsection{Cut norms and weighted graph counts}
The cut norm and step graphons are defined in
Section~\ref{row:sec:conditional-rows}.

The finite cut formula is

\begin{equation}\label{v8:proof:row:lem:macro:finite-cut-formula}
 \|\mathcal W_X-p\|_\square
 =\frac1{n^2}\max_{U,V\subseteq[n]}
       \left|\sum_{i\in U,\,j\in V}(X_{ij}-p)\right|.
\end{equation}

Indeed, put $I_i=((i-1)/n,i/n]$ and
$a_i=n|B_1\cap I_i|$, $b_j=n|B_2\cap I_j|$. The layer-cake representation gives
\[
 \int_{B_1\times B_2}(\mathcal W_X-p)
 =n^{-2}\int_0^1\!\!\int_0^1
       \sum_{\substack{i:a_i>s\\j:b_j>t}}(X_{ij}-p)\,\dd s\,\dd t .
\]
Taking absolute values inside the outer integrals gives the bound on the right-hand side of the finite cut
formula. Conversely, taking $B_1=\bigcup_{i\in U}I_i$ and
$B_2=\bigcup_{j\in V}I_j$ proves equality, including the diagonal terms
$-p$. The argument applies at any graph size.

For a finite simple graph $H$, set

\begin{equation}\label{row:eq:hom-density}
 t(H,W)=\int_{[0,1]^{v(H)}}
          \prod_{\{a,b\}\in E(H)}W(x_a,x_b)
          \prod_{j=1}^{v(H)}\dd x_j.
\end{equation}

Thus $t(H,\mathcal W_X)=n^{-v(H)}\Hom(H,X)$, with each isolated vertex
contributing a factor of one. Define

\begin{equation}\label{v8:nearby:app:graphon-energy}
 T(W)=\frac{h_*}{2}t(K_2,W)+\sum_{r=1}^k b_rt(F_r,W).
\end{equation}

For fixed coefficients, define the reference law by

\begin{equation}\label{v8:nearby:app:reference-graph-law}
 \begin{aligned}
 Z_n^{\mathrm{graph},0}&=\sum_x\exp\{n^2T(\mathcal W_x)\},\qquad
 \Pp_n^0(X=x)=\frac{\exp\{n^2T(\mathcal W_x)\}}{Z_n^{\mathrm{graph},0}},
 \end{aligned}
\end{equation}

Let $a_v:[0,1]\to[0,1]$ be measurable vertex weights. For a finite graph
$A$ on $m$ vertices, take weights $a_v(i)\in[0,1]$ and $p\in[0,1]$.

\begin{lemma}
\label{lem:weighted-cut-counting}
Let $H$ be a finite simple graph, let $U,W$ be graphons, and let the
measurable vertex weights be as specified above. Then
\begin{equation}\label{row:eq:weighted-cut-counting}
 \begin{aligned}
 &\left|\int_{[0,1]^{v(H)}}
    \left(\prod_{v\in V(H)}a_v(z_v)\right)
    \left\{\prod_{uv\in E(H)}U(z_u,z_v)
               -\prod_{uv\in E(H)}W(z_u,z_v)\right\}
          \prod_{v\in V(H)}\dd z_v\right|\\
 &\quad\le e(H)\|U-W\|_\square.
 \end{aligned}
\end{equation}
For the finite graph and weights specified above, the corresponding injective estimate is
\begin{equation}\label{row:eq:weighted-injective-counting}
 \begin{aligned}
 &\Biggl|\frac1{m^{v(H)}}
   \sum_{\varphi:V(H)\hookrightarrow V(A)}
       \prod_{v\in V(H)}a_v(\varphi(v))
       \prod_{uv\in E(H)}A_{\varphi(u)\varphi(v)}
\\
 &\quad-p^{e(H)}\prod_{v\in V(H)}
                         \left(\frac1m\sum_{i\in V(A)}a_v(i)\right)
            \Biggr|\\
 &\qquad\le e(H)\|\mathcal W_A-p\|_\square
                  +\frac{\binom{v(H)}2}{m}.
 \end{aligned}
\end{equation}
Both bounds are uniform in the vertex weights, and $H$ may have isolated
vertices. Empty products are one.
\end{lemma}

\begin{proof}
For $0\le\alpha,\beta\le1$, layer cake and Fubini give
\begin{equation}\label{audit:weighted-kernel-test}
 \int D(x,y)\alpha(x)\beta(y)\,\dd x\,\dd y
 =\int_0^1\!\!\int_0^1
     \int_{\{\alpha>s\}\times\{\beta>t\}}D\,\dd s\,\dd t,
 \qquad \left|\int D\alpha\beta\right|\le\|D\|_\square.
\end{equation}
Enumerate the $\ell=e(H)$ edges and telescope:
\[
 \prod_{s=1}^{\ell}U_s-\prod_{s=1}^{\ell}W_s
 =\sum_{s=1}^{\ell}(U_s-W_s)\prod_{t<s}U_t\prod_{t>s}W_t.
\]
For edge $uv$, fixing all other vertices leaves the integrand
$\gamma\alpha(z_u)\beta(z_v)(U-W)(z_u,z_v)$ with all three
weights in $[0,1]$, because $H$ is simple. Applying
\eqref{audit:weighted-kernel-test} to each term proves
\eqref{row:eq:weighted-cut-counting}.

Extend finite vertex weights as step functions. For the normalized
unrestricted and injective counts, with $h=v(H)$,
\[
 \begin{aligned}
 |T_A^{\rm all}-p^\ell\prod_v(m^{-1}\sum_i a_v(i))|
   &\le\ell\|\mathcal W_A-p\|_\square,\\
 0\le T_A^{\rm all}-T_A^{\rm inj}
   &\le m^{-h}\sum_{\{u,v\}}\#\{\varphi:\varphi(u)=\varphi(v)\}
     =\binom h2/m.
 \end{aligned}
\]
The second line uses bounded summands and a union bound over
colliding vertex pairs. Combining these estimates proves
\eqref{row:eq:weighted-injective-counting}.
Taking all weights equal to one also gives
\begin{equation}\label{row:eq:cut-energy-continuity}
 \begin{aligned}
 |t(H,U)-t(H,W)|&\le e(H)\|U-W\|_\square,\\
 |T(U)-T(W)|
 &\le\left(\frac{|h_*|}{2}+\sum_r|b_r|e_r\right)
                                      \|U-W\|_\square.
 \end{aligned}
\end{equation}
The formulas include isolated vertices and empty products.
\end{proof}

\subsection{Conditional energies and the row-potential gap}\label{sec:deterministic-row-bounds}
Use the conditional energies in \eqref{row:eq:fixed-rows-split}, retaining
the ambient normalization $n$. For multi-indices $\gamma$, write
$\partial^\gamma=\prod_j\partial_{x_j}^{\gamma_j}$ and
$|\gamma|=\sum_j\gamma_j$.
\begin{lemma}
\label{lem:row-deterministic-bounds}
For the energies of Section~\ref{sec:deterministic-row-bounds}, fix
$q_0\ge1$ independently of $n$ and assume
$|h_n|+\sum_r|b_{r,n}|$ is uniformly bounded. Then the following deterministic
bounds hold on the continuous cubes, uniformly in $A,\mathbf e$:
\begin{equation}\label{audit:row-energy-interaction-bounds}
 \begin{gathered}
 \|l_{i,A}\|_\infty+\|n\partial_j l_{i,A}\|_\infty\le C
       \quad(i\ne j),\\
 \|\partial^\gamma l_{i,A}\|_\infty+
 \|\partial^\gamma(n\partial_j l_{i,A})\|_\infty
       \le C_hn^{-h},\quad h=|\gamma|\ge1,\\
 \|V_{A,\mathbf e}\|_\infty\le C_V,\qquad
 \|\partial_uV_{A,\mathbf e}\|_\infty\le C/n,\\
 \|\partial_{(c,j)}\partial_{(b,i)}V_{A,\mathbf e}\|_\infty
 \le \frac Cn\one_{\{i=j,\ b\ne c\}}
       +\frac C{n^2}\one_{\{i\ne j\}}.
 \end{gathered}
\end{equation}
Moreover, the Hessian $H_A''=(\partial_{ij}H_A)$ satisfies
$\|H_A''(x)\|_{\infty\to\infty}+\|H_A''(x)\|_{\mathrm{op}}\le C$
uniformly over all exteriors $A$ and $x\in[0,1]^{V(A)}$.
More precisely, for a fixed maximum degree $D$,
\begin{equation}\label{audit:interrow-coefficient-budget}
 \partial_uV_{A,\mathbf e}(\mathbf x)
 =\sum_{\substack{S\subseteq\mathcal I\setminus\{u\}\\|S|\le D}}
       a_{u,S}\prod_{v\in S}\mathbf x_v,
 \qquad \sum_S|a_{u,S}|\le C/n.
\end{equation}
The coefficients may depend on $A,\mathbf e$, but the bound is uniform in both.
All interaction terms vanish for $q_0=1$. The constants and thresholds depend
only on the templates, coefficient bound, $q_0$, tolerances, and derivative
orders, and are uniform in $n,A,\mathbf e,s$.
\end{lemma}

\begin{proof}
Partition template injections by their distinguished images, retaining
all vertices. Fixing the root and a set $J$ of distinct row coordinates
leaves at most $C_{F_r,|J|}n^{v_r-|J|-1}$ maps. Thus
\begin{equation}\label{v8:proof:lem:row-deterministic-bounds:single-row-count}
 \begin{aligned}
 \sup_{A,x}|\partial_JH_A(x)|
 &\le |h_n|\one_{\{|J|=1\}}
   +\sum_{r:v_r\ge|J|+1}
     C_{F_r,|J|}|b_{r,n}|n^{2-v_r}n^{v_r-|J|-1}\\
 &\le C_{|J|}n^{1-|J|},\qquad |J|\ge1.
 \end{aligned}
\end{equation}
Repeated derivatives vanish by multi-affinity. Applying
\eqref{v8:proof:lem:row-deterministic-bounds:single-row-count} at orders $h+1,h+2$ gives
\begin{equation}\label{row:eq:dense-derivatives}
 \begin{gathered}
 \|l_{i,A}\|_\infty+\|n\partial_jl_{i,A}\|_\infty\le C,\\
 \|\partial_{j_1}\cdots\partial_{j_{h}}l_{i,A}\|_\infty
 +\|\partial_{j_1}\cdots\partial_{j_{h}}(n\partial_jl_{i,A})\|_\infty
       \le C_{h}n^{-h},\qquad h\ge1.
 \end{gathered}
\end{equation}
The symmetric Hessian has $O(n^{-1})$ entries, so its absolute
row and column sums, and hence its operator norm, are bounded.

An interrow term fixes at least two distinguished images.
Each additional distinct exterior image supplies a factor $n^{-1}$:
\[
 \begin{gathered}
 |V|\le C,\qquad |\partial_uV|\le C/n,\\
 |\partial_v\partial_uV|
 \le\begin{cases}
 C/n,&\text{same exterior site, distinct rows},\\
 C/n^2,&\text{distinct exterior sites}.
 \end{cases}
 \end{gathered}
\]
Repeated coordinates give zero. Before collecting equal monomials,
the absolute coefficient sum of $\partial_uV$ is at most
\[
 \sum_r C_{F_r,q_0}|b_{r,n}|n^{2-v_r}n^{v_r-3}\le C/n.
\]
Collecting equal monomials cannot increase this sum; multi-affinity excludes $u$,
and the fixed templates bound the degree. This proves
\eqref{audit:interrow-coefficient-budget}.
\end{proof}

The block features $\chi_{i,S}=\prod_{a\in S}x_i^a$,
$\varnothing\ne S\subseteq[q_0]$, have support size at most $q_0$,
overlap $2^{q_0-1}$, and $O(n)$ labels. The interrow energy expands as
\begin{equation}\label{ep:row:interrow-polynomial}
 V_{A,\mathbf e}=c+\sum_{d=1}^{D_V}\sum_{\beta_1,\ldots,\beta_d}
 v_{\beta_1\ldots\beta_d}\chi_{\beta_1}\cdots\chi_{\beta_d},
 \qquad |v_{\beta_1\ldots\beta_d}|\le Cn^{-d}.
\end{equation}
Here $d$ counts distinct exterior blocks. A template on $v$ vertices
meeting $h\ge2$ roots and $d$ prescribed exterior sites has at most
$Cn^{v-h-d}$ completions; multiplication by $n^{2-v}$ gives
$Cn^{2-h-d}\le Cn^{-d}$. Collecting over the fixed template family
preserves this bound.

Write $I(u)=u\log u+(1-u)\log(1-u)$, with $0\log0=0$, and define the row potential by

\begin{equation}\label{row:eq:row-potential}
 \begin{aligned}
 f_p(u)&=h_*u+\sum_r b_r\sum_{v\in V(F_r)}
                    p^{e_r-d_{F_r}(v)}u^{d_{F_r}(v)},\\
 J_p(u)&=f_p(u)-I(u),\qquad 0\le u\le1.
 \end{aligned}
\end{equation}

\begin{proposition}\label{row:prop:row-gap}
Assume $b_r\ge0$ and that the variational function $\Phi(x)$ in
\eqref{eq:scalar} has a unique global maximizer $p\in(0,1)$ satisfying \eqref{eq:critical} for a fixed
integer $\mc\ge2$. Then $p$ is the unique global
maximizer of $J_p(u)$, and there is a constant $c_1>0$, depending only on
the fixed templates and leading parameters, such that
\begin{equation}\label{row:eq:row-gap}
 \begin{aligned}
 J_p'(p)&=0,\qquad
 J_p''(p)=2t(p)-\frac1{\nu(p)}\le-\frac1{2\nu(p)},\\
 J_p(p)-J_p(u)&\ge c_1(u-p)^2
                  \qquad\text{for every }u\in[0,1].
 \end{aligned}
\end{equation}
The row potential $J_p$ and $\Phi$ satisfy
\begin{equation}\label{row:eq:row-Young}
 J_p(u)-J_p(p)\le2\{\Phi(u)-\Phi(p)\}
                  \qquad(0\le u\le1).
\end{equation}
No absolute interaction budget or global sign condition for $\Phi'(x)$
is required.
\end{proposition}

\begin{proof}
The identities $\sum_vd_F(v)=2e(F)$ and
$\sum_vd_F(v)(d_F(v)-1)=2a(F)$ give
\[
 J_p'(p)=h_*+2B'(p)-I'(p)=0,\qquad
 J_p''(p)=2t(p)-1/\nu(p)\le-1/(2\nu(p)),
\]
using $a(F)\le\binom{e(F)}2$, $b_r\ge0$, and $\Phi''(p)=0$ from
\eqref{eq:critical}.
Weighted AM--GM gives
$p^{e-d}u^d-p^e\le(d/e)(u^e-p^e)$ for $0\le d\le e$.
Summing over vertices yields
\[
 f_p(u)-f_p(p)\le h_*(u-p)+2\{B(u)-B(p)\},
\]
hence \eqref{row:eq:row-Young}. By this comparison and uniqueness of the global maximizer,
$2\{\Phi(u)-\Phi(p)\}<0$ for $u\ne p$. Near $p$, Taylor's formula gives
$J_p(p)-J_p(u)\ge(u-p)^2/(8\nu(p))$ by the curvature in
\eqref{row:eq:row-gap}; on the compact complement,
the positive continuous ratio
$(J_p(p)-J_p(u))/(u-p)^2$ has a positive minimum.
This proves \eqref{row:eq:row-gap} on $[0,1]$.
\end{proof}

\subsection{Graph localization and rooted-count concentration}\label{app:macro-proof}\label{row:sec:roots}
For undirected graphs, the entropy contribution to the variational
functional is $-\frac12\int I(W)$. The following specialization of
\cite[Theorem~3.2]{CD} gives exponential concentration when the variational
functional has a unique maximizer given by a constant graphon.

\begin{theorem}\label{row:thm:published-concentration}
Suppose that $p\in(0,1)$ satisfies the following variational gap:
for each fixed $\delta>0$ there is $\kappa_\delta>0$ such that
\begin{equation}\label{row:eq:variational-gap}
 \sup_{\|W-p\|_\square\ge\delta}
 \left\{T(W)-\frac12\int_{[0,1]^2}I(W(x,y))\,\dd x\,\dd y\right\}
 \le T(p)-\frac12 I(p)-\kappa_\delta.
\end{equation}
Here and below a supremum over an empty set is $-\infty$.
Then
\begin{equation}\label{row:eq:reference-concentration}
 \Pp_n^0(\|\mathcal W_X-p\|_\square>\delta)
   \le C_\delta e^{-c_\delta n^2}\qquad(n\ge1).
\end{equation}
The constants depend only on the fixed coefficients, templates, $p$
and $\delta$. No sign condition on the coefficients is needed once
\eqref{row:eq:variational-gap} is available.
\end{theorem}

By \eqref{row:eq:cut-energy-continuity}, $T$ is bounded, cut-continuous,
and invariant under measure-preserving relabeling. Condition
\eqref{row:eq:variational-gap} identifies the constant graphon with density $p$
as the unique maximizer of the variational functional. The cut distance to this graphon equals $\|W-p\|_\square$, and the entropy
in \cite[Theorem~3.2]{CD} is $\frac12\int I(W)$ in our notation.
That theorem therefore gives \eqref{row:eq:reference-concentration}.

We next transfer this concentration bound to the original graph law
and obtain a simultaneous bound for its induced exteriors. For $S\subseteq[n]$, let
$A_S=X_{[n]\setminus S}$, whose step graphon uses $n-|S|$ equal intervals.

\begin{proof}[Proof of Lemma~\ref{row:lem:macro}]
Taylor's formula near $p$ and uniqueness on the compact complement give
\begin{equation}\label{even:scalar-gap}
 \Phi(p)-\Phi(u)\ge c|u-p|^{2\mc},\qquad 0\le u\le1.
\end{equation}
H\"older's inequality, $t(F,W)\le\int W^{e(F)}$, $b_r\ge0$, and
\eqref{even:scalar-gap} yield
\begin{equation}\label{audit:dependency:app:graphon-gap-bound}
 T(W)-\tfrac12\int I(W)
 \le\int\Phi(W)\le\Phi(p)-c\int|W-p|^{2\mc}
 \le\Phi(p)-c\|W-p\|_\square^{2\mc}.
\end{equation}
The constant graphon attains $\Phi(p)$, and
\eqref{audit:dependency:app:graphon-gap-bound} verifies \eqref{row:eq:variational-gap}.
Theorem~\ref{row:thm:published-concentration} therefore gives
\eqref{row:eq:reference-concentration}.
The counting bound
$0\le\Hom(F_r,X)-\inj(F_r,X)\le\binom{v_r}2n^{v_r-1}$
and \eqref{eq:parameter-neighborhood} imply
\begin{equation}\label{audit:dependency:app:reference-law-comparison}
 \sup_X|\mathcal H_n-n^2T(\mathcal W_X)|\le Cn,\qquad
 \Pp_n(B)\le e^{2Cn}\Pp_n^0(B).
\end{equation}
Combining \eqref{audit:dependency:app:reference-law-comparison} with
\eqref{row:eq:reference-concentration} transfers the $e^{-cn^2}$ tail,
including signed perturbations.
Also $\int\mathcal W_X=(1-n^{-1})P$, which gives the density tail.
Finally, restricting the finite cuts gives, simultaneously for $|S|=d$,
\[
 \|\mathcal W_{A_S}-p\|_\square
 \le\left(\frac n{n-d}\right)^2\|\mathcal W_X-p\|_\square
 \le2\|\mathcal W_X-p\|_\square
\]
for large $n$. Apply the graph tail at $\delta/2$.
\end{proof}

\Needspace{7\baselineskip}
\subsubsection{Conditional concentration of empirical row types}
For distinct roots $s_1,\ldots,s_{q_0}$ put $S=\{s_1,\ldots,s_{q_0}\}$,
and condition on

\begin{equation}\label{row:eq:type-conditioning-data}
 X_S=(X_{s_as_b})_{1\le a<b\le q_0},\qquad
 A=(X_{ik})_{i,k\in[n]\setminus S},\qquad
 V(A)=[n]\setminus S.
\end{equation}

The random rows are $x_i^j=X_{s_ji}$, and their columns
$\boldsymbol x_i=(x_i^1,\ldots,x_i^{q_0})$ take values in
$\{0,1\}^{q_0}$. With $m=n-q_0$, define

\begin{equation}\label{row:eq:empirical-type-definition}
 \nu_{n-q_0}(w)
 =\frac1{n-q_0}\sum_{i\in V(A)}\one_{\{\boldsymbol{x}_i=w\}},
       \qquad w\in\{0,1\}^{q_0}.
\end{equation}

The normalized row degrees satisfy
\begin{equation}\label{row:eq:type-marginal-degree}
 u_j=\sum_w w_j\nu_m(w)=m^{-1}\sum_{i\in V(A)}x_i^j\in[0,1].
\end{equation}
Write $\pi_p=\operatorname{Ber}(p)^{\otimes q_0}$ and use the full
$\ell^1$ distance $\|\nu-\lambda\|_1=\sum_w|\nu(w)-\lambda(w)|$.

\begin{lemma}\label{row:lem:type}
Assume \hyperref[ass:PC]{\textup{(PC)}}.
Use \eqref{row:eq:type-conditioning-data} and
\eqref{row:eq:empirical-type-definition}. For every fixed $\epsilon>0$,
there is a fixed $\delta>0$ such that,
whenever $\|\mathcal W_A-p\|_\square\le\delta$,
\begin{equation}\label{row:eq:type-tail}
 \Pp_n(\|\nu_{n-q_0}-\pi_p\|_1>\epsilon\mid A,X_S)
 \le C_\epsilon e^{-c_\epsilon n}.
\end{equation}
The bound is uniform in the
internal edge configuration $X_S$ and over all exteriors satisfying this
cut condition.
\end{lemma}

\begin{proof}
Fix $(A,X_S)$ and write $m=n-q_0$. Applying
\eqref{row:eq:weighted-injective-counting}
to $F_r-\{v\}$, with neighbor weights $x_i^j$, gives an error
$O(\delta+m^{-1})$ in each normalized rooted count.
Since $n^{2-v_r}m^{v_r-1}=m+O(1)$, parameter errors are $O(n^{-1})$,
and $V_{A,X_S}=O(1)$,
\begin{equation}\label{audit:type-energy}
 \sum_jH_A(x^j)+V_{A,X_S}(\mathbf x)
       =m\sum_jf_p(u_j)+O(n\delta+1).
\end{equation}
For $\lambda=\bigotimes_j\operatorname{Ber}(u_j)$,
the marginal constraints give $\nu\ll\lambda$, even at $u_j=0,1$.
Pinsker's inequality \cite[Lemma~2.5(i)]{Tsybakov2009}
(see also \cite[Theorem~4.19]{BLM2013}) gives
\begin{equation}\label{audit:dependency:app:type-pinsker-bound}
 D(\nu\|\lambda)
 =\sum_w\nu(w)\log\nu(w)-\sum_jI(u_j)
 \ge\tfrac12\|\nu-\lambda\|_1^2.
\end{equation}
Product telescoping yields
$\|\lambda-\pi_p\|_1\le2\sum_j|u_j-p|$ and hence
$\|\nu-\pi_p\|_1^2
\le2\|\nu-\lambda\|_1^2+8q_0\sum_j(u_j-p)^2$.
Combining this comparison with \eqref{audit:dependency:app:type-pinsker-bound}
and \eqref{row:eq:row-gap} gives
\begin{equation}\label{audit:type-free-energy-gap}
 \begin{aligned}
 &q_0J_p(p)-\{\textstyle\sum_jf_p(u_j)-\sum_w\nu(w)\log\nu(w)\}\\
 &=\sum_j\{J_p(p)-J_p(u_j)\}+D(\nu\|\lambda)\\
 &\ge c_{q_0}\|\nu-\pi_p\|_1^2,\qquad
 c_{q_0}=\min\{1/4,c_1/(8q_0)\}>0.
 \end{aligned}
\end{equation}
The combinatorial type bound \cite[Theorems~11.1.1 and~11.1.3]{CoverThomas2006}
is $m!/\prod_w(m\nu(w))!\le e^{-m\sum_w\nu(w)\log\nu(w)}$.
For the denominator, Jensen's inequality under $\pi_p^{\otimes m}$ gives
\begin{equation}\label{audit:dependency:app:type-partition-lower-bound}
 \log Z_{A,X_S}\ge
 \E_{\pi_p^{\otimes m}}\!\left[
   \sum_jH_A(x^j)+V_{A,X_S}
            -\log\pi_p^{\otimes m}(\mathbf x)\right]
 \ge mq_0J_p(p)-C(n\delta+1),
\end{equation}
using $\E u_j=p$, $\E(u_j-p)^2=\nu(p)/m$, and bounded $f_p''(u)$ on $[0,1]$.
Equations~\eqref{audit:type-energy}, \eqref{audit:type-free-energy-gap}, and
\eqref{audit:dependency:app:type-partition-lower-bound}, with the type-count bound, give
\begin{equation}\label{audit:type-probability}
 \Pp_n(\nu_m=\nu\mid A,X_S)
 \le\exp\{-c_{q_0}m\|\nu-\pi_p\|_1^2+C(n\delta+1)\}.
\end{equation}
There are at most $(m+1)^{2^{q_0}}$ types. Choose
$0<\delta\le c_{q_0}\epsilon^2/(8C)$, and then take $n$ sufficiently large, with $m\ge n/2$.
Summing \eqref{audit:type-probability} proves \eqref{row:eq:type-tail},
uniformly in $(A,X_S)$.
\end{proof}

\subsubsection{Concentration of rooted graph counts}
A rooted template $(H,\mathbf u)$ consists of a fixed simple graph and
an ordered tuple $\mathbf u=(u_1,\ldots,u_{q_0})$ of distinct vertices,
where $q_0\ge0$ is fixed for each template and may vary across a finite family.
For ordered distinct images $\mathbf i=(i_1,\ldots,i_{q_0})$, define

\begin{equation}\label{row:eq:rooted-count}
 T_{H,\mathbf u;\mathbf i}(X)=n^{q_0-v(H)}
 \sum_{\substack{\varphi:V(H)\hookrightarrow[n]\\\varphi(u_j)=i_j}}
 \prod_{ab\in E(H)\setminus E(H[\{u_1,\ldots,u_{q_0}\}])}
                         X_{\varphi(a)\varphi(b)}.
\end{equation}

Thus $u_a\mapsto i_a$ and all root-to-root edge factors are omitted.
Put $e_0(H,\mathbf u)=e(H[\{u_1,\ldots,u_{q_0}\}])$. The full count retaining
these factors obeys

\begin{equation}\label{row:eq:rooted-full-factor}
 T^{\mathrm{full}}_{H,\mathbf u;\mathbf i}(X)
 =\left(\prod_{\substack{1\le a<b\le q_0\\u_au_b\in E(H)}}X_{i_ai_b}\right)
      T_{H,\mathbf u;\mathbf i}(X).
\end{equation}

For $q_0=2$, the full count is written $T^{\mathrm{full}}_{ij}(H;u,v;X)$,
as in \eqref{eq:full-two-root-count}.
Both ordered tuples are part of the data, since different template roots
can define different functions. On an induced graph, the maps are
restricted to its vertices while the ambient factor $n^{q_0-v(H)}$ is retained.
The unrooted case uses empty tuples. Let $\mathcal R$ be a fixed finite
family of such rooted pairs and write $D_i=\sum_{j\ne i}X_{ij}$.

\begin{theorem}\label{row:thm:rooted}
Assume \hyperref[ass:PC]{\textup{(PC)}}.
For the counts in \eqref{row:eq:rooted-count}, a fixed finite $\mathcal R$, and fixed
$\epsilon>0$, there are $C_\epsilon,c_\epsilon>0$ such that
\begin{equation}\label{row:eq:rooted-tail}
 \Pp_n\left(\max_{(H,\mathbf u)\in\mathcal R}
       \max_{\mathbf i\ \mathrm{distinct}}
 |T_{H,\mathbf u;\mathbf i}(X)-p^{e(H)-e_0(H,\mathbf u)}|>\epsilon\right)
 \le C_\epsilon e^{-c_\epsilon n}.
\end{equation}

In particular, the one-rooted edge $K_2$ gives the degree estimate
\begin{equation}\label{row:eq:rooted-degree-tail}
 \Pp_n\left(\max_{1\le i\le n}
            \left|\frac{D_i}{n}-p\right|>\epsilon\right)
       \le C_\epsilon e^{-c_\epsilon n}.
\end{equation}
The positive constants and the threshold may additionally depend on
$\mathcal R,\epsilon$; they are uniform over root images.
\end{theorem}

\begin{proof}
For $q_0=0$, the weighted counting bound and macro localization give the assertion.
For $q_0\ge1$, fix the root images $\mathbf i$ and put
$S=\{i_1,\ldots,i_{q_0}\}$,
$A=X_{[n]\setminus S}$, $m=n-q_0$, and
$x_k^j=X_{i_jk}$ for $j\in[q_0]$ and $k\in V(A)$.
Condition on $(A,X_S)$, where $X_S$ contains the edges within $S$, and set
$a_v(k)=\prod_{j:u_jv\in E(H)}x_k^j$ for $v\notin\mathbf u$.
Equation~\eqref{row:eq:weighted-injective-counting} and $m/n=1+O(n^{-1})$ give
\[
 \left|T_{H,\mathbf u;\mathbf i}
    -p^{e(H-\mathbf u)}\prod_{v\notin\mathbf u}
                    \left(m^{-1}\sum_i a_v(i)\right)\right|
 \le e(H-\mathbf u)\|\mathcal W_A-p\|_\square+C_H/n.
\]
Each average is a $[0,1]$-valued linear functional of $\nu_m$,
equal at $\pi_p$ to $p^{|\{j:u_jv\in E(H)\}|}$.
Product telescoping therefore gives
\begin{equation}\label{audit:rooted-type-error}
 |T_{H,\mathbf u;\mathbf i}-p^{e(H)-e_0(H,\mathbf u)}|
 \le C_H\{\|\mathcal W_A-p\|_\square+\|\nu_m-\pi_p\|_1+1/n\}.
\end{equation}
Choose a fixed type tolerance, then a fixed cut tolerance, and finally
take $n$ large enough that the right-hand side of \eqref{audit:rooted-type-error} is at most $\epsilon$.
Equations~\eqref{row:eq:type-tail} and \eqref{audit:macro-deletion-tail}
give a uniform $Ce^{-cn}$ tail. If all vertices are rooted, both
quantities equal one. A union bound over root images and the fixed finite family gives
$\sum_{(H,\mathbf u)\in\mathcal R}n^{|\mathbf u|}C_He^{-c_Hn}
\le C_\epsilon e^{-c_\epsilon n}$, proving \eqref{row:eq:rooted-tail}.
Taking $T_{K_2,(o);(i)}=D_i/n$ gives the degree bound.
\end{proof}

\subsection{Contractive row kernels and higher derivatives}\label{row:sec:extension}
Define the kernel constants

\begin{equation}\label{row:eq:ideal-kernel}
 \alpha=2t(p)-6p\tau(p)\ge0,\qquad \beta=6\tau(p)\ge0.
\end{equation}

Their nonnegativity follows from $a(F)\ge3\triangle(F)$; moreover

\begin{equation}\label{row:eq:weighted-row}
 \nu(p)(\alpha+\beta p)=2\nu(p)t(p)\le1/2.
\end{equation}

Fix the tolerances for the auxiliary kernels and the mean-map estimates
before choosing good exterior sets. Choose $\delta_0>0$ with
$[p-2\delta_0,p+2\delta_0]\subset J$, and a smooth
$\chi:\R\to[0,1]$ equal to one on $(-\infty,1/2]$ and zero on $[1,\infty)$.
Lemma~\ref{lem:row-deterministic-bounds} gives deterministic constants
$C_H,C_{\rm det}$ bounding
$n\|\partial_j(n\partial_kl_{i,A})\|_\infty$ and the full-cube sums
$\sum_k|\partial_kl_{i,A}|$, $\sum_k|n\partial_k\partial_jl_{i,A}|$,
respectively. For large $n$ these depend only on the leading model.
Enlarge $K$ to satisfy $K\ge2(1+C_{\rm det}+4/\nu(p))$, increasing the
threshold for $n$ if needed. By the gap in \eqref{row:eq:weighted-row},
choose fixed $\delta,\eta,\varepsilon_{\rm deg}>0$ and
$I=[p-\delta,p+\delta]$ such that
\begin{equation}\label{v11:eq:tolerance-contraction}
 \begin{gathered}
 (\nu(p)+K\delta)
   (\alpha+\beta p+\beta\varepsilon_{\rm deg}+\eta/4+K\delta)
       \le5/8,\qquad I\Subset(0,1),\\
 \delta\le\delta_0,\qquad 2\nu(u)t(u)\le3/4\quad(u\in I),\\
 \eta\le1,\quad \nu(p)\beta\varepsilon_{\rm deg}\le1/16,\quad
 \nu(p)\eta\le1/16 .
 \end{gathered}
\end{equation}
Next choose $0<a\le\min\{\delta/28,p/8,(1-p)/8\}$ so that
\begin{equation}\label{v11:eq:tolerance-saturation}
\begin{gathered}
 2a<\min\{p,1-p\}/2,\qquad
 2a(\alpha+\beta p+\beta\varepsilon_{\rm deg}+\eta)
       +2\|\chi'\|_\infty C_Ha/\eta\le1/8 .
\end{gathered}
\end{equation}
All these choices depend only on deterministic model constants and
are made independently of $n$ and of the exterior. The gap in
\eqref{row:eq:weighted-row} allows $\delta,\eta,\varepsilon_{\rm deg}$
to be chosen first and $a$ last.

For the reference law and its field, set

\begin{equation}\label{row:eq:Omega}
 \begin{aligned}
 &\Omega_A=\Bigl\{x\in\{0,1\}^{V(A)}:\\
 &\qquad\max_i|\pi_A(x_i=1\mid x_{-i})-p|\le a,\\
 &\qquad\max_{i\ne j}|n\partial_j l_{i,A}(x)-(\alpha+\beta A_{ij})|
                         \le\eta/4\Bigr\}.
 \end{aligned}
\end{equation}

Put $\Omega(A)=\Omega_A^{q_0}$ and
$d_{A,i}=\sum_{j\in V(A)\setminus\{i\}}A_{ij}$.
These sets depend only on the exterior and rows, not on $\mathbf e,s$.

\begin{lemma}\label{row:lem:good-exteriors}
Assume \hyperref[ass:PC]{\textup{(PC)}}.
Fix $q_0$ and sufficiently small tolerances as specified above.
Then there exist $c,C>0$ and deterministic sets $\mathcal A_n$ of exteriors such that
\begin{equation}\label{row:eq:good-properties}
 \begin{gathered}
 \Pp_n(A\notin\mathcal A_n)\le Ce^{-cn},\\
 |P_A-p|\le\delta,\qquad
 \max_i|d_{A,i}/n-p|\le \varepsilon_{\mathrm{deg}},\qquad
 \pi_A(\Omega_A^c)\le e^{-cn}\quad(A\in\mathcal A_n).
 \end{gathered}
\end{equation}
The sets $\mathcal A_n$ can be used simultaneously for all internal-edge
configurations and interpolation parameters, with
\[
 \sup_{A\in\mathcal A_n}\sup_{\mathbf e\in\{0,1\}^{\binom{q_0}{2}}}
 \sup_{0\le s\le1}
 \mu_{A,\mathbf e,s}(\Omega(A)^c)\le Ce^{-cn}.
\]
For $q_0=2$ this includes
$\mu_{A,(e),s}((\Omega_A\times\Omega_A)^c)\le Ce^{-cn}$ for both values
of $e$ and every $s\in[0,1]$; for $q_0=1$ it is the single-row bound.

\end{lemma}

\begin{proof}
Let $L_{bi}=\partial_{bi}\mathcal H_n$. The deterministic bounds give
\begin{equation}\label{audit:physical-reference-field-error}
 |L_{bi}-l_{i,A}(x^b)|+
 |n\partial_{bi}\partial_{bj}\mathcal H_n
                         -n\partial_jl_{i,A}(x^b)|\le C/n
 \qquad(i\ne j).
\end{equation}
Differentiating template injections gives
\begin{equation}\label{audit:dependency:app:physical-rooted-field-expansion}
 \begin{aligned}
 L_{bi}&=h_n+\sum_r b_{r,n}
       \sum_{uv\in E(F_r)\ {\rm ordered}}T_{F_r,(u,v);(b,i)},\\
 n\partial_{bi}\partial_{bj}\mathcal H_n
 &=\sum_r b_{r,n}
   \sum_{\substack{u,v,w\ {\rm distinct}\\uv,uw\in E(F_r)}}
       A_{ij}^{\one_{\{vw\in E(F_r)\}}}
                         T_{F_r,(u,v,w);(b,i,j)}.
 \end{aligned}
\end{equation}
The full templates are retained; $A_{ij}$ restores the third root edge.
There are $2e_r$ ordered edges,
$2a(F_r)-6\triangle(F_r)$ open ordered wedges, and
$6\triangle(F_r)$ closed ones. Applying \eqref{row:eq:rooted-tail} to
\eqref{audit:dependency:app:physical-rooted-field-expansion}, then using
\eqref{audit:physical-reference-field-error}, gives profiles $\ell(p)$ and $\alpha+\beta A_{ij}$ with simultaneous
error $C\epsilon_0+C/n$, except with probability $Ce^{-cn}$.
Choose a fixed $\epsilon_0$ small enough for the tolerances $a,\eta/4$.
Since $\sig(\ell(p))=p$, all rows then lie in $\Omega_A$.
The degree bound follows from
$d_{A,i}=D_i-\sum_bX_{bi}$ and \eqref{row:eq:rooted-degree-tail}.

The uniform interaction bound gives
\begin{equation}\label{audit:conditional-density-ratio}
 e^{-2C_V}\le\frac{\dd\mu_{A,\mathbf e,s}}
                       {\dd\pi_A^{\otimes q_0}}\le e^{2C_V}.
\end{equation}
Comparing with the true conditional law using
\eqref{audit:conditional-density-ratio}, and averaging over $A$, yields
\begin{equation}\label{audit:dependency:app:mean-good-row-defect}
 \E_A\pi_A^{\otimes q_0}(\Omega(A)^c)\le Ce^{-c_0n}.
\end{equation}
By Lemma~\ref{row:lem:macro} and
$|P_A-p|\le\|\mathcal W_A-p\|_\square+1/m$,
\[
 \Pp_n(|P_A-p|>\delta)\le C_\delta e^{-c_\delta n^2}.
\]
Define $\mathcal A_n$ by $|P_A-p|\le\delta$, the degree condition, and
$\pi_A^{\otimes q_0}(\Omega(A)^c)\le e^{-c_0n/2}$.
Markov's inequality applied to \eqref{audit:dependency:app:mean-good-row-defect}
and a union bound give
$\Pp(A\notin\mathcal A_n)\le Ce^{-c_0n/2}$, while for $A\in\mathcal A_n$,
\[
 \pi_A(\Omega_A^c)\le e^{-c_0n/2},\qquad
 \sup_{\mathbf e,s}\mu_{A,\mathbf e,s}(\Omega(A)^c)
       \le e^{2C_V-c_0n/2}.
\]
The construction depends on $A$ alone; relabeling and a polynomial
union bound also permit all choices of the fixed number of roots.
\end{proof}

\subsubsection{Smooth contractive extensions of row updates}
The next proposition extends the reference-row updates to the full cube while retaining a uniform contraction gap.
We use the fixed cutoff and tolerances above, with $m=n-q_0$.

\begin{proposition}
\label{row:prop:extension}
Assume \hyperref[ass:PC]{\textup{(PC)}}.
Use the good sets in \eqref{row:eq:Omega}, the fixed cutoff $\chi(u)$, and the tolerances of
Appendix~\ref{row:sec:extension}. For all
sufficiently large $n$ and every $A\in\mathcal A_n$, there exist
smooth functions
$\widetilde\pi_{A,i}:[0,1]^{V(A)}\to(0,1)$, independent of $x_i$, with
\begin{equation}\label{row:eq:agreement}
 \widetilde\pi_{A,i}(x)=\pi_A(x_i=1\mid x_{-i})\quad(x\in\Omega_A),\qquad
 |\widetilde\pi_{A,i}(x)-p|\le2a\quad(x\in[0,1]^{V(A)}).
\end{equation}
In particular,
$p/2\le\widetilde\pi_{A,i}(x)\le(1+p)/2$ on the whole continuous
cube. They admit a symmetric nonnegative majorant
$c^{(n,A)}=(c_{ij}^{(n,A)})_{i,j\in V(A)}$, with
$c_{ii}^{(n,A)}=0$, such that
\begin{equation}\label{row:eq:extension-majorant}
 |\partial_j\widetilde\pi_{A,i}(x)|\le c_{ij}^{(n,A)},\quad
 \max_i\sum_jc_{ij}^{(n,A)}\le3/4,
 \qquad c_{ij}^{(n,A)}\le C/n.
\end{equation}
The constants and threshold depend only on the fixed model, $K,q_0$,
cutoff, and tolerances in
\eqref{v11:eq:tolerance-contraction}--\eqref{v11:eq:tolerance-saturation}.
The kernels and majorant depend on $n,A$ but not on $\mathbf e,s$;
the bounds are uniform over good exteriors and coordinate indices.
The auxiliary kernels need not be Gibbs compatible.
\end{proposition}

\begin{proof}
Use the smooth-maximum construction of
\cite[Section~2 and Appendix~A.2]{CCK2013}
and \cite[Supplemental Material, Appendix~C]{CCKK2022}.
Set
$T_a(u)=\int_0^u\chi(t^2/(2a^2))\,\dd t$.
Then $T_a(u)=u$ on $|u|\le a$, $|T_a(u)|\le2a$,
$0\le T_a'(u)\le1$, $T_a'(u)=0$ on $|u|\ge2a$, and
$\|T_a^{(h)}\|_\infty\le C_ha^{1-h}$.
With $\lambda=\eta/[4\log(2m)]$, put
\begin{equation}\label{audit:dependency:app:soft-maximum-definition}
 \mathcal M_{A,i}(x)=\lambda\log
 \sum_{k\ne i}\sum_{\epsilon=\pm1}
 \exp\{\epsilon\bigl[n\partial_kl_{i,A}(x)-(\alpha+\beta A_{ik})\bigr]/\lambda\}.
\end{equation}
The maximum error is at most $\mathcal M_{A,i}(x)$ and differs from it
by at most $\eta/4$. Define
\begin{equation}\label{row:eq:cut-kernel}
 \widetilde\pi_{A,i}(x)
 =p+\chi(\mathcal M_{A,i}(x)/\eta)T_a(\sig(l_{i,A}(x))-p).
\end{equation}
Multi-affinity makes this independent of $x_i$. Its range lies in
$[p-2a,p+2a]\subset[p/2,(1+p)/2]$.
On $\Omega_A$, $\mathcal M_{A,i}(x)\le\eta/2$ and
$|\sig(l_{i,A}(x))-p|\le a$, proving agreement.

Differentiation gives the normalized weighted average
\begin{equation}\label{v8:proof:row:prop:extension:softmax-gradient}
\begin{gathered}
 \partial_j\mathcal M_{A,i}(x)=
 \frac{\displaystyle
  \sum_{k\ne i}\sum_{\epsilon=\pm1}
  e^{\epsilon\bigl[n\partial_kl_{i,A}(x)-(\alpha+\beta A_{ik})\bigr]/\lambda}
                      \epsilon n\partial_j\partial_kl_{i,A}(x)}
 {\displaystyle
  \sum_{k\ne i}\sum_{\epsilon=\pm1}
  e^{\epsilon\bigl[n\partial_kl_{i,A}(x)-(\alpha+\beta A_{ik})\bigr]/\lambda}},\\
 |\partial_j\mathcal M_{A,i}(x)|\le C_H/n.
\end{gathered}
\end{equation}
On the cutoff support,
$|\partial_jl_{i,A}(x)|\le(\alpha+\beta A_{ij}+\eta)/n$;
where $T_a'(\sig(l_{i,A}(x))-p)\ne0$, $\sig'(l_{i,A}(x))\le\nu(p)+2a$.
The product rule in \eqref{row:eq:cut-kernel}, together with
\eqref{v8:proof:row:prop:extension:softmax-gradient}, gives
\begin{equation}\label{audit:extension-majorant-profile}
 |\partial_j\widetilde\pi_{A,i}(x)|
 \le\frac{(\nu(p)+2a)(\alpha+\beta A_{ij}+\eta)
                  +2\|\chi'\|_\infty C_Ha/\eta}{n}
       =:c_{ij}^{(n,A)}\quad(i\ne j).
\end{equation}
Set $c_{ii}=0$. This majorant is symmetric, nonnegative, and $O(n^{-1})$.
Summing \eqref{audit:extension-majorant-profile}, using
\eqref{row:eq:good-properties}, \eqref{row:eq:weighted-row}, and
\eqref{v11:eq:tolerance-contraction}--\eqref{v11:eq:tolerance-saturation}, gives
\[
 \sum_jc_{ij}\le
 \tfrac12+\nu(p)\beta\varepsilon_{\rm deg}+\nu(p)\eta
 +2a(\alpha+\beta p+\beta\varepsilon_{\rm deg}+\eta)
 +2\|\chi'\|_\infty C_Ha/\eta\le\tfrac34.
\]
No further restriction of the exterior sets is needed.
\end{proof}

The log odds of the kernels in Proposition~\ref{row:prop:extension} are

\begin{equation}\label{v8:nearby:app:reference-log-odds}
 \log\frac{\widetilde\pi_{A,i}(x)}{1-\widetilde\pi_{A,i}(x)}.
\end{equation}

On binary states write $\nabla_jF(x)=F(x^{j,1})-F(x^{j,0})$ and
$\nabla_B=\prod_{j\in B}\nabla_j$.

\begin{lemma}
\label{row:lem:extension-derivatives}
Assume \hyperref[ass:PC]{\textup{(PC)}}.
For the kernels of Proposition~\ref{row:prop:extension}, let $\gamma$
be a fixed multi-index with $h=|\gamma|\ge1$. Then
\begin{equation}\label{row:eq:extension-continuous-higher}
 \begin{aligned}
 \sup_{i\in V(A)}
   \|\partial^\gamma\widetilde\pi_{A,i}\|_{\infty,[0,1]^{V(A)}}
     &\le C_{h}n^{-h}(1+\log n)^{h-1},\\
 \sup_{i\in V(A)}
   \|\partial^\gamma\!\left(\log\frac{\widetilde\pi_{A,i}}{1-\widetilde\pi_{A,i}}\right)\|_{\infty,[0,1]^{V(A)}}
     &\le C_{h}n^{-h}(1+\log n)^{h-1}.
 \end{aligned}
\end{equation}
Repeated derivatives in one coordinate are included. On binary states,
\begin{equation}\label{row:eq:extension-higher}
 \sup_{i\in V(A)}\sum_{\substack{B\subseteq V(A)\\|B|=h}}
       \|\nabla_B\widetilde\pi_{A,i}\|_\infty
       \le C_{h}(1+\log n)^{h-1}.
\end{equation}
Here the continuous supremum is over the displayed cube and the
finite-difference supremum is over its binary vertices. Constants depend on $h$ and the fixed data of
Proposition~\ref{row:prop:extension}, uniformly in $n$,
$A\in\mathcal A_n$ and all coordinate indices.
\end{lemma}

\begin{proof}
For a coordinate list $\boldsymbol j=(j_1,\ldots,j_{h})$ in which repeats are allowed,
\[
 \partial_{\boldsymbol j}e^{f/\lambda}
 =e^{f/\lambda}\sum_{\mathcal P\vdash[h]}\lambda^{-|\mathcal P|}
       \prod_{D\in\mathcal P}\partial_{\boldsymbol j_D}f.
\]
Apply this to \eqref{audit:dependency:app:soft-maximum-definition}
and isolate the one-block term. The normalized weights sum to one;
\eqref{row:eq:dense-derivatives} and induction give
\begin{equation}\label{audit:dependency:app:soft-maximum-derivatives}
 \|\partial_{\boldsymbol j}\mathcal M_{A,i}\|_\infty
       \le C_{h}n^{-h}\lambda^{1-h}.
\end{equation}
Indeed, each $k$-block term is
$\lambda^{1-k}\prod_D O(n^{-|D|}\lambda^{1-|D|})
=O(n^{-h}\lambda^{1-h})$; the base case is
\eqref{v8:proof:row:prop:extension:softmax-gradient}.
Combining \eqref{audit:dependency:app:soft-maximum-derivatives} with the chain
and product rules in \eqref{row:eq:cut-kernel} gives
$C_{h}n^{-h}(1+\log n)^{h-1}$.
Composition with the logit preserves this bound because the kernels
stay in a fixed compact subset of $(0,1)$.
Finally,
\[
 \nabla_BF(x)=\int_{[0,1]^B}\partial_BF(x^{B,t})\,\dd t,\qquad
 \binom mh\le n^{h}/h!
\]
and \eqref{row:eq:extension-continuous-higher} prove \eqref{row:eq:extension-higher}.
\end{proof}

\section{Resolvents and preliminary conditional cumulant bounds}\label{app:resolvent}
We study possibly nonreversible coordinate-update chains on binary
cubes of dimension $d=O(n)$, with bounds uniform in external parameters.
\subsection{Resolvents and conditional moment bounds}\label{row:sec:perturbation}\label{app:local-kernel-proof}\label{app:difference-moments}
Following the influence-matrix comparison approach of \cite{RvH2014},
we control resolvents and the higher differences needed for cumulants.
\subsubsection{Coordinate updates and resolvent assumptions}
\label{row:sec:auxiliary-hypotheses}
Let $\widetilde\pi_i(x)=\widetilde\pi_i(x_{-i})$ on $\{0,1\}^d$;
these updates need not form a Gibbs specification. For
$\nabla_i f=f(x^{i,1})-f(x^{i,0})$, define

\begin{equation}\label{v8:nearby:app:difference-seminorms}
 \begin{aligned}
 S_{h}(f)&:=\sum_{\substack{J\subseteq[d]\\|J|=h}}
                       \|\nabla_Jf\|_\infty,\qquad h\ge1,\\
 S_1(f)&=\sum_{i=1}^d\max_{x\in\{0,1\}^d}
                    |f(x^{i,1})-f(x^{i,0})|.
 \end{aligned}
\end{equation}

We impose the following assumptions uniformly in the dimension and
in any external parameters.
\begin{enumerate}[label=\textup{(R\arabic*)},ref=R\arabic*,leftmargin=3em]
\item\label{row:ass:aux-nondegenerate}
$\varepsilon_*\le\widetilde\pi_i\le1-\varepsilon_*$ for a fixed
$\varepsilon_*>0$.
\item\label{row:ass:aux-contraction}
There is a symmetric nonnegative matrix $\mathsf C=(c_{ij})$ with
$c_{ii}=0$, $\|\nabla_j\widetilde\pi_i\|_\infty\le c_{ij}$, and
$\max_i\sum_jc_{ij}\le\rho<1$.
\item\label{row:ass:aux-higher}
For every fixed $b\ge2$,
$\sup_i\sum_{|J|=b}\|\nabla_J\widetilde\pi_i\|_\infty
\le C_b(1+\log n)^{b-1}$. For estimates of order $h$, only the bounds with $b\le h$ are required.
\end{enumerate}
Define
\[
 \widetilde E_if=\widetilde\pi_if(x^{i,1})
                     +(1-\widetilde\pi_i)f(x^{i,0}),\qquad
 \widetilde{\mathcal L}=\sum_i(\widetilde E_i-I),\qquad
 P_t=e^{t\widetilde{\mathcal L}} .
\]
Each coordinate updates at rate one. Assumption
\ref{row:ass:aux-nondegenerate} makes the finite chain irreducible,
with unique stationary law $\widetilde\mu$. Define
\begin{equation}\label{v8:nearby:app:resolvent-definition}
 (\widetilde Rf)(x)
 :=\int_0^\infty[(P_tf)(x)-\E_{\widetilde\mu}f]\,\dd t.
\end{equation}

\begin{lemma}\label{row:lem:resolvent}
Assume \ref{row:ass:aux-nondegenerate} and \ref{row:ass:aux-contraction}
on the finite cube above. The auxiliary chain has a unique stationary
law $\widetilde\mu$. For every $f:\{0,1\}^d\to\R$, the defining integral
for $\widetilde Rf$ converges absolutely in the supremum norm, and
\begin{equation}\label{row:eq:resolvent-basic}
 \begin{aligned}
 -\widetilde{\mathcal L}\widetilde Rf
     &=f-\E_{\widetilde\mu}[f],\qquad
       \E_{\widetilde\mu}[\widetilde Rf]=0,\\
 \|\widetilde Rf\|_\infty
     &\le S_1(\widetilde Rf)\le\frac{S_1(f)}{1-\rho}.
 \end{aligned}
\end{equation}
The coordinatewise estimates are
\begin{equation}\label{row:eq:resolvent-gradient}
 \begin{aligned}
 \|\nabla_j\widetilde Rf\|_\infty
    &\le\sum_k[(I_d-\mathsf C)^{-1}]_{jk}\|\nabla_kf\|_\infty,\\
 [(I_d-\mathsf C)^{-1}]_{jk}&=[(I_d-\mathsf C)^{-1}]_{kj}\ge0,\qquad
 \max_j\sum_k[(I_d-\mathsf C)^{-1}]_{jk}\le(1-\rho)^{-1}.
 \end{aligned}
\end{equation}
If additionally $c_{ij}\le C_*/n$, then
\begin{equation}\label{row:eq:resolvent-dense-entries}
 [(I_d-\mathsf C)^{-1}]_{jk}\le\mathbf1_{\{j=k\}}+
                         \frac{C_*}{n(1-\rho)}.
\end{equation}
Finally, if \ref{row:ass:aux-higher} holds through a fixed order $h\ge2$, then
\begin{equation}\label{row:eq:R-bound}
 S_{h}(\widetilde Rf)\le C_{h}(1+\log n)^{h-1}\sum_{j=1}^{h}S_j(f).
\end{equation}
For $h=1$, the same bound follows from \eqref{row:eq:resolvent-basic}.
Here $C_{h}$ depends only on $h$, $1-\rho$, and the constants $C_b$,
$2\le b\le h$, in \ref{row:ass:aux-higher}. Under these uniform
hypotheses, the constants are independent of $d,n$ and of the particular
kernels. No reversibility, Gibbs representation, or good-set assumption
is needed.
\end{lemma}

\begin{proof}
\proofstep{First differences.}
Irreducibility follows from \ref{row:ass:aux-nondegenerate}
\cite[Theorem~20.1]{LP17}. Each coordinate update measure has mass one,
$\nabla_i\widetilde\pi_i=0$, and
$|\nabla_j\widetilde\pi_i|\le c_{ij}$. The Bernoulli Wasserstein
distance is the parameter difference, so
\cite[Theorem~3.1, \textup{(3.5)}]{Wu2006} gives
\begin{equation}\label{audit:semigroup-differences}
 \begin{aligned}
 \|\nabla_jP_tf\|_\infty
 &\le\sum_k[e^{(\mathsf C-I_d)t}]_{jk}\|\nabla_kf\|_\infty,\\
 S_1(P_tf)&\le e^{-(1-\rho)t}S_1(f).
 \end{aligned}
\end{equation}
The second line uses symmetry and
$e^{(\mathsf C-I_d)t}\one\le e^{-(1-\rho)t}\one$.
Since $P_t$ preserves $\widetilde\mu$-centering and
$\|u\|_\infty\le S_1(u)$ for centered $u$,
\eqref{audit:semigroup-differences} gives absolute convergence of the
integral defining $\widetilde Rf$ and
\[
 \E_{\widetilde\mu}\widetilde Rf=0,\qquad
 \|\widetilde Rf\|_\infty\le S_1(\widetilde Rf)
 \le\int_0^\infty S_1(P_tf)\,\dd t
 \le(1-\rho)^{-1}S_1(f).
\]
The Laplace resolvent theorem
\cite[Theorem~II.1.10\textup{(i)}]{EngelNagel2006}, applied on the
invariant space of $\widetilde\mu$-centered functions, gives
\begin{equation}\label{v8:proof:row:lem:resolvent:poisson}
 -\widetilde{\mathcal L}\widetilde Rf=f-\E_{\widetilde\mu}f.
\end{equation}
Thus \eqref{row:eq:resolvent-basic} holds. Integrating the first line of
\eqref{audit:semigroup-differences} gives \eqref{row:eq:resolvent-gradient}:
\begin{equation}\label{audit:influence-inverse-series}
 \int_0^\infty e^{(\mathsf C-I_d)t}\,\dd t
 =\sum_{m\ge0}\mathsf C^m=(I_d-\mathsf C)^{-1},\qquad
 \sum_k[(I_d-\mathsf C)^{-1}]_{jk}\le\sum_{m\ge0}\rho^m.
\end{equation}
The series is symmetric and nonnegative. If $c_{ij}\le C_*/n$, then
$ (\mathsf C^m)_{jk}\le C_*\rho^{m-1}/n$ for $m\ge1$.
Summing in \eqref{audit:influence-inverse-series} proves
\eqref{row:eq:resolvent-dense-entries}.

\proofstep{Higher differences.}
Write $u^{D=a}(x)=u(x^{D,a\mathbf1})$. The discrete product rule is
\begin{equation}\label{audit:discrete-leibniz}
 \nabla_J(uv)=\sum_{B\subseteq J}
 (\nabla_Bu)^{J\setminus B=0}(\nabla_{J\setminus B}v)^{B=1}.
\end{equation}
Fix $|J|=h\le d$ and let $\widetilde E_i^{J=0}$ use the probability
$\widetilde\pi_i(x^{J,\mathbf0})$. Since
$\widetilde E_i \widetilde Rf=(\widetilde Rf)^{i=0}+\widetilde\pi_i\nabla_i \widetilde Rf$,
\[
 \begin{aligned}
 \nabla_J\widetilde E_i\widetilde Rf
 &=\widetilde E_i^{J=0}\nabla_J\widetilde Rf+
 \sum_{\varnothing\ne B\subseteq J}
 (\nabla_B\widetilde\pi_i)^{J\setminus B=0}
 (\nabla_{(J\setminus B)\cup\{i\}}\widetilde Rf)^{B=1},&&i\notin J,\\
 \nabla_J\widetilde E_i \widetilde Rf&=0,&&i\in J.
 \end{aligned}
\]
For $Q_J=\sum_{i\notin J}(\widetilde E_i^{J=0}-I)$,
\eqref{v8:proof:row:lem:resolvent:poisson} becomes
\begin{equation}\label{audit:differenced-poisson}
 \begin{aligned}
 (hI-Q_J)\nabla_J\widetilde Rf
 &=\nabla_Jf\\
 &\quad+\sum_{i\notin J}\sum_{\varnothing\ne B\subseteq J}
 (\nabla_B\widetilde\pi_i)^{J\setminus B=0}
 (\nabla_{(J\setminus B)\cup\{i\}}\widetilde Rf)^{B=1}.
 \end{aligned}
\end{equation}
The finite Markov semigroup $e^{tQ_J}$ is a contraction in the supremum
norm. Since $h>0$, \cite[Theorem~II.1.10\textup{(ii)--(iii)}]{EngelNagel2006}
gives
\begin{equation}\label{audit:killed-resolvent}
 (hI-Q_J)^{-1}u=\int_0^\infty e^{-ht}e^{tQ_J}u\,\dd t,
 \qquad \|(hI-Q_J)^{-1}u\|_\infty\le h^{-1}\|u\|_\infty.
\end{equation}
Equations~\eqref{audit:differenced-poisson} and \eqref{audit:killed-resolvent} imply
\begin{equation}\label{audit:resolvent-difference-inequality}
 \begin{aligned}
 h\|\nabla_J\widetilde Rf\|_\infty
 &\le\|\nabla_Jf\|_\infty\\
 &\quad+\sum_{i\notin J}\sum_{\varnothing\ne B\subseteq J}
 \|\nabla_B\widetilde\pi_i\|_\infty
 \|\nabla_{(J\setminus B)\cup\{i\}}\widetilde Rf\|_\infty.
 \end{aligned}
\end{equation}
For $B=\{j\}$, reindex by $D=(J\setminus\{j\})\cup\{i\}$:
\begin{equation}\label{audit:resolvent-singleton-contribution}
 \begin{aligned}
 &\sum_{|J|=h}\sum_{i\notin J}\sum_{j\in J}
 \|\nabla_j\widetilde\pi_i\|_\infty
 \|\nabla_{(J\setminus\{j\})\cup\{i\}}\widetilde Rf\|_\infty\\
 &\quad\le\sum_{|D|=h}\|\nabla_D\widetilde Rf\|_\infty
       \left(\sum_{i\in D}\sum_{j\notin D}c_{ij}\right)
 \le h\rho S_{h}(\widetilde Rf).
 \end{aligned}
\end{equation}
For $|B|=b\ge2$, use $D=(J\setminus B)\cup\{i\}$ and
$J=(D\setminus\{i\})\cup B$:
\begin{equation}\label{audit:resolvent-higher-block-contribution}
 \begin{aligned}
 &\sum_{|J|=h}\sum_{i\notin J}\sum_{\substack{B\subseteq J\\|B|=b}}
 \|\nabla_B\widetilde\pi_i\|_\infty
 \|\nabla_{(J\setminus B)\cup\{i\}}\widetilde Rf\|_\infty\\
 &=\sum_{|D|=h-b+1}\|\nabla_D\widetilde Rf\|_\infty
   \left(\sum_{i\in D}\sum_{\substack{|B|=b\\B\cap D=\varnothing}}
   \|\nabla_B\widetilde\pi_i\|_\infty\right)\\
 &\le(h-b+1)C_b(1+\log n)^{b-1}S_{h-b+1}(\widetilde Rf).
 \end{aligned}
\end{equation}
Summing \eqref{audit:resolvent-difference-inequality}, using
\eqref{audit:resolvent-singleton-contribution} and
\eqref{audit:resolvent-higher-block-contribution}, gives
\begin{equation}\label{audit:resolvent-seminorm-recurrence}
 \begin{aligned}
 h(1-\rho)S_{h}(\widetilde Rf)
 &\le S_{h}(f)\\
 &\quad+\sum_{b=2}^{h}(h-b+1)C_b(1+\log n)^{b-1}
                            S_{h-b+1}(\widetilde Rf).
 \end{aligned}
\end{equation}
Induction in \eqref{audit:resolvent-seminorm-recurrence}, starting from
\eqref{row:eq:resolvent-basic}, proves
\eqref{row:eq:R-bound}, since $(b-1)+(h-b)=h-1$.
For $h=d$, $Q_J=0$; for $h>d$, $S_{h}(\widetilde Rf)=0$.
\end{proof}

\begin{proof}[Proof of Lemma~\ref{lem:local-input}]
Lemma~\ref{row:lem:good-exteriors} gives the good sets and tails,
uniformly in $\mathbf e,s$. Proposition~\ref{row:prop:extension} gives
single-row kernels in $[p/2,(1+p)/2]$, agreeing with
$\sigma(l_{i,A}(x))$ on $\Omega_A$, with symmetric influence majorant
$c_{ij}\le C/n$ and row sums at most $3/4$. Define
\begin{equation}\label{audit:dependency:app:multirow-update-kernel}
 \widetilde\pi_{(b,i)}^{A,\mathbf e,s}(\mathbf x)
 =\sig\left(\log\frac{\widetilde\pi_{A,i}(x^b)}
 {1-\widetilde\pi_{A,i}(x^b)}
              +s\partial_{(b,i)}V_{A,\mathbf e}(\mathbf x)\right).
\end{equation}
The log odds in \eqref{audit:dependency:app:multirow-update-kernel}
equal the true field on $\Omega_A^{q_0}$.
For distinct $u=(b,i),v=(c,j)$, Lemma~\ref{lem:row-deterministic-bounds}
gives
\[
 \|\partial_uV\|_\infty\le C/n,\qquad
 \|\partial_v\partial_uV\|_\infty
 \le Cn^{-1}\one_{\{i=j,b\ne c\}}+Cn^{-2}\one_{\{i\ne j\}}.
\]
The map $\sigma(\log(\theta/(1-\theta))+z)$ has derivative bounds
$e^{|z|}$ in $\theta$ and $1/4$ in $z$. Thus the symmetric matrix
\begin{equation}\label{audit:dependency:app:multirow-influence-majorant}
 T_{uv}=e^{C/n}\one_{\{b=c\}}c_{ij}
 +\tfrac14\max_{\mathbf e}\|\nabla_v\nabla_uV_{A,\mathbf e}\|_\infty,
 \qquad T_{uu}=0,
\end{equation}
majorizes the first differences and satisfies
$T_{uv}\le C/n$ and
$\max_u\sum_vT_{uv}\le(3/4)e^{C/n}+C/n<7/8$.
There are only $q_0-1$ other rows at the same site.
The update perturbation is at most $C/(4n)$, proving nondegeneracy.

For every fixed $h\ge2$, we also have
\begin{equation}\label{audit:local-higher-differences}
 \sup_{A\in\mathcal A_n,\,\mathbf e,\,s\in[0,1],\,u\in\mathcal I}
 \sum_{\substack{J\subseteq\mathcal I\\|J|=h}}
 \|\nabla_J\widetilde\pi_u^{A,\mathbf e,s}\|_\infty
 \le C_{h}(1+\log n)^{h-1}.
\end{equation}
Indeed, \eqref{audit:interrow-coefficient-budget} gives
\begin{equation}\label{eq:local-interrow-coefficients}
 \partial_uV_{A,\mathbf e}
 =\sum_{\substack{S\subseteq\mathcal I\setminus\{u\}\\|S|\le D_*}}
 a_{u,S}\prod_{v\in S}x_v,
 \qquad\sum_S|a_{u,S}|\le C/n.
\end{equation}
Hence $\sum_{|J|=h}\|\partial_J\partial_uV\|_\infty\le C_{h}/n$.
Equation~\eqref{row:eq:extension-continuous-higher} bounds the order-$h$
single-row log-odds derivatives by $Cn^{-h}(1+\log n)^{h-1}$.
In the partition chain rule for $\sigma$, a $t$-block partition
therefore contributes at most $C_{h}(1+\log n)^{h-t}$ after summation.
Summing the finitely many partitions and using
$\nabla_JF=\int_{[0,1]^J}\partial_JF(\mathbf x^{J,t})\,\dd t$
proves \eqref{audit:local-higher-differences}.
The construction in \eqref{audit:dependency:app:multirow-update-kernel},
the bound in \eqref{audit:dependency:app:multirow-influence-majorant}, and
the estimate \eqref{audit:local-higher-differences} verify
\ref{row:ass:aux-nondegenerate}--\ref{row:ass:aux-higher};
Lemma~\ref{row:lem:resolvent} supplies the stationary law,
\eqref{eq:local-resolvent-definition}, and
\eqref{eq:local-resolvent-bounds}.
\end{proof}

\begin{lemma}
\label{lem:conditional-difference-moments}
Assume that the hypotheses of Lemma~\ref{lem:cov-error} hold.
Then, for every real-valued function $f$ on the conditional cube
and every fixed real $q\ge1$,
\begin{equation}\label{eq:conditional-difference-moments}
 \|f-\E_\mu f\|_{L^{q}(\mu)}
 \le C_{q}\left(\sum_{i\in\mathcal I}
                       \|\nabla_i f\|_\infty^2\right)^{1/2}.
\end{equation}
The constants depend only on $q$ and the uniform hypotheses of that
lemma. In particular, they are uniform for all fixed conditional data
$(A,\mathbf e,s)$ of Lemma~\ref{lem:local-input}.
\end{lemma}
\begin{proof}
Set $a_i=\|\nabla_i f\|_\infty$ and $d=|\mathcal I|=O(n)$.
Symmetry of $(I_d-T)^{-1}$ gives
\[
 a^{\mathsf T}(I_d-T)^{-1}a\le(1-\rho)^{-1}\|a\|_2^2,\qquad
 \operatorname{osc}(f)\le\sqrt d\,\|a\|_2.
\]
Apply \eqref{eq:cov-error} with $g=(f-\E_\mu f)^{2j-1}$. Using
$|\nabla_i g|\le C_{j}a_i(|f-\E_\mu f|^{2j-2}+a_i^{2j-2})$ and
$\|\nabla_i\widetilde Rf\|_\infty\le((I_d-T)^{-1}a)_i$ yields
\begin{equation}\label{audit:conditional-moment-recursion}
 \E_\mu|f-\E_\mu f|^{2j}
 \le C_{j}\|a\|_2^2\E_\mu|f-\E_\mu f|^{2j-2}
 +(C_{j}+Ce^{-cn}d^{j})\|a\|_2^{2j}.
\end{equation}
The last term follows from \eqref{audit:covariance-discrepancy}
and the oscillation bound. Since $\sup_n e^{-cn}d^{j}<\infty$,
induction in \eqref{audit:conditional-moment-recursion} from the zeroth moment gives
$\|f-\E_\mu f\|_{2j}\le C_{j}\|a\|_2$. Choosing $2j\ge q$ proves
\eqref{eq:conditional-difference-moments}.
\end{proof}

\subsection{Preliminary conditional cumulant bounds}\label{row:sec:gibbs-comparison}
Let $\mathcal P(A)$ denote the set of partitions of the formal argument positions
$A$. Expanding $\exp(\log M)$ and $\log M$ for
$M(\mathbf t)=\E e^{\sum_jt_jZ_j}$ gives
\eqref{row:eq:moment-partition} and its inverse
\begin{equation}\label{row:eq:cumulant-partition}
 \cum(Z_j:j\in A)
 =\sum_{\pi\in\mathcal P(A)}(-1)^{|\pi|-1}(|\pi|-1)!
 \prod_{B\in\pi}\E\prod_{j\in B}Z_j.
\end{equation}
Cumulants are symmetric and multilinear; those of order at least
two vanish if an argument is constant. For
$\dd\mu_z/\dd\mu=e^{\sum_az_aZ_a}/\E_\mu e^{\sum_az_aZ_a}$,
differentiating the tilted expectation and applying the product rule gives
\begin{equation}\label{row:eq:tilted-expectation-derivative}
 \left.\partial_{z_{[m]}}\E_{\mu_z}u_z\right|_{z=0}
 =\sum_{I\subseteq[m]}
 \cum_\mu\!\left(\left.\partial_{z_I}u_z\right|_{z=0},
 Z_a:a\notin I\right).
\end{equation}
Likewise, differentiating
$\E_{\mu_z}UV=\Cov_{\mu_z}(U,V)+\E_{\mu_z}U\,\E_{\mu_z}V$ yields
\begin{equation}\label{row:eq:binary-product-cumulant}
 \begin{aligned}
 \cum_\mu(UV,Z_1,\ldots,Z_m)
 &=\cum_\mu(U,V,Z_1,\ldots,Z_m)\\
 &\quad+\sum_{S\subseteq[m]}
 \cum_\mu(U,Z_a:a\in S)\cum_\mu(V,Z_a:a\notin S).
 \end{aligned}
\end{equation}
Repeated variables occupy distinct formal positions, preserving all multiplicities.

We combine Lemma~\ref{lem:local-input}\textup{(iii)} with
\eqref{row:eq:moment-partition} and
\eqref{row:eq:cumulant-partition}--\eqref{row:eq:binary-product-cumulant}
to bound cumulants under the true conditional Gibbs laws.

\begin{theorem}\label{row:thm:main}
Assume \hyperref[ass:PC]{\textup{(PC)}}.
Fix $q_0\ge1$ independently of $n$, and let the good exterior sets be
those supplied by Lemma~\ref{lem:local-input}. They satisfy
\begin{equation}\label{row:eq:good-exterior-prob}
 \Pp_n(A\notin\mathcal A_n)\le Ce^{-cn}.
\end{equation}
For $A\in\mathcal A_n$, take $\mu=\pi_A$ or $\mu_{A,\mathbf e,s}$
from Lemma~\ref{lem:local-input}, on its own coordinate cube, with
prescribed $\mathbf e$ and $s\in[0,1]$. Then
\begin{align}
 \max_{i\in[d]}\sum_{j\in[d]}|\Cov_\mu(X_i,X_j)|&\le C,
                                                   \label{row:eq:cov-bound}\\
 \max_{i_1\in[d]}\sum_{i_2,\ldots,i_{\ell}\in[d]}
 |\cum_\mu(X_{i_1},\ldots,X_{i_{\ell}})|
 &\le C_{\ell}(1+\log n)^{(\ell-1)(\ell-2)/2},\qquad \ell\ge3.\label{row:eq:cum-bound}
\end{align}
The orders are fixed independently of $n$, and the sums include repeated
coordinate labels. More generally, the response estimate
\eqref{row:eq:response-seminorm} holds for every real-valued test function
$f$ on the chosen cube, without any locality assumption.
The constants are uniform in $n$, $A\in\mathcal A_n$, $\mathbf e$, $s$,
and the coordinate locations; they may depend on $q_0$, the indicated
orders, and the fixed data of Lemma~\ref{lem:local-input}. No sign
condition on the cumulants is required.
\end{theorem}

\begin{proof}[Proof of Theorem~\ref{row:thm:main}]
The exterior bound follows from \eqref{eq:local-good-exteriors}.
Fix the conditional data and write $\mu\propto e^H$ on
$\{0,1\}^d$, $d=O(n)$. Lemma~\ref{lem:local-input} and
\eqref{audit:local-higher-differences} supply the resolvent hypotheses and
\[
 \widetilde\pi_i=\sigma(\partial_iH)\text{ on }\Omega,
 \qquad\mu(\Omega^c)\le Ce^{-cn}.
\]
For every test function $f$, we prove
\begin{equation}\label{row:eq:response-seminorm}
 \sum_{j_1,\ldots,j_\ell}
 |\cum_\mu(f,X_{j_1},\ldots,X_{j_\ell})|
 \le C_\ell(1+\log n)^{\ell(\ell-1)/2}
                        \sum_{k=1}^{\ell}S_k(f),\qquad\ell\ge1.
\end{equation}

\proofstep{Tilt identity.}
Keep $\widetilde Rf$ and its generator fixed. Give each formal
position a separate parameter, even when coordinate labels repeat:
\[
 \mu_z\propto e^{H+\sum_{a=1}^{\ell}z_aX_{j_a}},\qquad
 \mathcal L_z\widetilde Rf=\sum_i\left[
 \sigma\!\left(\partial_iH+\sum_a z_a\one_{\{j_a=i\}}\right)-x_i
 \right]\nabla_i \widetilde Rf.
\]
Stationarity and the auxiliary Poisson identity give
\begin{equation}\label{audit:tilt-poisson-identity}
 \begin{aligned}
 \E_{\mu_z}\mathcal L_z\widetilde Rf&=0,\\
 \E_{\mu_z}f-\E_{\widetilde\mu}f
 &=-\E_{\mu_z}\widetilde{\mathcal L}\widetilde Rf\\
 &=\E_{\mu_z}(\mathcal L_z-\widetilde{\mathcal L})\widetilde Rf.
 \end{aligned}
\end{equation}
For $\varnothing\ne I\subseteq[\ell]$, set
\begin{equation}\label{audit:GI-definition}
 \begin{aligned}
 G_I(\widetilde Rf)
 &:=\left.\partial_{z_I}\mathcal L_z\widetilde Rf\right|_{z=0}\\
 &=\sum_i\sigma^{(|I|)}(\partial_iH)
          \prod_{a\in I}\one_{\{j_a=i\}}\nabla_i\widetilde Rf.
 \end{aligned}
\end{equation}
Differentiating \eqref{audit:tilt-poisson-identity} by
\eqref{row:eq:tilted-expectation-derivative} yields
\begin{equation}\label{audit:cumulant-response-identity}
 \begin{aligned}
 &\cum_\mu(f,X_{j_1},\ldots,X_{j_\ell})\\
 &=\cum_\mu((\mathcal L-\widetilde{\mathcal L})\widetilde Rf,
                                  X_{j_1},\ldots,X_{j_\ell})\\
 &\quad+\sum_{\varnothing\ne I\subseteq[\ell]}
       \cum_\mu(G_I(\widetilde Rf),X_{j_a}:a\notin I),\qquad\mathcal L=\mathcal L_0.
 \end{aligned}
\end{equation}
A one-argument cumulant is an expectation.

\proofstep{Difference bounds.}
The coefficient bounds \eqref{row:eq:dense-derivatives} and
\eqref{eq:local-interrow-coefficients} imply, for a fixed $D$,
\begin{equation}\label{audit:response-H-derivatives}
 \begin{gathered}
 H(x)=\sum_{|T|\le D}h_T\prod_{i\in T}x_i,\\
 \sup_i\sum_{T\ni i}|h_T|
 \le\sum_{h=1}^D\binom{n-q_0-1}{h-1}C_hn^{1-h}+C/n\le C,\\
 \sup_i\sum_{|J|=b}\|\partial_J\partial_iH\|_\infty
 \le\sup_i\sum_{T\ni i}\binom{|T|-1}{b}|h_T|\le C_b.
 \end{gathered}
\end{equation}
On the continuous cube, the partition chain rule gives
\[
 \partial_J\{\sigma^{(m)}(\partial_iH)\}
 =\sum_{\pi\in\mathcal P(J)}\sigma^{(m+|\pi|)}(\partial_iH)
                       \prod_{B\in\pi}\partial_B\partial_iH.
\]
Equation~\eqref{audit:response-H-derivatives}, bounded sigmoid derivatives,
and integration of mixed derivatives give
\begin{equation}\label{audit:sigmoid-difference-sums}
 \begin{aligned}
 \sup_i\sum_{|J|=b}\|\nabla_J\{\sigma^{(m)}(\partial_iH)\}\|_\infty
 &\le\sup_i\sum_{|J|=b}\|\partial_J\{\sigma^{(m)}(\partial_iH)\}\|_\infty\\
 &\le C_{b,m}\sum_{t=1}^b
       \sum_{\substack{b_1+\cdots+b_t=b\\b_v\ge1}}
       \prod_{v=1}^t C_{b_v}\le C_{b,m}.
 \end{aligned}
\end{equation}
The $b=0$ bound is immediate. The indicators in \eqref{audit:GI-definition}
force all $I$-labels to denote the same coordinate, so
\begin{equation}\label{audit:GI-supnorm}
 \sum_{(j_a)_{a\in I}}\|G_I(\widetilde Rf)\|_\infty
 =\sum_i\|\sigma^{(|I|)}(\partial_iH)\nabla_i \widetilde Rf\|_\infty
 \le C_{|I|}S_1(\widetilde Rf).
\end{equation}
The discrete product rule \eqref{audit:discrete-leibniz}, grouped by
the number $b$ of differences on $\sigma^{(m)}(\partial_iH)$, and
\eqref{audit:sigmoid-difference-sums} yield
\begin{equation}\label{audit:dependency:app:response-product-differences}
 \sum_iS_k\{\sigma^{(m)}(\partial_iH)\nabla_i \widetilde Rf\}
 \le\sum_{b=0}^k C_{b,m}
      \sum_i\sum_{|J|=k-b}\|\nabla_J\nabla_i \widetilde Rf\|_\infty.
\end{equation}
Repeated differences vanish, and each $(h+1)$-set has $h+1$
possible distinguished elements:
\begin{equation}\label{audit:difference-extension-count}
 \sum_i\sum_{|J|=h}\|\nabla_J\nabla_i \widetilde Rf\|_\infty
 =\sum_i\sum_{\substack{|J|=h\\i\notin J}}
       \|\nabla_{J\cup\{i\}}\widetilde Rf\|_\infty
 =(h+1)S_{h+1}(\widetilde Rf).
\end{equation}
Combining \eqref{audit:dependency:app:response-product-differences} with
\eqref{audit:difference-extension-count} gives
\begin{equation}\label{audit:GI-seminorm}
 \begin{aligned}
 \sum_{(j_a)_{a\in I}}S_k(G_I(\widetilde Rf))
 &=\sum_iS_k(\sigma^{(|I|)}(\partial_iH)\nabla_i\widetilde Rf)\\
 &\le C_{k,|I|}\sum_{v=1}^{k+1}S_v(\widetilde Rf),\qquad k\ge1.
 \end{aligned}
\end{equation}

\proofstep{Kernel discrepancy.}
Under the original law, agreement on $\Omega$ gives
\begin{equation}\label{audit:response-kernel-discrepancy}
 \begin{aligned}
 (\mathcal L-\widetilde{\mathcal L})\widetilde Rf
 &=\sum_i\{\sigma(\partial_iH)-\widetilde\pi_i\}\nabla_i\widetilde Rf,\\
 |(\mathcal L-\widetilde{\mathcal L})\widetilde Rf|
 &\le\one_{\Omega^c}S_1(\widetilde Rf).
 \end{aligned}
\end{equation}
Using \eqref{audit:response-kernel-discrepancy} in
\eqref{row:eq:cumulant-partition}, bound the block containing the
kernel-discrepancy term by its first absolute moment and all other
blocks by one. Therefore
\begin{equation}\label{audit:response-bad-set}
 \begin{aligned}
 &|\cum_\mu((\mathcal L-\widetilde{\mathcal L})\widetilde Rf,
                                  X_{j_1},\ldots,X_{j_\ell})|\\
 &\quad\le C_\ell\mu(\Omega^c)S_1(\widetilde Rf),\\
 &\sum_{j_1,\ldots,j_\ell}
 |\cum_\mu((\mathcal L-\widetilde{\mathcal L})\widetilde Rf,
                                  X_{j_1},\ldots,X_{j_\ell})|\\
 &\quad\le C_\ell d^\ell e^{-cn}S_1(\widetilde Rf)
 \le C_\ell n^\ell e^{-cn}S_1(f).
 \end{aligned}
\end{equation}

\proofstep{Order induction.}
For $\ell=1$, \eqref{audit:cumulant-response-identity}, with
\eqref{audit:response-bad-set} and \eqref{row:eq:resolvent-basic}, gives
\[
 \sum_j|\Cov_\mu(f,X_j)|
 \le Cne^{-cn}S_1(f)+\tfrac14S_1(\widetilde Rf)\le CS_1(f).
\]
Assume \eqref{row:eq:response-seminorm} for all lower orders and all
test functions. For $\varnothing\ne I\subsetneq[\ell]$, let
$t=\ell-|I|$. Apply the order-$t$ bound with $I$-labels fixed, then sum
over those labels using \eqref{audit:GI-seminorm} and \eqref{row:eq:R-bound}:
\begin{equation}\label{audit:response-induction-term}
 \begin{aligned}
 &\sum_{j_1,\ldots,j_\ell}
       |\cum_\mu(G_I(\widetilde Rf),X_{j_a}:a\notin I)|\\
 &\le C_t(1+\log n)^{t(t-1)/2}
       \sum_{(j_a)_{a\in I}}\sum_{k=1}^tS_k(G_I(\widetilde Rf))\\
 &\le C_\ell(1+\log n)^{t(t-1)/2}\sum_{v=1}^{t+1}S_v(\widetilde Rf)\\
 &\le C_\ell(1+\log n)^{t(t+1)/2}\sum_{v=1}^{\ell}S_v(f),\\
 &\quad t(t+1)/2\le\ell(\ell-1)/2.
 \end{aligned}
\end{equation}
For $I=[\ell]$, \eqref{audit:GI-supnorm} gives
\begin{equation}\label{audit:response-full-label-term}
 \sum_{j_1,\ldots,j_\ell}|\E_\mu G_{[\ell]}(\widetilde Rf)|
 \le C_\ell S_1(\widetilde Rf)\le C_\ell S_1(f).
\end{equation}
Combining \eqref{audit:response-bad-set},
\eqref{audit:response-induction-term}, and \eqref{audit:response-full-label-term}
in \eqref{audit:cumulant-response-identity} proves
\eqref{row:eq:response-seminorm}, since
$\sup_n n^\ell e^{-cn}<\infty$. Taking $f=X_i$, so
$S_1(f)=1$ and $S_k(f)=0$ for $k\ge2$, proves
\eqref{row:eq:cov-bound} and \eqref{row:eq:cum-bound}.
For $\pi_A$, take $s=0$ and restrict the product law and coordinate sums
to one row. All estimates are uniform in the fixed conditional data.
\end{proof}

\section{Proofs of uniform conditional cumulant bounds}\label{ep:row:section}
We establish the conditional mean equation and a common contraction
neighborhood, then remove the logarithmic factors in the preliminary
bounds of Theorem~\ref{row:thm:main}. Coordinate cumulants are treated
first, followed by local observables and normalized block polynomials.
Throughout, the conditional laws and seminorms are those of
Section~\ref{row:sec:conditional-rows}; all orders are fixed.

\subsection{Conditional mean equations and reference-row contraction}\label{app:mean-equation-proof}\label{app:common-neighborhood-proof}
\begin{proof}[Proof of Lemma~\ref{ep:row:marginal-mean}]
Work under the full conditional law and set $\bar x=\E X^b$.
The degree-$h$ coefficients are $O(n^{-h})$, so
\eqref{eq:row-polynomial-mean} and
\eqref{eq:conditional-difference-moments} give
$|\E\{l_{i,A}(X^b)-l_{i,A}(\bar x)\}|\le C/n$ and
$\|l_{i,A}(X^b)-l_{i,A}(\bar x)\|_{q}\le C_{q}n^{-1/2}$.
Taylor's formula, using \eqref{ep:row:field-moments} just established,
retains the signed linear bias:
\begin{equation}\label{ep:row:mean-field-error}
 \begin{aligned}
 &|\E\sigma(l_{i,A}(X^b))-\sigma(l_{i,A}(\bar x))|\\
 &\quad\le\|\sigma'\|_\infty
           |\E\{l_{i,A}(X^b)-l_{i,A}(\bar x)\}|\\
 &\qquad+\tfrac12\|\sigma''\|_\infty
           \E|l_{i,A}(X^b)-l_{i,A}(\bar x)|^2\le C/n.
 \end{aligned}
\end{equation}
The full field is $l_{i,A}(X^b)+s\partial_{(b,i)}V$, with
$\|s\partial_{(b,i)}V\|_\infty\le C/n$ by
\eqref{audit:row-energy-interaction-bounds}. Its conditional-mean identity,
together with \eqref{ep:row:mean-field-error}, gives
\begin{equation}\label{ep:row:joint-mean-error}
 |\bar x_i-\E\sigma(l_{i,A}(X^b))|\le C/n,\qquad
 |\bar x_i-\sigma(l_{i,A}(\bar x))|\le C/n.
\end{equation}
Equations~\eqref{row:eq:Omega} and \eqref{eq:local-kernel-agreement} give
$|\sigma(l_{i,A}(X^b))-p|\le a$ on $\Omega_A$ and a probability bound of $Ce^{-cn}$
for its complement. Together with \eqref{ep:row:joint-mean-error},
this yields $|\bar x_i-p|\le a+C/n$.
Finally, \eqref{ep:row:joint-mean-error}, $\sigma'(u)=\nu(\sigma(u))$,
and $\Lip(\nu)\le1$ give
$|\sigma'(l_{i,A}(\bar x))-\nu(\bar x_i)|\le C/n$.
These are the asserted estimates.
\end{proof}

\begin{proof}[Proof of Lemma~\ref{ep:row:mean-jacobian}]
Use the fixed choices
$K\ge2(1+C_{\rm det}+4/\nu(p))$,
$a\le\min\{\delta/28,p/8,(1-p)/8\}$ and
\eqref{v11:eq:tolerance-contraction}--\eqref{v11:eq:tolerance-saturation}.
Equation~\eqref{ep:row:mean-neighborhood} gives
\begin{equation}\label{audit:reference-mean-localization}
 \|\bar x^b-p\one\|_\infty\le a+C/n\le\delta,
 \qquad\bar x^b\in I^d\quad(n\ge n_0).
\end{equation}
This includes $\mu_A=\E_{\pi_A}X$ by taking $s=0$.
On $\Omega_A$, $\nu(u)\ge\nu(p)/4$ for $|u-p|\le a$ and
\eqref{row:eq:Omega} imply
\[
 |l_{i,A}(x)-\ell(p)|\le4a/\nu(p),\qquad
 |n\partial_jl_{i,A}(x)-(\alpha+\beta A_{ij})|\le\eta/4.
\]
Both polynomials satisfy \eqref{eq:row-polynomial-mean}.
With $\epsilon_n=C_*/n+C_*e^{-c_*n}$, this estimate and
\eqref{row:eq:good-properties} give
\begin{equation}\label{audit:reference-field-at-mean}
 \begin{aligned}
 |l_{i,A}(\mu_A)-\ell(p)|
 &\le C_*/n+4a/\nu(p)+C_*\pi_A(\Omega_A^c)
 \le4a/\nu(p)+\epsilon_n,\\
 |n\partial_jl_{i,A}(\mu_A)-(\alpha+\beta A_{ij})|
 &\le\eta/4+\epsilon_n.
 \end{aligned}
\end{equation}
For $y\in I^d$, \eqref{audit:reference-mean-localization} gives
$\|y-\mu_A\|_\infty\le2\delta$. Integrating the derivative bounds
\eqref{row:eq:dense-derivatives} and using \eqref{audit:reference-field-at-mean} yields
\begin{equation}\label{audit:reference-segment-transfer}
 \begin{aligned}
 |l_{i,A}(y)-\ell(p)|
 &\le4a/\nu(p)+2C_{\rm det}\delta+\epsilon_n
 \le K\delta+\epsilon_n,\\
 |n\partial_jl_{i,A}(y)-(\alpha+\beta A_{ij})|
 &\le\eta/4+2C_{\rm det}\delta+\epsilon_n
 \le\eta/4+K\delta+\epsilon_n.
 \end{aligned}
\end{equation}
Using \eqref{audit:reference-segment-transfer}, $\sigma'(\ell(p))=\nu(p)$,
$\|\sigma''\|_\infty\le1/4$,
and $\Lip(\nu)\le1$, we obtain, for $y,z\in I^d$,
\begin{equation}\label{audit:reference-logistic-factors}
 \begin{aligned}
 |\sigma'(l_{i,A}(y))-\nu(p)|&\le(K\delta+\epsilon_n)/4,\\
 \left|\sigma'\!\left(d^{-1}\sum_i l_{i,A}(y)\right)-\nu(p)\right|
 &\le(K\delta+\epsilon_n)/4,\\
 |\nu(z_i)-\nu(p)|&\le\delta.
 \end{aligned}
\end{equation}
Define
\begin{equation}\label{audit:reference-majorant-entries}
 C^{(n,A)}_{ij}=\frac{\one_{\{i\ne j\}}}{n}
 (\nu(p)+K\delta+\epsilon_n)
 (\alpha+\beta A_{ij}+\eta/4+K\delta+\epsilon_n).
\end{equation}
This matrix is symmetric, nonnegative, and independent of $\mathbf e,s,b$.
As $\partial_i l_{i,A}=0$, \eqref{ep:row:jacobian-definition},
\eqref{audit:reference-segment-transfer}, and \eqref{audit:reference-logistic-factors} give
\[
 |DM_A(y)|,\quad
 |\operatorname{diag}((\nu(z_i))_i)H_A''(y)|\le C^{(n,A)}.
\]
Using $d_{A,i}/n\le p+\varepsilon_{\rm deg}$ from \eqref{row:eq:good-properties}
in \eqref{audit:reference-majorant-entries}, together with
\eqref{v11:eq:tolerance-contraction}, we obtain
\begin{equation}\label{audit:reference-majorant-row-sum}
 \begin{aligned}
 \sum_jC^{(n,A)}_{ij}
 &=(\nu(p)+K\delta+\epsilon_n)
   \left[\frac{d-1}{n}(\alpha+\eta/4+K\delta+\epsilon_n)
                         +\beta\frac{d_{A,i}}n\right]\\
 &\le(\nu(p)+K\delta+\epsilon_n)
   (\alpha+\beta p+\beta\varepsilon_{\rm deg}
                         +\eta/4+K\delta+\epsilon_n)\\
 &\le5/8+1/8=\rho
 \end{aligned}
\end{equation}
for large $n$. Thus \eqref{ep:row:symmetric-majorant} holds.
For any $|B|\le C^{(n,A)}$, symmetry and
\eqref{audit:reference-majorant-row-sum} also control column sums:
\begin{equation}\label{audit:reference-dominated-norms}
 \|B\|_{\infty\to\infty},\ \|B\|_{1\to1}\le\rho,\qquad
 \|B\|_{\mathrm{op}}
 \le\sqrt{\|B\|_{1\to1}\|B\|_{\infty\to\infty}}\le\rho.
\end{equation}
Taking $B=DM_A(y)$ proves \eqref{ep:row:jacobian-gap}; taking
$B=\operatorname{diag}((\nu(z_i))_i)H_A''(y)$ and summing
$(I_d-B)^{-1}=\sum_{m\ge0}B^m$ proves \eqref{eq:shape-common-inverse}.
Finally, \eqref{audit:reference-logistic-factors},
\eqref{audit:reference-majorant-entries}, and \eqref{audit:reference-majorant-row-sum} give
\begin{equation}\label{audit:reference-averaged-gradient}
 \begin{aligned}
 \sum_j\left|\partial_j\sigma\!\left(d^{-1}\sum_i l_{i,A}(y)\right)\right|
 &\le d^{-1}\sum_{i,j}
 \sigma'\!\left(d^{-1}\sum_k l_{k,A}(y)\right)|\partial_jl_{i,A}(y)|\\
 &\le d^{-1}\sum_{i,j}C^{(n,A)}_{ij}\le\rho.
 \end{aligned}
\end{equation}
Integrating \eqref{audit:reference-averaged-gradient} along $[z,y]\subset I^d$ gives
\eqref{ep:row:average-field-contraction}.
\end{proof}

\subsection{Uniform coordinate cumulant bounds}
We will combine the preliminary bounds of Theorem~\ref{row:thm:main}
with the contraction estimate of Lemma~\ref{ep:row:mean-jacobian}
to obtain bounds without logarithmic loss. Fix a conditional law
$\mu=\mu_{A,\mathbf e,s}$ from Lemma~\ref{lem:local-input}, and put

\begin{equation}\label{v8:nearby:app:coordinate-cumulant-sum}
 \begin{aligned}
 K_{\ell}(\mu)&=\max_{i\in\mathcal I}\sum_{j_2,\ldots,j_{\ell}\in\mathcal I}
   |\cum_\mu(X_i,X_{j_2},\ldots,X_{j_{\ell}})|,\qquad \ell\ge2,\\
 \xi_i&=X_i-\E_\mu X_i.
 \end{aligned}
\end{equation}

\begin{lemma}\label{ep:row:coordinates}
Assume \hyperref[ass:PC]{\textup{(PC)}}.
For $\mu=\mu_{A,\mathbf e,s}$ in
\eqref{v8:def:conditional-multirow-law} and each fixed integer $\ell\ge2$,
$K_{\ell}(\mu)$ in \eqref{v8:nearby:app:coordinate-cumulant-sum} satisfies,
uniformly over
\[
 n\ge n_0,\qquad A\in\mathcal A_n,\qquad
 \mathbf e\in\{0,1\}^{\binom{q_0}{2}},\qquad s\in[0,1],
\]
\begin{equation}\label{ep:row:coordinate-bound}
 K_{\ell}(\mu)\le C_{\ell}.
\end{equation}
The ordered tuple sums include repeated labels. Constants depend only
on $\ell$ and the fixed data of Lemma~\ref{lem:local-input}.
In particular,
$\max_i\sum_j|\Cov_\mu(X_i,X_j)|\le C_2$ and
$\max_i\sum_{j,k}|\E_\mu\xi_i\xi_j\xi_k|\le C_3$.
\end{lemma}

\begin{proof}
We argue by induction on $\ell$. Theorem~\ref{row:thm:main} gives $K_2\le C$ and,
independently of the induction, $K_h\le C_h(1+\log n)^{(h-1)(h-2)/2}$
for every fixed $h$.
Fix $\ell\ge3$ and assume $K_h\le C_h$ for $2\le h<\ell$.
Equations~\eqref{ep:row:jacobian-definition} and \eqref{ep:row:jacobian-gap} give
\begin{equation}\label{audit:dependency:app:coordinate-jacobian-contraction}
 J_{(a,u),(b,v)}=\one_{\{a=b\}}\sigma'(l_{u,A}(\bar x^a))
                      \partial_vl_{u,A}(\bar x^a),\qquad
 \max_i\sum_k|J_{ik}|\le\rho=3/4.
\end{equation}
We prove $K_{\ell}\le\rho K_{\ell}+C_{\ell}$.

\proofstep{1. Exact response.}
For $i=(a,u)$, tilt by $\sum_{h=2}^{\ell}z_hX_{j_h}$, using separate
parameters for all formal positions. Differentiating the conditional-mean identity
\[
 \E_{\mu_z}X_i=\E_{\mu_z}\sigma\!\left(
 l_{u,A}(X^a)+s\partial_iV(X)
             +\sum_{h=2}^{\ell}z_h\one_{\{j_h=i\}}\right)
\]
using \eqref{row:eq:tilted-expectation-derivative} gives
\begin{equation}\label{audit:coordinate-response}
 \begin{aligned}
 &\cum(X_i,X_{j_2},\ldots,X_{j_{\ell}})\\
 &\quad=\sum_{S\subseteq\{2,\ldots,\ell\}}
   \left(\prod_{h\in S}\one_{\{j_h=i\}}\right)
   \cum(\sigma^{(|S|)}(l_{u,A}(X^a)+s\partial_iV(X)),X_{j_h}:h\notin S).
 \end{aligned}
\end{equation}
Repeated labels, including the anchor $i$, are allowed. All estimates
below use the original law and its fixed means.

\proofstep{2. Centered products.}
For $d\ge1$ and $t\ge0$, define, with unrestricted ordered tuples,
\begin{equation}\label{audit:normalized-product-definition}
 T_{d,t}
 :=n^{-d}\sum_{k_1,\ldots,k_d\in\mathcal I}
          \sum_{j_1,\ldots,j_t\in\mathcal I}
 \left|\cum\left(\prod_{b=1}^d\xi_{k_b},
                         X_{j_1},\ldots,X_{j_t}\right)\right|.
\end{equation}
When $t=0$, the external-label sum is absent and the one-argument
cumulant is an expectation.
Equation~\eqref{row:eq:binary-product-cumulant} and $\E\xi_k=0$ give,
for $d\ge2$,
\begin{equation}\label{audit:centered-product-split}
 \begin{aligned}
 &\cum\left(\prod_{b=1}^d\xi_{k_b},X_{j_1},\ldots,X_{j_t}\right)\\
 &=\cum\left(\prod_{b=2}^d\xi_{k_b},X_{k_1},X_{j_1},\ldots,X_{j_t}\right)\\
 &\quad+\sum_{\varnothing\ne S\subseteq[t]}
   \cum(X_{k_1},X_{j_h}:h\in S)
   \cum\left(\prod_{b=2}^d\xi_{k_b},X_{j_h}:h\notin S\right).
 \end{aligned}
\end{equation}
The subset sum is empty for $t=0$; subsets always index positions.
First, for every fixed $d,t$, the preliminary bounds imply
\begin{equation}\label{audit:preliminary-product-cumulants}
 \sup_{k_1,\ldots,k_d\in\mathcal I}\sum_{j_1,\ldots,j_t\in\mathcal I}
   |\cum(\prod_{b=1}^d\xi_{k_b},X_{j_1},\ldots,X_{j_t})|
       \le C_{d,t}(1+\log n)^{c_{d,t}} .
\end{equation}
Indeed, the case $t=0$ follows from $|\xi_k|\le1$, and the case $d=1,t\ge1$
follows from $K_{t+1}$. For $d\ge2,t\ge1$, summing
\eqref{audit:centered-product-split} gives
\[
 \begin{aligned}
 &\sup_{k_1,\ldots,k_d\in\mathcal I}\sum_{j_1,\ldots,j_t\in\mathcal I}
       |\cum(\prod_{b=1}^d\xi_{k_b},X_{j_1},\ldots,X_{j_t})|\\
 &\le\sup_{k_2,\ldots,k_d\in\mathcal I}\sum_{j_0,\ldots,j_t\in\mathcal I}
       |\cum(\prod_{b=2}^d\xi_{k_b},X_{j_0},\ldots,X_{j_t})|\\
 &\quad+\sum_{h=1}^t\binom th K_{h+1}
       \sup_{k_2,\ldots,k_d\in\mathcal I}\sum_{j_1,\ldots,j_{t-h}\in\mathcal I}
       |\cum(\prod_{b=2}^d\xi_{k_b},X_{j_1},\ldots,X_{j_{t-h}})|.
 \end{aligned}
\]
Induction on $d$, simultaneously for all fixed $t$, proves
\eqref{audit:preliminary-product-cumulants} using only the preliminary
logarithmic estimates. By \eqref{audit:preliminary-product-cumulants},
\begin{equation}\label{audit:normalized-product-preliminary}
 \begin{aligned}
 T_{1,0}&=0,\qquad
 T_{1,t}\le\frac{|\mathcal I|}{n}K_{t+1}
                                      \quad(t\ge1),\\
 T_{d,t}&\le\left(\frac{|\mathcal I|}{n}\right)^d
       \sup_{k_1,\ldots,k_d\in\mathcal I}\sum_{j_1,\ldots,j_t\in\mathcal I}
       \left|\cum\left(\prod_{b=1}^d\xi_{k_b},
                              X_{j_1},\ldots,X_{j_t}\right)\right|\\
 &\le C_{d,t}(1+\log n)^{c_{d,t}}
                                      \quad(d\ge1,\ t\ge0).
 \end{aligned}
\end{equation}
At $t=0$, \eqref{audit:centered-product-split} gives the exact identity
\begin{equation}\label{audit:centered-product-expectation}
 T_{d,0}=n^{-1}T_{d-1,1}
        \le C_dn^{-1}(1+\log n)^{c_d},\qquad d\ge2.
\end{equation}
Under the outer induction hypothesis, we next prove
\begin{equation}\label{audit:nonlinear-product-cumulants}
 T_{d,t}\le C_{\ell,d,t},\qquad d\ge2,\quad1\le t\le \ell-1.
\end{equation}
Sum \eqref{audit:centered-product-split}. Reclassifying $k_1$ as an
external label gives the first term below; bounding one factor
by $K_{h+1}$ before summing gives the second:
\begin{equation}\label{ep:row:normalized-product-recursion}
 T_{d,t}\le n^{-1}T_{d-1,t+1}
 +\frac{|\mathcal I|}{n}\sum_{h=1}^t\binom th
                        K_{h+1}T_{d-1,t-h}.
\end{equation}
For $d=2$, the endpoint $h=t$ vanishes because $T_{1,0}=0$:
\begin{equation}\label{audit:normalized-product-base}
 \begin{aligned}
 T_{2,t}
 &\le n^{-1}T_{1,t+1}
 +\frac{|\mathcal I|}{n}\sum_{h=1}^{t-1}\binom th K_{h+1}T_{1,t-h}\\
 &\le C_{\ell}n^{-1}(1+\log n)^{c_{\ell}}
 +\left(\frac{|\mathcal I|}{n}\right)^2
       \sum_{h=1}^{t-1}\binom th K_{h+1}K_{t-h+1}\le C_{\ell,2,t}.
 \end{aligned}
\end{equation}
All cumulant orders in this sum are at most $\ell-1$; for $t=1$ the
sum is empty. For $d\ge3$, \eqref{audit:normalized-product-preliminary} and
\eqref{audit:centered-product-expectation} give the endpoint bound
\begin{equation}\label{audit:normalized-product-endpoints}
 n^{-1}T_{d-1,t+1}+K_{t+1}T_{d-1,0}
 \le C_{\ell,d,t}n^{-1}(1+\log n)^{c_{\ell,d,t}}.
\end{equation}
This uses only preliminary bounds, including when $t+1=\ell$.
The terms $1\le h<t$ use $K_{h+1}\le C$ and the degree-$d-1$
induction bound for $T_{d-1,t-h}$. Since $|\mathcal I|/n\le q_0$,
\eqref{ep:row:normalized-product-recursion}, with
\eqref{audit:normalized-product-base} and \eqref{audit:normalized-product-endpoints},
proves \eqref{audit:nonlinear-product-cumulants}. Every cumulant outside the
outer induction range retains a factor $n^{-1}$.

\proofstep{3. Complete-field expansion.}
Fix $i=(a,u)$ and $0\le m\le \ell-1$. Set
$H_i=l_{u,A}(X^a)-l_{u,A}(\bar x^a)$.
The reference field has fixed degree $D$; \eqref{row:eq:dense-derivatives}
gives the coefficient bound in its exact expansion
\begin{equation}\label{audit:centered-reference-field}
 H_i=\sum_{d=1}^D\sum_{|J|=d}
       (\partial_Jl_{u,A})(\bar x^a)\prod_{j\in J}\xi_j,
 \qquad |(\partial_Jl_{u,A})(\bar x^a)|\le C_dn^{-d}.
\end{equation}
Derivatives outside row $a$ vanish, and
\begin{equation}\label{audit:reference-field-coefficient-sum}
 \sum_{d=1}^D\sum_{|J|=d}|(\partial_Jl_{u,A})(\bar x^a)|
 \le\sum_{d=1}^DC_dn^{-d}\binom{|\mathcal I|}{d}\le C.
\end{equation}
By \eqref{eq:local-interrow-coefficients},
$\partial_iV(X)=\sum_Sa_{i,S}\prod_{v\in S}X_v$, where
$|S|\le D_*$ and $\sum_S|a_{i,S}|\le C/n$.
Substituting $X_v=\bar x_v+\xi_v$ gives
$s\partial_iV(X)=\sum_Tu_{i,T}\prod_{v\in T}\xi_v$ with
$u_{i,T}=s\sum_{S\supseteq T}a_{i,S}\prod_{v\in S\setminus T}\bar x_v$.
Thus
\begin{equation}\label{audit:centered-interaction-coefficient-sum}
 \sum_T|u_{i,T}|
 \le\sum_S|a_{i,S}|\prod_{v\in S}(1+|\bar x_v|)
 \le2^{D_*}\sum_S|a_{i,S}|\le C/n.
\end{equation}
Taylor's formula through order $2\ell$ yields
\begin{equation}\label{audit:complete-field-taylor}
 \begin{aligned}
 &\sigma^{(m)}(l_{u,A}(X^a)+s\partial_iV(X))\\
 &\quad=\sum_{j=0}^{2\ell}\frac{\sigma^{(m+j)}(l_{u,A}(\bar x^a))}{j!}
                         (H_i+s\partial_iV(X))^{j}+R_{i,m},\\
 &|R_{i,m}|\le C_{\ell}|H_i+s\partial_iV(X)|^{2\ell+1}.
 \end{aligned}
\end{equation}
Separate the constant and reference linear terms from the pure-reference
terms of degrees $2\le d\le2\ell D$ and the terms containing an interaction
factor. Substituting \eqref{audit:centered-reference-field} and the interaction
expansion in \eqref{audit:centered-interaction-coefficient-sum} into
\eqref{audit:complete-field-taylor}, with ordered-tuple coefficients
$c_{i,m,\mathbf k}$ and $b_{i,m,\mathbf k}$, gives
\begin{equation}\label{audit:field-coefficient-budgets}
 \begin{aligned}
 \sigma^{(m)}(l_{u,A}(X^a)+s\partial_iV(X))
 &=\sigma^{(m)}(l_{u,A}(\bar x^a))
   +\sigma^{(m+1)}(l_{u,A}(\bar x^a))
                 \sum_k\partial_kl_{u,A}(\bar x^a)\xi_k\\
 &\quad+\sum_{d=2}^{2\ell D}\sum_{\mathbf k\in\mathcal I^d}
                    c_{i,m,\mathbf k}\prod_{b=1}^d\xi_{k_b}\\
 &\quad+\sum_{d=0}^{2\ell\max\{D,D_*\}}\sum_{\mathbf k\in\mathcal I^d}
                    b_{i,m,\mathbf k}\prod_{b=1}^d\xi_{k_b}+R_{i,m},\\
 |c_{i,m,\mathbf k}|&\le C_{\ell}n^{-d}\qquad(\mathbf k\in\mathcal I^d).
 \end{aligned}
\end{equation}
Distribute subset coefficients among their orderings. A pure degree-$d$
term has coefficient $O(n^{-d})$ and admits only finitely many allocations,
proving the displayed bound for $c$. By \eqref{audit:reference-field-coefficient-sum}
and \eqref{audit:centered-interaction-coefficient-sum}, submultiplicativity gives
\begin{equation}\label{audit:interaction-taylor-coefficient-budget}
 \sum_{d,\mathbf k}|b_{i,m,\mathbf k}|
 \le C_{\ell}\sum_{j=1}^{2\ell}\sum_{h=1}^{j}
       \binom jh C^{j-h}(C/n)^h\le C_{\ell}/n.
\end{equation}
The degree-zero term uses the empty product. Repeated labels remain
separate factors; no power reduction changes the formal degrees.

For $1\le t\le \ell-1$, \eqref{audit:field-coefficient-budgets},
\eqref{audit:normalized-product-definition}, and \eqref{audit:nonlinear-product-cumulants}
give the pure-reference error bound
\begin{equation}\label{audit:pure-reference-cumulant-error}
 \begin{aligned}
 &\sum_{j_1,\ldots,j_t}
 \left|\cum\left(
  \sum_{d=2}^{2\ell D}\sum_{\mathbf k\in\mathcal I^d}
        c_{i,m,\mathbf k}\prod_{b=1}^d\xi_{k_b},
                    X_{j_1},\ldots,X_{j_t}\right)\right|\\
 &\quad\le C_{\ell}\sum_{d=2}^{2\ell D}n^{-d}
  \sum_{\mathbf k\in\mathcal I^d}\sum_{j_1,\ldots,j_t}
  \left|\cum\left(\prod_{b=1}^d\xi_{k_b},X_{j_1},\ldots,X_{j_t}\right)\right|\\
 &\quad=C_{\ell}\sum_{d=2}^{2\ell D}T_{d,t}\le C_{\ell}.
 \end{aligned}
\end{equation}
Equations~\eqref{audit:preliminary-product-cumulants} and
\eqref{audit:interaction-taylor-coefficient-budget} give
\begin{equation}\label{audit:interaction-polynomial-cumulant-error}
 \begin{aligned}
 &\sum_{j_1,\ldots,j_t}
 \left|\cum\left(
  \sum_{d=0}^{2\ell\max\{D,D_*\}}\sum_{\mathbf k\in\mathcal I^d}
        b_{i,m,\mathbf k}\prod_{b=1}^d\xi_{k_b},
                    X_{j_1},\ldots,X_{j_t}\right)\right|\\
 &\quad\le\sum_{d=1}^{2\ell\max\{D,D_*\}}
       \sum_{\mathbf k\in\mathcal I^d}|b_{i,m,\mathbf k}|
                           C_{d,t}(1+\log n)^{c_{d,t}}\\
 &\quad\le C_{\ell}n^{-1}(1+\log n)^{c_{\ell}}\le C_{\ell}.
 \end{aligned}
\end{equation}
The constant term contributes zero. The field-moment bound
\eqref{ep:row:field-moments} and $\|s\partial_iV(X)\|_\infty\le C/n$ give
\begin{equation}\label{audit:complete-field-remainder-moment}
 \E|R_{i,m}|\le C_{\ell}\E|H_i+s\partial_iV(X)|^{2\ell+1}
 \le C_{\ell}\bigl(n^{-(2\ell+1)/2}+n^{-(2\ell+1)}\bigr)
 \le C_{\ell}n^{-(2\ell+1)/2}.
\end{equation}
In \eqref{row:eq:cumulant-partition}, the block containing $R_{i,m}$ is bounded
by $\E|R_{i,m}|$ and the others by one. Hence
\begin{equation}\label{audit:field-remainder-cumulants}
 \begin{aligned}
 |\cum(R_{i,m},X_{j_1},\ldots,X_{j_t})|&\le C_t\E|R_{i,m}|,\\
 \sum_{j_1,\ldots,j_t}|\cum(R_{i,m},X_{j_1},\ldots,X_{j_t})|
 &\le C_t|\mathcal I|^t\E|R_{i,m}|
 \le C_{\ell}n^{t-(2\ell+1)/2}\le C_{\ell}/n.
 \end{aligned}
\end{equation}
The reference constant also vanishes in these cumulants. Combining
\eqref{audit:pure-reference-cumulant-error},
\eqref{audit:interaction-polynomial-cumulant-error}, and
\eqref{audit:field-remainder-cumulants} in \eqref{audit:field-coefficient-budgets} gives
\begin{equation}\label{audit:nonlinear-field-cumulants}
 \begin{aligned}
 &\sum_{j_1,\ldots,j_t}\Bigl|
  \cum\bigl(\sigma^{(m)}(l_{u,A}(X^a)+s\partial_iV(X)),
                                      X_{j_1},\ldots,X_{j_t}\bigr)\\
 &\qquad-\sigma^{(m+1)}(l_{u,A}(\bar x^a))
     \sum_k\partial_kl_{u,A}(\bar x^a)
                  \cum(X_k,X_{j_1},\ldots,X_{j_t})\Bigr|\\
 &\quad\le C_{\ell}\sum_{d=2}^{2\ell D}T_{d,t}
       +C_{\ell}n^{-1}(1+\log n)^{c_{\ell}}+C_{\ell}/n\le C_{\ell},\\
 &\quad0\le m\le \ell-1,\quad1\le t\le \ell-1.
 \end{aligned}
\end{equation}

\proofstep{4. Absorb the current order.}
For $S=\varnothing$ in \eqref{audit:coordinate-response}, the linear
coefficient is $J_{ik}$. Using \eqref{audit:nonlinear-field-cumulants} at $m=0$,
\eqref{audit:dependency:app:coordinate-jacobian-contraction}, and centering invariance gives
\begin{equation}\label{audit:coordinate-empty-contraction}
 \begin{aligned}
 &\sum_{j_2,\ldots,j_{\ell}}
                |\cum(\sigma(l_{u,A}(X^a)+s\partial_iV(X)),X_{j_2},\ldots,X_{j_{\ell}})|\\
 &\le\sum_k|J_{ik}|\left(\sum_{j_2,\ldots,j_{\ell}}
                |\cum(X_k,X_{j_2},\ldots,X_{j_{\ell}})|\right)+C_{\ell}\\
 &\le\rho K_{\ell}+C_{\ell}.
 \end{aligned}
\end{equation}
For $\varnothing\ne S\subseteq\{2,\ldots,\ell\}$, the indicators fix
$j_h=i$ for $h\in S$, leaving $t=\ell-1-|S|$ free positions. If $t\ge1$,
\eqref{audit:nonlinear-field-cumulants} at $m=|S|$ gives
\begin{equation}\label{audit:coordinate-repeated-response}
 \begin{aligned}
 &\sum_{j_2,\ldots,j_{\ell}}
 \left(\prod_{h\in S}\one_{\{j_h=i\}}\right)
 |\cum(\sigma^{(|S|)}(l_{u,A}(X^a)+s\partial_iV(X)),X_{j_h}:h\notin S)|\\
 &=\sum_{j_h:h\notin S}
          |\cum(\sigma^{(|S|)}(l_{u,A}(X^a)+s\partial_iV(X)),X_{j_h}:h\notin S)|\\
 &\le\|\sigma^{(|S|+1)}\|_\infty
       \sum_k|\partial_kl_{u,A}(\bar x^a)|K_{t+1}+C_{\ell}\le C_{\ell}.
 \end{aligned}
\end{equation}
Here $t+1=\ell-|S|\le \ell-1$ and
$\sum_k|\partial_kl_{u,A}(\bar x^a)|\le C$; the free labels may still
equal $i$. If $t=0$, there is no free-label sum:
\begin{equation}\label{audit:coordinate-no-free-labels}
 \begin{aligned}
 &\sum_{j_2,\ldots,j_{\ell}}\left(\prod_{h=2}^{\ell}\one_{\{j_h=i\}}\right)
        |\E\sigma^{(\ell-1)}(l_{u,A}(X^a)+s\partial_iV(X))|\\
 &\quad=|\E\sigma^{(\ell-1)}(l_{u,A}(X^a)+s\partial_iV(X))|
        \le\|\sigma^{(\ell-1)}\|_\infty\le C_{\ell}.
 \end{aligned}
\end{equation}
Combining \eqref{audit:coordinate-empty-contraction}--\eqref{audit:coordinate-no-free-labels}
in \eqref{audit:coordinate-response} and maximizing over $i$ gives
\begin{equation}\label{audit:coordinate-final-absorption}
 K_{\ell}\le\rho K_{\ell}+C_{\ell},\qquad
 (1-\rho)K_{\ell}\le C_{\ell},\qquad K_{\ell}\le C_{\ell}/(1-\rho).
\end{equation}
$K_{\ell}$ is finite for each $n$, so absorption is valid. Every use of a
preliminary bound beyond the outer induction range carries $n^{-1}$,
and all constants are uniform in the conditional data.
The last assertions follow from
$\cum(X_i,X_j)=\Cov(X_i,X_j)$ and
$\cum(X_i,X_j,X_k)=\E\xi_i\xi_j\xi_k$.
\end{proof}

\subsection{Local and polynomial cumulant bounds}\label{app:sharp-proof}
The coordinate cumulant bounds extend to local observables by expanding
each observable into monomials on its support. For Boolean coordinates
under a fixed law $\mu$, let
$J_\alpha=\operatorname{supp}f_\alpha
=\{i:\nabla_if_\alpha\not\equiv0\}$.

\begin{lemma}\label{ep:row:features}
Suppose that, for every fixed order $m\ge2$,
\[
 \sup_{i_1\in[d]}\sum_{i_2,\ldots,i_m\in[d]}
   |\cum_\mu(X_{i_1},\ldots,X_{i_m})|\le C_m,
\]
as in \eqref{ep:row:coordinate-bound}. Let $(f_\alpha)$ have supports
\eqref{v8:def:local-support} and satisfy
\[
 \|f_\alpha\|_\infty\le1,\qquad |J_\alpha|\le d_0,
 \qquad
 \sup_{i\in[d]}\sum_\alpha\mathbf1_{\{i\in J_\alpha\}}\le d_1,
\]
where $d_0,d_1$ are fixed. Then, for every fixed $\ell\ge2$,
\[
 \sup_{\alpha_1}\sum_{\alpha_2,\ldots,\alpha_{\ell}}
 |\cum_\mu(f_{\alpha_1},\ldots,f_{\alpha_{\ell}})|\le C_{\ell}.
\]
More generally, if $F$ is any fixed
bounded local function, then, for every fixed $\ell\ge1$,
\[
 \sum_{\alpha_1,\ldots,\alpha_{\ell}}
 |\cum_\mu(F,f_{\alpha_1},\ldots,f_{\alpha_{\ell}})|\le C_{\ell,F}.
\]
The latter constant depends on $F$ only through its supremum norm and a
fixed bound on the size of its support. Repeated coordinate and feature
labels are allowed in all sums. For each asserted order, only finitely
many of the coordinate bounds are needed; for the first assertion,
orders at most $\ell d_0$ suffice.
\end{lemma}

\begin{proof}
Fix $\ell\ge2$ and discard constant functions, whose cumulants vanish.
If none remain, the claims are trivial; otherwise $d_0,d_1\ge1$.
Write $X_S=\prod_{i\in S}X_i$ and $X_\varnothing=1$.

\proofstep{1. Boolean expansion.}
Let $\mathcal A_s$ be the distinct $s$-element subsets of all supports
$J_\alpha$, for $1\le s\le d_0$. Inclusion--exclusion gives
\begin{equation}\label{audit:local-observable-expansion}
 f_\alpha(X)=\sum_{S\subseteq J_\alpha}a_{\alpha,S}X_S,
 \qquad
 a_{\alpha,S}=\sum_{T\subseteq S}(-1)^{|S|-|T|}f_\alpha(\mathbf1_T),
 \qquad |a_{\alpha,S}|\le2^{d_0}.
\end{equation}
Indeed, evaluation at $\mathbf1_R$ uses
$\sum_{T\subseteq S\subseteq R}(-1)^{|S|-|T|}=\one_{\{T=R\}}$.
For every nonempty $S$ and $i\in S$, bounded overlap gives
\begin{equation}\label{audit:local-support-counts}
 \begin{aligned}
 \#\{\alpha:S\subseteq J_\alpha\}&\le d_1,\\
 \sup_i\sum_{S\in\mathcal A_s}\one_{\{i\in S\}}
 &\le d_1\binom{d_0-1}{s-1}\le2^{d_0}d_1.
 \end{aligned}
\end{equation}

\proofstep{2. Monomial cumulants.}
For $k\ge1$, $1\le s_j\le d_0$, and $\sum_js_j\le \ell d_0$, we prove
\begin{equation}\label{audit:local-monomial-cumulants}
 \sup_{S_1\in\mathcal A_{s_1}}
 \sum_{S_j\in\mathcal A_{s_j},\,2\le j\le k}
 |\cum_\mu(X_{S_1},\ldots,X_{S_k})|\le C
\end{equation}
by strong induction on $\sum_j(s_j-1)$, simultaneously over these
choices. Repeated supports are allowed; a sum over an empty collection of index positions
has one summand. All support sums use $S_j\in\mathcal A_{s_j}$.
The case $k=1$ follows from $|\E X_{S_1}|\le1$.
For $\sum_j(s_j-1)=0$, all supports are singletons, so the coordinate cumulant
assumption proves the bound.

For $\sum_j(s_j-1)\ge1$, choose $s_a\ge2$ and split each $S\in\mathcal A_{s_a}$
as $S=\{i(S)\}\cup S^-$, where $i(S)=\min S$. Then
\begin{equation}\label{audit:local-split-support-counts}
 \begin{gathered}
 \{i(S)\}\in\mathcal A_1,\qquad S^-\in\mathcal A_{s_a-1},\\
 \#\{S:i(S)=i\}\le2^{d_0}d_1,\qquad
 \#\{S:S^-=T\}\le2^{d_0}d_1\quad(T\ne\varnothing).
 \end{gathered}
\end{equation}
The first count fixes a contained coordinate; for the second, fix any
coordinate of $T$. The joint split map is injective because its
components have union $S$. Equation~\eqref{row:eq:binary-product-cumulant} gives
\begin{equation}\label{audit:local-monomial-split}
 \begin{aligned}
 &\cum_\mu(X_{S_a},X_{S_j}:j\ne a)\\
 &=\cum_\mu(X_{i(S_a)},X_{S_a^-},X_{S_j}:j\ne a)\\
 &\quad+\sum_{I\subseteq[k]\setminus\{a\}}
 \cum_\mu(X_{i(S_a)},X_{S_j}:j\in I)
 \cum_\mu(X_{S_a^-},X_{S_j}:j\notin I\cup\{a\}).
 \end{aligned}
\end{equation}
The first term reduces $\sum_j(s_j-1)$ by one. The product factors have
nonnegative induction parameters summing to $\sum_j(s_j-1)-1$. Their total coordinate count
does not increase, so each is covered by induction with one support fixed.

If $a\ne1$, take $a=k$ by symmetry. Injectivity allows the first term
to be summed over independent split supports:
\begin{equation}\label{audit:dependency:app:free-first-split-sum}
 \begin{aligned}
 &\sum_{S_2,\ldots,S_k}
 |\cum_\mu(X_{S_1},\ldots,X_{S_{k-1}},X_{i(S_k)},X_{S_k^-})|\\
 &\quad\le\sum_{S_2,\ldots,S_{k-1}}
 \sum_{\substack{U\in\mathcal A_1\\V\in\mathcal A_{s_k-1}}}
 |\cum_\mu(X_{S_1},\ldots,X_{S_{k-1}},X_U,X_V)|\le C.
 \end{aligned}
\end{equation}
For a product term, let $I,J$ partition $\{2,\ldots,k-1\}$.
According to which split support occurs in the factor containing $S_1$,
the two projection bounds give
\begin{equation}\label{audit:local-free-split-sum}
 \begin{aligned}
 &\sum_{\substack{S_j\in\mathcal A_{s_j}\\2\le j\le k}}
 |\cum_\mu(X_{S_1},X_{i(S_k)},X_{S_j}:j\in I)|
 |\cum_\mu(X_{S_k^-},X_{S_j}:j\in J)|\\
 &\quad\le C\sum_{S_k\in\mathcal A_{s_k}}\sum_{\substack{S_j\in\mathcal A_{s_j}\\j\in I}}
 |\cum_\mu(X_{S_1},X_{i(S_k)},X_{S_j}:j\in I)|\\
 &\quad\le C\,2^{d_0}d_1\sum_{U\in\mathcal A_1}\sum_{\substack{S_j\in\mathcal A_{s_j}\\j\in I}}
 |\cum_\mu(X_{S_1},X_U,X_{S_j}:j\in I)|\le C,\\[2pt]
 &\sum_{\substack{S_j\in\mathcal A_{s_j}\\2\le j\le k}}
 |\cum_\mu(X_{S_1},X_{S_k^-},X_{S_j}:j\in I)|
 |\cum_\mu(X_{i(S_k)},X_{S_j}:j\in J)|\\
 &\quad\le C\sum_{S_k\in\mathcal A_{s_k}}\sum_{\substack{S_j\in\mathcal A_{s_j}\\j\in I}}
 |\cum_\mu(X_{S_1},X_{S_k^-},X_{S_j}:j\in I)|\\
 &\quad\le C\,2^{d_0}d_1\sum_{V\in\mathcal A_{s_k-1}}\sum_{\substack{S_j\in\mathcal A_{s_j}\\j\in I}}
 |\cum_\mu(X_{S_1},X_V,X_{S_j}:j\in I)|\le C.
 \end{aligned}
\end{equation}
Each first inequality sums the other factor with its split support fixed;
each second inequality uses \eqref{audit:local-split-support-counts}.
The induction bounds are uniform in the fixed supports, and a
one-argument factor is an expectation bounded by one.

If $a=1$, both split supports are fixed. For the first term, enlarge
the sum to all $T\in\mathcal A_{s_1-1}$:
\begin{equation}\label{audit:dependency:app:fixed-first-split-sum}
 \begin{aligned}
 &\sum_{S_2,\ldots,S_k}
 |\cum_\mu(X_{i(S_1)},X_{S_1^-},X_{S_2},\ldots,X_{S_k})|\\
 &\quad\le\sum_T\sum_{S_2,\ldots,S_k}
 |\cum_\mu(X_{i(S_1)},X_T,X_{S_2},\ldots,X_{S_k})|\le C.
 \end{aligned}
\end{equation}
For $I,J$ partitioning $\{2,\ldots,k\}$, the product sum factors:
\begin{equation}\label{audit:local-fixed-split-sum}
 \begin{aligned}
 &\sum_{S_2,\ldots,S_k}
 |\cum_\mu(X_{i(S_1)},X_{S_j}:j\in I)|
 |\cum_\mu(X_{S_1^-},X_{S_j}:j\in J)|\\
 &=\left[\sum_{S_j,\,j\in I}
 |\cum_\mu(X_{i(S_1)},X_{S_j}:j\in I)|\right]
 \left[\sum_{S_j,\,j\in J}
 |\cum_\mu(X_{S_1^-},X_{S_j}:j\in J)|\right]\le C.
 \end{aligned}
\end{equation}
Empty $I$ or $J$ again gives a bounded expectation.
There are at most $2^{k-1}$ product terms and $k\le \ell d_0$.
Equations~\eqref{audit:dependency:app:free-first-split-sum},
\eqref{audit:local-free-split-sum}, \eqref{audit:dependency:app:fixed-first-split-sum},
and \eqref{audit:local-fixed-split-sum}, applied in \eqref{audit:local-monomial-split},
complete the proof of \eqref{audit:local-monomial-cumulants}.

\proofstep{3. Recover the local functions.}
Multilinearity and \eqref{audit:local-observable-expansion} give
\[
 \cum_\mu(f_{\alpha_1},\ldots,f_{\alpha_{\ell}})
 =\sum_{\varnothing\ne S_j\subseteq J_{\alpha_j},\,1\le j\le \ell}
       \left(\prod_j a_{\alpha_j,S_j}\right)
       \cum_\mu(X_{S_1},\ldots,X_{S_{\ell}}).
\]
For each free label, regroup by its nonempty support, using
\eqref{audit:local-observable-expansion} and \eqref{audit:local-support-counts}:
\begin{equation}\label{audit:local-support-rearrangement}
 \begin{aligned}
 &\sum_{\alpha_j}\sum_{\substack{S_j\subseteq J_{\alpha_j}\\|S_j|=s_j}}
 |a_{\alpha_j,S_j}|\,|\cum_\mu(X_{S_1},\ldots,X_{S_{\ell}})|\\
 &=\sum_{S_j\in\mathcal A_{s_j}}
   \left(\sum_{\alpha_j:S_j\subseteq J_{\alpha_j}}|a_{\alpha_j,S_j}|\right)
   |\cum_\mu(X_{S_1},\ldots,X_{S_{\ell}})|\\
 &\le2^{d_0}d_1\sum_{S_j\in\mathcal A_{s_j}}
   |\cum_\mu(X_{S_1},\ldots,X_{S_{\ell}})|.
 \end{aligned}
\end{equation}
Fix $\alpha_1$, apply \eqref{audit:local-support-rearrangement} to $j=2,\ldots,\ell$,
and use
\eqref{audit:local-monomial-cumulants}. There are at most $2^{d_0}-1$
choices of $S_1$ and $d_0^{\ell-1}$ degree choices, giving the first
assertion. Splitting preserves total coordinate count, so only
coordinate cumulants through order $\ell d_0$ are used.

For a nonconstant $F$, adjoin $F/\|F\|_\infty$ to the family.
Its support-size and overlap bounds are
$\max\{d_0,|\operatorname{supp}F|\}$ and $d_1+1$.
Fixing this new function as the first argument gives
\begin{equation}\label{v8:proof:ep:row:features:fixed-local-anchor}
 \sum_{\alpha_1,\ldots,\alpha_{\ell}}
 |\cum_\mu(F,f_{\alpha_1},\ldots,f_{\alpha_{\ell}})|
 \le C_{\ell,d_0,d_1,|\operatorname{supp}F|}\|F\|_\infty.
\end{equation}
For constant $F$, these cumulants vanish.
\end{proof}

\begin{proof}[Proof of Theorem~\ref{ep:row:sharp}]
Lemmas~\ref{ep:row:coordinates} and \ref{ep:row:features} give
\eqref{ep:row:feature-bound} and its fixed-anchor version.
We suppress the conditional law from the notation. The block features
$\chi_\alpha=\chi_{(i,S)}$ form a bounded-overlap family with $O(n)$
labels. For positive degrees $s_1,\ldots,s_{\ell}$ and $L=\sum_bs_b$, we prove
\begin{equation}\label{audit:block-product-cumulant-bounds}
 \begin{aligned}
 \sum_{\boldsymbol\alpha}
 \left|\cum\left(\prod_{k=1}^{s_1}\chi_{\alpha_{1,k}},
              \ldots,\prod_{k=1}^{s_{\ell}}\chi_{\alpha_{\ell,k}}\right)\right|
 &\le Cn^{L-\ell+1},\\
 \sum_{\boldsymbol\alpha}
 \left|\cum\left(F,\prod_{k=1}^{s_1}\chi_{\alpha_{1,k}},
              \ldots,\prod_{k=1}^{s_{\ell}}\chi_{\alpha_{\ell,k}}\right)\right|
 &\le C_Fn^{L-\ell}.
 \end{aligned}
\end{equation}
Here all $L$ formal label positions are summed independently, allowing
repetitions. This enlarges the nonnegative sums from the polynomial
expansions. Empty collections of positions have one summand.
Equations~\eqref{ep:row:feature-bound} and
\eqref{v8:proof:ep:row:features:fixed-local-anchor}, with $|\chi_\alpha|\le1$, give the base cases
\begin{equation}\label{audit:block-product-base-cases}
 \begin{aligned}
 \sum_{\alpha_1,\ldots,\alpha_{\ell}}
 |\cum(\chi_{\alpha_1},\ldots,\chi_{\alpha_{\ell}})|
 &\le C_{\ell}n,&&\ell\ge2,\\
 \sum_{\alpha_1,\ldots,\alpha_{\ell}}
 |\cum(F,\chi_{\alpha_1},\ldots,\chi_{\alpha_{\ell}})|
 &\le C_{F,\ell},&&\ell\ge1,\\
 \sum_{\alpha_1,\ldots,\alpha_L}
 \left|\E\prod_{k=1}^L\chi_{\alpha_k}\right|
 &\le(Cn)^L,&&L\ge1.
 \end{aligned}
\end{equation}

Induct on $L-\ell$, simultaneously over degrees of bounded total size.
For the unanchored bound, the cases $L-\ell=0$ and $\ell=1$ are covered above.
Otherwise take $s_1\ge2$ by symmetry and put
$U=\chi_{\alpha_{1,1}}$,
$V=\prod_{k=2}^{s_1}\chi_{\alpha_{1,k}}$,
$Y_b=\prod_{k=1}^{s_b}\chi_{\alpha_{b,k}}$ for $b\ge2$.
Equation~\eqref{row:eq:binary-product-cumulant} gives
\begin{equation}\label{audit:block-product-split}
 \begin{aligned}
 \cum(UV,Y_2,\ldots,Y_{\ell})
 &=\cum(U,V,Y_2,\ldots,Y_{\ell})\\
 &\quad+\sum_{S\subseteq\{2,\ldots,\ell\}}
       \cum(U,Y_b:b\in S)\cum(V,Y_b:b\notin S).
 \end{aligned}
\end{equation}
Complements here and below are within $\{2,\ldots,\ell\}$.
The first term reduces $L-\ell$ by one and sums to $O(n^{L-\ell})$.
For a product term, set
\[
 L_1=1+\sum_{b\in S}s_b,\quad \ell_1=1+|S|,\qquad
 L_2=L-L_1,\quad \ell_2=\ell+1-\ell_1.
\]
The two nonnegative induction parameters sum to $L-\ell-1$.
Their formal positions are disjoint, so the finite sums factor:
\begin{equation}\label{audit:block-product-free-sum}
 \sum_{\boldsymbol\alpha}
 |\cum(U,Y_b:b\in S)\cum(V,Y_b:b\notin S)|
 \le Cn^{L_1-\ell_1+1}n^{L_2-\ell_2+1}=Cn^{L-\ell+1}.
\end{equation}
The bases include one-argument factors. Summing the finitely many
subsets in \eqref{audit:block-product-split}, using
\eqref{audit:block-product-free-sum}, proves the unanchored bound.

For the anchored bound, use the same induction on $L-\ell$, counting $F$ in neither
$L$ nor $\ell$. Its base case $L-\ell=0$ was established above, and the split is
\begin{equation}\label{audit:block-product-anchored-split}
 \begin{aligned}
 &\cum(F,UV,Y_2,\ldots,Y_{\ell}) =\cum(F,U,V,Y_2,\ldots,Y_{\ell})\\
 &\quad+\sum_{S\subseteq\{2,\ldots,\ell\}}
 \left\{\begin{aligned}[t]
 &\cum(F,U,Y_b:b\in S)\cum(V,Y_b:b\notin S)\\
 &\quad+\cum(U,Y_b:b\in S)\cum(F,V,Y_b:b\notin S)
 \end{aligned}\right\}.
 \end{aligned}
\end{equation}
The first term sums to $O(n^{L-\ell-1})$. Anchored induction and the
unanchored bound give
\begin{equation}\label{audit:block-product-anchored-sum}
 \sum_{\boldsymbol\alpha}
 |\cum(F,U,Y_b:b\in S)\cum(V,Y_b:b\notin S)|
 \le C_Fn^{L_1-\ell_1}n^{L_2-\ell_2+1}=C_Fn^{L-\ell}.
\end{equation}
Interchanging the anchored factor covers the other product, including
$\ell=1$, whose products are $\Cov(F,U)\E V$ and $\E U\Cov(F,V)$.
Combining \eqref{audit:block-product-anchored-split} with
\eqref{audit:block-product-anchored-sum} proves \eqref{audit:block-product-cumulant-bounds};
the factorization
follows from finite sums over formal positions and does not require probabilistic independence.

Expand $G_b$ into block monomials of degrees at most $d_b$.
Constants vanish in the relevant cumulants, and a degree choice with
total $L$ has a product of coefficients of order $O(n^{-L})$ by
\eqref{v8:def:block-coefficient-normalization}. Applying
\eqref{audit:block-product-cumulant-bounds} gives
\begin{equation}\label{audit:block-polynomial-cumulant-scaling}
 \begin{aligned}
 |\cum(G_1,\ldots,G_{\ell})|
 &\le C\sum_{1\le s_b\le d_b}n^{-L}n^{L-\ell+1}
 \le C_{\ell}n^{1-\ell},&&\ell\ge2,\\
 |\cum(F,G_1,\ldots,G_{\ell})|
 &\le C_F\sum_{1\le s_b\le d_b}n^{-L}n^{L-\ell}
 \le C_{F,\ell}n^{-\ell},&&\ell\ge1.
 \end{aligned}
\end{equation}
There are finitely many degree choices, and degrees count blocks.
These are \eqref{ep:row:polynomial-bound} and
\eqref{ep:row:anchored-polynomial}, uniformly in the conditional data.
\end{proof}

\section{Stein solution estimates and exchangeable-pair bounds}
\label{app:stein-proofs}
We prove the analytic estimates and the exchangeable-pair theorem stated
in Section~\ref{sec:stein}. The potentials, target distributions and compact
parameter set are those of \eqref{stein:quartic-potential}.

\subsection{Uniform bounds for Stein solutions}\label{app:stein-solutions}
The Stein equation for $Y$ is
\begin{equation}\label{stein:equation}
 f_h'(w)-V'(w)f_h(w)=h(w)-\E h(Y).
\end{equation}
For indicator tests $h_z(w)=\one_{\{w\le z\}}$, $1$-Lipschitz functions
with $h(0)=0$, and, when $\mc\ge2$, the quadratic test $h(w)=w^2$,
the bounded solution is
\begin{equation}\label{stein:solution-formula}
 \begin{aligned}
 f_h(w)&=\rho(w)^{-1}\int_{-\infty}^w
                    \{h(t)-\E h(Y)\}\rho(t)\,\dd t\\
       &=-\rho(w)^{-1}\int_w^\infty
                    \{h(t)-\E h(Y)\}\rho(t)\,\dd t.
 \end{aligned}
\end{equation}
For indicator tests $h_z(w)=\one_{\{w\le z\}}$, we set
$f_{h_z}'(z)=V'(z)f_{h_z}(z)+1-\Pp(Y\le z)$; elsewhere,
$f_h'$ denotes the ordinary derivative.

\begin{lemma}\label{stein:solutions}
For the potentials in \eqref{stein:quartic-potential} and the tests
specified above, the solutions in \eqref{stein:solution-formula} satisfy
\begin{equation}\label{stein:solution-norms}
 \|f_h\|_\infty+\Lip f_h+\|f_h'\|_\infty\le C_{\mathcal K}.
\end{equation}
For Lipschitz tests, and for the quadratic test when $\mc\ge2$,
$f_h'$ is locally absolutely continuous and
\begin{equation}\label{stein:second-derivative}
 |f_h''(w)|\le C_{\mathcal K}(1+|w|^{2\mc-1})\quad\text{for almost every }w.
\end{equation}
For indicators, $V'(w)f_{h_z}(w)$ is locally absolutely continuous and
\begin{equation}\label{stein:indicator-continuous-part}
 \left|\frac{\dd}{\dd w}\{V'(w)f_{h_z}(w)\}\right|\le C_{\mathcal K}(1+|w|^{2\mc-1})
                 \quad\text{for almost every }w.
\end{equation}
The constants for indicator tests are independent of $z$.
Each solution is locally absolutely continuous, with increments given by
integrals of $f_h'$.
\end{lemma}

\begin{proof}
Compactness of $\mathcal K$ gives, uniformly in $\gamma\in \mathcal K$,
\[
 V_\gamma(w)\ge c_{\mathcal K}|w|^{2\mc}-C_{\mathcal K},\qquad
 \sup_{\gamma\in \mathcal K,\,|w|\le R}|V_\gamma(w)|\le C_{\mathcal K,R}
 \quad(R<\infty).
\]
Integrating over $[-1,1]$ for the lower bound on $\mathcal Z_\gamma$
and using the first inequality for its upper bound yields
\begin{equation}\label{stein:quartic-common-bounds}
 \begin{gathered}
 2e^{-C_{\mathcal K,1}}\le\mathcal Z_\gamma\le C_{\mathcal K},\qquad
 \rho_\gamma(w)\le C_{\mathcal K}e^{-c_{\mathcal K}|w|^{2\mc}},\\
 \inf_{\gamma\in \mathcal K,\,|w|\le R}\rho_\gamma(w)\ge c_{\mathcal K,R}>0,
 \qquad \sup_{\gamma\in \mathcal K}\E|Y_\gamma|^j\le C_{\mathcal K,j}
 \quad(j\ge0\text{ fixed}).
 \end{gathered}
\end{equation}

The local lower bound in \eqref{stein:quartic-common-bounds} also gives a
uniform positive variance. Indeed, the independent-copy identity yields
\[
 \begin{aligned}
 \Var(Y_\gamma)
 &=\frac12\int_{\R^2}(u-v)^2\rho_\gamma(u)\rho_\gamma(v)\,\dd u\,\dd v\\
 &\ge\frac12\int_{-1}^{-1/2}\int_{1/2}^{1}
       \rho_\gamma(u)\rho_\gamma(v)\,\dd v\,\dd u
 \ge c_{\mathcal K}>0.
 \end{aligned}
\]
The second-moment bound gives the uniform upper bound on the variance.

The centered integral in \eqref{stein:solution-formula} gives a locally
absolutely continuous solution of \eqref{stein:equation}, with the
threshold convention stated above.

Suppose first that $\mc\ge2$. For some $R\ge1$ depending only on $\mathcal K$,
$V(w+s)-V(w)\ge c_{\mathcal K}w^{2\mc-1}s$ whenever $w\ge R$ and $s\ge0$.
Hence, for every fixed integer $j\ge0$,
\begin{equation}\label{stein:even-tail-integrals}
 \int_0^\infty s^j e^{V(w)-V(w+s)}\,\dd s
 \le C_{\mathcal K,j}w^{-(2\mc-1)(j+1)},\qquad w\ge R.
\end{equation}
For indicators, $|h_z-\E h_z(Y)|\le1$, so
\eqref{stein:solution-formula}, \eqref{stein:even-tail-integrals}, and
\eqref{stein:equation} give, uniformly in $z$,
\[
 |f_{h_z}(w)|\le C_{\mathcal K}w^{1-2\mc},\qquad
 |f_{h_z}'(w)|\le |V'(w)f_{h_z}(w)|+1\le C_{\mathcal K}.
\]
For the Lipschitz and quadratic tests,
$|h(t)-\E h(Y)|\le C_{\mathcal K}(1+t^2)$ and
$|h'(t)|\le C(1+|t|)$ almost everywhere. Thus
\eqref{stein:solution-formula} and \eqref{stein:even-tail-integrals} yield
\[
 |f_h(w)|\le C_{\mathcal K}\int_0^\infty(w+s)^2
               e^{V(w)-V(w+s)}\,\dd s\le C_{\mathcal K}w^{3-2\mc}.
\]
Integration by parts in the upper-tail expression, together with
\eqref{stein:equation}, gives
\begin{equation}\label{stein:even-tail-cancellation}
 \begin{aligned}
 f_h'(w)
 &=-\int_0^\infty h'(w+s)e^{V(w)-V(w+s)}\,\dd s\\
 &\quad+\int_0^\infty[h(w+s)-\E h(Y)]
       [V'(w+s)-V'(w)]e^{V(w)-V(w+s)}\,\dd s.
 \end{aligned}
\end{equation}
The boundary term at infinity vanishes by \eqref{stein:quartic-common-bounds}.
Using $|V'(w+s)-V'(w)|\le C_{\mathcal K}s(w+s)^{2\mc-2}$ and
\eqref{stein:even-tail-integrals}, we obtain
\[
 |f_h'(w)|\le C_{\mathcal K}\int_0^\infty
       [(w+s)+s(w+s)^{2\mc}]e^{V(w)-V(w+s)}\,\dd s
       \le C_{\mathcal K}w^{2-2\mc}.
\]
The reflected potential $V(-w)$ satisfies the same coefficient bounds,
so the same argument treats the left tail. On $[-R,R]$,
\eqref{stein:quartic-common-bounds}, \eqref{stein:solution-formula}, and
\eqref{stein:equation} bound $f_h$ and $f_h'$ uniformly. Local absolute
continuity then gives $\Lip f_h\le\|f_h'\|_\infty$, proving
\eqref{stein:solution-norms}.

For the continuous tests, $f_h'$ is locally absolutely continuous by
\eqref{stein:equation}; for indicators, $V'f_{h_z}$ is locally absolutely
continuous. Almost everywhere,
\[
 f_h''=V''f_h+V'f_h'+h',\qquad
 (V'f_{h_z})'=V''f_{h_z}+V'f_{h_z}'.
\]
Together with \eqref{stein:solution-norms} and the polynomial degrees,
these identities give \eqref{stein:second-derivative} and
\eqref{stein:indicator-continuous-part}. Local absolute continuity gives
the increment formula.

For $\mc=1$, the affine change of variable $\sqrt a(w+c_1/a)$ reduces
\eqref{stein:equation} to the standard normal Stein equation.
The classical bounds \cite[Lemmas~2.3--2.4]{CGS2011}, as used in
\cite[Section~4]{FLSZ}, give, uniformly over $\gamma\in \mathcal K$,
\[
 \begin{gathered}
 \|f_h\|_\infty+\|f_h'\|_\infty\le C_{\mathcal K},\qquad
 \Lip(f_h')\le C_{\mathcal K}\quad(\Lip h\le1),\\
 |(V'f_{h_z})'(w)|\le C_{\mathcal K}(1+|w|)\quad(w\ne z).
 \end{gathered}
\]
These imply \eqref{stein:solution-norms}--\eqref{stein:indicator-continuous-part}
for the admissible tests. Local absolute continuity gives the increment
formula, with the same threshold convention.
\end{proof}

For the target distributions in \eqref{eq:target}, condition
\eqref{even:calibration-system} places the coefficient vectors in a fixed
compact set of the form specified in \eqref{stein:quartic-potential}.
The bounds in \eqref{stein:quartic-common-bounds} and the variance estimate
in the preceding proof therefore give, for every fixed $q\ge1$,
\begin{equation}\label{even:target-basic-bounds}
 \sup_n\|\rho_n\|_\infty\le C,\qquad
 \sup_n\E|Y_n|^q\le C_q,\qquad
 0<c\le\Var(Y_n)\le C.
\end{equation}
The constants depend on the fixed leading model and the calibration
bound, and, for the moment estimate, on the moment order.

\subsection{Proof of the exchangeable-pair approximation bound}

\begin{proof}[Proof of Theorem~\ref{prop:Stein}]
Let $f_h$ be the solution of \eqref{stein:equation} given by
\eqref{stein:solution-formula}, with the threshold convention stated there.
For any bounded Lipschitz function $f$, exchangeability gives
\begin{equation}\label{stein:exchangeability-identity}
 \E[\Delta\{f(W')-f(W)\}]=-2\E[\Delta f(W)].
\end{equation}
Conditioning on $X$ also gives the equality in
\eqref{stein:weak-regression}, with
$\mathcal R\le\delta/\Lambda+\E|V'(W)|<\infty$.
The increment admits the expansion
\[
 f_h(W')-f_h(W)=\Delta f_h'(W)
       +\int_0^\Delta
          \{f_h'(W+t)-f_h'(W)\}\,\dd t.
\]
Substituting this into \eqref{stein:exchangeability-identity} and
conditioning on $X$ yields the Stein error decomposition
\begin{equation}\label{stein:exact-error-identity}
 \begin{aligned}
 \E h(W)-\E h(Y)
 &=\E[(1-D)f_h'(W)]\\
 &\quad-\E\left[
   \left\{\frac{\E[\Delta\mid X]}\Lambda+V'(W)\right\}
                   f_h(W)\right]\\
 &\quad-\frac1{2\Lambda}\E\left[\Delta\int_0^\Delta
        \{f_h'(W+t)-f_h'(W)\}\,\dd t\right].
 \end{aligned}
\end{equation}
By \eqref{stein:solution-norms} and \eqref{stein:weak-regression},
the first two terms are bounded in absolute value by
$C_{\mathcal K}\E|1-D|$ and $C_{\mathcal K}\mathcal R$, respectively.

For the Wasserstein bound, take $h$ to be $1$-Lipschitz with $h(0)=0$.
When $\mc\ge2$, also take $h(w)=w^2$ for the second-moment estimate.
Equation~\eqref{stein:second-derivative}
and $\delta\le1$ give
\[
 |f_h'(W+t)-f_h'(W)|
 \le C_{\mathcal K}|t|(1+|W|^{2\mc-1}),\qquad |t|\le\delta.
\]
Hence the remainder satisfies
\begin{equation}\label{audit:dependency:app:stein-smooth-remainder}
 \begin{aligned}
 \frac1{2\Lambda}\left|\E\left[\Delta\int_0^\Delta
       \{f_h'(W+t)-f_h'(W)\}\,\dd t\right]\right|
 &\le\frac{C_{\mathcal K}}{\Lambda}\E[|\Delta|^3(1+|W|^{2\mc-1})]\\
 &\le\frac{C_{\mathcal K}\delta^3}{\Lambda}\E(1+|W|^{2\mc-1})<\infty,
 \end{aligned}
\end{equation}
where the second inequality uses $|\Delta|\le\delta$ and the moment assumption.
Taking the supremum over the Lipschitz class proves the Wasserstein bound.
For $\mc\ge2$, choosing $h(w)=w,w^2$ also gives the moment estimates.

For the Kolmogorov bound, we first establish a concentration inequality.
For $u\le v$, define the piecewise linear function
\[
 f_{u,v}(w)=\int_{(u+v)/2}^w
               \one_{\{u-\delta\le t\le v+\delta\}}\,\dd t.
\]
It is nondecreasing and satisfies
\[
 \begin{gathered}
 \|f_{u,v}\|_\infty\le\frac{v-u+2\delta}{2},\qquad
 \Lip f_{u,v}\le1,\\
 \Delta\{f_{u,v}(W')-f_{u,v}(W)\}
       \ge\Delta^2\one_{\{u\le W\le v\}}.
 \end{gathered}
\]
Thus \eqref{stein:exchangeability-identity} and
\eqref{stein:weak-regression} imply
\begin{equation}\label{eq:short-interval}
 \begin{aligned}
 \E[D\one_{\{u\le W\le v\}}]
 &\le-\Lambda^{-1}\E[\Delta f_{u,v}(W)]\\
 &\le\frac{v-u+2\delta}{2}\E|V'(W)|
       +\left(1+\frac{v-u+2\delta}{2}\right)\mathcal R.
 \end{aligned}
\end{equation}

Now take $h=h_z$, where $h_z(w)=\one_{\{w\le z\}}$. By
\eqref{stein:equation},
\[
 \begin{aligned}
 f_{h_z}'(W+t)-f_{h_z}'(W)
 &=V'(W+t)f_{h_z}(W+t)-V'(W)f_{h_z}(W)\\
 &\quad+h_z(W+t)-h_z(W).
 \end{aligned}
\]
By \eqref{stein:indicator-continuous-part}, the contribution of the first
difference is bounded by the right-hand side of
\eqref{audit:dependency:app:stein-smooth-remainder}. For the indicator difference,
\begin{equation}\label{audit:dependency:app:stein-jump-remainder}
 \left|\Delta\int_0^\Delta
       \{h_z(W+t)-h_z(W)\}\,\dd t\right|
 \le\Delta^2\one_{\{z-\delta\le W\le z+\delta\}}.
\end{equation}
Applying \eqref{eq:short-interval} with $u=z-\delta$ and $v=z+\delta$ gives
\[
 \begin{aligned}
 &\frac1{2\Lambda}\left|\E\left[\Delta\int_0^\Delta
          \{h_z(W+t)-h_z(W)\}\,\dd t\right]\right|\\
 &\quad\le\E[D\one_{\{z-\delta\le W\le z+\delta\}}]
 \le2\delta\E|V'(W)|+(1+2\delta)\mathcal R.
 \end{aligned}
\]
These estimates hold for every $z$, including atoms of the law of $W$,
by the threshold convention and the closed intervals used above. Combining them with
\eqref{stein:exact-error-identity} and taking the supremum over $z$
proves the Kolmogorov estimate in \eqref{eq:Stein-bound}.

\end{proof}

\section{Solvability of the calibration equations}\label{app:calibration-proof}
We prove Proposition~\ref{even:bounded-calibration} under its stated
hypotheses, using the notation of \eqref{even:calibration-system}.

\begin{proof}[Proof of Proposition~\ref{even:bounded-calibration}]
We first prove that the coefficient matrix has full row rank, then use a fixed
right inverse to select a bounded solution. By \eqref{eq:critical},
$2B^{(j)}(p)=I^{(j)}(p)$ for $2\le j\le2\mc-1$.
Suppose $D=\{e_r:b_r>0\}$ has at most $\mc-2$ elements and set
$Q(z)=\prod_{d\in D}(z-d+2)$. Then
$\deg(Q^2)\le2\mc-4$, $Q(e_r-2)=0$ for $b_r>0$, and $Q(-1)\ne0$.
Near $t=1$,
\begin{equation}\label{audit:entropy-derivative-series}
 p(1-p)I''(pt)=\frac{1-p}{t(1-pt)}
      =(1-p)\sum_{\ell=0}^{\infty}p^\ell t^{\ell-1}.
\end{equation}
Its derivatives through order $2\mc-3$ at $t=1$ agree with those of
$2p(1-p)B''(pt)$. Write $(z)_d=\prod_{a=0}^{d-1}(z-a)$ for the
falling factorial, with $(z)_0=1$. Expanding $Q^2$ in this basis and
differentiating \eqref{audit:entropy-derivative-series} gives
\[
 \begin{aligned}
 0&=2p(1-p)\sum_r b_re_r(e_r-1)p^{e_r-2}Q(e_r-2)^2\\
  &=(1-p)\sum_{\ell\ge0}p^\ell Q(\ell-1)^2
    \ge(1-p)Q(-1)^2>0,
 \end{aligned}
\]
a contradiction. Hence, after relabeling, choose
$b_1,\ldots,b_{\mc-1}>0$ with $e_1<\cdots<e_{\mc-1}$.

Let $A\in\R^{\mc\times(k+1)}$ be the coefficient matrix of the linear
equations in \eqref{even:calibration-system}, with columns ordered as
$(\kappa_n,u_{1,n},\ldots,u_{k,n})$. Its entries are
\begin{equation}\label{even:calibration-determinant}
 A_{j,1}=\frac12\mathbf1_{\{j=1\}},\qquad
 A_{j,r+1}=\left.\frac{\dd^j}{\dd x^j}x^{e_r}\right|_{x=p},
 \qquad 1\le j\le\mc,\quad 1\le r\le k.
\end{equation}
For $2\le j\le\mc$ and $1\le i\le\mc-1$,
\[
 A_{j,i+1}
 =e_i(e_i-1)p^{e_i-2}\,p^{2-j}(e_i-2)_{j-2}.
\]
Expanding the determinant of the square submatrix formed by the first $\mc$ columns
along the field column and factoring the displayed column and row terms gives
\[
 \frac12\left\{\prod_{i=1}^{\mc-1}e_i(e_i-1)p^{e_i-2}\right\}
   p^{-(\mc-1)(\mc-2)/2}
   \prod_{1\le i<r\le\mc-1}(e_r-e_i)>0.
\]
Here the last product is the Vandermonde determinant of the monic
polynomials $(z)_{j-2}$ evaluated at $e_i-2$. Thus
$\operatorname{rank}A=\mc$, and $AA^{\mathsf T}$ is positive definite.

Writing indexed vectors as column vectors, take
\begin{equation}\label{even:calibration-explicit}
 (\kappa_n,u_{1,n},\ldots,u_{k,n})^{\mathsf T}
 =A^{\mathsf T}(AA^{\mathsf T})^{-1}
   \bigl(n^{-1+j/\mc}\beta_{j,n}
                  -\mathcal F_1^{(j)}(p)\bigr)_{j=1}^{\mc}.
\end{equation}
Multiplication by $A$ verifies the equations in
\eqref{even:calibration-system}. Since the right inverse is fixed,
\eqref{even:calibration-explicit} and the bound in
\eqref{even:calibration-system} give
\[
 |\kappa_n|+\sum_r|u_{r,n}|
 \le C\left(1+\sum_{j=1}^{\mc}n^{-1+j/\mc}|\beta_{j,n}|\right)
 \le C(1+L),
\]
which proves \eqref{even:inverse-calibration-bound}.
All solutions are obtained by adding a sequence of vectors in $\ker A$
to \eqref{even:calibration-explicit}; uniformly bounded additions preserve
boundedness. For $\boldsymbol\beta_n=0$, the solution
\eqref{even:calibration-explicit} is independent of $n$.
\end{proof}

\medskip
\noindent\textbf{Concurrent work.}
While this work was in preparation, Liu, Shao, Shi and Xu~\cite{LSSX2026}
posted their preprint on arXiv on 30 October 2026. Their results are specific
to the edge--triangle model. They obtain quantitative Gaussian and quartic
approximations for the edge density using a Dobrushin contraction bound
for the conditional row laws, uniform over all exterior graphs for parameters
near the critical point~\cite[Lemma~4.2]{LSSX2026}.
Our results cover general fixed simple interaction graphs at any fixed even
critical order $2\mc\ge4$ under the stated single-phase assumptions.
For a set of exterior graphs of probability $1-O(e^{-cn})$, we construct
auxiliary row-update kernels that contract uniformly on the full row state
space and agree with the Gibbs updates on typical row configurations
(Lemma~\ref{lem:local-input}).

\bookmarksetup{depth=-1}
\begingroup
\renewcommand{\addcontentsline}[3]{}
\section*{Acknowledgments}
The authors thank Xiao Fang for helpful discussions.

\section*{AI Statement}
OpenAI's ChatGPT and Codex assisted with English revision, LaTeX typesetting,
and the presentation and structure of proofs. The tools additionally supported
checks of bibliographic details, notation, internal cross-references, and
derivations. The authors subsequently examined and amended the manuscript,
checked the mathematical reasoning, and accept responsibility for all of its
contents.

\par
\endgroup

\end{document}